\documentclass[reqno, 10pt]{amsart}

\usepackage[usenames,dvipsnames]{color}
\usepackage[pdfauthor={Andrea Appel, Valerio Toledano Laredo}, pdftitle={Uniqueness of Coxeter structures on Kac-Moody algebras}]{hyperref} 
\hypersetup{
    colorlinks=true,       			
    linkcolor=blue,          			
    citecolor=blue,		 		
    filecolor=blue,      				
    urlcolor=blue           			
}
\usepackage{amsmath, amssymb, empheq, stmaryrd}
\usepackage{tikz}
\usetikzlibrary{shapes,backgrounds}
\usepackage[all,cmtip,knot]{xy}
\xyoption{arc}

\newcommand{\Omit}[1]{}

\newtheorem*{theorem}{Theorem}
\newtheorem*{proposition}{Proposition}
\newtheorem*{corollary}{Corollary}
\newtheorem*{lemma}{Lemma}
\newtheorem*{conjecture}{Conjecture}

\theoremstyle{definition}
\newtheorem*{definition}{Definition}

\newtheorem*{example}{Example}

\numberwithin{equation}{section}

\newenvironment{pf}{\paragraph{{\sc Proof}}}{\qed\par\medskip}

\newcommand{\remark}{{\bf Remark.\ }}

\newcommand {\IC}{\mathbb{C}}
\newcommand {\IN}{\mathbb{N}}

\newcommand {\IZ}{\mathbb{Z}}

\newcommand {\bfI}{\mathbf I}

\newcommand {\A}{\mathcal A}

\newcommand {\C}{\mathcal C}
\newcommand {\D}{\mathcal D}
\newcommand {\E}{\mathcal E}
\newcommand {\F}{\mathcal F}
\newcommand {\G}{\mathcal G}
\renewcommand {\H}{\mathcal H}
\newcommand {\I}{\mathcal I}

\renewcommand {\L}{\mathcal L}

\newcommand {\N}{\mathcal N}
\renewcommand {\O}{\mathcal O}
\newcommand {\compl}{_{\scriptstyle{\O}}}
\renewcommand {\P}{\mathcal P}
\newcommand {\Q}{\mathcal Q}

\newcommand {\T}{\mathcal T}
\newcommand {\U}{\mathcal U}

\newcommand {\X}{\mathcal X}

\renewcommand {\a}{\mathfrak a}
\renewcommand {\b}{\mathfrak b}
\renewcommand {\c}{\mathfrak c}

\newcommand {\g}{\mathfrak{g}}
\newcommand {\h}{\mathfrak h}
\renewcommand {\ll}{\mathfrak l}
\newcommand {\m}{\mathfrak m}
\newcommand {\n}{\mathfrak n}

\renewcommand {\SS}{\mathfrak{S}}

\newcommand{\olh}{\ol{\h}}
\newcommand{\olb}{\ol{\b}}

\newcommand {\sfA}{\mathsf{A}}

\newcommand {\sfC}{\mathsf{C}}

\newcommand {\sfS}{\mathsf{S}}

\newcommand {\ie}{{\it i.e., }}
\newcommand {\eg}{{\it e.g.}, }

\newcommand {\fd}{finite--dimensional }

\newcommand {\wrt}{with respect to }

\newcommand {\ol}{\overline}
\newcommand {\wt}{\widetilde}
\newcommand {\wh}{\widehat}

\newcommand {\mns}{maximal nested set }
\newcommand {\mnss}{maximal nested sets }

\newcommand {\DCP}{De Concini--Procesi }

\newcommand {\DYt}{Drinfeld--Yetter }

\newcommand {\nEK}{Etingof--Kazhdan }

\newcommand {\KM}{Kac--Moody }
\newcommand {\KMA}{Kac--Moody algebra }

\newcommand {\EEK}{^{\scriptscriptstyle{\operatorname{EK}}}}

\def\iip#1#2{\langle{#1},{#2}\rangle}

\def\slantfrac#1#2{\kern.2em\raise.2em\hbox{$#1$}\kern-.2em\left/\lower.4em\hbox{$#2$}\kern-.2em\right.}
\def\backslantfrac#1#2{\kern.2em\lower.4em\hbox{$#1$}\kern-.3em\left\backslash\raise.4em\hbox{$#2$}\right.}

\DeclareMathOperator{\Hom}{Hom}
\DeclareMathOperator{\End}{End}
\DeclareMathOperator{\Aut}{Aut}

\renewcommand{\Im}{\operatorname{Im}}

\DeclareMathOperator{\Res}{Res}
\DeclareMathOperator{\id}{id}
\DeclareMathOperator{\rk}{rk}

\DeclareMathOperator{\Rep}{Rep}

\DeclareMathOperator{\vect}{Vect}

\newcommand {\fml}{[\nts[\hbar]\nts]}

\newcommand {\tens}[1]{^{\otimes #1}}

\newcommand{\ctp}{\hat{\otimes}}
\newcommand{\ten}{\otimes}

\newcommand {\ul}[1]{\underline{#1}}
\newcommand {\Alt}{\operatorname{Alt}}

\newcommand{\ind}{\operatorname{Ind}}

\newcommand{\intm}{\mathsf{int}}

\renewcommand {\sl}[1]{\mathfrak{sl}_{#1}}

\newcommand {\Ug}{U\g}

\newcommand {\Uhg}{U_{\hbar}\g}

\newcommand {\aand}{\qquad\text{and}\qquad}

\newcommand{\DCPA}[2]{\Upsilon_{#1#2}}

\newcommand{\Mns}[1]{\sfMns(#1)}
\newcommand{\sfMns}{\mathsf{Mns}}

\renewcommand{\1}{\mathbf{1}}

\newcommand{\EK}{U_{\hbar}} 

\newcommand{\Eq}[1]{\E_{#1}}

\newcommand{\Oint}{\O^{\operatorname{int}}}

\newcommand{\DrY}[1]{\mathsf{DY}_{#1}}
\newcommand{\ff}{\mathsf{f}}

\newcommand{\DYC}[1]{\mathcal{U}_{#1}}

\newcommand{\gb}{\g_{\b}}
\newcommand{\ga}{\g_{\a}}

\newcommand {\op}{^{\scriptscriptstyle{\operatorname{op}}}}

\newcommand {\<}{\langle}
\renewcommand {\>}{\rangle}

\newcommand{\PROP}{{\sf PROP}}

\newcommand{\Sym}{\operatorname{Sym}}
\newcommand{\SA}[1]{S{#1}}

\newcommand{\LA}{\ul{\sf LA}}
\newcommand{\LCA}{\ul{\sf LCA}}

\newcommand{\preLA}{{\sf LA}}

\newcommand{\HA}{\ul{\sf HA}}
\newcommand{\LBA}{\ul{\sf LBA}}
\newcommand{\PLBA}{\ul{\sf{PLBA}}}
\newcommand{\DY}{\ul{\sf DY}}

\newcommand{\PLCA}{\ul{\sf PLCA}}
\newcommand{\PLA}{\ul{\sf PLA}}
\newcommand{\PDY}{\ul{\sf PDY}}

\newcommand{\Idem}[1]{\sfk[\I^{#1}]}

\newcommand{\sfEnd}[1]{\mathsf{End}\left(#1\right)}

\newcommand{\sfAd}[1]{\mathsf{Ad}(#1)}
\newcommand{\sfAut}[1]{\mathsf{Aut}(#1)}

\newcommand{\Cat}{\mathsf{Cat}}

\newcommand{\sfAlt}{\mathsf{Alt}}

\newcommand{\sfad}{\mathsf{ad}}

\newcommand{\GCM}[1]{\mathsf{#1}}

\newcommand{\fcw}{\lambda^{\vee}}
\newcommand {\hcor}[1]{\alpha^\vee_{#1}}

\newcommand{\euniv}{_{{\operatorname{\mathsf{LBA}}}}}
\newcommand{\univ}{_{{\mathsf{DY}}}}
\newcommand{\Uguniv}{\mathfrak{U}\univ}
\newcommand{\TenUguniv}[1]{\Uguniv^{#1}}

\newcommand{\UU}{\wh{\mathfrak{U}}}

\newcommand{\osplit}{_{{{\mathsf{PDY}}}}}
\newcommand{\Ugsplit}{\mathfrak{U}\osplit}

\newcommand{\TenUgsplit}[1]{\Ugsplit^{#1}}

\renewcommand{\HA}{\ul{\sf HA}}
\renewcommand{\LBA}{\ul{\sf LBA}}
\renewcommand{\PROP}{\sf PROP}
\renewcommand{\PLBA}{\ul{\sf{PLBA}}}
\renewcommand{\DY}[1]{\ul{\sf DY}^{#1}}

\newcommand{\uPLCA}{\mathsf{PLCA}}
\newcommand{\uPLA}{\mathsf{PLA}}
\newcommand{\uPLBA}{\mathsf{PLBA}}
\renewcommand{\PDY}[1]{\ul{\mathsf{PDY}}^{#1}}

\newcommand{\Gsplit}{\Gamma_N}

\newcommand{\gauge}[2]{{#1}_1\cdot{#1}_2\cdot #2\cdot{#1}_{12}^{-1}}
\newcommand{\gaugebis}[2]{#1_1\cdot#1_2\cdot #2\cdot{#1_{12}}^{-1}}
\newcommand{\gaugetris}[2]{(#1)_1\cdot(#1)_2\cdot #2\cdot(#1)_{12}^{-1}}

\newcommand{\VDY}[1]{[\mathsf{V}_{#1}]} 
\newcommand{\ADY}[1]{[#1]} 
\newcommand{\VPDY}[1]{\VDY{#1}} 
\newcommand{\APDY}[1]{[#1]} 
\newcommand{\VCDY}[1]{\VDY{#1}} 
\newcommand{\ACDY}[1]{[#1]} 

\newcommand{\pEnd}[2]{{\mathsf{End}}_{#1}\left(#2\right)}

\newcommand {\fin}{^{\mathsf{fin}}}

\newcommand{\LBATfin}{\LBA(\whT{[1]}^{\ten n}, T[1]^{\ten n})}
\newcommand{\LBASfin}{\LBA(\whS{[1]}^{\ten n}, S[1]^{\ten n})}
\newcommand{\PLBASfin}{\PLBA(\whS{[1]}^{\ten n}, S[1]^{\ten n})}
\newcommand{\PLBATfin}{\PLBA(\wh{T}[1]^{\ten n},T[1]^{\ten n})}

\newcommand{\ulp}{{\ul{p}}}
\newcommand{\ulq}{{\ul{q}}}
\newcommand {\sfQ}{\mathsf Q}
\renewcommand{\int}{^{\scriptscriptstyle{\operatorname{int}}}}

\newcommand {\Euniv}{_{\scriptscriptstyle{\operatorname{univ}}}}

\newcommand {\sfa}{\mathsf a}
\newcommand {\sfsym}{\mathsf{Sym}}
\newcommand {\sfP}{{\mathsf P}}
\newcommand {\sfuP}{\underline{{\mathsf P}}}
\newcommand {\Sch}{{\mathsf{Sch}}}

\newcommand {\rel}{^{\operatorname{\scriptstyle{rel}}}}

\newcommand{\sarch}[2]{r^{\;#2}_{#1,#1}}
\newcommand{\rarch}[3]{r^{\;#2}_{#1,#3}}
\newcommand{\colarch}[4]{r^{\;#2,#3}_{#1,#4}}

\newcommand{\action}{
\xy
(-6,-3)*{};(6,-3)*{}**\dir{-};
(-6,-3.1);(6,-3.1)**\dir{-};
(-6,-3.2);(6,-3.2)**\dir{-};
(-6,-3.3);(6,-3.3)**\dir{-};
(-6,4);(4,-3)**\crv{(0,5)}?(.5)*\dir{};
\endxy
}

\newcommand{\coaction}{
\xy
(-6,-3);(6,-3)**\dir{-};
(-6,-3.1);(6,-3.1)**\dir{-};
(-6,-3.2);(6,-3.2)**\dir{-};
(-6,-3.3);(6,-3.3)**\dir{-};
(-4,-3);(6,4)**\crv{(0,5)}?(.5)*\dir{};
\endxy
}

\newcommand{\bracket}{
\xy
(0,4);(4,0)**\dir{-}?(.35)*\dir{};
(0,-4);(4,0)**\dir{-}?(.35)*\dir{};
(4,0);(8,0)**\dir{-}?(.5)*\dir{};
\endxy
}

\newcommand{\cobracket}{
\xy
(0,0);(4,0)**\dir{-}?(.5)*\dir{};
(8,4);(4,0)**\dir{-}?(.25)*\dir{};
(8,-4);(4,0)**\dir{-}?(.25)*\dir{};
\endxy
}

\newcommand{\diagramrmatrix}{
\xy
(0,-5)*{};
(0,2)*{
	\xy
	(-5,-5)*{};
	(-5,-10)*{};
	(0,-5);(20,-5)**\dir{-};
	(0,-5.1);(20,-5.1)**\dir{-};
	(0,-5.2);(20,-5.2)**\dir{-};
	(0,-5.3);(20,-5.3)**\dir{-};
	(0,-10);(20,-10)**\dir{-};
	(0,-10.1);(20,-10.1)**\dir{-};
	(0,-10.2);(20,-10.2)**\dir{-};
	(0,-10.3);(20,-10.3)**\dir{-};
	(10,0)*{}="A";
	(5,-10);(15,-5)**\crv{(10,5)}?(.65)*\dir{};
	\endxy
};
\endxy
}

\newcommand{\diagrammaticCasimir}{
\xy
(0,2)*{
\xy
(0,-5);(20,-5)**\dir{-};
(0,-5.1);(20,-5.1)**\dir{-};
(0,-5.2);(20,-5.2)**\dir{-};
(0,-5.3);(20,-5.3)**\dir{-};
(5,-5);(15,-5)**\crv{(10,5)}?(.5)*\dir{};
\endxy
};
\endxy
}

\newcommand{\anticasimirtensor}{
	\xy
	(0,-5)*{};
	(0,2)*{
		\xy
		(-5,-5)*{};
		(-5,-10)*{};
		(0,-5);(20,-5)**\dir{-};
		(0,-5.1);(20,-5.1)**\dir{-};
		(0,-5.2);(20,-5.2)**\dir{-};
		(0,-5.3);(20,-5.3)**\dir{-};
		(0,-10);(20,-10)**\dir{-};
		(0,-10.1);(20,-10.1)**\dir{-};
		(0,-10.2);(20,-10.2)**\dir{-};
		(0,-10.3);(20,-10.3)**\dir{-};
		(10,0)*{}="A";
		(5,-10);(15,-5)**\crv{(10,5)}?(.65)*\dir{};
		\endxy
	};
	(17.5,-0.5)*{-};
	(30,2)*{
		\xy(-5,-5)*{};
		(-5,-10)*{};
		(0,-5);(20,-5)**\dir{-};
		(0,-5.1);(20,-5.1)**\dir{-};
		(0,-5.2);(20,-5.2)**\dir{-};
		(0,-5.3);(20,-5.3)**\dir{-};
		(0,-10);(20,-10)**\dir{-};
		(0,-10.1);(20,-10.1)**\dir{-};
		(0,-10.2);(20,-10.2)**\dir{-};
		(0,-10.3);(20,-10.3)**\dir{-};
		(10,0)*{}="A";
		(5,-5);(15,-10)**\crv{(10,5)}?(.45)*\dir{};
		\endxy
	};
	\endxy
}

\newcommand{\polarizedexample}{
	\xy
	(0,2)*{
		\xy
		(-5,-5)*{};
		(-5,-10)*{};
		(0,-5);(20,-5)**\dir{-};
		(0,-5.1);(20,-5.1)**\dir{-};
		(0,-5.2);(20,-5.2)**\dir{-};
		(0,-5.3);(20,-5.3)**\dir{-};
		(0,-10);(20,-10)**\dir{-};
		(0,-10.1);(20,-10.1)**\dir{-};
		(0,-10.2);(20,-10.2)**\dir{-};
		(0,-10.3);(20,-10.3)**\dir{-};
		(10,0)*{}="A";
		(5,-10);(15,-5)**\crv{(10,5)}?(.65)*\dir{};
		\endxy
	};
	(21,2)*{^{m_2^1}};
   (15,-0.5);(27.5,-0.5)**\dir{-}*\dir{>};
	(40,2)*{
	\xy
	(-5,-5)*{};
	(-5,-10)*{};
	(0,-7.5);(20,-7.5)**\dir{-};
	(0,-7.6);(20,-7.6)**\dir{-};
	(0,-7.7);(20,-7.7)**\dir{-};
	(0,-7.8);(20,-7.8)**\dir{-};
	(5,-7.5);(15,-7.5)**\crv{(10,5)}?(.45)*\dir{};
	\endxy
	};
	(64,2)*{^{m_2^1}};
	(70,-0.5);(57.5,-0.5)**\dir{-}*\dir{>};
	(80,2)*{
		\xy(-5,-5)*{};
		(-5,-10)*{};
		(0,-5);(20,-5)**\dir{-};
		(0,-5.1);(20,-5.1)**\dir{-};
		(0,-5.2);(20,-5.2)**\dir{-};
		(0,-5.3);(20,-5.3)**\dir{-};
		(0,-10);(20,-10)**\dir{-};
		(0,-10.1);(20,-10.1)**\dir{-};
		(0,-10.2);(20,-10.2)**\dir{-};
		(0,-10.3);(20,-10.3)**\dir{-};
		(10,0)*{}="A";
		(5,-5);(15,-10)**\crv{(10,5)}?(.45)*\dir{};
		\endxy
	};
	\endxy
}

\newcommand{\arch}[3]
{
\xy
0;/r.#3pc/:
(-22,0)*{}="A";(22,0)*{}="B" **\dir{-};
(-22,-0.1);(22,-0.1)**\dir{-};
(-22,-0.2);(22,-0.2)**\dir{-};
(-22,-0.3);(22,-0.3)**\dir{-};
(-22,-0.4);(22,-0.4)**\dir{-};
(-15,-3)*{{}_{#2}};
(15,-3)*{{}_{#2}};
(0,10)*{#1};(0,10)*\xycircle(4,4){-}="s";
"s";(-12,0)**\crv{(-6,10)}?(.65)*\dir{};
"s";(12,0)**\crv{(6,10)}?(.65)*\dir{};
"s";(-14,0)**\crv{(-7,13)}?(.65)*\dir{};
"s";(14,0)**\crv{(7,13)}?(.65)*\dir{};
"s";(-16,0)**\crv{(-8,16)}?(.65)*\dir{};
"s";(16,0)**\crv{(8,16)}?(.65)*\dir{};
"s";(-18,0)**\crv{(-9,19)}?(.65)*\dir{};
"s";(18,0)**\crv{(9,19)}?(.65)*\dir{};
\endxy
}

\newcommand{\actcoact}{
\xy
(-6,0);(6,0)**\dir{-};
(-6,-0.1);(6,-0.1)**\dir{-};
(-6,-0.2);(6,-0.2)**\dir{-};
(-6,-0.3);(6,-0.3)**\dir{-};
(-6,6)*{};(-2,0)*{}**\crv{(-4,4)} ?(.5)*\dir{};
(2,0)*{};(6,6)*{}**\crv{(4,4)} ?(.6)*\dir{};
(9,3)*{=};
(12,0);(24,0)**\dir{-};
(12,-0.1);(24,-0.1)**\dir{-};
(12,-0.2);(24,-0.2)**\dir{-};
(12,-0.3);(24,-0.3)**\dir{-};
(12,6)*{};(20,0)*{}**\crv{(16,4)} ?(.3)*\dir{};
(16,0)*{};(24,6)*{}**\crv{(20,4)} ?(.8)*\dir{};
(27,3)*{+};
(30,0);(42,0)**\dir{-};
(30,-0.1);(42,-0.1)**\dir{-};
(30,-0.2);(42,-0.2)**\dir{-};
(30,-0.3);(42,-0.3)**\dir{-};
(34,0)*{};(42,6)*{}**\crv{(38,4)} ?(.8)*\dir{} ?(.2)*\dir{};
(30,6)*{};(37,3)*{}**\dir{-}?(.5)*\dir{};
(45,3)*{-};
(48,0);(60,0)**\dir{-};
(48,-0.1);(60,-0.1)**\dir{-};
(48,-0.2);(60,-0.2)**\dir{-};
(48,-0.3);(60,-0.3)**\dir{-};
(48,6)*{};(56,0)*{}**\crv{(52,4)} ?(.3)*\dir{} ?(.9)*\dir{};
(53,3)*{};(60,6)*{}**\dir{-} ?(.5)*\dir{};
\endxy
}
\newcommand{\comodule}{
\xy
(-6,0);(6,0)**\dir{-};
(-6,-0.1);(6,-0.1)**\dir{-};
(-6,-0.2);(6,-0.2)**\dir{-};
(-6,-0.3);(6,-0.3)**\dir{-};
(-3,0);(2,5)**\crv{(0,5)}?(.5)*\dir{};
(2,5);(6,6)**\dir{-};
(2,5);(6,4)**\dir{-};
(9,3)*{=};
(12,0);(24,0)**\dir{-};
(12,-0.1);(24,-0.1)**\dir{-};
(12,-0.2);(24,-0.2)**\dir{-};
(12,-0.3);(24,-0.3)**\dir{-};
(14,0);(20,6)**\crv{(17,6)} ?(.5)*\dir{};
(16,0);(20,4)**\crv{(18,4)} ?(.5)*\dir{};
(20,6);(24,4)**\dir{-};
(20,4);(24,6)**\dir{-};
(27,3)*{-};
(30,0);(42,0)**\dir{-};
(30,-0.1);(42,-0.1)**\dir{-};
(30,-0.2);(42,-0.2)**\dir{-};
(30,-0.3);(42,-0.3)**\dir{-};
(32,0);(42,6)**\crv{(36,6)}?(.5)*\dir{};
(34,0);(42,4)**\crv{(37,5)}?(.5)*\dir{};
\endxy
}
\newcommand{\module}{
\xy
(-6,0);(6,0)**\dir{-};
(-6,-0.1);(6,-0.1)**\dir{-};
(-6,-0.2);(6,-0.2)**\dir{-};
(-6,-0.3);(6,-0.3)**\dir{-};
(-2,5);(4,0)**\crv{(0,5)}?(.5)*\dir{};
(-6,6);(-2,5)**\dir{-};
(-6,4);(-2,5)**\dir{-};
(9,3)*{=};
(12,0);(24,0)**\dir{-};
(12,-0.1);(24,-0.1)**\dir{-};
(12,-0.2);(24,-0.2)**\dir{-};
(12,-0.3);(24,-0.3)**\dir{-};
(12,6);(23,0)**\crv{(19,6)} ?(.5)*\dir{};
(12,4);(20,0)**\crv{(17,4)} ?(.5)*\dir{};
(27,3)*{-};
(30,0);(42,0)**\dir{-};
(30,-0.1);(42,-0.1)**\dir{-};
(30,-0.2);(42,-0.2)**\dir{-};
(30,-0.3);(42,-0.3)**\dir{-};
(34,6);(41,0)**\crv{(37,6)} ?(.5)*\dir{};
(34,4);(38,0)**\crv{(35,4)} ?(.5)*\dir{};
(30,4);(34,6)**\dir{-};
(30,6);(34,4)**\dir{-};
\endxy
}

\newcommand{\diagA}{
\xy
(-10,0)*{}="B-IN"; (50,0)*{}="B-OUT";
{\ar@{-} "B-IN";"B-OUT"};
{\ar@{-}@<-.5pt> "B-IN";"B-OUT"};
{\ar@{-}@<-1pt> "B-IN";"B-OUT"};
{\ar@{-}@<-1.5pt> "B-IN";"B-OUT"};
(5,0)*{}="A1";
(7.5,0)*{}="A2";
(10,0)*{}="A3";
(15,7.5)*{}="B1";
(15,10)*{}="B2";
(15,12.5)*{}="B3";
(20,7.5)*{}="C1";
(20,10)*{}="C2";
(20,12.5)*{}="C3";
(25,7.5)*{}="D1";
(25,10)*{}="D2";
(25,12.5)*{}="D3";
(30,7.5)*{}="E1";
(30,10)*{}="E2";
(30,12.5)*{}="E3";
(35,0)*{}="F1";
(37.5,0)*{}="F2";
(40,0)*{}="F3";
(17.5,10)*=<.5cm,.75cm>{\ul{\alpha}}*\frm{-}="L1";
(27.5,10)*=<.5cm,.75cm>{\sigma}*\frm{-}="L2";
{\ar@/^10pt/@{-} "A1";"B3"};
{\ar@/^7.5pt/@{-} "A2";"B2"};
{\ar@/^5pt/@{-} "A3";"B1"};
{\ar@{-} "C1";"D1"};
{\ar@{-} "C2";"D2"};
{\ar@{-} "C3";"D3"};
{\ar@/_5pt/@{-} "F1";"E1"};
{\ar@/_7.5pt/@{-} "F2";"E2"};
{\ar@/_10pt/@{-} "F3";"E3"};
(-10,15)*{}="A-IN";
(-7.5,15)*\xycircle(2,2){-}="A";
"A"*{\beta};
(0,0)*{}="A-OUT";
(-10,25)*{}="A-AUX";
{\ar@{-} "A-IN";"A"};
"A";"A-OUT"**\crv{(0,12.5)};
(7.5,-4.5)*{\scriptstyle N};
(37.5,-4.5)*{\scriptstyle N};
\endxy
}

\newcommand{\diagB}
{
\xy
(0,0)*{}="B-IN"; (45,0)*{}="B-OUT";
{\ar@{-} "B-IN";"B-OUT"};
{\ar@{-}@<-.5pt> "B-IN";"B-OUT"};
{\ar@{-}@<-1pt> "B-IN";"B-OUT"};
{\ar@{-}@<-1.5pt> "B-IN";"B-OUT"};
(5,0)*{}="A1";
(7.5,0)*{}="A2";
(10,0)*{}="A3";
(15,7.5)*{}="B1";
(15,10)*{}="B2";
(15,12.5)*{}="B3";
(20,7.5)*{}="C1";
(20,10)*{}="C2";
(20,12.5)*{}="C3";
(25,7.5)*{}="D1";
(25,10)*{}="D2";
(25,12.5)*{}="D3";
(30,7.5)*{}="E1";
(30,10)*{}="E2";
(30,12.5)*{}="E3";
(35,0)*{}="F1";
(37.5,0)*{}="F2";
(40,0)*{}="F3";
(17.5,10)*=<.5cm,.75cm>{\ul{\alpha}}*\frm{-}="L1";
(27.5,10)*=<.5cm,.75cm>{\sigma}*\frm{-}="L2";
{\ar@/^10pt/@{-} "A1";"B3"};
{\ar@/^7.5pt/@{-} "A2";"B2"};
{\ar@/^5pt/@{-} "A3";"B1"};
{\ar@{-} "C1";"D1"};
{\ar@{-} "C2";"D2"};
{\ar@{-} "C3";"D3"};
{\ar@/_5pt/@{-} "F1";"E1"};
{\ar@/_7.5pt/@{-} "F2";"E2"};
{\ar@/_10pt/@{-} "F3";"E3"};
(0,15)*{}="A-IN";
(2.5,15)*\xycircle(2,2){-}="A";
"A"*{\beta};
(15,0)*{}="A-OUT";
(-10,25)*{}="A-AUX";
{\ar@{-} "A-IN";"A"};
"A";"A-OUT"**\crv{(15,10)};
(7.5,-4.5)*{\scriptstyle N};
(37.5,-4.5)*{\scriptstyle N};
\endxy
}

\newcommand{\diagC}
{
\xy
(0,0)*{}="B-IN"; (45,0)*{}="B-OUT";
{\ar@{-} "B-IN";"B-OUT"};
{\ar@{-}@<-.5pt> "B-IN";"B-OUT"};
{\ar@{-}@<-1pt> "B-IN";"B-OUT"};
{\ar@{-}@<-1.5pt> "B-IN";"B-OUT"};
(5,0)*{}="A1";
(7.5,0)*{}="A2";
(10,0)*{}="A3";
(15,7.5)*{}="B1";
(15,10)*{}="B2";
(15,12.5)*{}="B3";
(20,7.5)*{}="C1";
(20,10)*{}="C2";
(20,12.5)*{}="C3";
(25,7.5)*{}="D1";
(25,10)*{}="D2";
(25,12.5)*{}="D3";
(30,7.5)*{}="E1";
(30,10)*{}="E2";
(30,12.5)*{}="E3";
(35,0)*{}="F1";
(37.5,0)*{}="F2";
(40,0)*{}="F3";
(17.5,10)*=<.5cm,.75cm>{\ul{\alpha}}*\frm{-}="L1";
(27.5,10)*=<.5cm,.75cm>{\sigma}*\frm{-}="L2";
{\ar@/^10pt/@{-} "A1";"B3"};
{\ar@/^7.5pt/@{-} "A2";"B2"};
{\ar@/^5pt/@{-} "A3";"B1"};
{\ar@{-} "C1";"D1"};
{\ar@{-} "C2";"D2"};
{\ar@{-} "C3";"D3"};
{\ar@/_5pt/@{-} "F1";"E1"};
{\ar@/_7.5pt/@{-} "F2";"E2"};
{\ar@/_10pt/@{-} "F3";"E3"};
(0,15)*{}="A-IN";
(2.5,15)*\xycircle(2,2){-}="A";
"A"*{\beta};
(10,7.5)*{}="A-OUT";
(-10,25)*{}="A-AUX";
{\ar@{-} "A-IN";"A"};
"A";"A-OUT"**\crv{(10,10)};
(7.5,-4.5)*{\scriptstyle N};
(37.5,-4.5)*{\scriptstyle N};
\endxy
}

\newcommand{\diagD}
{
\xy
(0,0)*{}="B-IN"; (45,0)*{}="B-OUT";
{\ar@{-} "B-IN";"B-OUT"};
{\ar@{-}@<-.5pt> "B-IN";"B-OUT"};
{\ar@{-}@<-1pt> "B-IN";"B-OUT"};
{\ar@{-}@<-1.5pt> "B-IN";"B-OUT"};
(5,0)*{}="A1";
(7.5,0)*{}="A2";
(10,0)*{}="A3";
(15,7.5)*{}="B1";
(15,10)*{}="B2";
(15,12.5)*{}="B3";
(20,7.5)*{}="C1";
(20,10)*{}="C2";
(20,12.5)*{}="C3";
(25,7.5)*{}="D1";
(25,10)*{}="D2";
(25,12.5)*{}="D3";
(30,7.5)*{}="E1";
(30,10)*{}="E2";
(30,12.5)*{}="E3";
(35,0)*{}="F1";
(37.5,0)*{}="F2";
(40,0)*{}="F3";
(17.5,10)*=<.5cm,.75cm>{\ul{\alpha}}*\frm{-}="L1";
(27.5,10)*=<.5cm,.75cm>{\sigma}*\frm{-}="L2";
{\ar@/^10pt/@{-} "A1";"B3"};
{\ar@/^7.5pt/@{-} "A2";"B2"};
{\ar@/^5pt/@{-} "A3";"B1"};
{\ar@{-} "C1";"D1"};
{\ar@{-} "C2";"D2"};
{\ar@{-} "C3";"D3"};
{\ar@/_5pt/@{-} "F1";"E1"};
{\ar@/_7.5pt/@{-} "F2";"E2"};
{\ar@/_10pt/@{-} "F3";"E3"};
(0,15)*{}="A-IN";
(2.5,15)*\xycircle(2,2){-}="A";
"A"*{\beta};
(42.5,0)*{}="A-OUT";
(-10,25)*{}="A-AUX";
{\ar@{-} "A-IN";"A"};
"A";"A-OUT"**\crv{(42.5,22.5)};
(7.5,-4.5)*{\scriptstyle N};
(37.5,-4.5)*{\scriptstyle N};
\endxy
}

\newcommand{\scs}{\scriptscriptstyle}

\newcommand{\sfR}{{\mathsf R}}

\newcommand{\UgM}{\mathfrak{U}_\mon}
\newcommand{\UgS}{\mathfrak{U}_\dsg}

\newcommand{\Kar}[1]{\mathsf{Kar}(#1)}
\newcommand{\cKar}[1]{\ul{#1}}
\newcommand{\hext}[1]{{#1}{\fml}}
\newcommand{\hextsub}[1]{{#1}_\hbar}
\newcommand{\hextsup}[1]{{#1}^\hbar}

\newcommand{\mon}{\mathsf{S}}
\newcommand{\dmon}{\mathsf{S}^{(2)}}

\newcommand{\wmLA}[1]{\wt{\preLA}_{#1}}
\newcommand{\mLA}[1]{{\preLA}_{#1}}
\newcommand{\MLA}{\LA_{\mon}}
\newcommand{\preLCA}{\mathsf{LCA}}

\newcommand{\mLCA}[1]{{\LCA}_{#1}}
\newcommand{\MLCA}{\LCA_{\mon}}
\newcommand{\preLBA}{\mathsf{LBA}}

\newcommand{\mLBA}[1]{{\LBA}_{#1}}
\newcommand{\MLBA}{\LBA_{\mon}}

\newcommand{\hatS}[1]{\wh{\SS}_{#1}}

\newcommand{\pb}{[\b]} 
\newcommand{\Erho}[2]{\rho^{#2}_{\g_{#1}}} 
\newcommand{\DYrho}[2]{\rho^{#2}_{#1}} 

\newcommand{\inv}{\scs\mathsf{inv}}

\newcommand{\DYUAinv}[1]{(\DYUA{#1})^{\inv}}
\newcommand{\DYAinv}[2]{(\DYA{#1}{#2})^{\inv}}

\newcommand{\whS}[1]{\wh{S}#1}
\newcommand{\whT}[1]{\wh{T}#1}

\newcommand{\MDY}[2]{\ul{\mathsf{DY}}_{#1}^{#2}}
\newcommand{\MUA}[2]{\mathfrak{U}_{#1}^{#2}}

\newcommand{\SymCat}{\mathsf{SymCat}}
\newcommand{\Sf}{\mathsf{f}}
\newcommand{\diag}{\triangle}
\newcommand{\kS}{\ol{\sfk\SS}}

\newcommand{\pa}{[\a]} 
\newcommand{\pmm}{[\m]} 
\newcommand{\PDYrho}[3]{\rho^{#3}_{#1,#2}} 

\newcommand{\PDYUA}[1]{\mathfrak{U}_{\mathsf{PDY}}^{#1}}
\newcommand{\PDYUAinv}[1]{(\PDYUA{#1})^{\pa}}
\newcommand{\PDYAinv}[3]{\mathcal{U}_{#1,#2}^{#3}}
\newcommand{\CPDYUAinv}[1]{(\wh{\mathfrak{U}}_{\mathsf{PDY}}^{#1})^{\pa}}
\newcommand{\PDYext}[1]{\wedge_{\mathsf{PDY}}^{#1}}
\newcommand{\PDYextinv}[1]{(\PDYext{#1})^{\pa}}

\newcommand{\DYext}[1]{\wedge_{\mathsf{DY}}^{#1}}
\newcommand{\DYextinv}[1]{({\DYext{#1}})^{[1]}}

\newcommand{\Shcp}{\Delta^{\mathsf{sh}}}
\newcommand{\sfU}{\mathsf{U}}
\newcommand{\Sh}{\mathsf{Sh}}

\newcommand{\dmp}[1]{\theta_{#1}}
\newcommand{\cdmp}[1]{\ol{\theta}_{#1}}

\newcommand{\homo}[1]{(#1)_{\delta_N}}

\newcommand{\drc}[1]{\delta_{#1}}

\newcommand{\MDYrho}[2]{\rho^{#2}_{\mon,#1}}

\newcommand{\MLBASfin}{\MLBA({\wh{S}[1]}^{\otimes n},S[1]^{\otimes n})}
\newcommand{\MLBATfin}{\MLBA({\wh{T}[1]}^{\otimes n},T[1]^{\otimes n})}

\newcommand{\RDY}[2]{\ul{\mathsf{DY}}^{#2}_{#1}}
\newcommand{\RDYUA}[2]{\mathfrak{U}_{#1}^{#2}}
\newcommand{\CRDYUA}[2]{\wh{\mathfrak{U}}_{#1}^{#2}}

\newcommand{\mDY}[2]{\DY{#2}_{#1}}
\newcommand{\noc}[1]{\kappa_{#1}} 
\newcommand{\FA}[1]{\A_{#1}}
\newcommand{\FL}[1]{\L_{#1}}
\newcommand{\abFL}[3]{\FL{#1,#2}^{#3}}

\newcommand{\nts}{\negthinspace}
\newcommand{\bb}[1]{\b_{#1}}
\newcommand{\vrtx}[1]{\{i\}}

\newcommand {\sfk}{\mathsf{k}}

\newcommand {\khvect}{\mathsf{Vect}_{\hext{\sfk}}}
\newcommand{\kvect}{\mathsf{Vect_{\sfk}}}

\newcommand{\SDY}[1]{\ul{\mathsf{SDY}}^{#1}} 
\newcommand{\gr}{\mathsf{gr}}
\newcommand{\olg}{\ekm{\g}}

\newcommand{\ekm}[1]{\ol{#1}}

\newcommand {\ulm}{\ul{m}}

\newcommand {\BDm}{B_D^{\ulm}}

\newcommand{\DYA}[2]{\U_{#1}^{#2}}  
\newcommand{\CDYA}[2]{{\wh{\U}_{#1}}^{#2}}  
\newcommand{\DYUA}[1]{\mathfrak{U}_{\mathsf{DY}}^{#1}}  
\newcommand{\whUA}[2]{\wh{U#1^{#2}}} 
\newcommand{\CDYUA}[1]{{\wh{\mathfrak{U}}}_{\mathsf{DY}}^{#1}}  

\newcommand{\hDrY}[2]{\mathsf{DY}_{{#1}}^{#2}}
\newcommand{\hDYA}[2]{\wh{\U}_{#1}^{#2}}

\newcommand{\nqcs}{Coxeter }
\newcommand{\sfBr}{\mathsf{Br}}
\newcommand{\sfT}{\mathsf{T}}
\newcommand{\sft}{\mathsf{t}}

\newcommand{\partialsg}{}
\newcommand{\commutativesg}{}

\newcommand{\sfRo}{\mathsf{R}_0}
\newcommand{\Fun}{\mathsf{Fun}}
\newcommand{\bfb}{\mathbf{b}}

\newcommand{\dsg}{\mathbb{S}}

\newcommand{\dLBA}[1]{{\LBA}_{#1}}

\newcommand{\Ns}[1]{\mathsf{Ns}(#1)}
\newcommand{\Nsr}[2]{\mathsf{Ns}_{#2}(#1)}

\newcommand{\multifibered}{}
\newcommand {\iH}[2]{\ul{#1}_{#2}} 

\newcommand{\hcmp}[2]{h_{#1#2}}

\newcommand{\cc}[1]{\mathsf{conn}(#1)}

\newcommand{\defDY}[1]{\hDrY{#1}{\hbar}}
\newcommand{\intDY}[1]{\hDrY{#1}{\hbar,\mathsf{int},0}}
\newcommand{\sfCox}{\C}
\newcommand{\Cox}{\sfCox}
\newcommand{\Coxb}{\sfCox{(\b)}}

\newcommand{\adm}{\scs\operatorname{adm}}

\newcommand{\mUA}[2]{\mathfrak{U}_{#1}^{#2}}
\newcommand{\Diag}{\operatorname{Diag}}
\newcommand{\sint}{{\mathsf{int}}}

\newcommand {\enrique}[1]{U\mathfrak{G}^{#1}\Euniv} 
\newcommand{\sfp}{\mathsf{p}}
\newcommand{\lba}{_{\preLBA}}

\newcommand{\zolh}{\olh^{c}}
\newcommand{\olgtwo}{\olg^{(2)}}
\newcommand{\olbtwo}{\olb^{(2)}}

\newcommand {\res}{^{\scs\operatorname{res}}}
\newcommand{\grb}{\g_{\b}\res}

\title[Uniqueness of Coxeter structures]
{Uniqueness of Coxeter structures \\ on Kac--Moody algebras}
\author[A. Appel]{Andrea Appel}
\address{School of Mathematics,
University of Edinburgh,
James Clerk Maxwell Building, 
The King's Buildings, 
Peter Guthrie Tait Road,
Edinburgh, EH9 3FD, UK}
\email{andrea.appel@ed.ac.uk}
\author[V. Toledano Laredo]{Valerio Toledano Laredo}
\address{Department of Mathematics,
Northeastern University,
360 Huntington Avenue, Boston MA 02115, USA}
\email{V.ToledanoLaredo@neu.edu}
\thanks{The first author was partially supported by the ERC
Grant 637618 and the NSF grant DMS--1255334.
The second author was partially supported by the NSF
grants DMS--1206305, DMS--1441467 and DMS--1802412.}


\begin{document}

\begin{abstract}
Let $\g$ be a symmetrisable Kac--Moody algebra, and $\Uhg$
the corresponding quantum group. We showed in \cite{ATL1,ATL1-2}
that the braided \nqcs structure on integrable, category
$\O$ representations of $\Uhg$ which underlies the $R$--matrix
actions arising from the Levi subalgebras of $\Uhg$ and the
quantum Weyl group action of the generalised braid group
$B_{\g}$ can be transferred to integrable, category $\O$
representations of $\g$. We prove in this paper that, up to
unique equivalence, there is a unique such structure on the
latter category with prescribed restriction functors, $R$--matrices,
and local monodromies. This extends, simplifies and strengthens a
similar result of the second author valid when $\g$ is semisimple,
and is used in \cite{ATL3} to describe the monodromy
of the rational Casimir connection of 
$\g$ in terms of the quantum Weyl group
operators of $\Uhg$. Our main tool is a refinement
of Enriquez's universal algebras, which is adapted to the
$\PROP$ describing a Lie bialgebra graded by the
non--negative roots of $\g$.
\end{abstract}

\maketitle
\setcounter{tocdepth}{1} 
\tableofcontents

\section{Introduction}\label{s:intro}

\subsection{} 

This is the second of three papers whose goal is to extend the description
of the monodromy of the rational Casimir connection of a semisimple Lie
algebra in terms of quantum Weyl groups given in \cite{vtl-2,vtl-3,vtl-4,vtl-6}
to the case of an arbitrary symmetrisable Kac--Moody algebra $\g$.

In \cite{ATL1-2}, we introduced the notion of braided Coxeter category,
which is informally a tensor category carrying commuting actions of Artin's
braid groups and a given generalised braid group on the tensor product
of its objects. We showed that such a structure arises from the quantum
group $\Uhg$, specifically on the
category $\O\int_\hbar$ of integrable, highest weight representations of
$\Uhg$. The corresponding Artin group actions are given by the universal
$R$--matrices of the Levi subalgebras of $\Uhg$, and the action of the
generalised braid group of $\g$ by the quantum Weyl group operators
of $\Uhg$. The main result of \cite{ATL1-2} is that this structure can be
transferred to the category $\O\int$ of integrable, highest weight modules
for $\g$. The transfer relies on a 2--categorical, relative version of Etingof--Kazhdan
quantisation, which takes as input a split inclusion of Lie bialgebras
$\a\subset\b$, and allows 
to construct an equivalence
$\O\int_\hbar\cong\O\int$ which is compatible with a given chain of
Levi subalgebras of $\g$ \cite{ATL1}.

\subsection{} 

The goal of the present paper is to prove that $\O\int$ possesses,
up to unique equivalence, a unique braided \nqcs structure with
prescribed restrictions functors, $R$--matrices and local monodromies.
This is used in \cite{ATL3} to prove that the monodromy of the rational
Casimir connection of $\g$ is described by the quantum Weyl group
operators of $\Uhg$, by showing that the monodromy of the rational
KZ and Casimir connections arise from a braided quasi--Coxeter
structure on $\O\int$.\footnote{More precisely, for an arbitrary
symmetrisable \KM algebra only the normally ordered version of the Casimir
connection introduced in \cite{FMTV} can be defined. We show
in \cite{ATL3}, however, that its monodromy can be modified by
a cocyle so as to become equivariant under the Weyl group, and
that the resulting action of the braid group action is described by
the quantum Weyl group operators of $\Uhg$.}

\subsection{} 

The uniqueness of braided \nqcs structures on $\O\int$
is obtained from a cohomological rigidity result, as is the case for
a semisimple Lie algebra. 
The proof of this result,
however, differs significantly from that given in \cite{vtl-3,vtl-4}.
Indeed, the latter relies on the well--known computation of the
Hochschild (coalgebra) cohomology of the enveloping algebra
$U\g$ in terms of the exterior algebra of $\g$. For an arbitrary
\KM algebra, the tensor powers of $U\g$ need to be replaced
by their completion $U\g^{\otimes n}\compl$ \wrt category $\O$.
Indeed, $U\g$ and $U\g^{\otimes 2}$ do not contain the Casimir operator
$C$ of $\g$ and the invariant tensor $2\Omega=\Delta(C)-C\otimes
1-1\otimes C$ respectively, and are therefore not appropriate
receptacles for the coefficients of the Casimir and KZ connections.
While the computation of the Hochschild cohomology of $U\g$
holds for an arbitrary Lie algebra, it is not known to do so, and
may in fact fail, for the topological coalgebra $U\g\compl$, which
seems to have a rather unwieldy cohomology.

\subsection{} 

Rather than using the completions $U\g\tens{n}\compl$, we rely
on a refinement of Enriquez's universal algebras $\enrique{n}$
\cite{e1}. These arise from the $\PROP$ of Lie bialgebras,
and were used by Enriquez to give a cohomological construction of
quantisation functors for Lie bialgebras \cite{e3}. They are universal
in the following sense: for any Lie bialgebra $\b$ with Drinfeld double
$\gb=\b\oplus\b^*$ and any $n\geq 1$, $\enrique{n}$ maps to a
completion $\wh{U\gb\tens{n}}$ of the $n$--fold tensor product
of the enveloping algebra of $\gb$. For $n=1$, the image of $\enrique
{n}$ in $\wh{U\gb}$ is the subalgebra spanned by the interlaced powers
of the normally ordered Casimir operator of $\gb$, \ie the elements
\begin{equation}\label{eq:interlaced-casimir}
\kappa_N^\sigma=
\sum_{i_1,\ldots,i_N}
b_{i_1}b_{i_2}\cdots b_{i_N}\cdot
b^{i_{\sigma(N)}}b^{i_{\sigma(N-1)}}\cdots b^{i_{\sigma(1)}}
\end{equation}
where $\{b_i\},\{b^i\}$ are dual bases of $\b$ and $\b^*$, $N$ is an
arbitrary integer, and $\sigma$ a permutation in ${\mathfrak S}_N$.

The completion $\wh{U\gb\tens{n}}$ is \wrt the category $\E_{\gb}$
of {\it equicontinuous} $\gb$--modules, which are those on which the
action of $\b^*$ (and therefore the sum \eqref{eq:interlaced-casimir})
is locally finite. If $\g$ is a symmetrisable \KM algebra with negative
Borel subalgebra $\b$, the realisation of $\g$ as a quotient of the
Drinfeld double of $\b$ gives rise to an embedding of category
$\O$ for $\g$ as a full subcategory of $\E_{\gb}$.

The coproduct on $U\gb$ gives rise to a cosimplicial structure
on the tower of algebras $\{\wh{U\gb\tens{n}}\}_{n\geqslant  0}$.
The latter can be lifted to $\{\enrique{n}\}_{n\geqslant  0}$, and gives rise to a
Hochschild complex. Enriquez's crucial insight is that this complex
contains enough elements to allow for the construction of quantisation
functors, yet has a manageable cohomology,
which is given by a universal version
of the exterior algebra of $\gb$. 

\subsection{}\label{ss:DYs and some elements}

We give in this paper an alternative, and perhaps more natural
construction of $\enrique{n}$ by using \DYt modules over a Lie
bialgebra $\b$. Such a module is a triple $(V,\pi,\pi^*)$ where
$\pi:\b\otimes V\to V$ gives $V$ the structure of a left $\b$--module,
$\pi^*:V\to\b\otimes V$ that of a right $\b$--comodule, and $\pi,\pi^*$
satisfy a compatibility condition \cite{ek-2}. The latter is designed
so as to give rise to an action of the Drinfeld double of $\b$, with
$\phi\in\b^*\subset\gb$ acting on $V$ by $\phi\otimes\id_V\circ\pi^*$.

The symmetric tensor category $\hDrY{\b}{}$ of such modules
coincides with that of equicontinuous $\gb$--modules, with the
coaction of $\b$ on $V\in\E_{\gb}$ given by $\pi^*(v)=\sum_i
b_i\otimes b^i\,v$ \cite{ek-1}. Under this correspondence, the
action of the normally ordered Casimir $\kappa=\sum_i b_i
b^i$ of $\gb$ on $V\in\hDrY{\b}{}$ is simply given by $\pi\circ\pi^*$. More generally,
the interlaced Casimir $\kappa_N^\sigma$ \eqref{eq:interlaced-casimir}
acts on $V$ by the composition of the iterated coaction $(\pi^*)
^{(N)}:V\to\b^{\otimes N}\otimes V$ followed by the permutation
$\sigma^{-1}\otimes\id_V$ 
and the iterated action $\pi^{(N)}:\b^{\otimes N}\otimes V\to V$.
Similarly, the $r
$--matrix of $\gb$ given by $r=\sum_i b_i\otimes b^i\in\b\wh
{\otimes}\b^*$ acts on a tensor product $V\otimes W$ as the
composition
\[
r_{VW}=
\pi_V\otimes\id_W\circ(1\,2)\circ\id_V\otimes\pi^*_W
\]


For any $n\geq 1$, we introduce a colored $\PROP$ $\DY{n}$
which describes a Lie bialgebra $\b$, together with $n$ \DYt 
modules $V_1,\ldots,V_n$ over $\b$. We then consider the algebra $\TenUguniv{n}
=\End_{\DY{n}}(V_1\otimes\cdots\otimes V_n)$, and show it to be
isomorphic to Enriquez's algebra $\enrique{n}$. This alternative
construction makes the algebra structure on $\enrique{n}$, and
its action on equicontinuous $\gb$--modules far more transparent. 

\subsection{} 

We then introduce three refinements of the algebras $\TenUguniv{n}$.
The first one, $\PDYUA{n}$, is obtained from the colored $\PROP$ describing
a split inclusion of Lie bialgebras $\a\subset\b$, together with $n$ \DYt
modules over $\b$. The image of $\TenUgsplit{1}$ in $\wh{U\g_\b}$
is spanned by the interlaced products of the normally ordered Casimir
operators of the doubles of $\a$ and $\b$. The Hochschild cohomology 
of the tower $\TenUgsplit{n}$ can be computed via the calculus of Schur
functors developed in \cite{e3}, and shows in particular that the relative
quantisation functor constructed in \cite{ATL1} is unique up to unique
isomorphism.\footnote{The uniqueness of the isomorphism follows from
the fact that the first Hochschild cohomology of $\TenUgsplit{\bullet}$ is
zero, as is the case for $\enrique{\bullet}$. Thus, figuratively speaking,
$\enrique{\bullet}$ and $\TenUgsplit{\bullet}$ behave like the tensor
powers of an enveloping algebra without primitive elements.}

The second refinement, $\UgM^n$, is obtained in a similar way from
a $\PROP$ describing $n$ \DYt modules over a Lie bialgebra $\b$
which is graded by a partial abelian semigroup $\mon$.
The image  of $\UgM^n$ in $\wh{U\g_\b}$ is then spanned by the
interlaced products of the normally ordered Casimir operators of
the subspaces $\b_\alpha\oplus\b_\alpha^*\subset\gb$, $\alpha\in\mon$.
When $\mon$ is the partial semigroup $\sfR_+\sqcup\{0\}$ consisting
of the positive roots of a symmetrisable \KM algebra $\g$ together
with zero, this makes $\UgM$ an appropriate
receptacle for the coefficients of the Casimir connection of $\g$. 

The third refinement, $\UgS^n$, is prompted by the following. We
show in \cite{ATL1-2} that the braided (pre--)Coxeter structure transferred
from $\Oint_{\hbar}$ to $\Oint$  is \emph{diagrammatic}, \ie compatible
with the Lie subalgebras $\g_B$ generated by the root vectors $\{e_i,f_i\}_
{i\in B}$ corresponding to a subdiagram $B$ of the Dynkin diagram
of $\g$. In particular, this structure cannot be lifted to a braided (pre--)Coxeter
structure on $\UgM^{\bullet}$, since the latter only accounts for the
Cartan subalgebra of $\g$ and not its subspaces $\h_B$ spanned 
by $\{\alpha^\vee_i\}_{i\in B}$.

The definition of $\UgS^n$ relies on a \emph{diagrammatic semigroup} $\dsg$,
and accounts for both the root space decomposition of $\g$ as well as for its
diagrammatic subalgebras $\g_B$. The braided (pre--)Coxeter structures coming from the quantum group and 
the Casimir connection can then both be realized in $\UgS^{\bullet}$, as we show in
\cite{ATL1-2} and \cite{ATL3}, respectively. The computation of the Hochschild 
cohomology of $\UgS^n$ yields the required rigidity result, thus allowing 
to prove they are isomorphic.

\subsection{} 

The use of the algebras $\UgS^n$ leads to far stronger uniqueness
results than had been obtained in \cite{vtl-3,vtl-4} for a semisimple
Lie algebra $\g$. Indeed, as is the case for the universal algebras
$\enrique{n}$, the tower $\{\UgS^n\}_{n\geq 0}$ has trivial first Hochschild
cohomology, which implies that the isomorphism of two braided,
\nqcs structures 
is unique up to a
{\it unique} gauge. This raises the hope that the equivalences
we construct may be convergent as series in the deformation
parameter $\hbar$, and could in particular be specialised to
numerical, non--rational, values of $\hbar$.
It is also worth pointing out that the vanishing of the first
Hochschild cohomology removes the need for the use of
Dynkin diagram cohomology developed in \cite{vtl-4} to deal
with secondary obstructions, thereby simplifying the proof
of rigidity even for a semisimple Lie algebra.

\subsection{} 

We now review our results in more detail. A $\PROP$ is a
categorical realisation of an algebraic structure. More precisely, 
given a field $\sfk$ of characteristic zero, a ${\sf PRO}$duct--$
{\sf P}$ermutation category is a $\sfk$--linear, symmetric monoidal
category with objects the non--negative numbers and tensor product
$[n]\ten[m]=[n+m]$. For example, the $\PROP$ $\LA$ of Lie algebras
is generated by an anti--symmetric morphism $\mu:[2]\to[1]$ satisfying
the Jacobi identity. One can then think of Lie algebras over $\sfk$
as symmetric monoidal functors from $\LA$ to $\sfk$--vector spaces,
and morphisms of Lie algebras as natural transformations of the
corresponding realisation functors. 

\subsection{} 

Richer structures can be described by \emph{colored} $\PROP$s, whose
objects are sequences in a given set of colors $\sfA$. A key example for
us is the $\PROP$ $\DY{n}$ on $n+1$ colors which we introduce in Section
\ref{s:univDY}. In this case, the category of symmetric monoidal functors
$\DY{n}\to\vect_{\sfk}$ is isomorphic to that of tuples $(\b;V_1,\ldots,V_n)$ 
consisting of a Lie bialgebra $\b$ over $\sfk$, and $n$ \DYt modules $V_1,
\ldots,V_n$.

A natural transformation of these functors amounts to a tuple $(\phi;f_1,
\ldots,f_n)$, where $\phi:\b\to\c$ is a morphism of Lie bialgebras, and each
$f_i:V_i\to W_i$ is both a morphism of $\b$--modules $V_i\to\phi^*W_i$ and
of $\c$--comodules $\phi_*V_i\to W_i$. In particular, choosing $\phi=\id$
shows that any endomorphism of $V_1\ten \cdots\ten V_n$ in $\DY{n}$
commutes with morphisms of \DYt modules over any Lie bialgebra $\b$.
Thus, if $\ff:\DrY{\b}\to\vect_\sfk$ is the forgetful functor, the algebra
\[\DYUA{n}=\mathsf{End}_{\DY{n}}(V_1\ten\cdots\ten V_n)\]
maps to the endomorphisms of $\ff^{\otimes n}$, and therefore to
the completion of $U\gb^{\ten n}$ \wrt equicontinuous modules.

\subsection{}

The category $\DY{n}$ is best described diagrammatically. The identity on the
universal Lie bialgebra (resp. \DYt module) is represented by a thin (resp. bold)
horizontal line, and the bracket, cobracket, action and coaction by the diagrams
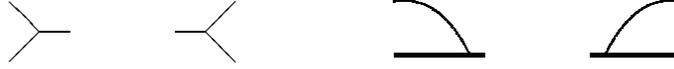
\[\bracket\qquad\qquad\cobracket\qquad\qquad\qquad\action\qquad\qquad\coaction\]
which are read from left to right. By \ref{ss:DYs and some elements}, the normally
ordered Casimir $\kappa\in\DYUA{1}$ and $r$--matrix $r\in\DYUA{2}$ are therefore
represented, respectively, by
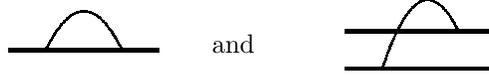
\[\diagrammaticCasimir
\aand
\diagramrmatrix\]

In Sections \ref{s:morlba} and \ref{s:univDY}, we explicitly describe the morphisms
in $\DY{n}$, and construct an integral basis for $\DYUA{n}$ which, for $n=1$ is
given by the diagrams
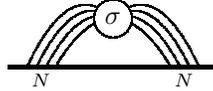
\[\arch{\sigma}{N}{15}\]
where $N\geqslant 0$ and $\sigma\in\SS_N$. By \ref{ss:DYs and some elements},
these correspond to the interlaced powers of the normally ordered Casimir \eqref
{eq:interlaced-casimir}. This description leads to a PBW theorem for $\DYUA{n}$,
and the explicit computation of its Hochschild cohomology, which is analogous to
the fact that $H^n(U\gb)=\wedge^n\gb$.

Pictorially, an element in $H^{n}(\DYUA
{\bullet})$ is a linear combination of anti--symmetric diagrams with $n$ bold lines,
and exactly one action or one coaction on each of these. For example, $H^{1}
(\DYUA{\bullet})=0$ (\ie $\DYUA{1}$ has no primitive elements), and the simplest
non--trivial element in $H^{2}(\DYUA{\bullet})$ is the anti--symmetric $r$--matrix
\[\frac{1}{2}\left(\anticasimirtensor\quad\right)\] 

\subsection{}\label{ss:receptacle}

The algebra $\DYUA{n}$ is a universal receptacle for the coefficient of the
KZ connection on $n$ points for a Drinfeld double $\gb$ since, for $n=2$,
it contains the invariant tensor $\Omega=\sum_i b_i\otimes b^i+b^i\otimes
b_i=r+r_{21}$. However, $\DYUA{1}$ is too small to contain the coefficients
of the Casimir connection of a symmetrisable \KMA $\g$, since it does not
account for the root space decomposition of $\g$, and in particular the
diagrammatic and Levi subalgebras
\[\g_B=\<e_i,f_i\>_{i\in B}\subset\ll_B=\g_B+\h\subset\g\]
corresponding to a subdiagram $B$ of the Dynkin diagram of $\g$.

To this end, we first introduce and study in Section \ref{s:plba}
the $\PROP$ $\PDY{n}$ obtained by adding to $\DY{n}$ an idempotent
endomorphism of the universal Lie bialgebra. Its image is then a split Lie
subbialgebra. The results of Section \ref{s:univDY}, in particular the PBW Theorem
and computation of Hochschild cohomology extend easily to the algebras
$\PDYUA{n}=\mathsf{End}_{\PDY{n}}(V_1\ten\cdots\ten V_n)$.

From here, the $\PROP$ic construction of Levi subalgebras is fairly straightforward.
We first observe that the negative Borel of a symmetrisable \KM algebra $\g$ is
graded, as a Lie bialgebra, by the partial semigroup $\sfR_0$ consisting of the
positive roots of $\g$ and zero. The $\PROP$ encoding this structure is denoted
$\mDY{\mon}{n}$, where $\mon$ is any partial abelian semigroup. It is obtained
from $\DY{n}$ by adding a complete family of orthogonal idempotents labelled by
the elements of $\mon$. In the case of the semigroup $\sfR_0$, by considering
the sum of the idempotents corresponding to zero or a root associated to a 
subdiagram $B$ of the Dynkin diagram,  one obtains a universal analogue of the
Levi subalgebra $\ll_B$.

The universal algebra $\RDYUA{\mon}{n}=\mathsf{End}_{\mDY{\mon}{n}}(V_1\ten
\cdots\ten V_n)$ is generated by arc diagrams, in which each thin line is now
labeled by $\mon$. In particular, $\RDYUA{\mon}{1}$ contains the elements
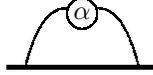
\[
\xy
(0,0)*{
\xy
(0,-7)*{};
(0,-5);(20,-5)**\dir{-};
(0,-5.1);(20,-5.1)**\dir{-};
(0,-5.2);(20,-5.2)**\dir{-};
(0,-5.3);(20,-5.3)**\dir{-};
(10,2)*{\alpha};(10,2)*\xycircle(2,2){-}="A";
(2.5,-5);"A"**\crv{(5,5)}?(.25)*\dir{};
"A";(17.5,-5)**\crv{(15,5)}?(.75)*\dir{};
\endxy
};
\endxy
\]
\noindent
for any $\alpha\in\mon$. In the case of the semigroup of non--negative roots
$\sfR_0$, these diagrams correspond precisely to the normally ordered Casimir
elements of the $\sl{2}$--triple of the root $\alpha$, and make $\RDYUA{\mon}
{1}$ a universal receptacle for the coefficients of the Casimir connection. 

The universal algebras $\RDYUA{\mon}{n}$, however, do not provide a universal
realization of the diagrammatic subalgebras $\g_B$,
which are necessary to describe the braided (pre--)Coxeter structure transferred
from $\Oint_{\hbar}$. We therefore introduce a refinement $\mDY{\dsg}{n}$ of
the $\PROP$ $\mDY{\mon}{n}$, where we further decompose the idempotent
corresponding to the zero element of $\mon$, so as to reproduce the subspaces
$\h_B=\<\alpha_i^\vee\>_{i\in B}\subset\h$.
The corresponding universal algebras $\RDYUA{\dsg}{n}$ account for both 
the root space decomposition of $\g$, and therefore the coefficients of the Casimir
connection, as well as for its diagrammatic subalgebras. 

\subsection{}

We now sketch the definition of a braided Coxeter category. 
We refer to an unoriented graph $D$ with no multiple edges or loops as a
diagram, and to its full subgraphs $B\subseteq D$ as subdiagrams. A braided
{\it pre--}Coxeter category $\Q$ of type $D$ consists of the following three
pieces of data
\begin{itemize}
\item {\bf Diagrammatic categories.} For any subdiagram $B\subseteq D$,
a braided tensor category $\Q_B$.
\item {\bf Restriction functors.} For any pair of subdiagrams $B'\subseteq B$,
a (not necessarily braided) monoidal functor $F^\F_{B'B}:\Q_B\to\Q_{B'}$
depending upon the choice of a maximal chain of subdiagrams
$B=B_1 \supsetneq B_2 \supsetneq\cdots \supsetneq B_m=B'$.
\item {\bf Associators.} For any $B'\subseteq B$, and pair of maximal chains
$\G,\F$ from $B$ to $B'$, an isomorphism of monoidal functors $\Upsilon
^{\G\F}:F_{B'B}^\F\rightarrow F_{B'B}^\G$.\footnote{The data labeling the restriction
functors $F^\F_{B'B}$ and isomorphisms $\Upsilon^{\G\F}$ actually consists
of a {\it maximal nested set} on $B$ relative to $B'$ (see Section \ref{s:diagrams}
for the definition). The collection of such nested sets is a quotient of the set of
set of maximal chains, and for simplicity we identify the two in the introduction.}
\end{itemize}
The above data satisfies various requirements. In particular, the restriction functors
and associators are compatible with the composition of chains corresponding to
triple inclusions $B''\subset B'\subset B$ and, for any $B'\subset B$ and maximal chains $\F,\G,\H$
from $B$ to $B'$, one has $\Upsilon^{\H\G}\circ\Upsilon^{\G\F}=\Upsilon^{\H\F}$.
 
$\Q$ is a braided Coxeter category if it is further endowed with distinguished 
elements $S^{\Q}_i\in\sfAut{F_i}$ where $i$ ranges over the vertices of $D$,
which satisfy the following version of the braid relations determined labeling
the edges of $D$ by multiplicities $m_{ij}=m_{ji}\in\{2,3,\ldots,\infty\}$. For
any $i\neq j$ such that $m_{ij}<\infty$, and maximal chains $\F,\G$ from $D$
to the empty subdiagram such that $\F$ (resp. $\G$) contains $i$ (resp. $j$)
among the connected components of its elements, the following holds in
$\Aut(F^\F_{\emptyset D})$
\[\underbrace{\mathsf{Ad}\left(\Upsilon^{\G\F}\right)(S_j^{\Q})\cdot S_i^{\Q}
	\cdot \mathsf{Ad}\left(\Upsilon^{\G\F}\right)(S_j^{\Q})\cdots}_{m_{ij}}=
\underbrace{S_i^{\Q}\cdot\mathsf{Ad}\left(\Upsilon^{\G\F}\right)(S_j^{\Q})\cdot S_i^{\Q}\cdots}_{m_{ij}}\]
This gives rise to an action $\rho^\F:B_D\to\Aut(F^\F_{\emptyset D})$ of the Artin
braid group $B_D$ determined by the labeling of $D$, 
which is intertwined by the associators $\Upsilon^{\G\F}$.

\subsection{}

Let $\g$ be the symmetrisable Kac--Moody algebra with Dynkin 
diagram $D$. For any subdiagram $B\subseteq D$, we denote by 
$\g_B\subseteq\g$ the diagrammatic subalgebra and by $\b_B\subseteq\b$ 
its negative Borel subalgebra. To study braided pre--Coxeter structure 
on \DYt modules over $\{\b_B\}_{B\subseteq D}$, $\b$ has to be
\emph{diagrammatic}, \ie for any $B'\subseteq B$, $\b_{B'}\subseteq\b_B$
and, for any $B'\perp B$, $[\b_{B'}, \b_{B}]=0$. \footnote{Recall that $B'\perp B$
if $B\cap B'=\emptyset$ and no vertex of $B$ is connected with one of $B'$.}
Although Kac--Moody algebras of finite, affine, or hyperbolic type are diagrammatic,
not all are diagrammatic with counterexamples already in rank $4$ (cf.\cite{ATL1-2}). 
To remedy this, in Section \ref{s:km-rigidity}, we consider certain split central extensions,
referred to as \emph{extended Kac--Moody algebras}, whose Borel subalgebras are canonically 
endowed with a split diagrammatic structure.

Let $\olg$ be a diagrammatic or extended symmetrisable Kac--Moody algebras with
diagrammatic semigroup of positive roots $\dsg$ and universal algebras $\RDYUA{\dsg}{n}$, 
$n\geqslant 1$. For any subdiagram $B\subseteq D$, there is a diagrammatic subalgebra 
$\olg_B\subseteq\olg$ with negative Borel subalgebra $\olb_B\subseteq\olb$. The corresponding
root subsystem defines a universal subalgebra $\RDYUA{\dsg, B}{n}\subseteq\RDYUA{\dsg}{n}$. 

The definition of a braided pre--Coxeter structure on \DYt modules
over $\{\olb_B\}_{B\subseteq D}$ can be lifted to an algebraic datum
on $\RDYUA{\dsg}{\bullet}$, which we call a \emph{universal} braided
pre--Coxeter structure. This consists of a collection of associators
$\Phi_B\in\RDYUA{\dsg, B}{3}$, $B\subseteq D$, twists $J^{\F}_{B'B}
\in\RDYUA{\dsg, B}{2}$, $\F\in\Mns{B,B'}$,  and gauge transformations
$\Upsilon^{\G\F}\in\RDYUA{\dsg, B}{1}$, $\F,\G\in\Mns{B,B'}$. 

\subsection{}

By relying on the computation of the Hochschild cohomology of 
$\RDYUA{\dsg}{\bullet}$, in particular the vanishing of $H^1(\RDYUA
{\dsg}{\bullet})$, we prove 
the uniqueness of universal braided pre--Coxeter structures on
$\RDYUA{\dsg}{\bullet}$ with prescribed braiding.
We also show 
that a universal braided pre--Coxeter structure with diagrammatic categories 
$\{\intDY{{\olb_B}}\}_{B\subseteq D}$  extends in at most one way to a
braided Coxeter one. This gives our main result.

\begin{theorem}
Let $\sfCox_1,\sfCox_2$ be two universal braided \nqcs structures with 
diagrammatic categories $\{\intDY{{\olb_B}}\}_{B\subseteq D}$.
Then,
\begin{enumerate}
\item $\sfCox_1$ and $\sfCox_2$ are twist equivalent.
\item The twist relating them is unique up to a unique gauge transformation.
\end{enumerate}
\end{theorem}

The categories $\intDY{{\olb_B}}$ naturally contains a generalisation $\O_{\infty,\olg_B}
^{\hbar,\sint}$ of category $\O$, where the weight spaces are allowed to be
infinite--dimensional. The theorem above readily restricts to the diagrammatic
categories $\{\O_{\infty,\olg_B}^{\hbar,\sint}\}_{B\subseteq D}$ and yields the following.

\begin{corollary}
There is, up to a unique universal equivalence, a unique universal braided
Coxeter structure with diagrammatic categories $\{\O_{\infty,\olg_B}^{\hbar,\sint}\}_{B\subseteq D}$
\end{corollary}

\subsection{Outline of the paper}

In Section \ref{s:props}, we review the \nEK quantisation of Lie
bialgebras and its description in terms of the $\PROP$ $\LBA$
of Lie bialgebras. In Section \ref{s:sch}, we review the theory of
Schur functors and their cohomology following \cite{baez} and
\cite{e3} respectively. In Section \ref{s:morlba},
we describe the factorised structure of morphisms in $\LBA$, and
their relation to free Lie algebras.
In Section \ref{s:univDY}, we introduce the $\PROP$ $\DY{n}$,
the algebra $\TenUguniv{n}$, and we study its properties and its
Hochschild cohomology. 
In Section \ref{s:plba}, we introduce the refined $\PROP$
$\PLBA$, describing a split inclusion $\a\subset\b$ of Lie bialgebras,
and the corresponding universal algebra $\PDYUA{n}$, for which we prove a number
of results analogous to those obtained for $\TenUguniv{n}$.
In Section \ref{s:twists} these are used to prove the uniqueness, up
to a unique gauge transformation, of the relative quantisation functor
constructed in \cite{ATL1}.
Section \ref{se:sgp} contains some background material on partial
semigroups and Lie bialgebras graded over them.
In Section \ref{s:sgp-ext}, we study the further refined $\PROP$ $\LBA_
\mon$, for a partial abelian semigroup $\mon$, and its universal
algebra $\UgM^n$. In particular, we compute its Hochschild cohomology. 
In Section \ref{se:sub U_S}, we study the subalgebras of $\UgM^n$
defined by the saturated subsemigroups of $\mon$. 
Section \ref{s:diagrams} reviews the combinatorial definitions of
diagrams and maximal nested sets.
In Section \ref{s:gr-prop}, we define diagrammatic (partial) semigroups,
and define for these an extension of $\LBA_\mon$ which allows in
particular to simultaneously account for both the diagrammatic structure
of the Borel subalgebra of a complex semisimple Lie algebra as well
as its root space decomposition.
In Section \ref{s:univ-w-qC} we define a universal braided pre--Coxeter
structure associated to a diagrammatic semigroup, and prove
its rigidity.
In Section \ref{se:braided Cox}, we review the definition of braided
(pre-)Coxeter categories following \cite{ATL1-2}. We then show that
the braided pre--Coxeter structures introduced in Section \ref{s:univ-w-qC} 
give rise to braided (pre--)Coxeter category structures on \DYt
modules over Lie bialgebras graded by a diagrammatic semigroup.
In the final Section \ref{s:km-rigidity}, we use these results to prove the
uniqueness of braided pre--Coxeter structures on the category of
integrable, \DYt modules over the Borel subalgebra of an arbitrary
symmetrisable diagrammatic or extended \KM algebra $\olg$ and on category 
$\O$--modules over $\olg$.

\subsection{Acknowledgments}

We are grateful to Pavel Etingof for a number of useful discussions. A substantial
portion of this paper was written while the authors attended the 2015 Park City
Mathematics Institute research program on the Geometry of Moduli Spaces and
Representation Theory. We are indebted to Roman Bezrukavnikov, Alexander
Braverman and Zhiwei Yun for organising the program and for their invitation,
and to Rafe Mazzeo and the staff at PCMI for the marvelous atmosphere,
flawless organisation, and financial support through the NSF grant DMS--1441467.


\section{Universal quantisation of Lie bialgebras}\label{s:props}

In this section, we review the Etingof--Kazhdan quantisation of Lie
bialgebras, and its description in terms of \emph{product--permutation
categories} ($\PROP$s). For more details, we refer the reader to
\cite{ek-1,ek-2}.

\subsection{Drinfeld double}
A Lie bialgebra over a field $\sfk$ is a triple $(\b,[\cdot,\cdot]_{\b},\delta_{\b})$ 
where $(\b,[\cdot,\cdot]_{\b})$ is a Lie algebra (\ie 
$[\cdot,\cdot]_{\b}:\b\ten\b\to\b$ is antisymmetric and satisfies 
the Jacobi identity), $(\b,\delta_{\b})$ is a Lie coalgebra (\ie 
$\delta_{\b}:\b\to\b\ten\b$ is antisymmetric and satisfies the co-Jacobi identity),
and $[\cdot,\cdot]_{\b}$, $\delta_{\b}$ satisfy the cocycle condition
\begin{equation}
\delta_{\b}([x,y]_{\b})=[x\ten1+1\ten x,\delta_{\b}(y)]-[y\ten1+1\ten y,\delta_{\b}(x)]
\end{equation}

The Drinfeld double $\ga$ of $\b$ is the Lie algebra defined
as follows. As a vector space, $\gb=\b\oplus\b^*$. The pairing
$\iip{\cdot}{\cdot}:\b\ten\b^*\to\sfk$ extends uniquely to a symmetric,
non--degenerate bilinear form on $\gb$, such that $\b,\b^*$
are isotropic subspaces. The Lie bracket on $\gb$ is then
defined as the unique bracket compatible with $\iip{\cdot}
{\cdot}$, \ie such that
\[
\iip{[x,y]}{z}=\iip{x}{[y,z]}
\]                                    
for all $x,y,z\in\gb$. It coincides with $[\cdot,\cdot]_{\b}$ on $\b$, and
with the bracket induced by $\delta_{\b}$ on $\b^*$. The
mixed bracket for $b\in\b,\phi\in\b^*$ is then equal to
\[
[b,\phi]=
\sfad^*(b)(\phi)-\sfad^*(\phi)(b)
=
\sfad^*(b)(\phi)
+\phi\otimes\id_\b\circ\delta(b)
\]
where $\sfad^*$ denotes the coadjoint action of $\b$ on $\b^*$ and of $\b^*$ on $\b$,
respectively. 

The Lie algebra $\gb$ is a (topological) quasitriangular Lie bialgebra, with cobracket 
$\delta=\delta_{\b}\oplus(-\delta_{\b^*})$, where $\delta_{\b^*}$ is the (topological)
cobracket on $\b^*$ induced by $[\cdot,\cdot]_{\b}$, and $r$--matrix $r\in\gb\ctp\gb$
corresponding to the identity in $\End(\b)\simeq\b\ctp\b^*\subset\gb\ctp\gb$. Explicitly,
if $\{b_i\}_{i\in I},\{b^i\}_{i\in I}$ are dual bases of $\b$ and $\b^*$ respectively, then
$r=\sum_{i\in I}b_i\ten b^i\in\b\ctp\b^*$.

\subsection{Drinfeld--Yetter modules}\label{ss:pre-DY}

A triple $(V,\pi,\pi^*)$ is a Drinfeld--Yetter module over a Lie bialgebra
$(\b,[\cdot,\cdot]_{\b}, \delta_{\b})$ if $(V,\pi)$ is a $\b$--module, that
is the map $\pi:\b\otimes V\to V$ satisfies
\begin{equation}\label{eq:action}
\pi\circ[\cdot,\cdot]_{\b}=\pi\circ(\id\ten\pi)-\pi\circ(\id\ten\pi)\circ(21)
\end{equation}
$(V,\pi^*)$ is a $\b$--comodule, that is the map $\pi^*:V\to \b\otimes V$
satisfies
\begin{equation}\label{eq:coaction}
\delta\circ\pi^*=(21)\circ(\id\ten\pi^*)\circ\pi^*-(\id\ten\pi^*)\circ\pi^*
\end{equation}
and the maps $\pi,\pi^*$ satisfy the following compatibility condition
in $\End(\b\ten V)$:
\begin{equation}\label{eq:actcoact}
\pi^*\circ\pi-\id\ten\pi\circ(12)\circ\id\ten\pi^*=[\cdot,\cdot]_{\b}\ten\id\circ\id\ten\pi^*-
\id\ten\pi\circ\delta_{\b}\ten\id
\end{equation}
The category $\DrY{\b}$ is a symmetric tensor category. 

In terms of representations
of the Drinfeld double, $\DrY{\b}$ is equivalent to the category $\Eq{\gb}$ 
of equicontinuous $\gb$--modules \cite{ek-1}. 
Roughly speaking, a $\gb$--module is equicontinuous if the action 
of $\b^*$ is locally finite.
In particular, there is a functor $\Eq{\gb}\to\DrY{\b}$ which assign to any 
$(V,\pi)\in\Eq{\gb}$, the Drinfeld--Yetter $\b$--module $(V,\pi,\pi^*)$ where $\pi$ is 
restricted to $\b\subset\gb$, and the coaction $\pi^*$ is given by
\begin{equation}
\pi^*(v)=\sum_{i}b_i\ten b^i\cdot v\in\b\ten V
\end{equation}
The equicontinuity condition ensures that the sum is finite and the coaction well--defined. 
Conversely, given a Drinfeld--Yetter $\b$--module $(V,\pi,\pi^*)$, the action of $\phi\in\b^*$ on $V$ 
is defined by the formula
\begin{equation}\label{eq:dy-to-eq}
\phi\cdot v=\phi\ten\id_V\circ\pi^*(v)
\end{equation}
The compatibility condition \eqref{eq:actcoact} guarantees that this
lifts to an equicontinuous action of the Drinfeld double $\gb$.
One can prove that this is an equivalence of symmetric tensor categories.

\subsection{Restricted Drinfeld double}\label{ss:res double}

Let $\b=\bigoplus_{n\in\IN}\b_n$ be an $\IN$--graded Lie bialgebra with \fd homogeneous 
components. Its restricted dual $\b^{\star}=\bigoplus_{n\in\IN}\b_n^*$ and its
restricted Drinfeld double $\grb=\b\oplus\b^{\star}$ are also Lie bialgebras with 
cobrackets $\delta_{\b^\star}=[\,,\,]_\b^t$ and $\delta_{\grb}=\delta_\b-\delta_{\b^{\star}}$,
respectively. Moreover, since $\b^\star$ is dense in $\b^*$, the Lie algebra $\grb$ is dense 
in $\gb$. Therefore, any continuous  action of $\grb$ extends automatically to one of $\gb$.
One can show easily that this induces a canonical isomorphism $\Eq{\gb}\simeq\Eq{\grb}$.
In particular, one has $\DrY{\b}\simeq\Eq{\grb}$.

\subsection{Completions}\label{ss:completions}

Let $\Sf_{\b}:\DrY{\b}\to\vect$, $\Sf_{\gb}:\Eq{\gb}\to\vect$ be the forgetful
functors and $\DYA{\b}{}=\sfEnd{\Sf_{\b}}$, $\whUA{\gb}{}=\sfEnd{\Sf_{\gb}}$
the corresponding algebras of endomorphisms. Since the equivalence $\Eq
{\gb}\simeq\DrY{\b}$ preserves the underlying vector space and commutes
with the forgetful functors, there is a canonical isomorphism $\DYA{\b}{}
\simeq\whUA{\gb}{}$. In particular, we can think of $U\gb$ as a subalgebra
in $\DYA{\b}{}$.

Since the equivalence preserves the tensor structure, the same identification
holds for the $n$--folds forgetful functor $\Sf^{\boxtimes n}(V_1,\dots, V_n)=
V_1\ten\cdots\ten V_n$, \ie 
\begin{equation}\label{eq:identification}
\DYA{\b}{n}=
\sfEnd{\Sf_{\b}^{\boxtimes n}}\simeq\sfEnd{\Sf_{\gb}^{\boxtimes n}}=
\whUA{\gb}{\ten n}
\end{equation}
and we can consider $U\gb^{\ten n}$ as a subalgebra in $\DYA{\b}{n}$.

Under the identification \eqref{eq:identification}, the $r$--matrix of $\gb$,
$r_{\b}=\sum_i b_i\otimes b^i\in\b\wh{\otimes}\b^*\subset\gb\wh{\otimes}
\gb$, where $\{b_i\}$ and $\{b^i\}$ are dual bases of $\b$ and $\b^*$, 
corresponds to the element of $\DYA{\b}{2}$ given by the 
maps $r_{VW}\in\End_\sfk(V\ten W)$, $V,W\in\DrY{\b}$, defined by
\begin{equation}\label{eq:rmatrix}
r_{VW}=\pi_V\ten\id\circ\,(12)\,\circ\id\ten\pi^*_W
\end{equation} 

\subsection{Etingof--Kazhdan quantisation}\label{ss:ek-quant}

In \cite{ek-1}, Etingof and Kazhdan give an explicit procedure to construct
a quantisation of $\b$, that is a Hopf algebra $\EK\b$ over $\hext{\sfk}$
endowed with an isomorphism
\[
\EK\b/\hbar\EK\b\simeq U\b
\]
of Hopf algebras, which induces the cobracket $\delta_{\b}$ on $\b$.

The construction proceeds as follows. One considers the Drinfeld category
$\hDrY{\b}{\Phi}$ of deformation Drinfeld--Yetter $\b$--modules,\ie
topologically free $\hext{\sfk}$--modules with a Drinfeld--Yetter structure over $\b$,
with associativity and commutativity constraints given by
\[
\Phi_{UVW}=\Phi(\hbar\Omega_{12},\hbar\Omega_{23})
\aand
\beta_{VW}=
(12)\circ e^{\hbar\Omega/2}
\]
where $U,V,W\in\hDrY{\b}{}$, $\Omega=r+r^{21}$, and $\Phi$ is a fixed Lie
associator. Let $\ff:\hDrY{\b}{\Phi}\to{\khvect}$ be the forgetful
functor. Etingof and Kazhdan construct an explicit tensor structure on $\ff$,
\ie a collection of natural isomorphisms
\[
J\EEK_{VW}:\ff(V)\ten\ff(W)\to\ff(V\ten W)
\]
which are the identity modulo $\hbar$ and satisfy the relation 
\begin{equation}\label{eq:kill Phi}
\ff(\Phi)\circ J\EEK_{U\ten V, W}\circ(J\EEK_{U, V} \ten\id)=
J\EEK_{U, V\ten W}\circ(\id\ten J\EEK_{V, W})
\end{equation}
in $\Hom(\ff(U)\ten\ff(V)\ten\ff(W), \ff(U\ten(V\ten W)))$.

The algebra $\hDYA{\b}{}=\sfEnd{\ff}$ is a topological Hopf algebra, with
coproduct induced by the tensor product in $\hDrY{\b}{}$. Twisting $\hDYA{\b}{}$ 
by $J\EEK$ produces a new Hopf algebra, with a coassociative
deformation coproduct $\Delta_J$. In order to produce a quantisation
of $\b$, one considers the Drinfeld--Yetter module corresponding to the 
Verma module 
\[
M_{\b}=\ind_{\b^*}^{\gb}\IC\simeq U\b
\]
and shows that there is a natural embedding $\ff(M_{\b})\subset\sfEnd{\ff}$.
The coproduct $\Delta_J$ induces a coproduct on $\ff(M_{\b})$ which
can explicitely computed as the composition
\[
\xymatrix{\ff(M_{\b})\ar[r]^(.4){\ff(\Delta_0)} & \ff(M_{\b}\ten M_{\b})\ar[rr]^{(J\EEK_{M_\b,M_\b})^{-1}} && 
\ff(M_{\b})\ten\ff(M_{\b})}
\]
This induces a Hopf algebra structure on the vector space $\ff(M_{\b})\simeq
\hext{U\b}$, which quantizes the Lie bialgebra $\b$. 
In \cite{ek-2}, Etingof and Kazhdan showed that the construction
of the quantum enveloping algebra $\ff(M_{\b})$ is universal. 
In \ref{ss:prop-intro}--\ref{ss:EK-univ}, we explain the precise meaning of
this statement.

\subsection{$\PROP$s \cite{La,mac,ee,ATL1}}\label{ss:prop-intro} %

A $\PROP$ is a $\sfk$--linear, strict, symmetric monoidal category
$\sfP$ whose objects are the non--negative integers, and such that
$[n]\ten[m]=[n+m]$. In particular $[0]$ is the unit object, and $[1]^
{\ten n}=[n]$. A morphism of $\PROP$s is a symmetric monoidal
functor $\G:\sfP\to\sfQ$ which is the identity on objects,
and is endowed with the trivial tensor structure
\[\id:\G[m]_{\C}\otimes\G[n]_{\C}=[m]_{\D}\otimes[n]_{\D}=[m+n]_{\D}=\G([m+n]_{\C})\]

Fix henceforth a complete bracketing $b_n$ on $n$ letters for
any $n\geqslant   2$, and set $\bfb=\{b_n\}_{n\geqslant   2}$. A {\it module}
over $\sfP$ in a symmetric monoidal category $\N$ is a symmetric
monoidal functor $(\G,J):\sfP\to\N$ such that\footnote{In a monoidal
category $(\C,\otimes)$, $V^{\otimes n}_{b_n}$ denotes the $n$--fold
tensor product of $V\in\C$ bracketed according to $b_n$. For example
$V^{\otimes 3}_{(\bullet\bullet)\bullet}=(V\otimes V)\otimes V$.}
\begin{equation*}
\G([n])=\G([1])^{\otimes n}_{b_n}
\end{equation*}
and the following diagram is commutative
\begin{equation*}
\xymatrix{
\G([m])\otimes \G([n]) \ar@{=}[d]\ar[rr]^{J_{[m],[n]}} && \G([m+n])\ar@{=}[d]\\
\G([1])^{\otimes m}_{b_m}\otimes \G([1])^{\otimes n}_{b_n}\ar[rr]_{\Phi}&& \G([1])^{\otimes(m+n)}_{b_{m+n}}
}
\end{equation*}
where $\Phi$ is the 
associativity constraint in $\N$. A {\it morphism} of modules over $\sfP$ is a natural
transformation of functors. The category of $\sfP$--modules is
denoted by $\Fun_{\bfb}^\otimes(\sfP,\N)$.

\subsection{The Karoubi envelope}\label{ss:karoubi}

Recall that the Karoubi envelope of a category $\C$ is the category $\Kar
{\C}$ whose objects are pairs $(X,\pi)$, where $X\in\C$ and $\pi:X\to X$ is
an idempotent. The morphisms in $\Kar{\C}$ are defined as
\[
{\Kar{\C}}((X,\pi), (Y,\rho))=\{f\in\C(X,Y)\;|\; \rho\circ f=f=f\circ\pi\}
\]
with $\id_{(X,\pi)}=\pi$. In particular, 
$
{\Kar{\C}}((X,\id),(Y,\id))={\C}(X,Y)
$, 
and the functor $\C\to\Kar{\C}$, mapping $X\mapsto(X,\id)$, $f\mapsto f$, 
is fully faithful.

Every idempotent in $\Kar{\C}$ splits canonically. Namely,
if $q\in\Kar{\C}((X,\pi),(X,\pi))$ satisfies $q^2=q$, the maps
\[i=q:(X,q)\to(X,\pi)\aand p=q:(X,\pi)\to(X,q)\]
satisfy $i\circ p=q$ and $p\circ i=\id_{(X,q)}$. 
 
We denote by $\cKar{\sfP}$ the closure under 
infinite direct sums of the Karoubi completion of $\sfP$. 
It is then clear that, if $\N$ is Karoubi complete, 
there is an essentially unique equivalence
 $\Fun_{\bfb}^\otimes(\cKar{\sfP},\N)\simeq\Fun_{\bfb}^\otimes(\sfP,\N)$.
 
\subsection{Example}\label{ss:prop-la}
Let $\preLA$ be the $\PROP$ generated by a morphism
$\mu:[2]\to[1]$ 
subject to the relations
\begin{equation}\label{eq:bracket}
\mu\circ(\id_{[2]}+(1\,2))=0
\aand \mu\circ(\mu\ten\id_{[1]})\circ(\id_{[3]}+(1\,2\,3)+(3\,1\,2))=0
\end{equation}
as morphisms $[2]\to[1]$ and $[3]\to[1]$ respectively.
Then, 
there is 
a canonical isomorphism of categories $\Fun_{\bfb}(\preLA,\vect_{\sfk})\simeq\preLA(\sfk)$,
where $\preLA(\sfk)$ is the category of Lie algebras over $\sfk$.


\subsection{The $\PROP$s $\LCA$ and $\LBA$}\label{ss:exa}

The $\PROP$ of Lie coalgebras $\LCA$ is generated by 
a morphism $\delta:[1]\to[2]$ satisfying
\begin{equation}\label{eq:cobracket}
(\id_{[2]}+(1\,2))\circ\delta=0
\aand
(\id_{[3]}+(123)+(312))\circ(\delta\ten\id_{[1]})\circ
\delta=0
\end{equation}

There is a natural identification of $\PROP$s 
\begin{equation}
\Theta:\preLCA\to\preLA\op
\end{equation}
defined by $\Theta(\delta)=\mu$. The relation between the functor $\Theta$ and
the standard duality between Lie algebras and Lie coalgebras is easily described.
Let $\b$ be a Lie algebra and $\c$ a Lie coalgebra, with a compatible pairing 
$\iip{}{}:\b\ten\c\to\sfk$, \ie
such that $\iip{[b_1,b_2]_{\b}}{c}=\iip{b_1\ten b_2}{\delta_{\c}(c)}$ for any $b_1,b_2\in\b$
and $c\in\c$. Let $\G_\b:\LA\to\vect_\sfk,\G_\c:\LCA\to\vect_\sfk$ be the corresponding
realisation functors, then for any $T\in\LA([N],[n])$, 
$\ul{b}_{N}\coloneqq b_1\ten\dots\ten b_N\in\b^{\ten N}$ and $\ul{c}_n\coloneqq c_1\ten\dots\ten c_n\in\c^{\ten n}$,
one has 
\begin{equation}\label{eq:contravariant}
\iip{\G_{\b}(T)(\ul{b}_N)}{\ul{c}_n}=
\iip{\ul{b}_N}{\G_{\c}(\Theta(T))(\ul{c}_n)}
\end{equation}

Finally, the $\PROP$ of Lie bialgebras $\LBA$ is generated by $\mu:[2]\to[1]$ 
and $\delta:[1]\to[2]$ satisfying \eqref{eq:bracket}, \eqref{eq:cobracket}, and 
the \emph{cocycle condition}
\begin{equation}\label{eq:cocycon}
\delta\circ\mu=(\id_{[2]}-(21))\circ\id\ten\mu\circ\delta\ten\id
\circ(\id_{[2]}-(21))
\end{equation}

\subsection{Etingof--Kazhdan quantisation in $\LBA$}\label{ss:EK-univ}

In \cite{ek-2}, Etingof and Kazhdan showed that the construction of $J=J\EEK_
{M_\b,M_\b}$ is universal, \ie that it can be realised in the $\PROP$ $\LBA$. To
this end, one first replaces the module $M_{\b}$ with a universal Drinfeld--Yetter
module in $\LBA$, by constructing an action and a coaction of the Lie bialgebra
$[1]\in\LBA$ on $M\coloneqq S[1]$. The twist $J$ is then defined, using the same formulae
as in \cite{ek-1}, as an element 
\[
J\in\hext{\LBA(M\ten M,M\ten M)}
\]
It induces a universal quantisation functor, that is a functor $\mathsf{Q}$ from
the $\PROP$ of Hopf algebras\footnote{The $\PROP$ $\HA$ is generated by the
morphisms $m:[2]\to[1]$, $\iota:[0]\to[1]$, $\Delta:[2]\to[1]$, $\epsilon:[1]\to[0]$, 
$S,S^{-1}:[1]\to[1]$ with the relations coming from the Hopf algebra axioms.}
$\HA$ to $\LBA$, mapping $[1]_{\scriptscriptstyle{\mathsf{HA}}}$ to $S[1]\euniv$.
 
A universal interpretation of the fiber functor $(\ff,J\EEK)$, rather than of the Hopf
algebra $(\ff(M_\a),J^{-1}\ff(\Delta_0))$ alone, will be given in the Section \ref{s:univDY}, 
by using the $\PROP$ of universal Drinfeld--Yetter modules.

\section{Schur functors}\label{s:sch}

We review in this section some basic facts about the cohomology of Schur functors
which are due to Enriquez \cite[Sec. 1]{e3}, and will be used repeatedly. The exposition
follows the approach to the theory of Schur functors of Baez and Trimble \cite{baez}.

\subsection{Schur functors}\label{ss:sch}

Let $\Cat$ be the $2$--category of categories and
$\SymCat$ the $2$--category of $\sfk$--linear, additive, 
Karoubi closed, symmetric monoidal categories.

\begin{definition}
A \emph{Schur functor} is an endomorphism of the
forgetful $2$--functor $\Sf:\SymCat\to\Cat$. 
That is, a collection of endofunctors $F_{\sfC}:\sfC\to\sfC$ in $\Cat$, indexed by
objects in $\SymCat$, and invertible natural transformations $F_{\G}$ 
in $\Cat$,
\begin{equation}
\xymatrix{
\sfC_1 \ar[r]^{\G} \ar[d]_{F_{\sfC_1}} & \sfC_2 \ar[d]^{F_{\sfC_2}} \ar@{=>}|{F_{\G}}[dl]\\
\sfC_1 \ar[r]_{\G}& \sfC_2
}
\end{equation}
indexed by functors $\G\in\SymCat(\sfC_1,\sfC_2)$, and such that
$F_{\id_{\sfC}}=\id_{F_{\sfC}}$ and $F_{\G_2\circ\G_1}=F_{\G_1}\circ_{\mathsf{h}}F_{\G_2}$,
\ie
\begin{equation}
\xymatrix@C=0.7in{
\sfC_1 \ar[r]^{\G_1} \ar[d]_{F_{\sfC_1}} & 
\sfC_2 \ar|{F_{\sfC_2}}[d] \ar@{=>}|{F_{\G_1}}[dl] \ar[r]^{\G_2}&
\sfC_3 \ar[d]^{F_{\sfC_3}} \ar@{=>}|{F_{\G_2}}[dl] \ar@{}|{=}[dr]&
\sfC_1 \ar[r]^{\G_2\circ\G_1} \ar[d]_{F_{\sfC_1}} & \sfC_3 \ar[d]^{F_{\sfC_3}} \ar@{=>}|{F_{\G_2\circ\G_1}}[dl]
\\
\sfC_1 \ar[r]_{\G_1}& 
\sfC_2 \ar[r]_{\G_2} &
\sfC_3 &
\sfC_1 \ar[r]_{\G_2\circ\G_1}& \sfC_3
}
\end{equation} 
A morphism of Schur functors $\phi:F^1\to F^2$ is a collection of natural
transformation $\phi_{\sfC}:F^1_{\sfC}\to F^2_{\sfC}$, indexed by $\sfC\in\SymCat$,
such that, for any functor $\G\in\SymCat(\sfC_1,\sfC_2)$,
\begin{equation}
\xymatrix{
\sfC_1\ar@{=}[r] \ar[d]_{F_{\sfC_1}^2}&
\sfC_1 \ar@{=>}|{\phi_{\sfC_1}}[dl] \ar[r]^{\G} \ar|{F_{\sfC_1}^1}[d] & \sfC_2 \ar[d]^{F^1_{\sfC_2}} \ar@{=>}|{F^1_{\G}}[dl] & &\ar@{}|{=}[dll] 
\sfC_1\ar[r]^{\G} \ar[d]_{F^2_{\sfC_1}}&
\sfC_2 \ar@{=>}|{F^2_{\G}}[dl] \ar@{=}[r] \ar|{F_{\sfC_1}^2}[d] & \sfC_2 \ar[d]^{F^1_{\sfC_2}} \ar@{=>}|{\phi_{\sfC_2}}[dl]
\\
\sfC_1 \ar@{=}[r] &
\sfC_1 \ar[r]_{\G}& \sfC_2 & &
\sfC_1 \ar[r]_{\G} &
\sfC_2 \ar@{=}[r]& \sfC_2
}
\end{equation} 
\end{definition}

The category of Schur functors 
$\Sch=\sfEnd{\Sf}$
is endowed with the following operations:
\begin{itemize}
\item {\bf Direct sum.} For any $F_1, F_2\in\Sch$, we set
\begin{equation}
F^1\oplus F^2 = \oplus \circ F^1\times F^2\circ \diag
\end{equation}
where $\oplus:\Sf\times\Sf\to\Sf$ is thought of as a morphism of $2$--functors, 
and $\diag: \Sf\to\Sf\times\Sf$ is the diagonal. The neutral element is the zero functor
$\Sigma_0\in\Sch$, which assigns to each object in $\sfC$ the zero object in $\sfC$.
\item {\bf Tensor product.} For any $F_1, F_2\in\Sch$, we set
\begin{equation}
F^1\ten F^2 = \ten \circ F^1\times F^2\circ \diag
\end{equation}
where, as before, $\ten:\Sf\times\Sf\to\Sf$.
The neutral element is the unit functor $T^0\in\Sch$,
which assigns to each object in $\sfC$ the unit object in $\sfC$.
\end{itemize}
Both assignments extend to natural transformations and give rise
to functors $\oplus,\ten:\Sch\times\Sch\to\Sch$, which endow $\Sch$
with a natural structure of additive symmetric monoidal category. 

\subsection{Representability}\label{ss:rep-sch}

Let $\sfk\SS$ denotes the permutation algebroid (\ie the free $\PROP$
generated by permutations) and $\kS$ be its additive and Karoubian envelope. 

\begin{theorem}\cite{baez}
The forgetful $2$--functor $\Sf:\SymCat\to\Cat$ is 
represented by $\kS$, \ie there is an equivalence
of $2$--functors
\begin{equation}
\Sf\simeq\SymCat(\kS,-)
\end{equation}
In particular, $\Sch\simeq\kS$ in $\Cat$.
\end{theorem}

\begin{pf}
The proof is essentially a $2$--categorical version
of the representability of the forgetful functor from the category of representations of
an associative algebra to vector spaces. Namely, for any $\sfC\in\SymCat$, 
there is a canonical functor $\sfU_{\sfC}:\SymCat(\kS,\sfC)\to\sfC$ defined by
$\sfU_{\sfC}(\G)=\G[1]$, for any $\G\in\SymCat(\kS,\sfC)$, and 
$\sfU_{\sfC}(\phi)=\phi_{[1]}:\G[1]\to\G'[1]$,
for any natural transformation $\phi:\G\Rightarrow\G'$.

For any functor $\F:\sfC_1\to\sfC_2$, the natural transformation $\sfU_{\F}$
\begin{equation}
\xymatrix@C=0.7in{
\SymCat(\kS,\sfC_1) \ar[d]_{\F_*}  \ar[r]^(.6){\sfU_{\sfC_1}}  & \sfC_1\ar[d]^{\F}
\ar@{=>}|(.4){\sfU_{\F}}[dl]\\
\SymCat(\kS,\sfC_2)  \ar[r]^(.6){\sfU_{\sfC_2}}  & \sfC_2
}
\end{equation}
is given by the identity on $\F\circ\G[1]$.

It is easy to see that this defines an essentially unique equivalence 
of $2$--functors $\Sf\simeq\SymCat(\kS,-)$. It then follows from Yoneda lemma
$\Sch=\sfEnd{\Sf}\simeq\kS$.
\end{pf}

\subsection{Abelianity}\label{ss:ab-sch}

We will use of the following general fact.

\begin{proposition}
Let $A$ be an associative $\sfk$--algebra, $\sfC_A$ the 
corresponding algebroid (\ie $\sfC_A$ is the category with 
one object $\bullet$ and $\mathsf{End}_{\sfC_A}(\bullet)=A$), and 
$\ol{\sfC}_A$ its additive and Karoubi envelope. Then $\ol{\sfC}_A$ 
is equivalent to the category $\mathsf{Proj}(A\op)$ of projective 
$A\op$--modules. 
\end{proposition}

\begin{pf}
Let $\sfC_A^{\oplus}$ be the additive envelope of $\sfC_A$. Then
the functor $\sfC_A^{\oplus}\to\Rep A\op$ mapping the generating object $[1]$
to $A$ induces an equivalence of categories $\sfC_A^{\oplus}\simeq\mathsf{Free}(A\op)$.
It follows
\begin{equation}
\ol{\sfC}_A=\mathsf{Kar}(\sfC_A^{\oplus})\simeq\mathsf{Kar}(\mathsf{Free}(A\op))\simeq\mathsf{Proj}(A\op)
\end{equation}
\end{pf}

\begin{corollary}
If $A$ is hereditary (resp. semisimple), $\ol{\sfC}_A$ is abelian (resp. semisimple).
In particular, the category $\kS$, and therefore $\Sch$, is a semisimple category.
\end{corollary}

\subsection{Representations of $\SS_n$}\label{ss:repSn}

Recall that the set of irreducible representations $\hatS{n}$ of the
symmetric group $\SS_n$ is in bijection with minimal idempotents
in $k\SS_n$, modulo the equivalence relation $p\sim upu^{-1}$,
$u\in k\SS_n^\times$. We henceforth regard $\hatS{n}$ as a subset
of $k\SS_n$ by choosing a representative for each class, and set
$\hatS{}=\bigsqcup_{n\geqslant   0}\hatS{n}$. If $\pi\in\hatS{n}$, we set
$|\pi|=n$.
We proved in \ref{ss:ab-sch} that the category of Schur functor 
is semisimple and equivalent to category $\Rep\sfk\SS$ of 
representations of $\sfk\SS=\bigoplus_N\sfk\SS_N$.
It follows that, up to isomorphism, any Schur functor has the form:
\begin{equation}\label{eq:sch-explicit}
F_{\sfC}(X)=\bigoplus_{\pi\in\hatS{}}\pi\left(X^{\ten |\pi|}\right)^{\oplus m_{\pi}}
\end{equation} 
for some $m_{\pi}\in\IN\cup\{\infty\}$.

\subsection{Schur bifunctors}\label{ss:bisch}

A \emph{Schur bifunctor} is a morphism of $2$--functors 
from $\Sf\times\Sf$ to $\Sf$, where $\Sf\times\Sf(\sfC)=\sfC\times\sfC$.
The category of Schur bifunctors is denoted 
$\Sch_2=\mathsf{Hom}(\Sf\times\Sf,\Sf)$.
In particular, $\oplus$ and $\otimes$ are Schur bifunctors.

Schur bifunctors can be obtained from Schur functors by 
using the following operations.
\begin{itemize}
\item {\bf External tensor product.} For any $F^1,F^2\in\Sch$,
set
\[F^1\boxtimes F^2=\ten\circ F^1\times F^2\]
\item {\bf Coproduct.} For any $F\in\Sch$, set
\[\Delta(F)=F\circ\oplus\]
\end{itemize}

\begin{example}
If $S=\bigoplus_{n\geqslant 0}S^n$ and $\wedge=\bigoplus_{n\geqslant 0}
\wedge^n$ are the symmetric and exterior algebra functors, then
\[\Delta(S)\simeq S\boxtimes S\aand
\Delta(\wedge)\simeq\wedge\boxtimes\wedge\]
\end{example}

The results from \ref{ss:rep-sch} and \ref{ss:ab-sch} readily extends to
$\Sch_2$.

\begin{theorem}\hfill
\begin{enumerate}
\item The $2$--functor $\Sf\times\Sf:\SymCat\to\Cat$ is represented by
the category $\kS\times\kS$, \ie there is an equivalence
\begin{equation}
\Sf\times\Sf\simeq\SymCat(\kS\times\kS,-)
\end{equation}
\item $\Sch_2\simeq\kS\times\kS$ in $\Cat$.
\item $\Sch_2$ is a semisimple abelian category. 
\end{enumerate}
\end{theorem}

\begin{pf}
The representability of $\Sf\times\Sf$ is straightforward. Then, by Yoneda lemma,
one gets
\[
\Sch_2=\mathsf{Hom}(\Sf\times\Sf,\Sf)\simeq\kS\times\kS
\]
From Corollary \ref{ss:ab-sch}, we conclude that $\Sch_2$ is semisimple and abelian.
\end{pf}

\subsection{Cohomology of Schur (bi)functors}\label{ss:co}

Since the category of Schur (bi)functors is abelian, we can consider
the cohomology of complexes in $\Sch$ or $\Sch_2$. 

\begin{proposition}\cite[Prop. 1.3]{e3}
Let $(F^n,d^n)_{n\geqslant 0}$ be a complex in $\Sch$. Then
\[H^i(\Delta(F^{\bullet}),\Delta(d^{\bullet}))\simeq\Delta(H^i(F^{\bullet},d^{\bullet}))\]
\end{proposition}

\begin{pf}
It is enough to observe that the functor $\Delta=-\circ\oplus:\Sch\to\Sch_2$ is 
additive, and therefore exact, due to the fact that $\Sch$ and $\Sch_2$ are
semisimple.
\end{pf}

\subsection{The Hochschild complex}\label{ss:hosch}

The Hochschild complex $(SV^{\otimes\bullet},d_H)$ of a symmetric
coalgebra $SV$ can be interpreted as a complex of Schur functors as
follows. 

Let $\Sigma_1$ be the Schur functor $\id_{\Sf}\in\Sch$ and, for any $n
\geqslant 1$, set $\Sigma_n=\Sigma_1^{\oplus n}=\Sigma_1\oplus\cdots
\oplus\Sigma_1$. Let $i_0:\Sigma_0\to\Sigma_1$ be the inclusion of the
zero object $\Sigma_0\in\Sch$, and $\delta:\Sigma_1\to\Sigma_2$ the
diagonal morphism, \ie for any $\sfC\in\SymCat$ and $X\in\sfC$,
$\delta_{\sfC, X}=(\id_X,\id_X):X\to X\oplus X$.
There are natural transformations $\{\delta_n^i\}_{i=0}^{n+1}:\Sigma
_n\to\Sigma_{n+1}$ defined as follows
\footnote{
In the category $\vect$, the Schur functor $\Sigma_n:\vect\to\vect$ 
is given by $V\to V^{\oplus
n}=V\otimes\sfk^n$, and the natural transformations $\{\delta^n_i\}_{i=0}^{n+1}:
\Sigma_n\to\Sigma_{n+1}$ are induced by the maps $\sfk^n\to\sfk^{n+1}$ given
by
\[(x_1,\ldots,x_n)\to
\left\{\begin{array}{cc}
(0,x_1,\ldots,x_n)&i=0\\
(x_1,\ldots,x_{i-1},x_i,x_i,x_{i+1},\ldots,x_n)&1\leqslant i\leqslant n\\
(x_1,\ldots,x_n,0)&i=n+1
\end{array}\right.\]
} 
\[\begin{array}{ccc}
(i=0) & &
\xymatrix@C=1in{\Sigma_n=\Sigma_0\oplus\Sigma_n\ar[r]^{i_0\oplus\id_{\Sigma_n}} 
& \Sigma_1\oplus\Sigma_n=\Sigma_{n+1}}
 \\
(i=n+1) & &
\xymatrix@C=1in{\Sigma_n=\Sigma_n\oplus\Sigma_0\ar[r]^{\id_{\Sigma_n}\oplus i_0} 
& \Sigma_n\oplus\Sigma_1=\Sigma_{n+1}}
\end{array}
\]
and, for $1\leqslant i\leqslant n$,
\[
\xymatrix@C=1in{\Sigma_n=\Sigma_{i-1}\oplus\Sigma_1\oplus\Sigma_{n-i}\ar[r]^{\id_{\Sigma_{i-1}}\oplus\delta\oplus\id_{\Sigma_{n-i}}} 
& \Sigma_{i-1}\oplus\Sigma_2\oplus\Sigma_{n-i}=\Sigma_{n+1}}
\]

The natural transformations $\{\delta_i^n\}$ give rise to a cosimplicial structure on the tower of Schur functors
$S^{\otimes n}=S\circ\Sigma_n$, whose associated differential is the Hochschild differential $d_H$. The latter restricts to zero on $T^\bullet
\subset S^{\otimes\bullet}$, where $T^n=\Sigma_1\ten\cdots\ten\Sigma_1$, and gives rise to a quasi--isomorphism
\[\iota:\left(\wedge^\bullet,0\right)\to (S^{\otimes\bullet}, d_H)\]

\subsection{The diagonal Hochschild complex}

\newcommand {\copr}[1]{\Delta(#1)}

By Proposition \ref{ss:co}, the map $\Delta(\iota)$ induces
a quasi--isomorphism\footnote{Here, and in the sequel, the
notation $S^{\otimes\bullet}\boxtimes S^{\otimes\bullet}$
refers to the complex whose $n$th term is $S^{\otimes n}
\boxtimes S^{\otimes n}$, not to the total complex underlying
the exterior product of the complex $S^{\otimes\bullet}$ with
itself.}

\begin{equation}
H^{n}(S^{\otimes\bullet}\boxtimes S^{\otimes\bullet}, d_H\boxtimes d_H)
=
H^{n}(\Delta(S^{\otimes\bullet}),\Delta(d_H))
\simeq
\Delta(\wedge^n)
=
\bigoplus_{j=0}^n\wedge^j\boxtimes\wedge^{n-j}
\label{eq:S S}
\end{equation}
To work out $\Delta(\imath)$ explicitly, note that it arises from the restriction
to $\copr{\wedge^n}$ of the inclusion $d_n=\copr{\iota_1^{\otimes n}}:\copr
{T^n}\subset\copr{S^{\otimes n}}\cong S^{\otimes n}\boxtimes S^{\otimes n}$,
where $\iota_1:T\to S$ is the inclusion. The restriction of $d_n$ to
$T^j\boxtimes T^{n-j}\subset \copr{T^n}$ is readily seen to be
\[\tau_j=
({\iota_1}^{\ten j}\ten {\iota_0}^{\ten n-j})\boxtimes ({\iota_0}^{\ten j}
\ten\iota_1^{\ten n-j})\]
where $\iota_0:T^0\to S$ is the inclusion of the unit. Since $d_n$ is equivariant
under $\SS_n$, its restriction to $\wedge^j\boxtimes\wedge^{n-j}\subset
\copr{\wedge^n}\subset\copr{T^n}$ is given by
\begin{equation}\label{eq:qiso}
\Alt_n^{(2)}\circ\tau_j\circ\Alt_j\ten\Alt_{n-j}
\end{equation}
where $\Alt_n=\frac{1}{n!}\sum_{\sigma\in\SS_n}(-1)^{\sigma}\sigma$
and $\Alt_n^{(2)}=\frac{1}{n!}\sum_{\sigma\in\SS_n}(-1)^{\sigma}\sigma\boxtimes\sigma$.

\begin{example}
For $n=2$, the restriction of $\Delta(\iota)$ to $\wedge^1\boxtimes\wedge^1$
is given by
\[
\Alt_2^{(2)}\circ\tau_1=
\frac{1}{2}
\left[
(\iota_1\ten\iota_0)\boxtimes(\iota_0\ten\iota_1)
-(\iota_0\ten\iota_1)\boxtimes(\iota_1\ten\iota_0)
\right]
\] 
In $\vect$, this reads: for any $V, W\in\vect$, $\Delta(\iota)$ is given on
$V\ten W\subset\wedge^2(V\oplus W)$ by
\[
v\ten w\mapsto \frac{1}{2}\left[
(v\ten 1)\ten(1\ten w)-(1\ten v)\ten(w\ten 1)
\right]
\]
\end{example}
\subsection{Tensor algebra}\label{ss:T-alg}

Set $T^0=1$, \ie the unit object in $\Sch$, and $T^n=\Sigma_1^{\ten n}=\Sigma_1\ten\cdots\ten\Sigma_1$ for any $n\geqslant 1$.
By Theorem \ref{ss:rep-sch}, $\Sch(T^n,T^n)=\sfk\SS_n$. The tensor algebra functor
$T=\bigoplus_{n\geqslant0}T^n$ is also endowed with a cosimplicial structure on $\{T^{\ten\bullet}\}$.
Namely, let $\Shcp:T\to T\ten T$ be the \emph{shuffle coproduct} defined on 
$T^n$ by 
\begin{equation}
\Shcp_n=\sum_{\substack{n_1+n_2=n\\\sigma\in\Sh(n_1,n_2)}}\psi_{n_1,n_2}\circ\sigma
: T^n\to\bigoplus_{n_1+n_2=n}T^{n_1}\ten T^{n_2}\subset T\ten T
\end{equation}
where $\Sh(n_1,n_2)\subset\SS_{n}$ is the set of $(n_1,n_2)$--shuffles, \ie
permutations $\sigma\in\SS_{n}$ such that $\sigma(i)<\sigma(j)$ whenever
$1\leqslant i<j\leqslant n_1$ or $n_1+1\leqslant i<j\leqslant n$, and $\psi_{n_1,n_2}$
is the \emph{deconcatenation} $T^n=T^{n_1}\ten T^{n_2}$. 
Then the Hochschild differential on $T$ is $d_H^n=\sum_{i=0}^{n+1}d_i^n$, with
face maps $\{d_i^n\}_{i=0}^{n+1}:T^{\otimes n}\to T^{\otimes(n+1)}$ given by
\[d_i^n 
=\left\{\begin{array}{ll}
1\otimes\id_{T^{\ten n}} &i=0\\
\id_{T^{\ten (i-1)}}\otimes\Shcp\otimes \id_{T^{\ten (n-i)}} &i=1,\ldots,n-1\\
\id_{T^{\ten n}}\otimes 1 &i=n
\end{array}\right.\]
The canonical inclusion $\Sym: S\hookrightarrow T$ preserves the differential defined on \ref{ss:hosch}.
and induces a morphism of complexes $(S^{\ten\bullet},d_H)\to(T^{\ten\bullet}, d_H)$.

\subsection{Duality in $\Sch$}\label{ss:dual-sch}
Let $\Sch^{\Pi}$ be the completion of $\Sch$ under infinite direct products.
That is, we formally add to $\Sch$ the objects $\prod_i F_i\in\Sch^{\IN}$ with morphisms
\begin{eqnarray}
\Sch^{\Pi}\left(\prod_iF_i, \prod_jF'_j\right)&=&\prod_j\bigoplus_{i}\Sch^{\Pi}\left(F_i,F'_j\right)
\end{eqnarray}
Objects in $\Sch^{\Pi}$ are well--defined, and universal, on symmetric monoidal categories
closed under infinite direct products. The equivalence $\Sch\simeq\Rep\sfk\SS$ extends the 
duality in $\Rep\sfk\SS$ to a contravariant functor 
$\Sch\to\Sch^{\Pi}$,
where
\[
\left(\bigoplus_iF_i\right)^*=\prod_iF_i^*
\]

\subsection{$\PROP$s and Schur bifunctors}\label{ss:prop}\label{ss:cohoprop}

A $\PROP$ $\sfP$ gives rise to a functor $\sfP_{\Sch}:\Sch_2\to\vect$
which is defined on a bifuntor $F=\bigoplus_i F_i\boxtimes G_i$ by
\footnote{The definition of $\sfP_{\Sch}$ requires to consider the
completion of $\sfP$ with respect to infinite direct limits, which,
by abuse of notation, we still denote by $\sfuP$.}
\begin{equation}\label{eq:PSch and Hom}
\sfP_{\Sch}(F)=\bigoplus_{i}\sfuP(F_i^*[1],G_i[1])
\end{equation}
Then for any morphism
\[
f=\sum_{i,j}f_{ij}\boxtimes g_{ij}: 
F=\bigoplus_i F_i\boxtimes G_i\to F'=\bigoplus_j F'_j\boxtimes G'_j
\]
we define $\sfP_{\Sch}(f):\sfP_{\Sch}(F)\to\sfP_{\Sch}(F')$ as follows.
For any $\phi=\sum_i\phi_i$ in $\sfP_{\Sch}(F)=\bigoplus_{i}\sfuP (F_i^*[1],G_i[1])$, 
we set
\[
\sfP_{\Sch}(f)(\phi)=\sum_j \left(\sum_i g_{ij}\circ\phi_i\circ f_{ij}^*\right)
\in\bigoplus_j\sfP((F_j')^*[1],G_j[1])=\sfP_{\Sch}(F')
\]

\begin{proposition}\label{pr:cohoprop}\cite[Prop. 1.2]{e3}
For any complex $(F^{\bullet},d^{\bullet})$ in $\Sch_2$, $(\sfP_{\Sch}
(F^{\bullet}), \sfP_{\Sch}(d^{\bullet}))$ is a complex of vector spaces, and
\[
H^i(\sfP_{\Sch}(F^{\bullet}), \sfP_{\Sch}(d^{\bullet}))
\simeq\sfP_{\Sch}(H^{i}(F^{\bullet}, d^{\bullet}))
\]
\end{proposition}

\begin{pf}
It is enough to observe that the functor $\sfP_{\Sch}:\Sch_2\to\vect$ is 
additive, and therefore exact.
\end{pf}

\subsection{Hochschild cohomology}\label{ss:coho}

The differential $\sfP_{\Sch}(d_H\boxtimes d_H)$ of the complex
$\sfP_{\Sch}(S^{\ten\bullet}\boxtimes S^{\ten\bullet})$
can be described more explicitely. From \ref{ss:dual-sch}, we have
\begin{equation}\label{ss:T-S-dual}
S^*=\prod_{n\geqslant0}S^n=:\wh{S}
\end{equation}
and the cosimplicial structure on $S^{\ten\bullet}$ described in \ref{ss:hosch}
induces a simplicial structure on $\wh{S}^{\ten \bullet}$ with associated differential $\partial_H$. 
Therefore, for any $\phi\in\sfP_{\Sch}(\wh{S}^{\ten n}[1],S^{\ten n}[1])$, one has
\[
\sfP_{\Sch}(d_H\boxtimes d_H)(\phi)=d_H\circ\phi\circ\partial_H
\]
Analogous considerations hold for the complex of Schur bifunctors 
$(T^{\ten\bullet}\boxtimes T^{\ten\bullet}, d_H\boxtimes d_H)$.

\begin{proposition}\label{pr:coho}
The following holds for any $\PROP$ $\sfP$.
\begin{enumerate}
\item The inclusion
\[
\sfuP(\wh{S}[1]^{\otimes \bullet},S[1]^{\otimes\bullet})\to\sfuP(\wh{T}[1]^{\otimes\bullet},T[1]^{\otimes\bullet})
\]
induced by the natural inclusion $\Sym: S\to T$ and projection $\Sym:\wh{T}\to\wh{S}$,
is a morphism of cosimplicial spaces.
\item The inclusion
\[
\sfP_{\Sch}(\wt{\iota}):
\left(\bigoplus_{j=0}^{\bullet}\sfuP\left(\wedge^j[1],\wedge^{\bullet-j}[1]\right),0\right)\to
\left(\sfuP(\wh{S}[1]^{\otimes\bullet},S[1]^{\otimes\bullet}),d_H\circ(-)\circ\partial_H\right)\]
obtained by \eqref{eq:qiso} is a quasi--isomorphism.
\end{enumerate}
\end{proposition}

\begin{pf}
$(1)$ It is enough to observe that, by duality, $\Sym:\wh{T}\to\wh{S}$ induces a morphism of
simplicial objects. 
 
$(2)$ We have
\[\begin{split}
H^i\left(\sfuP(\wh{S}[1]^{\otimes\bullet},S[1]^{\otimes\bullet}),d_H\circ(-)\circ\partial_H\right)
&=
H^i\left(\sfP_\Sch(S^{\otimes\bullet}\boxtimes S^{\otimes\bullet}),\sfP_\Sch(d_H\boxtimes d_H\right))\\
&=
\sfP_\Sch\left(H^i\left(S^{\otimes\bullet}\boxtimes S^{\otimes\bullet},d_H\boxtimes d_H\right)\right)\\
&=
\sfP_\Sch\left(\bigoplus_{j=0}^i\wedge^j\boxtimes\wedge^{i-j}\right)\\
&=
\bigoplus_{j=0}^{i}\sfuP\left(\wedge^j[1],\wedge^{i-j}[1]\right)
\end{split}\]
where the first and last equalities hold by definition of the functor $\sfP_\Sch$,
the second one by Proposition \ref{pr:cohoprop}, and the third one by \eqref{eq:S S}.
\end{pf}

The quasi--isomorphism $\sfP_{\Sch}(\wt{\iota})$ is described as follows. 
Let $\iota_0: T^0\to S$, $\iota_1: T^1\to S$ be the canonical inclusions,
and $\iota_0^*:\wh{S}\to T^0$, $\iota_1^*:\wh{S}\to T^1$ the corresponding projections.
Set $\tau_j': T^j\to S^{\ten i}$ by $\tau_j'=\iota_1^{\ten j}\ten\iota_0^{\ten i-j}$,
$\tau_{i-j}'':\wh{S}\to T^{i-j}$ by $\tau''_{i-j}=\iota_0^{\ten j}\ten\iota_1^{\ten i-j}$, and let
$\wt{\tau}_j':\wedge^j\to S^{\ten i}$, $\wt{\tau}_{i-j}'':\wh{S}^{\ten i}\to\wedge^{n-j}$
be the compositions with $\Alt_j$ and $\Alt_{n-j}$, respectively.
Then, for any $\phi\in\sfuP(\wedge^j[1],\wedge^{i-j}[1])$, 
$\sfP_{\Sch}(\wt{\iota})(\phi)\in\sfuP(\wh{S}[1]^{\ten i}, S[1]^{\ten i})$ is given by
\begin{equation}\label{eq:qiso-P}
\sfP_{\Sch}(\phi)=\frac{1}{i!}\sum_{\sigma\in\SS_i}(-1)^{\sigma}
\sigma\circ\wt{\tau}''_{n-j}\circ\phi\circ\wt{\tau}_j'\circ\sigma^{-1}
\end{equation}

\section{Factorisation of morphisms in \ul{$\mathsf{LBA}$}}\label{s:morlba}

In \ref{ss:morlba}--\ref{ss:simplex-T-S}, we review the polarised structure 
of morphisms in the $\PROP$ $\LBA$, and their relation to free Lie algebras 
obtained in \cite{e1,pol}. We include proofs for the reader's convenience, and
because they readily carry over to the refinements of $\LBA$ introduced in
Sections \ref{s:plba} and \ref{s:sgp-ext}.

\subsection{Factorisation of morphisms in \ul{$\mathsf{LBA}$}}\label{ss:morlba}


The inclusions $\LCA,\LA\subset\LBA$ induce maps
\[
i^N_{p,q}:\LCA([p],[N])\ten\LA([N], [q])\to\LBA([p],[N])\ten\LBA([N],[q])\to\LBA([p],[q])
\]
given by the composition of morphisms in $\LBA$. 

\begin{proposition}\label{pr:LBA factorisation}
The maps $\{i^N_{p,q}\}_{N\geqslant0}$ induce an isomorphism
\begin{equation*}
\LBA([p],[q])\simeq\bigoplus_{N\geqslant0}
\LCA([p],[N])\ten_{\SS_N}\LA([N],[q])
\end{equation*}
\end{proposition}
\begin{pf}
Morphisms in $\LBA$ can be represented as linear combinations of oriented graphs
with no loops or multiple edges, obtained by (horizontal) composition
\begin{align*}
\LBA([p],[q])\ten\LBA([q],[s])\stackrel{\circ^{op}}{\to}\LBA([p],[s])
\end{align*}
or tensor product (vertical composition)
\begin{align*}
\LBA([p],[q])\ten\LBA([p'],[q'])\stackrel{\ten}{\to}\LBA([p+p'],[q+q'])
\end{align*}
The cocycle condition \eqref{eq:cocycon}
allows to reorder every morphism as a linear combination of diagrams where the
cobrakets horizontally precede the brackets. Finally, all permutations can be moved after the
cobrackets and before the brackets, and identified with elements in $\SS_N$.

The decomposition in terms of the morphisms in the $\PROP$ $\LA$ and $\LCA$ 
follows, and the tensor product in the proposition should be interpreted as horizontal
composition of graphs. The natural map to $\LBA$ factors through the simultaneous action 
of $\SS_N$, and provides a surjective map.

The injectivity follows by the evaluation of the morphism in $\LBA$  on the 
Lie bialgebra $F(\c)=T\c$, obtained from a Lie coalgebra $(\c,\delta)$ with the free Lie 
algebra structure. 
\end{pf}

\subsection{Morphisms in $\LA$ and $\LCA$}\label{ss:morla}

Let $\FL{N}$ be the free Lie algebra over $\sfk$ with
generators $x_1,\dots, x_N$. The relation between $\FL{N}$
and morphisms in the $\PROP$s $\LA,\LCA$ is easily explained
by considering the following description of $\FL{N}$ in terms
of binary trees (see, \eg \cite[\S 0.2]{reut}). 

Let $\mathsf{T}(N)$ denote the set of binary trees over $X$,
recursively defined as follows: $x_1,\dots, x_N\in\mathsf{T}(N)$
and, for any $t_1,t_2\in\mathsf{T}(N)$, $(t_1,t_2)\in\mathsf{T}(N)$. 
Let $\T_N$ denote the $\sfk$--vector space with basis $\mathsf{T}(N)$.
The composition law $(\cdot,\cdot)$ extends to a bilinear mapping
$(\cdot,\cdot):\T_N\ten\T_N\to\T_N$. Let $J\subset\T_N$ be the
ideal generated by all elements of the form $(t,t)$, $t\in\T_N$,
and $(t_1,(t_2,t_3))+(t_2,(t_3,t_1))+(t_3,(t_1,t_2))$, $t_1,t_2,
t_3\in\T_N$, and set $\FL{N}=\T_N/J$. It is easy to see that
$\FL{N}$ is the free Lie algebra over $X$. We consider on
$\FL{N}$ the natural $\IN^N$--grading given by $\deg(x_i)=
e_i$.

\begin{lemma}\label{le:morla}
There are natural isomorphisms, compatible with the actions of 
$\SS_N$ and $\SS_n$
\[
\LA([N],[n])\simeq\homo{\FL{N}^{\ten n}}\simeq\LCA([n],[N])
\]
where $\delta_N=e_1+\cdots+e_N$, and $\homo{\FL{N}^{\ten n}}$
is the subspace of its $n$--fold tensor product spanned by homogeneous
elements of degree one in each variable.
\end{lemma}
\begin{pf}
The identification with $\LA([N],[n])$ is straightforward. We first observe that
\begin{equation}\label{eq:la-operad}
\LA([N],[n])\simeq
\bigoplus_{(I_1,\dots, I_n)\in \mathsf{P(N,n)}}\LA([|I_1|],[1])\ten\cdots\ten\LA([|I_n|],[1])
\end{equation}
where $\mathsf{P(N,n)}$ is the set of partitions of $\{1,\dots, N\}$ by $n$
unordered sets. 

Assume now that $n=1$.
Every morphism in $\LA([N],[1])$ is represented by a linear combination of 
trees with $N$ leaves precomposed with a permutation $\sigma\in\SS_N$.
The permutation $\sigma$ determines uniquely a labeling by $\{x_1,\dots, x_N\}$,
where the $i$th leaf is labeled by $x_{\sigma^{-1}(i)}$. 
This provides a surjective map from $\homo{\FL{N}}$ to $\LA([N],[1])$.
Conversely, every morphism $f\in\LA([N],[1])$ determines an element in $\homo{\FL{N}}$
by evaluating $f$ on $x_1\otimes\cdots\otimes x_N\in\FL{N}^{\otimes N}$, 
and the two maps are inverses of each other. Combined with \eqref{eq:la-operad}, this extends to
a canonical isomorphism $\LA([N],[n])\simeq\homo{\FL{N}^{\ten n}}$.
The identification with $\LCA([n],[N])$ follows by the equivalence $\preLCA\simeq\preLA\op$.
\end{pf}

\subsection{Morphisms in $\LBA$ and free Lie algebras}\label{ss:LBApoly}

\begin{proposition}\label{pr:LBA FLN}\hfill
\begin{enumerate}
\item There is an isomorphism of $(\SS_q,\SS_p)$--bimodules
\[\LBA([p],[q])\simeq\bigoplus_{N\geqslant1}
\left(\homo{\FL{N}^{\ten p}}\ten\homo{\FL{N}^{\ten q}}\right)_{\SS_N}\]
\item Let $F\in k\SS_p$ and $G\in k\SS_q$ be idempotents, and
$F[p]=([p],F)$, $G[q]=([q],G)$ the corresponding objects in $\LBA$.
Then one has
\[
\LBA(F[p], G[q])
\simeq\bigoplus_{N\geqslant0}
\left(F(\FL{N}^{\ten p})_{\delta_N}\ten G(\FL{N}^{\ten q})_{\delta_N}\right)_{\SS_N}
\]
\end{enumerate}
\end{proposition}

\begin{pf}
(1) follows from Proposition \ref{ss:morlba} and Lemma \ref{ss:morla}.
(2) Normal ordering in $\LBA$ gives
\[
\LBA(F[p], G[q])\simeq\bigoplus_{N\geqslant0}\LCA(F[p], [N])\ten_{\SS_N}\LA([N], G[q])
\]
By \ref{ss:morla} and \ref{ss:karoubi}, $\LCA(F[p],[N])\simeq F(\FL{N}^{\ten p})_{\delta_N}$
and $\LA([N], G[q])\simeq G(\FL{N}^{\ten q})_{\delta_N}$.
\end{pf}

\subsection{The tensor and symmetric algebras in $\LBA$}\label{ss:simplex-T-S}

The objects $T[1], S[1]$ in $\LBA$ play an important role in understanding
the structure of the universal algebras which will be introduced in Section
\ref{s:univDY},

Let $A=\bigoplus_{p\geqslant 0}A^p\in\LBA$ be either $T[1]$ or $S[1]$.
If follows from \ref{ss:coho} that the
tower $\LBA(\wh{A}^{\otimes n},A^{\otimes n})$ has a cosimplicial
structure, and the map
\begin{equation}\label{eq:sym LBA}
\Sym:
\LBASfin
\to\LBATfin
\end{equation}
obtained by combining the natural projection $\wh{T}[1]\to\wh{S}[1]$ and
injection $S[1]\to T[1]$ is a morphism of cosimplicial spaces.

The following result relates this structure to the standard
cosimplicial structure on the tensor and symmetric algebras
of the free Lie algebras $\FL{N}$ via the identifications provided
by Proposition \ref{pr:LBA FLN}.

\begin{lemma}\label{le:cosi and Lie}
Let $\sfsym:S\FL{N}\to T\FL{N}$ be the symmetrisation map.
The following is a commutative diagram of cosimplicial spaces
\[\xymatrix{
\LBA(\whT{[1]}^{\ten n},T[1]^{\ten n})
\ar[r]
&
\bigoplus_{N\geqslant0}
\left((T\FL{N}^{\ten n})_{\delta_N}\ten (T\FL{N}^{\ten n})_{\delta_N}\right)_{\SS_N} 
\\
\LBA(\whS{[1]}^{\ten n},S[1]^{\ten n})\ar[u]^{\sfsym} 
\ar[r]
&
\bigoplus_{N\geqslant0}
\left((S\FL{N}^{\ten n})_{\delta_N}\ten (S\FL{N}^{\ten n})_{\delta_N}\right)_{\SS_N}
\ar[u]_{\sfsym\ten\sfsym}
}\]
where the horizontal maps are those defined in Proposition \ref{ss:LBApoly}.
\end{lemma}

\section{Universal Drinfeld--Yetter modules}\label{s:univDY}

We introduce in this section the $\PROP$ $\DY{n}$ describing $n$ 
\DYt modules $\VDY{1},\ldots,\VDY{n}$ over a Lie bialgebra.
The algebra \[\DYUA{n}=\End_{\DY{n}}(\VDY{1}\otimes\cdots
\otimes\VDY{n})\] is universal in that, for any Lie bialgebra $\b$
with Drinfeld double $\gb$, it is endowed with a canonical morphism
$\DYUA{n}\longrightarrow\wh{\Ug_\b^{\otimes n}}$ to the completion
of the $n$--fold tensor product of the enveloping algebra of $\gb$
considered in \ref{ss:completions}. We show that the tower $\{\DYUA
{n}\}_{n\geqslant   0}$ shares many properties of $\{\Ug_\b^{\otimes
n}\}_{n\geqslant   0}$, namely that it has a cosimplicial structure,
satisfies the PBW theorem, and that its Hochschild cohomology is
given by a universal version of the exterior algebra of $\gb$.

\subsection{Colored $\PROP$s}

A \emph{colored} $\PROP$ $\mathsf{P}$ is a $\sfk$--linear, strict,
symmetric monoidal category whose objects are finite sequences
over a set $\sfA$, \ie
\[\mathsf{Obj}(\mathsf{P})=\coprod_{n\geq0}\sfA^n\]
with tensor product given by concatenation of sequences, and
tensor unit given by the empty sequence. 
Modules over a colored $\PROP$ are defined as in \ref{ss:prop-intro}.

\subsection{The $\PROP$ $\DY{n}$ and the algebra $\DYUA{n}$}\label{ss:univdy}

\begin{definition}
Let $n\geq 1$.\hfill
\begin{enumerate}
\item $\DY{n}$ is the colored $\PROP$ generated by $n+1$ objects
$\ADY{1}$ and $\{\VDY{k}\}_{k=1}^n$, and morphisms
\begin{gather*}
\mu:\ADY{2}\to\ADY{1}\qquad\delta:\ADY{1}\to\ADY{2}\\
\pi_k:\ADY{1}\ten \VDY{k}\to \VDY{k}\qquad
\pi_k^*:\VDY{k}\to\ADY{1}\ten \VDY{k}
\end{gather*}
such that $(\ADY{1},\mu,\delta)$ is a Lie bialgebra in $\DY{n}$, and
every $(\VDY{k}, \pi_k,\pi_k^*)$ is a Drinfeld--Yetter module over
$\ADY{1}$.
\item $\DYUA{n}$ is the algebra given by
\[\DYUA{n}
=\pEnd{\DY{n}}{\VDY{1}\otimes\cdots\otimes\VDY{n}}
\]
\end{enumerate}
\end{definition}

If $\N$ is a $k$--linear symmetric monoidal category, the category of
$\DY{n}$--modules in $\N$ is isomorphic to the category whose objects
are tuples $(\b;V_1,\ldots,V_n)$ consisting of a Lie bialgebra $\b$ in
$\N$, and $n$ \DYt modules $V_1,\ldots,V_n\in\N$ over $\b$. A morphism
$(\b;V_1,\ldots,V_n)\mapsto(\c;W_1,\ldots,W_n)$ is a tuple
$(\phi;f_1,\ldots,f_n)$, where $\phi:\b\to\c$ is a morphism of Lie bialgebras,
and $f_i:V_i\to W_i$ are such that the following diagrams are commutative
\[
\xymatrix{
\b\otimes V_i\ar[r]^{\pi_{V_i}}\ar[d]_{\phi\otimes f_i}	&V_i\ar[d]^{f_i}
&&V_i\ar[r]^{\pi^*_{V_i}}\ar[d]_{f_i}	&\b\otimes V_i\ar[d]^{\phi\otimes f_i}\\
\c\otimes W_i\ar[r]_{\pi_{W_i}}					&W_i
&&W_i\ar[r]_{\pi^*_{W_i}}						&\c\otimes W_i\\
}
\]
so that $f_i$ is a morphism of $\b$--modules $V_i\to\phi^*W_i$ as well
as a morphism of $\c$--comodules $\phi_*V_i\to W_i$.

\subsection{Action of $\DYUA{n}$ on \DYt modules}\label{ss:DY action}

Let $(\b,[\cdot,\cdot],\delta)\in{\kvect}$ be a Lie bialgebra, $\{V_k,\pi_k,
\pi_k^*\}_{k=1}^n$ $n$ \DYt modules over $\b$, and
\[\G_{(\b,V_1,\dots, V_n)}:\DY{n}\longrightarrow{\kvect}\]
the corresponding (symmetric, monoidal) realisation functor such that
$\ADY{1}\mapsto\b$ and $\VDY{k}\mapsto V_k$.

\begin{proposition}\label{pr:universal action}
Let $\ff\hspace{-0.1cm}:\DrY{\b}\to{\kvect}$ be the forgetful functor,
and $\DYA{\b}{n}\coloneqq \sfEnd{\ff^{\boxtimes n}}$. Then, there is an algebra
homomorphism 
\[\DYrho{\b}{n}:\DYUA{n}\to\DYA{\b}{n}\]
which assigns to any $T\in\DYUA{n}$, and any $V_1,\ldots,V_n\in
\DrY{\b}$ the endomorphism $\G_{(\b,V_1,\dots, V_n)}(T)\in\End_
\sfk(V_1\otimes\cdots\otimes V_n)$.
\end{proposition}
\begin{pf} 
We need to prove that, for any $\{V_i,W_i\}_{i=1}^n\subset\DrY{\b}$
and $f_i\in\Hom_{\DrY{\b}}(V_i,W_i)$, one has 
\[ f_1\otimes\cdots\otimes f_n\circ \G_{(\b,V_1,\dots, V_n)}(T)=
\G_{(\b,W_1,\dots, W_n)}(T)\circ f_1\otimes\cdots\otimes f_n\]
This follows from the fact that $(\id_\b;f_1,\ldots,f_n)$ induces a natural
transformation $\G_{(\b,V_1,\dots, V_n)}\Rightarrow\G_{(\b,W_1,\dots,
W_n)}$.
\end{pf}

\subsection{Distinguished elements in $\DYUA{}$}

The algebra $\DYUA{2}$ has a distinguished element, which is given by
\[r_{\VDY{1},\VDY{2}}=
\pi_{\VDY{1}}\ten\id_{\VDY{2}}\circ\;(12)\circ\id_{\VDY{1}}\ten\pi^*_{\VDY{2}}\]
and is easily seen to be a solution of the classical Yang--Baxter equation
in $\DYUA{3}$
\[[r_{\VDY{1},\VDY{2}},r_{\VDY{1},\VDY{3}}]+
[r_{\VDY{1},\VDY{2}},r_{\VDY{2},\VDY{3}}]+[r_{\VDY{1},\VDY{3}},r_{\VDY{1},\VDY{3}}]=0\]
Under the homomorphism $\DYrho{\b}{2}:\DYUA{2}\to\DYA{\b}{2}$, $r_{\VDY{1},
\VDY{2}}$ corresponds to the action of the $r$--matrix $r_\b=\sum_i b_i\otimes
b^i$ of $\gb$ defined in \eqref{eq:rmatrix}.

The algebra $\DYUA{}$ contains the element $\kappa=\pi_{\VDY{1}}\circ\pi^*_{\VDY{1}}$,
which corresponds to the normally ordered Casimir operator $\kappa_\b
=\sum_i b_ib^i=m(r_\b)$ of $\gb$. We note further that while one can
consider the following elements in $\DYUA{2}$
\begin{gather*}
r_{21}=\id_{\VDY{1}}\otimes\pi_{\VDY{2}}\circ\;(12)\circ \pi^*_{\VDY{1}}\otimes\id_{[V_2]}\\
\kappa_1=\left(\pi_{\VDY{1}}\circ\pi_{\VDY{1}}^*\right)\otimes\id_{[V_2]}
\aand
\kappa_2=\id_{\VDY{1}}\otimes\left(\pi_{[V_2]}\circ\pi_{[V_2]}^*\right)
\end{gather*}
which correspond to $r_\b^{21}$, $\kappa_\b\otimes 1$ and
$1\otimes\kappa_\b$ in $\DYA{\b}{2}$ respectively, there is
no analogue in $\DYUA{}$ of the non--normally ordered Casimir
operator $\sum_i b^ib_i=m(r_\b^{21})$, which does not converge
in $\DYA{\b}{2}$ if $\dim\b=+\infty$.

\subsection{Universal invariants}\label{ss:LBAinv}

\begin{definition}
An element $\phi\in\DYUA{n}$ is {\it invariant} if it commutes with the action
and coaction of the Lie bialgebra $[1]$ on $\VDY{1}\otimes\cdots\otimes\VDY{n}$,
that is satisfies 
\[\pi_{\VDY{1}\otimes\cdots\otimes\VDY{n}}\circ\id_{[1]}\otimes\phi
=
\phi\circ\pi_{\VDY{1}\otimes\cdots\otimes\VDY{n}}\]\\
as maps $[1]\otimes\VDY{1}\otimes\cdots\otimes\VDY{n}\to\VDY{1}\otimes\cdots\otimes\VDY{n}$,
and 
\[\pi^*_{\VDY{1}\otimes\cdots\otimes\VDY{n}}\circ\phi
=
\id_{[1]}\otimes\phi\circ \pi^*_{\VDY{1}\otimes\cdots\otimes\VDY{n}}\]\\
as maps $\VDY{1}\otimes\cdots\otimes\VDY{n}\to[1]\otimes\VDY{1}\otimes\cdots\otimes\VDY{n}$.
\end{definition}

\noindent 
Let $\DYUAinv{n}\subset\DYUA{n}$ be the subalgebra of
invariant elements. The following is clear.

\begin{proposition}
The map $\DYrho{\b}{n}:\DYUA{n}\to\DYA{\b}{n}$
defined in \ref{ss:DY action} restricts to a homomorphism
\[\DYrho{\b}{n}:\DYUAinv{n}\to\DYAinv{\b}{n}\coloneqq\sfEnd{\sf id^{\boxtimes n}}\] 
\end{proposition}

\subsection{Cosimplicial structure of $\DYA{\b}{}$}\label{ss:cosim}

The monoidal structure on $\DrY{\b}$ endows the tower
$\{\DYA{\b}{n}\}$ with the structure of a cosimplicial complex
of algebras
\[\xymatrix@C=.5cm{\sfk\ar@<-2pt>[r]\ar@<2pt>[r] & 
\DYA{\b}{} \ar@<3pt>[r] \ar@<0pt>[r] \ar@<-3pt>[r] & 
\DYA{\b}{2} \ar@<-6pt>[r]\ar@<-2pt>[r]\ar@<2pt>[r]\ar@<6pt>[r] & \DYA{\b}{3}
\quad\cdots}\]
The corresponding face morphisms $\{d_i^n\}_{i=0}^{n+1}:\DYA{\b}{n}\to\DYA{\b}{n+1}$,
are given by $(d_0^0 \varphi)_V=(d_1^0 \varphi)_V=\varphi\cdot\id_V$,
for $\varphi\in\sfk$ and $V\in\DrY{\b}$, and, for $n\geqslant   1$, $\varphi\in\DYA{\b}{n}$, 
and $\{V_i\}_{i=1}^{n+1}\subset\in\DrY{\b}$
\[
(d_i^n \varphi)_{V_1,\dots, V_{n+1}}=
\left\{
\begin{array}{cc}
\id_{V_1}\ten \varphi_{V_2,\dots, V_{n+1}} & i=0\\[1.1ex]
\varphi_{V_1,\dots, V_i\ten V_{i+1},\dots,V_{n+1}} & 1\leqslant i \leqslant n\\[1.1ex]
\varphi_{V_1,\dots, V_n}\ten\id_{V_{n+1}} & i=n+1
\end{array}
\right.
\]
The degeneration homomorphisms $\varepsilon_n^i:\DYA{\b}{n}\to
\DYA{\b}{n-1}$, for $i=1,\dots, n$, are
\[(\varepsilon_n^i \varphi)_{X_1, \dots, X_{n-1}}=\varphi_{X_1,\dots, X_{i-1}, {\bf 1}, X_{i}, \dots, X_{n-1} }\]
The morphisms $\{\varepsilon_n^i\}, \{d_i^n\}$ satisfy the standard relations 
\[
\begin{array}{lr}
d_{n+1}^jd^i_n=d^i_{n+1}d_n^{j-1} & i<j \\[1.1ex]
\varepsilon_n^j\varepsilon_{n+1}^i=\varepsilon_n^i\varepsilon_{n+1}^{j+1} & i\leqslant j
\end{array}
\]
and
\[
\varepsilon_{n+1}^jd_i^n=\left\{
\begin{array}{cl}
d_{n-1}^i\varepsilon_n^{j-1} & i<j\\[1.1ex]
\id & i=j,j+1\\[1.1ex]
d_{n-1}^{i-1}\varepsilon_n^j & i>j+1
\end{array}
\right.
\]
These give rise to the Hochschild differential 
\[
d^n=\sum_{i=0}^{n+1}(-1)^{i} d_i^n: \DYA{\b}{n}\to\DYA{\b}{n+1}
\]

\subsection{Cosimplicial structure of $\DYUA{n}$}\label{ss:cosimp-DY}

The above structures can be lifted to the $\PROP$s $\DY{n}$. For
every $n\geqslant1$ and $i=0,1,\dots, n+1$, there are faithful functors
\[
\D_{i}^n:\DY{n}\to\DY{n+1}
\]
mapping $[1]$ to $[1]$, and given by
\[\D^n_0(\VDY{k})=\VDY{k+1}
\aand
\D^n_{n+1}(\VDY{k})=\VDY{k}\]
for $1\leqslant k\leqslant n$, and, for $1\leqslant i\leqslant n$,
\[\D^n_i(\VDY{k})=
\left\{\begin{array}{cc}
\VDY{k} 				& 1\leqslant k\leqslant i-1	\\[1.1ex]
\VDY{i}\otimes \VDY{i+1} 	& k=i					\\[1.1ex]
\VDY{k+1}				& i+1\leqslant k\leqslant n	
\end{array}\right.\]
{
and $\E_n^{i}:\DY{n}\to\DY{n-1}$
\[
\E_n^{i}=\G_{([1], \VDY{1},\dots,\VDY{i-1},\1,\VDY{i+1},\dots, \VDY{n-1})}
\]
where $\1$ is the trivial representation in $\DY{n}$.
}
These induce algebra homomorphisms
\[
\Delta_i^n:\DYUA{n}\to\DYUA{n+1}
\]
which are universal analogues of the insertion/coproduct maps on $U\gb^{\otimes n}$.
They give the tower of algebras $\{\DYUA{n}\}_{n\geqslant0}$ the structure of a cosimplicial
complex, with Hochschild differential $d^n=\sum_{i=0}^{n+1}(-1)^i\Delta_i^n:\DYUA
{n}\to\DYUA{n+1}$. The morphisms $\DYrho{\b}{n}:\DYUA{n}\to\DYA{\b}{n}$ defined
in \ref{ss:univdy} are compatible with the face and degeneration morphisms, and
therefore with the differentials $d^n$.

\subsection{The algebra $T[1]$ and the coalgebra $\whT{[1]}$}

Regard $T[1]=\bigoplus_{p\geqslant   0}[p]$ as a graded algebra, with
concatenation product given by the identification $[p_1]\otimes[p_2]=
[p_1+p_2]$ and unit given by the embedding $\imath:[0]\hookrightarrow
T[1]$. If $([V],\pi_{[V]})$ is a module over the Lie algebra $[1]$ in $\DY
{1}$, the iterated action maps
\[\pi_{[V]}^{(p)}:[p]\otimes[V]\to[V]\]
endow $[V]$ with the structure of a module over $T[1]$.

Dually, regard $\whT{[1]}=\prod_{p\geqslant   0}[p]$ as a topological
graded coalgebra, with deconcatenation coproduct given by the direct
sum of identifications
\[[p]\to\bigoplus_{p_1+p_2=p}[p_1]\otimes [p_2]=\bigoplus_{p_1+p_2=p}[p]\]
and counit given by the projection $\epsilon:\whT{[1]}\to[0]$. Then, if
$([V],\pi^*_{[V]})$ is a comodule over the Lie coalgebra $[1]$ in $\DY
{1}$, the iterated coaction maps
\[{\pi^*}^{(p)}:[V]\to[p]\otimes[V]\]
endow $[V]$ with the structure of a comodule over $\whT{[1]}$.

\subsection{Action of morphisms in $\LBA$ on $\DYUA{n}$}\label{ss:LBA on DY}

Consider now the vector space
\[\LBA(\whT{[1]},T[1])=
\bigoplus_{p,q\geqslant   0}\LBA([p],[q])\]
with the convolution product $\phi_1\star\phi_2=m_{T[1]}\circ \phi_1\otimes
\phi_2\circ\Delta_{\whT{[1]}}$, and unit $1_{T[1]}\circ\epsilon_{\whT{[1]}}$. 
Then, regarding $\VDY{1}\in\DY{1}$ as a module over $T[1]$ and a
comodule over $\whT{[1]}$ yields a convolution action of $\LBA
(\whT{[1]},T[1])$ on $\DY{1}(\VDY{1},\VDY{1})$ given by
\[\phi\cdot \,X=
\pi_{T[1]}\circ \phi\otimes X\circ\pi_{\whT{[1]}}^*\]
In particular, specialising to $X=\id_{\VDY{1}}$ yields a map
\[\sfa^1:\LBA(\whT{[1]},T[1])\longrightarrow\DYUA{1}=\DY{1}(\VDY{1},\VDY{1})\]
mapping $\phi$ to $\phi\cdot\id_{\VDY{1}}$.

More generally, for any $n\geqslant   1$, the algebra structure on $T[1]^
{\otimes n}$ and the coalgebra structure on $\whT{[1]}^{\otimes n}$
yield a map
\[
\sfa^n:
\LBA(\whT{[1]}^{\otimes n},T[1]^{\otimes n})
\longrightarrow
\DYUA{n}=\DY{n}(\ten_{k=1}^n\VDY{k},\ten_{k=1}^n\VDY{k}),
\]
which maps $\phi$ to $\phi\cdot\id_{\ten_{k=1}^n\VDY{k}}$.

Recall from \ref{ss:simplex-T-S} that the tower $\LBA(\whT{[1]}
^{\ten n},T[1]^{\ten n})$ is cosimplicial.

\begin{proposition}\label{pr:a cosimplicial}
The collection of maps $\{\sfa^n\}$ is a morphism of cosimplicial spaces.
\end{proposition}
\begin{pf}
It suffices to prove that
	\[\sfa^n:\LBA(\whT{[1]}^{\ten n},T[1]^{\ten n})\to\DY{n}(\ten_{k=1}^n
	\VDY{k},\ten_{k=1}^n\VDY{k})\]
	is compatible with the face maps $\{d_i^n\}_{i=0}^{n+1}$. The case
	$i=0,n+1$ is easily checked. To check the compatibility with $d^n_1$,
	it suffices to consider the case $n=1$. Let $\phi\in\LBA([p],[q])$. We
	need to check the equality of
	\[d_1^1\circ\sfa^1(\phi)=
	\pi_{\VDY{1}\otimes\VDY{2}}^{(q)}\circ\,\phi\otimes\id_{\VDY{1}\otimes\VDY{2}}\circ\,{{\pi^*}^{(p)}_{\VDY{1}\otimes\VDY{2}}}\]
	and
	\[\sfa^2\circ d_1^1(\phi)=
	\sum_{\substack{\ulp,\ulq\in\IN^2\\|\ulp|=p,|\ulq|=q}}
	\pi^{(\ulq)}_{\VDY{1}\otimes\VDY{2}}\circ\,d_1^1(\phi)_{\ulp,\ulq}\otimes\id_{\VDY{1}\otimes\VDY{2}}\circ\,{\pi^*}^{(\ulp)}_{\VDY{1}\otimes\VDY{2}}\]
	where $d_1^1(\phi)_{\ulp,\ulq}\in\LBA(T^\ulp[1],T^\ulq[1])$ are the
	homogeneous components of $d_1^1(\phi)=\Delta\circ\phi\circ m$,
	and $m,\Delta$ are the multiplication and comultiplication of $T[1]$.
	
	The equality now follows from the identities
	\begin{align*}
	\bigoplus_{\ulp:|\ulp|=p}
	m\otimes\id_{\VDY{1}\otimes\VDY{2}}\circ\,{\pi^*}^{(\ulp)}_{\VDY{1}\otimes\VDY{2}}
	&=
	{{\pi^*}^{(p)}_{\VDY{1}\otimes\VDY{2}}}\\
	\bigoplus_{\ulq:|\ulq|=q}
	\pi^{(\ulq)}_{\VDY{1}\otimes\VDY{2}}\circ\,\Delta\otimes\id_{\VDY{1}\otimes\VDY{2}}
	&=
	\pi_{\VDY{1}\otimes\VDY{2}}^{(q)}
	\end{align*}
	of maps $\VDY{1}\otimes\VDY{2}\to[p]\otimes\VDY{1}\otimes\VDY{2}$
	and $[q]\otimes\VDY{1}\otimes\VDY{2}\to\VDY{1}\otimes\VDY{2}$
	respectively. The first (resp. second) one holds because both sides are
	the components of the coaction (resp. action) of $T[1]$ on $\VDY{1}
	\otimes\VDY{2}$.
\end{pf}

\subsection{A basis of $\DYUA{1}$}\label{ss:mordy}

In the following paragraphs, we describe the vector space underlying
$\DYUA{n}$, and the convolution action of $\LBA(\whT{[1]}^{\otimes n},
T[1]^{\otimes n})$ on it in terms of free algebras and Lie algebras, in
analogy with \ref{pr:LBA FLN}--\ref{le:cosi and Lie}. We then use this
description to prove an analogue of the PBW theorem for $\DYUA{n}$
in \ref{ss:LBAPBW}.

Let $\pi^{(N)}:[N]\ten\VDY{1}\to\VDY{1}$ (resp. ${\pi^*}^{(N)}:\VDY{1}
\to[N]\ten\VDY{1}$) be the $N$th iterated action (resp. coaction) on
$\VDY{1}$.

\begin{proposition}
The endomorphisms of $\VDY{1}\in\DY{1}$ given by
\[
\sarch{N}{\sigma}=\pi^{(N)}\circ\sigma\ten\id_{[V_1]}\circ\;{\pi^*}^{(N)}=\sigma\cdot\id_{[V_1]}
\]
for $N\geqslant 0$ and $\sigma\in\SS_N$, form a basis of $\DYUA{}=\pEnd{\DY{1}}
{\VDY{1}}$.
\end{proposition}
\begin{pf}
We represent $\id_{[1]}$ with a line and $\id_{\VDY{1}}$ with a bold line.
The morphisms $\mu,\delta,\pi,\pi^*$ in $\DY{1}$ are then represented by the diagrams
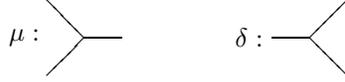
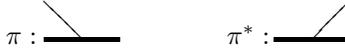
\[
\xy
(0,0)*{
\xy
(-3,0)*{\mu:};
(5,0);(10,0)**\dir{-}?(.5)*\dir{};
(0,5);(5,0)**\dir{-}?(.5)*\dir{};
(0,-5);(5,0)**\dir{-}?(.5)*\dir{};
\endxy
};
(30,0)*{
\xy
(-3,0)*{\delta:};
(0,0);(5,0)**\dir{-}?(.5)*\dir{};
(5,0);(10,5)**\dir{-}?(.5)*\dir{};
(5,0);(10,-5)**\dir{-}?(.5)*\dir{};
\endxy
}
\endxy
\]
and
\[
\xy
(30,5)*{
\xy
(-3.5,.5)*{\pi^*:};
(0,0);(10,0)**\dir{-};
(0,-0.1);(10,-0.1)**\dir{-};
(0,-0.2);(10,-0.2)**\dir{-};
(0,-0.3);(10,-0.3)**\dir{-};
(5,0);(10,5)**\dir{-}?(.5)*\dir{};
\endxy
};
(0,5)*{
\xy
(-3,0)*{\pi:};
(0,0);(10,0)**\dir{-};
(0,-0.1);(10,-0.1)**\dir{-};
(0,-0.2);(10,-0.2)**\dir{-};
(0,-0.3);(10,-0.3)**\dir{-};
(0,5);(5,0)**\dir{-}?(.5)*\dir{};
\endxy
}
\endxy
\]
which are read from left to right. A non--trivial endomorphism of $\VDY{1}$
is represented as a linear combination of oriented diagrams, necessarily
starting with a coaction and ending with an action. The compatibility relation
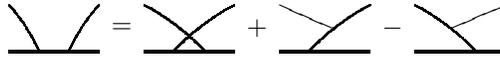
\eqref{eq:actcoact} 
\begin{align*}
&\\
&\actcoact\\
\end{align*}
allows to reorder $\pi$ and $\pi^*$, moving every coaction before any action. 
The cocycle condition \eqref{eq:cocycon} allows to reorder brackets and cobrackets 
as in $\LBA$. Finally, the relations \eqref{eq:action}, \eqref{eq:coaction}\\
\begin{align*} 
&\module\\
&\comodule\\
\end{align*}
allow to remove from the graph every $\mu$ and every $\delta$ involved.
It follows that every endomorphism of $\VDY{1}$ is a linear combination
of the elements $\sarch{N}{\sigma}$ given by
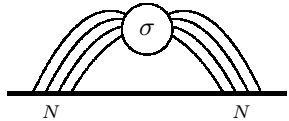
\[
\arch{\sigma}{N}{20}
\]
where $N\geqslant 0$ and $\sigma\in\SS_N$. These morphisms are linearly
independent in $\DY{1}$, since they are on the free Drinfeld--Yetter module
constructed over the comodule $\VDY{1}$, following an argument similar to
\ref{ss:morlba}. 
\end{pf}

\noindent
\begin{remark}
Under the map $\DYrho{\b}{}:\DYUA{}\to\DYA{\b}{}$, the basis element
$\sarch{N}{\sigma}$ maps to the interlaced $N$th power of the normally
ordered Casimir operator of $\gb$ given by
\[\kappa_N^\sigma=
\sum_{i_1,\ldots,i_N}
b_{i_{\sigma(1)}}b_{i_{\sigma(2)}}\cdots b_{i_{\sigma(N)}}\cdot 
b^{i_N}\cdots b^{i_2}b^{i_1}\]
\end{remark}

\noindent
\begin{remark}
Proposition \ref{ss:mordy} yields in particular an isomorphism of vector spaces
\begin{equation}\label{eq:iso-DY-Sym}
\pEnd{\DY{1}}{\VDY{1}}\simeq\bigoplus_{N\geqslant0}\sfk\SS_N
\end{equation}
mapping $\sarch{N}{\sigma}$ to $\sigma\in\SS_N$. It is clear from the
description above that the multiplication in $\pEnd{\DY{1}}{\VDY{1}}$ is
$\IN$--graded, as the normal ordering on the product of two elements of
the basis preserves the total number of strings. Namely, for any $N,M>0$,
$\sigma\in\SS_N$, $\tau\in\SS_M$, one gets
\[
\xy
(-28,0)*{\arch{\sigma}{N}{15}};
(0,0)*{\arch{\tau}{M}{15}};
(30,0)*{={\sum_{\substack{{\rho\in\SS_{N+M}}}}}c_{\sigma,\tau}^{\rho}};
(60,0)*{\arch{\rho}{{N+M}}{15}};
\endxy
\]
for some $c_{\sigma,\tau}^{\rho}\in\IZ$. It seems an interesting problem
to determine the structure constants $c_{\sigma,\tau}^{\rho}$ explicitly.
\end{remark}

\subsection{Convolution product on $\DYUA{}$}

Under the isomorphism \eqref{eq:iso-DY-Sym}, the exterior product
of permutations $\ten:\SS_N\times\SS_M\to\SS_{N+M}$ gives rise
to a convolution product $\star$ on $\DYUA{1}$, defined on the basis
elements by
\[
\sarch{N}{\sigma}\star\sarch{M}{\tau}=\sarch{N+M}{\sigma\ten\tau}
\]

Pictorially, the product $\star$ corresponds to the encapsulation of 
$\sarch{M}{\tau}$ inside $\sarch{N}{\sigma}$. In particular, the action
of $\LBA(\whT{[1]}^{\otimes n},T[1]^{\otimes n})$ commutes with
convolution in $\DYUA{}$ on the right, \ie 
\[\phi\cdot (X\star Y)=(\phi\cdot X)\star Y\]
for any $\phi\in\LBA(\whT{[1]}^{\otimes n},T[1]^{\otimes n})$ and $X,Y\in\DYUA{}$. 
It follows that the action is given by left convolution with $\sfa^1(\phi)
=\phi\cdot\id_{[V]}$, and that 
$\sfa^1$ is a morphism \wrt convolution.

In terms of the interlaced Casimir operators of $\DYA{\b}{}$, one has
\begin{multline*}
\rho_{\b}(\sarch{N}{\sigma}\star\sarch{M}{\tau})=\rho_{\b}(\sarch{N+M}{\sigma\ten\tau})\\[1.1ex]
=\sum_{\substack{i_1,\ldots,i_N\\j_1,\dots, j_M}}
\left(b_{i_{\sigma(1)}}\cdots b_{i_{\sigma(N)}}\right)\cdot 
\left(b_{j_{\tau(1)}}\cdots b_{j_{\tau(M)}}\right)\cdot
\left(b^{j_M}\cdots b^{j_1}\right)\cdot
\left(b^{i_N}\cdots b^{i_1}\right)
\end{multline*}
The product $\star$ can therefore be thought of as a $\PROP$ic analogue
of the polarised multiplication on $U\gb\simeq U\b\ten U\b^*$ given by $(x
\ten y)\star(x'\ten y')=(x\cdot x')\ten (y'\cdot y)$. Thus, $\DYUA{}$ is endowed
with two distinct products: the canonical one coming from its definition as
$\End([V_1])$, which corresponds to the usual product on $U\gb$, and
the convolution product $\star$, which is defined in terms of the basis
$\sarch{N} {\sigma}$ and corresponds to the polarised product on $U\gb$.

\subsection{A basis of $\DYUA{n}$, $n>1$}\label{ss:univdyn}

The description of the morphisms in $\DY{n}$ is similar to the case $n=1$.
For any $N\in\IN$ and $\ul{N}=(N_1,\dots, N_n)\in\IN^n$ such that $|\ul{N}|
=N$, let
\begin{equation*}
\pi^{(\ul{N})}: [N]\ten\bigotimes_{k=1}^n \VDY{k}\to\bigotimes_{k=1}^n \VDY{k}
\aand
{\pi^*}^{(\ul{N})}: \bigotimes_{k=1}^n \VDY{k}\to[N]\ten\bigotimes_{k=1}^n \VDY{k}
\end{equation*}
be the ordered composition of $N_i$ actions (resp. coactions) on $\VDY{i}$.

\begin{proposition}\label{pr:basis of U_DY^n}
The endomorphisms of $\VDY{1}\ten\cdots\ten\VDY{n}$ given by
\[
\rarch{\ul{N}}{\sigma}{\ul{N}^\prime}=\pi^{(\ul{N})}\circ\sigma\ten\id\circ{\pi^*}^{(\ul{N}^\prime)}
\]
where $N\geqslant 0$, $\ul{N},\ul{N}^\prime\in\IN^n$ are such that $|\ul{N}|=
N=|\ul{N}^\prime|$, and $\sigma\in\SS_N$, form a basis of $\DYUA{n}
=\pEnd{\DY{n}}{\ten_{k=1}^n\VDY{k}}$.
\end{proposition}

As in the case of $n=1$, the basis $\rarch{\ul{N}}{\sigma}{\ul{N}^\prime}$
gives rise to a convolution product on $\DYUA{n}$ given by
\[\sarch{\ul{N}}{\sigma}\star\sarch{\ul{N'}}{\tau}=\sarch{\ul{N}+\ul{N'}}{\sigma\ten\tau}\]
providing a $\PROP$ic analogue of the polarised multiplication in $U\gb^{\ten n}$. 
One checks easily that, with respect to $\star$, the map $\sfa^n$ defined in \ref{ss:LBA on DY}
is a morphism of algebras.

\subsection{Cosimplicial structure and basis elements}\label{ss:cosimplicial}
The cosimplicial structure of $\DYUA{n}$ introduced in \ref{ss:cosimp-DY} is 
defined on the elements $\rarch{\ul{N}}{\sigma}{\ul{N}^\prime}$ as follows.
The degeneration map $\E_n^{i}:\DYUA{n}\to\DYUA{n-1}$ is given by
\[
\E_n^i(\rarch{\ul{N}}{\sigma}{\ul{N}^\prime})=
\begin{cases}
\rarch{\ul{N}_{\wh{i}}}{\sigma}{\ul{N}_{\wh{i}}^\prime} & \text{if } N_i=0=M_i\\
0 & \text{otherwise}
\end{cases}
\]
where $\ul{N}_{\wh{i}}$ is obtained from $\ul{N}$ by removing $N_i$.
The face map $\Delta_i^n:\DYUA{n}\to\DYUA{n+1}$ is given by
\[
\Delta_i^n(\rarch{\ul{N}}{\sigma}{\ul{N}^\prime})=
\sum_{\substack{p=0,\dots, N_i\\ q=0, \dots N_i^\prime}}
\sum_{\substack{\tau\in\mathfrak{G}(N_i,p)\cup\{\id\}\\
		\tau^\prime\in\mathfrak{G}(N_i^\prime, N_i^\prime-q)\cup\{\id\}}}
\rarch{\ul{N}_p}{(\tau^\prime)^{-1}\circ\sigma\circ\tau}{\ul{N}_q^\prime}
\] 
where $\ul{N}_p=(N_1,\dots, N_{i-1}, p, N_i-p, N_{i+1},\dots, N_n)$ and 
$\mathfrak{G}(N_i,p)\subset\SS_N$ is the set of permutations $\tau$ acting on 
$(1,\dots, N_{i-1}, N_{i-1}+N_i+1,\dots N)$ as the identity and on $(N_{i-1}+1, \dots, N_{i-1}+N_i)$  
as a Grassmannian permutations with a unique descent at $N_{i-1}+p$\footnote{Recall that a Grassmannian permutation is a permutation $\tau\in\SS_N$ with a unique descent. 
In other words there exists $k\in \{1,\ldots,N-1\}$ such that $\tau(i) < \tau(i+1)$ if $i \neq k$ and
$\tau(k) > \tau(k+1)$.}. Similarly for $\ul{N}_q$ and $\mathfrak{G}(N_i^\prime-q, N_i^\prime)$. 
Note that the appearance of the corrections $\tau, \tau^\prime$ are due to the prescribed 
order of actions and coactions on the basis elements $\rarch{\ul{N}}{\sigma}{\ul{N}^\prime}$.

Moreover, one can verify by direct inspection that the face and degenerations maps are morphisms 
of convolution algebras.

\subsection{Identification with free algebras}\label{ss:Un Fn}

Let $\FA{N}$ be the free algebra in $N$ variables, and,
for any $n\geq 1$, denote by $\homo{\FA{N}^{\ten n}}\subset\FA{N}$
the subspace spanned by elements of degree one in each variable.
The symmetric group $\SS_N$ acts diagonally on $\homo{\FA{N}^{\ten n}}\ten\homo{\FA{N}^{\ten n}}$
by simultaneous permutation of the variables. The corresponding
space of coinvariants $\left(\homo{\FA{N}^{\ten n}}\ten\homo{\FA{N}
^{\ten n}}\right)_{\SS_N}$ has the following basis. For any $\ul{N},
\ul{N}^\prime\in\IN^n$ such that $|\ul{N}|=N=|\ul{N}^\prime|$ and
$\sigma\in\SS_N$, define $x_{\ul{N}},y_{\sigma(\ul{N}')}\in\homo
{\FA{N}^{\ten n}}$ by
\begin{align*}
x_{\ul{N}}
&=
x_{1}\cdots x_{N_1} \otimes
x_{N_1+1}\cdots x_{N_1+N_2}\otimes\cdots\otimes
x_{N_1+\cdots+N_{n-1}+1}\cdots x_N\\
y_{\sigma(\ul{N}')}
&=
y_{\sigma(1)}\cdots y_{\sigma(N'_1)} \otimes
\cdots\otimes 
y_{\sigma(N'_1+\cdots+N'_{n-1}+1)}\cdots y_{\sigma(N)}
\end{align*}
Then, $\{x_{\ul{N}}\otimes y_{\sigma(\ul{N}')}\}_{\ul{N},\ul{N'},\sigma}$
is a basis of $\left(\homo{\FA{N}^{\ten n}}\ten\homo{\FA{N}^{\ten n}}
\right)_{\SS_N}$.

The following is an immediate consequence of Proposition \ref{pr:basis of U_DY^n}.

\begin{corollary}
The linear map
\[
\xi^n\univ:\DYUA{n} \to
\bigoplus_{N\geqslant0}\left(\homo{\FA{N}^{\ten n}}\ten\homo{\FA{N}^{\ten n}}\right)_{\SS_N}
\]
given by
\begin{equation}\label{eq:xi n DY}
\xi^n\univ(\rarch{\ul{N}}{\sigma}{\ul{N}'})=x_{\ul{N}}\ten y_{\wt{\sigma}(\ul{N}')}
\end{equation}
where $\wt{\sigma}=\sigma^{-1}\circ\tau$ and $\tau\in\SS_N$, such that
$\tau(i)=N-i$, is an isomorphism of vector spaces.\footnote{The involution $\tau_N$ is required because of the
contravariance of the expression \eqref{eq:contravariant} with respect to
the Lie polynomial $Q$, and to ensure the commutativity of the diagram
in Theorem \ref{ss:LBAPBW}.}
\end{corollary}

\subsection{Module structure on coinvariants}\label{ss:free-lie-pbw-conv}

The space of coinvariants 
\[\A^n\coloneqq\bigoplus_{N\geqslant0}\left((\FA{N}^{\ten n})_{\delta_N}
\ten (\FA{N}^{\ten n})_{\delta_N}\right)_{\SS_N}\] is an associative algebra,
with product map in degree $(M,N)$
\[
\left(\FA{M}^{\ten n}\ten\FA{M}^{\ten n}\right)_{\SS_M}
\ten 
\left(\FA{N}^{\ten n}\ten\FA{N}^{\ten n}\right)_{\SS_N}
\to
\left(\FA{M+N}^{\ten n}\ten\FA{M+N}^{\ten n}\right)_{\SS_{M+N}}
\]
given by the formula
\begin{equation}\label{eq:pol-conv-act}
(x_{\ul{M}}\ten y_{\wt{\sigma}(\ul{M}')})\star (x_{\ul{N}}\ten y_{\wt{\tau}(\ul{N}')})=
(x_{\ul{M}}\cdot x_{\ul{N}})\ten(y_{\wt{\tau}(\ul{N}')}\cdot y_{\wt{\sigma}(\ul{M}')})
\end{equation}
where $x_{\ul{M}}\cdot x_{\ul{N}}$ and $y_{\wt{\tau}(\ul{N}')}\cdot y_{\wt{\sigma}
(\ul{M}')}$ are identified with elements in $\FA{M+N}^{\ten n}$. Note that, under the 
identification provided by $\xi^n\univ$, \eqref{eq:pol-conv-act} reads
\[\xi^n\univ(\sarch{M}{{\sigma}})\star\xi^n\univ(\sarch{N}{{\tau}})
=\xi^n\univ(\sarch{M+N}{{\sigma}\ten{\tau}})=\xi^n\univ(\sarch{M}{{\sigma}}\star\sarch{N}{{\tau}})\]

The formula \eqref{eq:pol-conv-act} is easily adapted to define an algebra
structure on
\[\T^n\coloneqq\bigoplus_{N\geqslant0}
\left((T\FL{N}^{\ten n})_{\delta_N}\ten (T\FL{N}^{\ten n})_{\delta_N}\right)_{\SS_N}\]
In particular, the linear surjection $\sfp^n:\T^n\to\A^n$,
defined componentwise by the quotient map
$T\FL{N}\to U\FL{N}=\FA{N}$ for the free Lie algebra $\FL{N}$,
is an algebra map and induces on $\A^n$ a natural structure of $\T^n$--module.

\subsection{Identification with free Lie algebras} \label{ss:lba-pbw-conv}
Let 
\[\xi^n\lba:
\L^n\coloneqq\LBA(\whT{[1]}^{\ten n},T[1]^{\ten n})
\to
\T^n\]
be the isomorphism of vector spaces given by Proposition \ref{ss:LBApoly}. One checks 
by direct inspection that $\xi^n\lba$ is a morphism of algebras, with respect to the convolution
product on $\L^n$ and the multiplication
on $\T^n$ defined in \ref{ss:free-lie-pbw-conv}, \ie 
$\xi^n\lba(\phi\star\psi)=\xi^n\lba(\phi)\star\xi^n\lba(\psi)$
for any $\phi,\psi\in\L^n$.
Therefore, through $\xi^n\lba$, we obtain a convolution action
of $\L^n$ on $\A^n$, \ie we set
\[
\phi\cdot (x_{\ul{N}}\ten y_{\sigma(\ul{N'})})=(\sfp^n\circ\xi^n\lba)(\phi)\cdot (x_{\ul{N}}\ten y_{\sigma(\ul{N'})})
\]
for any $\phi\in\L^n$.

\begin{proposition}\label{pr:PBW DY}
The isomorphism $\xi^n\univ:\DYUA{n}\to\A^n$ given by \eqref
{eq:xi n DY} intertwines the convolution actions of $\L^n$, that
is satisfies
\[
\xi^n\univ(\phi\cdot X)=(\sfp^n\circ\xi^n\lba)(\phi)\cdot\xi^n\univ(X)
\]
for any $\phi\in\L^n$ and $X\in\DYUA{n}$.
\end{proposition}

\begin{pf}
Assume for simplicity that $n=1$. The proof for $n>1$ is identical.

Let $P_1\ten\cdots\ten P_p\in \FL{N}^{\ten p}$ be an element of degree
$\delta_N$, and
$\mu_{P_1\ten\cdots\ten P_p}\in\LA([N],[p])$
the element corresponding to it by Lemma \ref{ss:morla}. In $\FA{N}=U\FL{N}$,
the product $P_1\cdots P_p$ corresponds to an element
$\sigma_{P_1\cdots P_p}\in(\FA{N})_{\delta_N}\simeq\sfk\SS_N$, which, 
by \eqref{eq:action}, satisfies the following relation in
$\DY{1}([N]\ten\VDY{1},\VDY{1})$
\begin{equation}\label{eq:highaction}
\pi^{(p)}\circ\mu_{P_1\ten\cdots\ten P_p}\ten X=\pi^{(N)}\circ\sigma_{P_1\cdots P_p}\ten X
\end{equation}
For example, if $p=1$, $N=2$, and $P\in\FL{2}$ is the element $[x_1,x_2]$,
then $\mu_P:[2]\to[1]$ is the Lie bracket and, by \eqref{eq:action}
\[\pi\circ\mu_P=\pi\circ(\id\ten\pi)-\pi\circ(\id\ten\pi)\circ(1\,2)=\pi^{(2)}\circ\sigma_P\]
with $\sigma_P=\id-(1\,2)$. Dually, for any $Q_1\ten\cdots\ten Q_q\in \FL{N}^{\ten q}$ of degree $\delta_N$,
there are elements
\begin{align*}
\delta_{Q_1\otimes\cdots\otimes Q_q}\in \LCA([q],[N])\qquad\mbox{and}\qquad
\wt{\sigma}_{Q_1\cdots Q_q}\in\sfk\SS_N
\end{align*}
such that the following holds in $\DY{1}(\VDY{1}, [N]\ten\VDY{1})$
\begin{equation}\label{eq:highcoaction}
\delta_{Q_1\otimes\cdots\otimes Q_q}\otimes X \circ {\pi^*}^{(q)}=
\wt{\sigma}_{Q_1\cdots Q_q}\otimes X\circ {\pi^*}^{(N)}
\end{equation}

The commutativity of the diagram then follows easily.
Namely, assume that $X=\sarch{M}{\sigma}$ for some $M>0$ and $\sigma\in\SS_M$. In particular,
we have $\xi\univ(X)=(x_1\cdots x_M)\ten (y_{\wt{\sigma}(1)}\cdots y_{\wt{\sigma}(M)})=: Q_X\ten P_X$ and
\[(\sfp\circ\xi\lba)(\mu_{P_1\ten\cdots\ten P_p}\circ \delta_{Q_1\otimes\cdots\otimes Q_q})=(Q_1\cdots Q_q)\otimes (P_1\cdots P_p)\]
Then,  by \eqref{eq:highaction} and \eqref{eq:highcoaction}, 
\[\begin{split}
(\mu_{P_1\ten\cdots\ten P_p}\circ \delta_{Q_1\otimes\cdots\otimes Q_q})\cdot \sarch{M}{\sigma}
&=
\pi^{(p)}\circ \mu_{P_1\ten\cdots\ten P_p}
\circ
\delta_{Q_1\otimes\cdots\otimes Q_q} \otimes \sarch{M}{\sigma}\circ\, {\pi^*}^{(q)}\\
&=
\pi^{(N)}\circ \sigma_{P_1\cdots Pp}\circ(\sigma_{Q_1\cdots Q_q}\circ\tau_N)
\otimes \sarch{M}{\sigma}\circ\,{\pi^*}^{(N)}\\
&=
\pi^{(N+M)}\circ \sigma_{P_1\cdots Pp}\circ\wt{\sigma}_{Q_1\cdots Q_q}
\otimes \sigma\circ\,{\pi^*}^{(N+M)}\\
\end{split}\]
which, under $\xi\univ$, corresponds precisely to the element 
$(Q_1\cdots Q_q)\cdot Q_X\otimes P_X\cdot(P_1\cdots P_p)$ in 
$((\FA{N+M})_{\delta_{N+M}}\ten(\FA{N+M})_{\delta_{N+M}})_{\SS_{N+M}}$.
Therefore
\begin{align*}
\xi\univ((\mu_{P_1\ten\cdots\ten P_p}\circ \delta_{Q_1\otimes\cdots\otimes Q_q})\cdot \sarch{M}{\sigma})&=\\
=(\sfp\circ\xi\lba)(\mu_{P_1\ten\cdots\ten P_p}&\circ \delta_{Q_1\otimes\cdots\otimes Q_q})\cdot\xi\univ(\sarch{M}{\sigma})
\end{align*}
and the result follows.
\end{pf}

Applying the result to $X=\id_{\ten_{k=1}^n\VDY{k}}$, yields the following

\begin{corollary}\label{cor:PBW DY}
	The following is a commutative diagram of convolution algebras.
	\[\xymatrix{
		\mbox{$\displaystyle
			\DY{n}(\bigotimes_{k=1}^n\VDY{k},\bigotimes_{k=1}^n\VDY{k}) 
			$}
		\ar[r]^(.4){\xi^n\univ}
		& 
		\mbox{$\displaystyle
			{\bigoplus_{N\geqslant0}}
			\left(\homo{\FA{N}^{\ten n}}\ten\homo{\FA{N}^{\ten n}}\right)_{\SS_N}
			$}
		\\
		\LBA(\whT{[1]}^{\ten n},T[1]^{\ten n})
		\ar[r]_(.4){\xi^n\lba} \ar[u]^{\sfa^n}
		&
		\mbox{$\displaystyle
			{\bigoplus_{N\geqslant0}}
			\left((T\FL{N}^{\ten n})_{\delta_N}\ten (T\FL{N}^{\ten n})_{\delta_N}\right)_{\SS_N}
			$} 
		\ar[u]_{\sfp^n}
	}
	\]
\end{corollary}

\subsection{PBW theorem for $\DYUA{n}$}\label{ss:LBAPBW}

Let 
\[\sfsym:\LBA(\whS{[1]}^{\ten n},S[1]^{\ten n})\to\LBA(\whT{[1]}^{\ten n},T[1]^{\ten n})\]
be the map \eqref{eq:sym LBA}. The following result shows that the
composition $\sfa\circ\sfsym$ can be thought of as the symmetrisation
map $S \ll\to U\ll$ for a Lie algebra $\ll$.

\begin{theorem}\label{th:PBW}
The following is a commutative diagram
\[\xymatrix{
\mbox{$\displaystyle
\DY{n}(\bigotimes_{k=1}^n\VDY{k},\bigotimes_{k=1}^n\VDY{k}) 
$}
\ar[r]
& 
\mbox{$\displaystyle
{\bigoplus_{N\geqslant0}}
\left(\homo{\FA{N}^{\ten n}}\ten\homo{\FA{N}^{\ten n}}\right)_{\SS_N}
$}
\\
\LBA(\whT{[1]}^{\ten n},T[1]^{\ten n})
\ar[r] \ar[u]^{\sfa^n}
&
\mbox{$\displaystyle
{\bigoplus_{N\geqslant0}}
\left((T\FL{N}^{\ten n})_{\delta_N}\ten (T\FL{N}^{\ten n})_{\delta_N}\right)_{\SS_N}
$} 
\ar[u]_{\sfp^n}
\\
\LBA(\whS{[1]}^{\ten n},S[1]^{\ten n})\ar[u]^{\sfsym} 
\ar[r]
&
\mbox{$\displaystyle
{\bigoplus_{N\geqslant0}}
\left((S\FL{N}^{\ten n})_{\delta_N}\ten (S\FL{N}^{\ten n})_{\delta_N}\right)_{\SS_N}
$} 
\ar[u]_{\sfsym\ten\sfsym}
}
\]
where the right vertical arrows are the symmetrisation map $S\FL{N}\to
T\FL{N}$ and quotient map $T\FL{N}\to U\FL{N}=\FA{N}$ for the Lie algebra
$\FL{N}$. 

Moreover, the map $\sfa\circ\sfsym$
is an isomorphism of cosimplicial spaces.
\end{theorem}
\begin{pf}
The commutativity of the diagram follows from Lemma \ref{le:cosi and Lie}
and Proposition \ref{cor:PBW DY}.
The fact that $\sfa\circ\sfsym$ is an isomorphism then follows
from the PBW Theorem for the Lie algebra $\FL{N}$, and the
fact that it is compatible with the cosimplicial structure from
Proposition \ref{pr:a cosimplicial}.
\end{pf}

\subsection{PBW conjecture for $\DY{n}$}\label{ss:PBW-DY}

Let $\SDY{n}$ be the colored $\PROP$ generated by an
$\LBA$--module $([1],\mu,\delta)$ and objects
$\VDY{k}$ endowed with maps 
$\pi_{\VDY{k}}:[1]\ten\VDY{k}\to\VDY{k}$, $\pi_{\VDY{k}}^*:\VDY{k}\to[1]\ten\VDY{k}$,
$k=1,\dots, n$, satisfying the relations
\begin{gather*}
\pi_{\VDY{k}}\circ\id_{[1]}\ten\pi_{\VDY{k}}\circ(\id_{[2]}-(1\,2))\ten\id_{\VDY{k}}=0\\[1.1ex]
(\id_{[2]}-(1\,2))\ten\id_{\VDY{k}}\circ\id\ten\pi_{\VDY{k}}^*\circ\pi_{\VDY{k}}^*=0\\[1.1ex]
\pi_{\VDY{k}}^*\circ\pi_{\VDY{k}}=\id\ten\pi_{\VDY{k}}\circ(1\,2)\circ\id\ten\pi_{\VDY{k}}^*
\end{gather*}
Thus, $\SDY{n}$ encodes a Lie bialgebra $\b$, together with $n$ \DYt
modules over the underlying vector space of $\b$ endowed with trivial
bracket and cobracket.

The $\PROP$ $\DY{n}$ (resp. $\SDY{n}$) is $\IN$--filtered (resp. graded) by
$\deg(\delta)=0=\deg(\mu)$, and $\deg{\pi_{\VDY{k}}}=1=\deg{\pi_{\VDY{k}}^*}$.
Moreover, there is a canonical filtered functor $\SDY{n}\to\gr(\DY{n})$
which is the identity on objects and is easily seen to be full. It is natural to
conjecture the following results which, together, extend Theorem \ref{ss:LBAPBW}.

\begin{conjecture}\hfill
\begin{enumerate}
\item The functor $\SDY{n}\to\gr(\DY{n})$ is faithful, and therefore
an isomorphism of $\PROP$s.
\item The map $\sfa\circ\sfsym:\LBA(\whS{[1]}^{\ten n},S[1]^{\ten n})\to
\SDY{n}(\bigotimes_{k=1}^n\VDY{k},\bigotimes_{k=1}^n\VDY{k})$ is an
isomorphism.
\end{enumerate}
\end{conjecture}

\subsection{Cohomology of $\DYUA{}$}

\begin{theorem}\label{th:Hochschild}
The map $\sfa\circ\sfsym$ induces an isomorphism
\[H^n(\TenUguniv{\bullet},d_H)\simeq
\bigoplus_{j=0}^{n}\LBA\left(\wedge^j[1],\wedge^{n-j}[1]\right)\]
In particular, $H^{0}(\TenUguniv{\bullet}, d_H)=\sfk$ and $H^{1}(\TenUguniv{\bullet}, d_H)=0$.
\end{theorem}
\begin{pf}
By Theorem \ref{th:PBW}, $\sfa\circ\sfsym$ is an isomorphism
of cosimplicial spaces. The result then follows from Proposition
\ref{pr:coho} applied to the $\PROP$ $\sfuP=\LBA$. Namely,
we have
\[\LBA(\whS{[1]}^{\otimes\bullet},S[1]^{\otimes\bullet})
=
\LBA(({S^*})^{\otimes\bullet}[1],S^{\otimes\bullet}[1])\\
=
\LBA_\Sch(S^{\otimes\bullet}\boxtimes S^{\otimes\bullet})\\
\]
where the first equality relies on the equality of
Schur functors $\displaystyle{\wh{S}=S^*}$, and the fact that the
cosimplicial structure on $S[1]
^{\otimes\bullet}$ (resp. the simplicial structure on 
$\wh{S}[1]^{\ten \bullet}$) is induced by that on the Schur
functors $S^{\otimes\bullet}$ (resp. ${S^*}^{\otimes\bullet}$), 
and the second one from \eqref
{eq:PSch and Hom}. The result now follows from
Proposition \ref{pr:coho}, and the equality of
Schur functors $(\wedge^n)^*=\wedge^n$ for any $n\geqslant0$. 
\end{pf}

\noindent
\begin{remark}
Theorem \ref{th:Hochschild} can also be obtained via Lemma \ref
{le:cosi and Lie} from the fact that the diagram in \ref{ss:LBAPBW}
is one of cosimplicial spaces, and the standard computation of the Hochschild
cohomology of a symmetric algebra. The proof via Schur bifunctors given
above yields a more uniform answer for the refinements of the $\PROP$
$\LBA$ introduced in Sections \ref{s:plba} and \ref{s:sgp-ext}. Note, however
that it still depends on the PBW Theorem for $\TenUguniv{\bullet}$, which is
obtained from the identification of $\DY{n}(\bigotimes_{k=1}^n\VDY{k},
\bigotimes_{k=1}^n\VDY{k})$ (resp. $\LBA(\whS{[1]}^{\ten n},S[1]^{\ten n})$)
with free algebras (resp. Lie algebras).
\end{remark}

\subsection{Explicit description of $H^n(\TenUguniv{\bullet},d_H)$}

The cohomology of $\DYUA{\bullet}$ can be described more explicitly, 
in the spirit of \ref{ss:coho}. Denote by $\wt{\iota}_{n,j}$ the inclusion 
\[
\LBA(\wedge^j[1],\wedge^{n-j}[1])\to
\LBA(\wh{S}[1]^{\ten n}, S[1]^{\ten n})
\]
defined by \eqref{eq:qiso-P}.
We first observe that in $\DY{1}$
\[
\pi_{T[1]}\circ\iota_0\ten\id_{\VDY{1}}=\id_{\VDY{1}}
\aand
\pi_{T[1]}\circ\iota_1\ten\id_{\VDY{1}}=\pi_{\VDY{1}}
\]
and dually
\[
{\iota_0}^*\ten\id_{\VDY{1}}\circ\pi^*_{\wh{T}[1]}=\id_{\VDY{1}}
\aand
{\iota_1}^*\ten\id_{\VDY{1}}\circ\pi^*_{\wh{T}[1]}=\pi^*_{\VDY{1}}
\]
It follows that the inclusion 
\[
\iota_{n,j}:
\LBA(\wedge^j[1],\wedge^{n-j}[1])\to
\DY{n}(\ten_{k=1}^n\VDY{k},\ten_{k=1}^n\VDY{k})
\]
where $\iota_{n,j}=\sfa\circ\Sym\circ\wt{\iota}_{n,j}$, sends 
a morphism $\phi\in\LBA(\wedge^j[1],\wedge^{n-j}[1])$ to
\[
\iota_{n,j}(\phi)=\frac{1}{n!}\sum_{\sigma\in\SS_n}(-1)^{\sigma}\pi_{\ol{J_{\sigma}}}\circ\phi\circ\pi^*_{J_{\sigma}}
\]
where $J_{\sigma}=\{\sigma(1),\dots, \sigma(j)\}$, $\ol{J}_{\sigma}$ is its complement,
$\pi_{\ol{J}_{\sigma}}$ denotes the ordered action of $[n-j]=[1]^{\ten n-j}$ on the components 
$\VDY{k}$, $k\in\ol{J}_{\sigma}$, 
and $\pi^*_{J_{\sigma}}$ the ordered coaction of $[j]=[1]^{\ten j}$ on the components 
$\VDY{k}$, $k\in{J}_{\sigma}$.
For example, for $n=2$, we have
\[
\iota_{2,1}(\id_{[1]})=
\frac{1}{2}\left(
\pi_{\VDY{2}}\circ(1\,2)\circ\pi_{\VDY{1}}^*
-
\pi_{\VDY{1}}\circ(1\,2)\circ\pi_{\VDY{2}}^*
\right)
=
\frac{1}{2}\left(
r_{\VDY{1},\VDY{2}}-r_{\VDY{2},\VDY{1}}
\right)
\]
\ie the antisymmetric $r$--matrix corresponds to the identity in $\LBA([1],[1])$.

Thus, the image of $\LBA(\wedge^j[1],\wedge^{n-j}[1])$ inside $\DY{n}(\ten_
{k=1}^n\VDY{k},\ten_{k=1}^n\VDY{k})$ consists of linear combinations of arc
diagrams with exactly one coaction or one action on each bold line, which are
antisymmetric under permutation of the bold lines.

In terms of the identification with free Lie algebras given by Proposition
\ref{pr:LBA FLN}, the above isomorphism yields
\[H^i(\DYUA{\bullet},d_H)
\cong
\bigoplus_{N\geqslant 0}
\bigoplus_{j=0}^i
\left[
\left(\wedge^j \FL{N}\right)_{\delta_N}
\otimes\left(\wedge^{i-j} \FL{N}\right)_{\delta_N}\right]_{\SS_N}\]
Then $H^i(\DYUA{\bullet},d_H)$ embeds in 
$\DYUA{i}\simeq{\bigoplus_{N\geqslant0}}
\left(\homo{\FA{N}^{\ten i}}\ten\homo{\FA{N}^{\ten i}}\right)_{\SS_N}$ via \eqref{eq:qiso}.

\subsection{Enriquez's universal algebras}\label{ss:enri}

In \cite{e1,e2,e3}, Enriquez introduced the universal algebras $\{\enrique{n}\}
_{n\geq 1}$ associated to the $\PROP$ $\LBA$. As a $\sfk$--vector space,
$\enrique{n}$ is defined as the space of coinvariants
\[
\enrique{n}=
\bigoplus_{N\geqslant0}\big(\homo{\FA{N}^{\ten n}}\ten\homo{\FA{N}^{\ten n}}\big)_{\SS_N}
\]
introduced in \ref{ss:Un Fn}--\ref{ss:free-lie-pbw-conv}. The multiplication on
$\enrique{n}$ is defined by an explicit formula in the basis $\{x_{\ul{N}}\otimes
y_{\sigma(\ul{N}')}\}_{\ul{N},\ul{N'},\sigma}$ given in \ref{ss:Un Fn}, and proved
to be associative by a lengthy calculation \cite[B.1--B.2]{e2}.\footnote{The multiplication on $\enrique{n}$ differs from the convolution product discussed
in \ref{ss:free-lie-pbw-conv}.} 

$\enrique{}$ is universal in the following sense. For any Lie bialgebra $\b$ with
Drinfeld double $\gb=\b\oplus\b^*$ and $r$--matrix $r_{\b}=\sum_{i\in I} b_i\ten
b^i\in\b\ctp\b^*$, the linear map $\Erho{\b}{}:\enrique{}\to\whUA{\gb}{}$ given by
\[
\Erho{\b}{}(x_1\cdots x_N\ten y_{\sigma(1)}\cdots y_{\sigma(N)})=
\sum_{\ul{i}\in I^N}b_{i_1}\cdots b_{i_N}b^{i_{\sigma(1)}}\cdots b^{i_{\sigma(N)}}
\]
is an algebra homomorphism. Similarly, for any $n\geq 2$, there is a map $\Erho{\b}
{n}:\enrique{n}\to\whUA{\gb}{\ten n}$ given by
\[\Erho{\b}{n}(x_{\ul{N}}\ten y_{\sigma(\ul{N}')})=
\sum_{\ul{i}\in I^N}b_{\ul{N}(\ul{i})}\cdot b^{{\sigma}(\ul{N}')(\ul{i})}\]
is an algebra homomorphism, where 
\begin{align*}
b_{\ul{N}(\ul{i})}
&=\bigotimes_{k=1}^n b_{i_{N_1+\cdots+N_{k-1}+1}}\cdots b_{i_{N+1+\cdots+N_k}}\\
b^{\sigma(\ul{N}')(\ul{i})}
&=
\bigotimes_{k=1}^n b^{i_{\sigma(N'_1+\cdots+N'_{k-1}+1)}}\cdots b^{i_{\sigma(N'_1+\cdots+N'_k)}}
\end{align*}

\subsection{The isomorphism $\DYUA{n}\simeq \enrique{n}$}\label{ss:identify}

The following result identifies the algebra $\enrique{n}$ with $\DYUA{n}$,
thereby considerably simplifying the proof of the existence of an algebra structure
on $\enrique{n}$ given in \cite[Appendices B and C]{e2}.

Let $\xi^n\univ:\DYUA{n}\to \enrique{n}=\bigoplus_{N\geqslant0}\big(\homo
{\FA{N}^{\ten n}}\ten\homo{\FA{N}^{\ten n}}\big)_{\SS_N}$ be the map defined in \ref
{ss:Un Fn}.

\begin{proposition}\hfill
 \begin{enumerate}
 \item $\xi^n\univ$ is an isomorphism of cosimplicial spaces.
 \item There is a commutative diagram
 \[
\xymatrix{
\DYUA{n}\ar[r]^{\rho^n_{\b}} \ar[d]_{\xi^n\univ} & \DYA{\b}{n} \ar@{=}[d]\\
\enrique{n} \ar[r]_{\rho^n_{\gb}} & \whUA{\gb}{\ten n}
}
\]
\end{enumerate}
\end{proposition}
\begin{pf}
$(1)$ The fact that $\xi^n\univ$ is an isomorphism was proved in \ref{ss:Un Fn},
and its compatibility with the cosimplicial structure in Lemma \ref{le:cosi and Lie}.
$(2)$ The commutativity of the diagram follows by direct inspection.
\end{pf}

\noindent
\begin{remark}
It seems very likely that the map $\xi^n\univ$ is an algebra homomorphism.
This would follow from a detailed inspection of the algebra structure on
$\enrique{n}$, or from the commutativity of the above diagram if the collection
of maps $\rho^n_{\gb}$ were known to be be injective. In any event, the above
proposition shows that $\DYUA{n}$ is an isomorphic replacement of
$\enrique{n}$ with a more naturally defined multiplication.
\end{remark}

\section{The universal algebra of a split pair}\label{s:plba}

In this section, we give a relative version of the results of Sections \ref{s:univDY}
by adapting them to case of a split pair of Lie bialgebras, as defined in \cite{ATL1}.

\subsection{The $\PROP$ $\PLBA$}\label{ss:plba}

Let $(\b,\a)$ be a split pair of Lie bialgebras, \ie $\b$ and $\a$ are Lie bialgebras 
endowed with Lie bialgebra morphisms
\[\a\xrightarrow{i_{\a}}\b\xrightarrow{p_{\a}}\a\]
such that $p_{\a}\circ i_{\a}=\id_{\a}$. These maps induce an isometric inclusion of the
corresponding Drinfeld doubles $\ga\subset\gb$, and a restriction functor 
$\Res_{\b,\a}: \DrY{\b}\to\DrY{\a}$.

\begin{definition}
Let $\PLBA$ be the colored $\PROP$ generated by two Lie bialgebra objects
$\pb,\pa$ related by Lie bialgebra morphisms $i_{\pa}:\pa\to\pb$, $p_{\pa}:
\pb\to\pa$ such that $p_{\pa}\circ i_{\pa}=\id_{\pa}$. 
\end{definition}

The kernel ${\pmm}$ of the projection $p_{\pa}$, 
is an object of $\PLBA$, and $\pb$ decomposes as $\pb=\pa\oplus{\pmm}$. 
${\pmm}$ is an ideal in $\pb$, and has a Lie algebra structure. 
It is also a coideal, but has no natural Lie coalgebra structure.

\subsection{Universal property of $\PLBA$} 
The following is clear.
\begin{proposition}\hspace{0.1cm}
\begin{enumerate}
\item
The $\PROP$ $\PLBA$ is endowed with a pair of functors
$\beta,\alpha:\LBA\to\PLBA$ given by
\[\beta[1]=\pb\aand\alpha[1]=\pa\]
The maps $i_{\pa},p_{\pa}$ in $\PLBA$ induce
two natural transformations $i_{\alpha}:\alpha\to\beta, {p_{\alpha}}:\beta\to\alpha$
such that ${p_{\alpha}}\circ{i_{\alpha}}=\id_{\alpha}$
\[
\xymatrix{\LBA\ar@/^15pt/[rr]^{\beta}="\upsilon"\ar@/_15pt/[rr]_{\alpha}="\lambda"  
& & \PLBA \ar@{=>}@<1ex> "\upsilon";"\lambda"^{p_{\alpha}} \ar@{=>}@<1ex> 
"\lambda";"\upsilon"^{i_{\alpha}}}
\]
\item $\PLBA$ is universal \wrt property (i): for any tensor category $\C$
for which it holds, there is a unique tensor functor $F:\PLBA\to\C$ such that
the following diagram commutes
\[
\xymatrix{\LBA\ar@/^15pt/[rr]^{\beta}="\upsilon"\ar@/_15pt/[rr]_{\alpha}="\lambda"  
\ar@/^40pt/[rrrrr]^(.55){\beta_{\C}}="UC" \ar@/_40pt/[rrrrr]_(.55){\alpha_{\C}}="LC" & & \PLBA   \ar@{-->}[rrr]
\ar@{=>}@<1ex> 
"UC";"LC"^(.2){p_{\C}} \ar@{=>}@<1ex> "LC";"UC"^(.8){i_{\C}}
\ar@{=>}@<1ex> 
"\upsilon";"\lambda"^{p_{\alpha}} \ar@{=>}@<1ex> "\lambda";"\upsilon"^{i_{\alpha}} 
& & & \C}
\]
\end{enumerate}
\end{proposition}

\subsection{Alternative presentation of $\PLBA$}\label{ss:plba2}

The following presentation of $\PLBA$ is more convenient for computations.
Let $\uPLA$ be the $\PROP$ generated by $\mu:[2]\to[1]$ and
$\dmp{0}:[1]\to[1]$ satisfying the relations \eqref{eq:bracket}, 
\begin{equation}\label{eq:projbra}
\dmp{0}^2=\dmp{0}\qquad\mbox{and}\qquad\mu\circ(\dmp{0}\ten\dmp{0})=\dmp{0}\circ\mu
\end{equation}
Let $\uPLCA$ be the $\PROP$ generated by
$\delta:[1]\to[2]$ and $\dmp{0}:[1]\to[1]$ satisfying the relations \eqref{eq:cobracket},
\begin{equation}\label{eq:projcobra}
\dmp{0}^2=\dmp{0}\qquad\mbox{and}\qquad(\dmp{0}\ten\dmp{0})\circ\delta=\delta\circ\dmp{0}
\end{equation} 
Let $\uPLBA$ be the $\PROP$ generated by $\mu:[2]\to[1]$,
$\delta:[1]\to[2]$, and $\dmp{0}:[1]\to[1]$, satisfying the relations \eqref{eq:bracket},
\eqref{eq:cobracket}, \eqref{eq:cocycon}, \eqref{eq:projbra}, \eqref{eq:projcobra}.
Finally, let $\PLA,\PLCA,\PLBA$ be their corresponding completions. 

The two presentations of $\PLBA$ are canonically equivalent by sending $[1]$
to $\pb$ and the idempotent $\dmp{0}$ to the composition $i_{\pa}\circ p_{\pa}:\pb\to\pb$. 

\begin{corollary}\hfill
\begin{enumerate}
\item There is a forgetful functor $\PLBA\to\LBA$, mapping $[1]_{\PLBA}$ to $[1]_{\LBA}$
and $\dmp{0}$ to $\id_{[1]_{\LBA}}$.
\item There is a forgetful functor $\PLBA\to\LBA$, mapping $[1]_{\PLBA}$ to $[1]_{\LBA}$
and $\dmp{0}$ to $0$.
\end{enumerate}
\end{corollary}

\subsection{Factorisation of morphisms in $\PLBA$}\label{ss:normplba}

Set $\dmp{1}=\id-\dmp{0}$, and $\I=\{0,1\}$. The projections
\[\dmp{\ul{i}}=\dmp{i_1}\ten\cdots\ten\dmp{i_N},\quad \ul{i}=(i_1,\ldots,i_N)\in \I^N\]
are a complete family of idempotents in $\PLA([N],[N])$
and $\PLCA([N],[N])$, \ie
\[\dmp{\ul{i}}\circ\dmp{\ul{i}\prime}=\drc{\ul{i}\ul{i}\prime}\dmp{\ul{i}}
\aand
\sum_{\ul{i}\in\I^N}\dmp{\ul{i}}=\id_{[N]}\]
There is a natural right (resp. left) action of $\Idem{N}$ on $\PLA([N],[q])$
and $\PLCA([p],[N])$ given by
\[
\phi\cdot f=\sum_{\ul{i}\in\I^N}f(\ul{i})\,\phi\circ\dmp{\ul{i}}
\aand
f\cdot\psi=\sum_{\ul{i}\in\I^N}f(\ul{i})\,\dmp{\ul{i}}\circ\psi
\]
Set $\Gsplit=\Idem{N}\rtimes\SS_N$, where $\sigma\in\SS_N$ acts on
$f\in\Idem{N}$ by $\sigma\cdot f=f\circ\sigma^{-1}$.

\begin{proposition}
The embeddings $\PLA,\PLCA\rightarrow\PLBA$ induce an isomorphism
of $(\SS_q,\SS_p)$--bimodules
\[
\PLBA([p], [q])\simeq\bigoplus_{N\geqslant0}
\PLCA([p], [N])\ten_{\Gsplit}\PLA([N],[q])
\]
\end{proposition}

\begin{pf}
The proof is similar to that of Proposition \ref{pr:LBA factorisation}.
The computation can be carried out with the $\PROP$s $\uPLA,\uPLCA,
\uPLBA$ introduced in \ref{ss:plba2} since these contain the objects
$[p],[N],[q]$. A morphism in $\uPLBA$ can be represented as an
oriented graph obtained from the composition of brackets, cobrackets,
permutations, and idempotents. The compatibility \eqref{eq:cocycon}
between $\delta$ and $\mu$, and the relations
\begin{equation*}
\delta\circ\dmp{0}=(\dmp{0}\ten\dmp{0})\circ\delta\qquad
\dmp{0}\circ \mu=\mu\circ(\dmp{0}\ten\dmp{0})
\end{equation*}
allow to reorder the morphisms so that the cobrackets precede the
brackets, and the idempotent $\dmp{0}$ occur in between. This yields
a surjective map
\[\bigoplus_{N\geqslant0}
\left(\PLCA([p],[N])\ten \PLA([N],[q])\right)\to\PLBA([p],[q])\]
which factors through the action of $\Idem{N}\rtimes\SS_N$. 
The injectivity follows as in \ref{ss:morlba}.
\end{pf}

\subsection{Morphisms in $\PLA$ and $\PLCA$}\label{ss:morpla}

\begin{lemma}
There are isomorphisms of left (resp. right) $\sfk[\I^N]\rtimes\SS_N$--modules
\begin{equation*}
\PLA([N],[p])\simeq\Idem{N}\ten\LA([N],[p])
\end{equation*}
\begin{equation*}
\PLCA([p],[N])\simeq\LCA([p],[N])\ten\Idem{N}
\end{equation*}
where $\SS_N$ acts diagonally on the right--hand side, which are compatible
with the action of $\SS_p$.
\end{lemma}
\begin{pf}
We only explain the isomorphism in $\PLA$. The result for $\PLCA$
follows by observing that $\mathsf{PLCA}\simeq\mathsf{PLA}\op$.
Every morphism in $\uPLA([N],[p])$ is represented by a linear
combination of oriented graphs from $N$ sources to $p$ targets. Since
$\id_{[1]}=\dmp{0}+\dmp{1}$, all the edges of these graphs can be assumed
to be decorated by the idempotents $\dmp{0}$ or $\dmp{1}$. The relations
\begin{align*}
\dmp{0}\circ\mu&=\mu\circ(\dmp{0}\ten\dmp{0})\\
\dmp{1}\circ\mu&=\mu\circ\big(\dmp{1}\ten\dmp{1}+\dmp{1}\ten\dmp{0}+\dmp{0}\ten\dmp{1}\big)
\end{align*}
allow to move all idempotents to the $N$ sources and yield the surjectivity
of the map
\[
\Idem{N}\ten\LA([N],[p])\to\PLA([N],[p])\qquad\qquad 
f\ten P\mapsto f\cdot P
\]
Its injectivity follows from the canonical embedding $\LA\to\PLA$ and the isomorphism
\[
\bigoplus_{\ul{i}\in\I^N}\PLA([N],[p])\circ\dmp{\ul{i}}\simeq\PLA([N],[p])
\]
\end{pf}

\subsection{$\PLBA$ and free Lie algebras}\label{ss:morplba}

The following is a direct consequence of \ref{ss:morpla}, \ref{ss:normplba} and
Lemma \ref{le:morla}.
\begin{proposition}\label{pr:morplba}\hfill
\begin{enumerate}
\item There is an isomorphism of $(\SS_q,\SS_p)$--bimodules
\[
\PLBA([p],[q])\simeq
\bigoplus_{\substack{N\geqslant0}}
\left(\homo{\FL{N}^{\ten p}}\ten\Idem{N}\ten\homo{\FL{N}^{\ten q}}\right)_{\SS_N}
\]
\item Let $F\in k\SS_p$ and $G\in k\SS_q$ be idempotents, and
$F[p]=([p],F)$, $G[q]=([q],G)$ the corresponding objects in $\PLBA$.
Then one has
\[
\PLBA(F[p], G[q])
\simeq\bigoplus_{N\geqslant0}
\left(F(\FL{N}^{\ten p})_{\delta_N}\ten\Idem{N}\ten G(\FL{N}^{\ten q})_{\delta_N}\right)_{\SS_N}
\]
In particular,
\begin{align*}
\PLBA(\wh{T}[1]^{\otimes n},T[1]^{\otimes n})
&\simeq
\bigoplus_{\substack{N\geqslant0}}
\left((T\FL{N}^{\otimes n})_{\delta_N}\ten\Idem{N}\ten(T\FL{N}^{\otimes n})_{\delta_N}\right)_{\SS_N}\\
\PLBA(\wh{S}[1]^{\otimes n},S[1]^{\otimes n})
&\simeq
\bigoplus_{\substack{N\geqslant0}}
\left((S\FL{N}^{\otimes n})_{\delta_N}\ten\Idem{N}\ten(S\FL{N}^{\otimes n})_{\delta_N}\right)_{\SS_N}
\end{align*}
\end{enumerate}
\end{proposition}

\subsection{Universal Drinfeld--Yetter modules and $\PLBA$}\label{ss:dyplba}\label{ss:PLBAcompletions}

\begin{definition}
The category $\PDY{n}$, $n\geqslant 1$, is the colored $\PROP$ generated
by $n+1$ objects $\APDY{1}$ and $\{\VPDY{k}\}_{k=1,\dots,n}$, and morphisms
\[
\mu:\APDY{2}\to\APDY{1}\qquad\delta:\APDY{1}\to\APDY{2}\qquad
\dmp{}:\APDY{1}\to\APDY{1}
\]
\[
\pi_k:\APDY{1}\ten \VPDY{k}\to \VPDY{k}\qquad
\pi_k^*:\VPDY{k}\to\APDY{1}\ten \VPDY{k}
\]
such that $(\APDY{1},\mu,\delta,\dmp{})$ is a  $\PLBA$--module in $\PDY{n}$, and, for
every $k=1,\dots, n$, $(\VPDY{k},\pi_k,\pi_k^*)$ is a Drinfeld--Yetter module 
over $\APDY{1}$.
\end{definition}

Set
\begin{equation}\label{eq:univ PLBA}
\PDYUA{n}=\pEnd{\PDY{n}}{\bigotimes_{k=1}^n\VPDY{k}}
\end{equation}

The algebras $\PDYUA{n}$ are universal in the following sense. Let $(\b,\a)$ be a
split pair of bialgebras over $\sfk$. Then, for any $n$--tuple  
$\{V_k,\pi_k,\pi_k^*\}_{k=1}^n$ of \DYt modules over $\b$, there is a realisation functor 
\[\G_{(\b,\a,V_1,\dots, V_n)}:\PDY{n}\longrightarrow{\kvect}\]
such that $\pb\mapsto\b$, $\pa\mapsto\a$, and 
$\VDY{k}\mapsto V_k$, $k=1,\ldots,n$. 

\begin{proposition}
Let $\ff\hspace{-0.1cm}:\DrY{\b}\to{\kvect}$ be the forgetful functor,
and $\DYA{\b}{n}=\sfEnd{\ff^{\boxtimes n}}$. The functors $\G_{(\b,\a,V_1,
\dots, V_n)}$ induce an algebra homomorphism
\[
\PDYrho{\b}{\a}{n}:\PDYUA{n}\to\DYA{\b}{n}
\]
\end{proposition}

\begin{pf}
The proof is identical to that of Proposition \ref{ss:DY action}.
\end{pf}

The following is a corollary of the proposition above and \ref{ss:plba}.

\begin{corollary}\label{pr:up-down}
Let $\beta,\alpha:\DYUA{n}\to\PDYUA{n}$ be the two algebra homomorphisms
defined by the functors $\beta,\alpha:\LBA\to\PLBA$.
For any split pair $(\b,\a)$, there are commutative diagrams
\[
\xymatrix{
\PDYUA{n} \ar[r]^{\PDYrho{\b}{\a}{n}} & \DYA{\b}{n} &
& \PDYUA{n} \ar[r]^{\PDYrho{\b}{\a}{n}} & \DYA{\b}{n}\\
\DYUA{n} \ar[u]^{\beta} \ar[ur]_{\DYrho{\b}{n}} & &
& \DYUA{n} \ar[u]^{\alpha} \ar[r]_{\DYrho{\a}{n}} & \DYA{\a}{n} \ar[u]_{\Res_{\a}^*}
}
\]
where $\Res_{\a}^*$ is the morphism induced by the restriction 
$\DrY{\b}\to\DrY{\a}$, and $\DYrho{\b}{n},\DYrho{\a}{n}$ are the
homomorphisms defined in \ref{ss:DY action}.
\end{corollary}

\subsection{Universal invariants}\label{ss:univ-inv}

In $\PLBA$ we can introduce the notion of invariants
with respect to the Lie bialgebra $\pa$.

\begin{definition}
The subalgebra of $\pa$--invariants $\PDYUAinv{n}\subset\PDYUA{n}$ is the
subspace of all $\phi\in\PDYUA{n}$ which commute with the 
action and the coaction of the Lie bialgebra $\pa$ on $\VDY{1},\dots,\VDY{n}$,
that is satisfy 
\begin{align}
\pi_{\pa,\VDY{1}\otimes\cdots\otimes\VDY{n}}\circ\id_{\pb}\otimes\phi
&=
\phi\circ\pi_{\pa,\VDY{1}\otimes\cdots\otimes\VDY{n}}\label{eq:a-inv-1}\\[1.1ex]
\pi^*_{\pa,\VDY{1}\otimes\cdots\otimes\VDY{n}}\circ\phi
&=
\id_{\pb}\otimes\phi\circ {\pi^*}_{\pa,\VDY{1}\otimes\cdots\otimes\VDY{n}}
\label{eq:a-inv-2}
\end{align}
where $\pi_{\pa,\VDY{1}\otimes\cdots\otimes\VDY{n}}
=\pi_{\VDY{1}\otimes\cdots\otimes\VDY{n}}\circ\dmp{}\ten\id_{\VDY{1}\otimes\cdots\otimes\VDY{n}}$
and
$\pi_{\pa,\VDY{1}\otimes\cdots\otimes\VDY{n}}^*
=\dmp{}\ten\id_{\VDY{1}\otimes\cdots\otimes\VDY{n}}\circ\pi^*_{\VDY{1}\otimes\cdots\otimes\VDY{n}}$.
\end{definition}

The following is clear.

\begin{proposition}
The algebra homomorphism $\PDYrho{\b}{\a}{n}:\PDYUA{n}\to\DYA{\b}{n}$
restricts to an algebra homomorphism
\[\PDYrho{\b}{\a}{n}:\PDYUAinv{n}\to\PDYAinv{\b}{\a}{n}\coloneqq 
\sfEnd{\Res_{\a}^{\boxtimes n}:\DrY{\b}^{n}\to\DrY{\a} }\] 
\end{proposition}

\subsection{A basis for $\PDYUA{n}$}
The description of the algebras $\PDYUA{n}$ is obtained along the
same lines of Propositions \ref{ss:univdy} and \ref{ss:univdyn}.

\begin{proposition}
The endomorphisms
\[
\colarch{\ul{N}}{\ul{i}}{\sigma}{\ul{N}^\prime}
=\pi^{(\ul{N})}\circ\dmp{\ul{i}}\ten\id\circ\;\sigma\ten\id\circ\;{\pi^*}^{(\ul{N}^\prime)}
\]
where $N\geqslant0$, $\ul{N},\ul{N}^\prime\in\IN^n$ are such that $|\ul{N}|=N=|\ul{N}
^\prime|$, $\ul{i}\in\I^N$, and $\sigma\in\SS_N$, are a basis of $\PDYUA{n}$.
In particular, the map
\[\xi^n\osplit:
\PDYUA{n}
\longrightarrow
\bigoplus_{N\geqslant0}
\left(\homo{FA^{\ten n}_N}\ten\Idem{N}\ten\homo{FA^{\ten n}_N}\right)_{\SS_N}\]
given by
$\xi^n\osplit(\colarch{\ul{N}}{\ul{i}}{\sigma}{\ul{N}'})=x_{\ul{N}}\ten\drc{\ul{i}}\ten 
y_{\wt{\sigma}(\ul{N}^\prime)}$
is a linear isomorphism.
\end{proposition}

\subsection{PBW theorem for $\PDYUA{n}$}\label{ss:PLBAPBW}

As in the case of $\Uguniv$, the tower of algebras $\{\PDYUA{n}\}_{n\geqslant1}$ is
endowed with face maps $\Delta^n_i:\PDYUA{n}\to\PDYUA{n+1}$ and degeneratation
maps $\E_n^{i}:\PDYUA{n}\to\PDYUA{n-1}$ defining a cosimplicial structure.

Let
\[
\sfa^n: \PLBATfin\to
\PDY{n}(\ten_{k=1}^n\VDY{k},\ten_{k=1}^n\VDY{k})
\]
be the map given on $\phi_{\ul{p},\ul{q}}\in\PLBA(T^\ulp[1],T^\ulq[1])$,
by
$\sfa(\phi_{\ul{p},\ul{q}})=\pi^{(\ul{p})}\circ\phi_{\ul{p},\ul{q}}\circ{\pi^*}^{(\ul{q})}$.

\begin{theorem}\hfill
\begin{enumerate}
\item The following diagram is commutative
\[
\xymatrix@C=0.2in{
\mbox{$\displaystyle
\PDY{n}(\bigotimes_{k=1}^n\VDY{k},\bigotimes_{k=1}^n\VDY{k})
$}
\ar[r] & 
\mbox{$\displaystyle
{\bigoplus_{\substack{N\geqslant0}}}
\left(\homo{\FA{N}^{\ten n}}
\ten\sfk[\I^N]\ten
\homo{\FA{N}^{\ten n}}\right)_{\SS_N}$}\\
\PLBATfin\ar[u]^{\sfa^n} \ar[r] &
\mbox{$\displaystyle
{\bigoplus_{\substack{N\geqslant0}}}
\left(\homo{T\FL{N}^{\ten n}}
\ten\sfk[\I^N]\ten
\homo{T\FL{N}^{\ten n}} \right)_{\SS_N}$}
\ar[u]\\
\PLBASfin \ar[u]^\sfsym \ar[r] &
\mbox{$\displaystyle
{\bigoplus_{\substack{N\geqslant0}}}
\left(\homo{S\FL{N}^{\ten n}}
\ten\sfk[\I^N]\ten
\homo{S\FL{N}^{\ten n}} \right)_{\SS_N}$}
\ar[u]_{\sfsym\ten\id\ten\sfsym}
}
\]
\item The map $\sfa^n\circ\sfsym$ is an isomorphism of cosimplicial spaces.
\end{enumerate}
\end{theorem}

\subsection{Hochschild cohomology}\label{ss:PLBAcoho}

The cosimplicial structure on $\{\PDYUA{n}\}_{n\geqslant1}$ gives rise to the 
\emph{relative} universal Hochschild complex with differential 
$d^n=\sum_{i=0}^{n+1}(-1)^i\Delta^n_i:\PDYUA{n}\to\PDYUA{n+1}$
The morphisms $\{\PDYrho{\b}{\a}{n}\}_{n\geqslant1}$ defined in \ref{ss:PLBAcompletions}
define a chain map between the corresponding Hochschild complexes.

\begin{theorem}\hfill
\begin{enumerate}
\item The map $\sfa^n\circ\sfsym$ induces an isomorphism
\[H^i(\PDYUA{\bullet},d_H)\cong
\bigoplus_{j=0}^{i}\PLBA\left(\wedge^j[1],\wedge^{i-j}[1]\right)\]
In particular, $H^{0}(\PDYUA{\bullet}, d_H)=\sfk$ and $H^{1}(\PDYUA{\bullet}, d_H)=0$.
\item The identification in terms of free Lie algebras of Proposition \ref{pr:morplba} 
yields
\[H^i(\PDYUA{\bullet},d_H)
\cong
\bigoplus_{N\geqslant 0}
\bigoplus_{j=0}^i
\left[
\left(\wedge^j \FL{N}\right)_{\delta_N}\otimes\Idem{N}\otimes\left(\wedge^{i-j} \FL{N}\right)_{\delta_N}\right]_{\SS_N}\]
\end{enumerate}
\end{theorem}

\begin{pf}
The proof is identical to that of Theorem \ref{th:Hochschild} where the
$\PROP$ $\LBA$ is replaced by $\PLBA$.
\end{pf}

\subsection{Hochschild cohomology and invariants}\label{ss:cohoPLBA}

\begin{lemma}
$(\PDYUAinv{n}, d^n)$ is a subcomplex of $(\PDYUA{n},d^n)$.
\end{lemma}

\begin{pf}
It is enough to observe that, if $\phi\in\PDYUA{n}$ satisfies \eqref{eq:a-inv-1}, \eqref{eq:a-inv-2},
then so does $d^n_i(\phi)\in\PDYUA{n+1}$. Namely, let $\D^n_i:\PDY{n}\to\PDY{n+1}$ be 
as in \ref{ss:cosimp-DY}. Then, for any $u\in\PDYUA{n}$, we have
\begin{align*}
\pi_{\pa,\VDY{1}\otimes\cdots\otimes\VDY{n+1}}&\circ\id_{\pb}\otimes d_i^n(u)-d_i^n(u)\circ\pi_{\pa,\VDY{1}\otimes\cdots\otimes\VDY{n+1}}\\
&=\D_i^n\left(\pi_{\pa,\VDY{1}\otimes\cdots\otimes\VDY{n}}\circ\id_{\pb}\otimes u-u\circ\pi_{\pa,\VDY{1}\otimes\cdots\otimes\VDY{n}}\right)
\end{align*}
Set $\D(\phi)=\sum_i\D^n_i(\phi)$. Then, in particular,
\begin{align*}
\pi_{\pa,\VDY{1}\otimes\cdots\otimes\VDY{n+1}}&\circ\id_{\pb}\otimes d(u)-d(u)\circ\pi_{\pa,\VDY{1}\otimes\cdots\otimes\VDY{n+1}}\\
&=\D\left(\pi_{\pa,\VDY{1}\otimes\cdots\otimes\VDY{n}}\circ\id_{\pb}\otimes u-u\circ\pi_{\pa,\VDY{1}\otimes\cdots\otimes\VDY{n}}\right)
\end{align*}
\end{pf}

Let $\PDYext{n}\subset\PDYUA{n}$ be the image of the injective map
\begin{equation}
\xymatrix@C=0.3in{
\PLBA(\wedge^n[1],\wedge^n[1])\ar@{^(->}[r] &
\PLBATfin \ar[r]^(.75){\sfa^n} &
\PDYUA{n}
}
\end{equation}
and set $\PDYextinv{n}=\PDYext{n}\cap\PDYUAinv{n}$.

\begin{proposition}\label{prop:cohoinv} 
$H^n(\PDYUAinv{\bullet}, d^n)\simeq\PDYextinv{n}$.
\end{proposition}

\begin{pf}
Let $f\in\PDYUA{n}$ such that $d(f)=0$. Then there are unique $d(u)\in\PDYUA{n}$ 
and $v\in\PDYext{n}$ such that $f=v+d(u)$. Namely, let $f=v^\prime+d(u^\prime)$
for some $v'\in\PDYext{n}$, $u'\in\PDYUA{n-1}$. It follows
from \ref{ss:PLBAcoho} that
\[
v-v^\prime=d(u-u^\prime)
\quad\Longrightarrow\quad 
v=v'\aand d(u)=d(u')
\]
Assume now $f\in\PDYUAinv{n}$ and $d(f)=0$. Since $f$ satisfies \eqref{eq:a-inv-1}, one has
\begin{align}
\pi_{\pa,\VDY{1}\otimes\cdots\otimes\VDY{n}}\circ\id_{\pb}\otimes v
&=
v\circ\pi_{\pa,\VDY{1}\otimes\cdots\otimes\VDY{n}}\\
\pi_{\pa,\VDY{1}\otimes\cdots\otimes\VDY{n}}\circ\id_{\pb}\otimes d(u)
&=
d(u)\circ\pi_{\pa,\VDY{1}\otimes\cdots\otimes\VDY{n}}\label{eq:inv-2}
\end{align}
and therefore
\begin{align}
\pi_{\pa,\VDY{1}\otimes\cdots\otimes\VDY{n-1}}\circ\id_{\pb}\otimes u
&=
u\circ\pi_{\pa,\VDY{1}\otimes\cdots\otimes\VDY{n-1}}
\end{align}
Similarly for \eqref{eq:a-inv-2}. It follows that $v\in\PDYextinv{n}$ and $u\in\PDYUAinv{n-1}$.
\end{pf}

\begin{corollary}
$(\DYUAinv{n}, d^n)$ is a subcomplex of $(\DYUA{n},d^n)$, and
$H^n(\DYUAinv{\bullet}, d^n)\simeq\DYextinv{n}$.
\end{corollary}

\section{Universal relative twists}\label{s:twists}

In this section, we discuss the existence and uniqueness of invertible elements
in the graded completion of $\PDYUAinv{2}$ satisfying the relative twist equation
\eqref{eq:reltwisteq}. This implies the uniqueness of the tensor structure on the
restriction functor of \DYt modules corresponding to a split pair of Lie
bialgebras.

\subsection{Gradings}\label{ss:grading}

The $\PROP$ $\DY{n}$ has a natural $\IN$--bigrading given by $\deg(\sigma)=
(0,0)$ for any $\sigma\in\SS_N$,
\[\deg(\mu)=(1,0)=\deg(\pi_{\VDY{k}})
\aand
\deg(\delta)=(0,1)=\deg(\pi_{\VDY{k}}^*)\]
for any $1\leq k\leq n$.
The algebra $\DYUA{n}$ inherits this bigrading, and $\deg(r^{\sigma}_{{\ul{N}},
{\ul{N'}}})=(N,N)$, where $r^{\sigma}_{{\ul{N}},{\ul{N'}}}$ is the basis element
defined in \ref{ss:univdyn}, and $|\ul{N}|=N=|\ul{N'}|$.

For any $a,b\in\IN$, the corresponding $\IN$--grading determined by mapping
$(1,0),(0,1)$ to $a,b$ respectively yields the same graded completion $\CDYUA
{n}$ of $\DYUA{n}$, so long as $a+b>0$. For definiteness, we set $a=0$ and
$b=1$.

\subsection{Notation}\label{ss:notation}

There is a natural action of $\SS_n$ on $\DYUA{n}$ 
given by permutations of $\VDY{1}\ten\cdots\ten\VDY{n}$. Specifically, 
for any $\sigma\in\SS_n$, there is an endofunctor $P_\sigma$ of $\DY{n}$
which is the identity on $([1],\mu,\delta)$ and maps each
$([V_k],\pi_k,\pi^*_k)$ to $([V_{\sigma(k)}],\pi_{\sigma(k)},\pi^*_{\sigma(k)})$.
The action of $\sigma\in\SS_n$ on $\DYUA{n}$ is then defined by
$X^{\sigma}\coloneqq \sfAd{\sigma}P_{\sigma}(X)$ for any $X\in\DYUA{n}$.
This is a propic version of the action of $\SS_n$ on $U\gb^{\ten n}$.

The generalisation of the insertion/coproduct maps introduced in \ref{ss:cosimp-DY} is 
defined as follows.
For any $m\geqslant n$, $1\leqslant   i\leqslant   m-n+1$, and $X\in\CDYUA{n}$, 
we define $X_{(i,\dots, i+n-1)}\in\CDYUA{m}$ by 
\[X_{(i,\dots, i+n-1)}\coloneqq \id_{\VDY{1}\ten\cdots\ten\VDY{i-1}}\ten X_{\VDY{i}\ten\cdots\ten\VDY{i+n-1}}
\ten\id_{\VDY{i+n}\ten\cdots\ten\VDY{m}}\]
Then, for any $\sigma\in\SS_m$, we set
\[
X_{(\sigma(i),\dots,\sigma(i+n-1))}\coloneqq (X_{(i,\dots, i+n-1)})^{\sigma}
\]
For any $p_1,\dots, p_n$ with $p_1+\dots+p_n=p\leqslant m$, $p_k\neq0$, 
and $1\leqslant i\leqslant m-p+1$, set $i_k=i+p_1+\cdots+p_{k-1}$, and
$I_k=(i_k, \cdots, i_{k+1}-1)$, $k=1,\dots, n$. Then, we define
$X_{(I_1,\dots, I_n)}\in\CDYUA{m}$ by
\[X_{(I_1,\dots, I_n)}\coloneqq \id_{\VDY{1}\ten\cdots\ten\VDY{i_1-1}}\ten X_{\VDY{[i_1,i_2-1]}\ten\cdots\ten\VDY{[i_{n-1},i_n]}}
\ten\id_{\VDY{i_n+1}\ten\cdots\ten\VDY{m}}\]
where $\VDY{[i_1,i_2-1]}$ denotes the Drinfeld--Yetter module $(\VDY{i_1}\ten\cdots\ten\VDY{i_2-1})$.
As before, for any $\sigma\in\SS_m$, we set $X_{\sigma(I_1,\dots, I_n)}=(X_{I_1,\dots,I_n})^{\sigma}$.

\subsection{Associators}\label{ss:assoc}

Define the $r$--matrix $r=r_{\VDY{1},\VDY{2}}\in\pEnd{\DY{2}}{\VDY{1}
\ten\VDY{2}}$ by
\eqref{eq:rmatrix}, and set $\Omega=r_{12}+r_{21}$.

\begin{definition}
An invertible element $\Phi\in\CDYUA{3}$ is called an \emph{associator}
if the following relations are satisfied (in $\CDYUA{4}$ and $\CDYUA{3}$ respectively).
\begin{itemize}
\item {\bf Pentagon relation}
\[
\Phi_{1,2,34}\Phi_{12,3,4}=\Phi_{2,3,4}\Phi_{1,23,4}\Phi_{1,2,3}
\]
\item {\bf Hexagon relations}
\begin{align*}
&e^{\Omega_{12,3}/2}=\Phi_{3,1,2}e^{\Omega_{13}/2}\Phi_{1,3,2}^{-1}e^{\Omega_{23}/
2}\Phi_{1,2,3}\\
&e^{\Omega_{1,23}/2}=\Phi_{2,3,1}^{-1}e^{\Omega_{13}/2}\Phi_{2,1,3}e^{\Omega_{12}/
2}\Phi_{1,2,3}^{-1}
\end{align*}
\item {\bf Duality}
\[\Phi_{3,2,1}=\Phi_{1,2,3}^{-1}\]
\item {\bf $2$--jet}
\[\Phi=1+\frac{1}{24}[\Omega_{12},\Omega_{23}]\qquad\mod(\DYUA{3})_{\geqslant 3}\]
\end{itemize}
\end{definition}
We denote by $\mathsf{Assoc}$ the set of associators. 

 \subsection{Deformation Drinfeld--Yetter modules}\label{ss:defDY}

Let $\b$ be a Lie bialgebra with Drinfeld double $\gb$. As we explain
below, the algebra $\CDYUA{n}$ introduced in \ref{ss:grading} is a
universal analogue of the topological algebra $\hext{U\gb^{\otimes n}}$.

Let for this purpose $\defDY{\b}$ be the category of \DYt $\b$--modules
in topologically free $\sfk{\fml}$--modules. $\defDY{\b}$ is isomorphic
to the category  $\hDrY{\hextsub{\b}}{\adm}$ of \DYt modules over the
Lie bialgebra $\hextsub{\b}=(\hext{\b},[\cdot,\cdot],\hbar\delta)$, whose
coaction is divisible by $\hbar$. We denote by $\hDYA{\b}{n}$ the algebra of endomorphisms of 
the $n$--fold tensor power of the forgetful functor $\ff:\defDY{\b}\to{\khvect}$.
$\hDYA{\b}{n}$ identifies canonically with the analogous completion
defined for $\hDrY{\hextsub{\b}}{\adm}$. Moreoever, the
realisation functors 
\[\G_{(\hextsub{\b},V_1,\dots, V_n)}:\DY{n}\longrightarrow{\vect_{\hext{\sfk}}}\]
induce a homomorphism $\wh{\rho}_{\b}^{n}:\DYUA{n}\to\hDYA{\b}{n}$
which naturally extends to $\CDYUA{n}$. In particular,\footnote{Note
that $\defDY{\b}$ can also be identified with the category of \DYt modules
over the Lie bialgebra $\hextsup{\b}=(\hext{\b},\hbar[\cdot,\cdot],\delta)$
whose action is divisible by $\hbar$. The corresponding realisation
functors for $\hextsup{\b}$ yield the same homomorphism 
$\wh{\rho}_{\b}^{n}:\DYUA{n}\to\hDYA{\b}{n}$.}
\[\wh{\rho}_{\b}^{1}(\pi_{\VDY{1}}\circ\pi^*_{\VDY{1}})=\hbar\sum_i b_ib^i
\aand
\wh{\rho}_{\b}^{2}(r_{\VDY{1},\VDY{2}})=\hbar\sum_i b_i\otimes b^i\]

In Section \ref{se:braided Cox}, we shall make use of the following
standard construction due to Drinfeld. Let $\Phi\in\CDYUA{3}$ be
an associator. Then, $\hDrY{\b}{\Phi}$ is the braided monoidal category 
with the same objects of $\defDY{\b}$ and commutativity, and 
associativity constraints given respectively by
\[\beta_{\b}=(1\,2)\circ \wh{\rho}_{\b}^{2}(e^{\Omega/2})
\aand
\Phi_{\b}=\wh{\rho}_{\b}^3(\Phi).
\]

\subsection{Universal twists in $\DY{2}$}

The associativity relation \eqref{eq:kill Phi} admits a natural lift to
the $\PROP$s $\DY{n}$. 
\begin{proposition}
Let $\Phi\in\mathsf{Assoc}$, and $J\in\CDYUA{2}$
be such that
\[J_{23}\cdot J_{1,23}\cdot \Phi=J_{12}\cdot J_{12,3}\]
Then, for any Lie bialgebra $\b$, the element $\rho^2_{\b}(J)\in\DYC{\b}^2$
defines a tensor structure on the forgetful functor $\ff:\hDrY{\b}{\Phi}\to{\khvect}$.
\end{proposition}

A simple argument in \cite[Section 6.11]{ATL1} shows that the Etingof--Kazhdan
tensor structure $J\EEK_{VW}$ can be lifted to the $\PROP$ $\DY{2}$.

\subsection{Existence of a universal relative twist}\label{ss:relative-twist}
Let $\Phi\in\mathsf{Assoc}$ and let $(\b,\a)$ be a split pair with corresponding 
Drinfeld doubles $(\gb,\ga)$. Let $\Phi_{\b}$, $\Phi_{\a}$ be the images of $\Phi$
in $\hDYA{\b}{3}$ and $\hDYA{\a}{3}$ respectively.

In \cite[Prop.\;3.17]{ATL1}, we constructed an element $J_{\Phi}\in\hDYA{\b}{2}$, 
which is invariant under $\a$, $J_{\Phi}=1\mod\hbar$, and satisfies the \emph{relative twist
equation}
\begin{equation*}\label{eq:rel-twist-eq}
	J_{\Phi}^{23}J_{\Phi}^{1,23}\Phi_{\b} = \Phi_{\a}J_{\Phi}^{12}J_{\Phi}^{12,3}
\end{equation*}
We also showed \cite[Sec.\;7.7]{ATL1} that the construction of $J_{\Phi}$ is universal \ie 
that it can be realised as an $\pa$--invariant element
\[
J_{\Phi}\in\PDYUAinv{2}\subset\pEnd{\PDY{2}}{\VPDY{1}\ten\VPDY{2}}
\]
We summarize this in the following

\begin{theorem}\label{th:ATL1}
There is a map $\mathsf{Assoc}\to\CPDYUAinv{2}$, $\Phi\to J_{\Phi}$ such that 
$(J_{\Phi})_0=1$ and
\begin{equation}\label{eq:reltwisteq}
(\Phi_{\beta})_{J_{\Phi}}=\Phi_{\alpha}
\end{equation}
where $\Phi_{\beta}$, $\Phi_{\alpha}$ are the images of $\Phi$ in $\CPDYUAinv{3}$
via $\alpha$ and $\beta$, and
\begin{equation}\label{eq:twisted-associator}
\Phi_{J_{\Phi}}\coloneqq J_{\Phi}^{23}J_{\Phi}^{1,23}\Phi (J_{\Phi}^{12,3})^{-1}
(J_{\Phi}^{12})^{-1}
\end{equation}
\end{theorem}

\subsection{Uniqueness of universal relative twists}\label{ss:uniqueness}
We now show the uniqueness of the twist $J_{\Phi}$ up to a unique gauge transformation.
\begin{theorem}\label{thm:unique-univ-twist}
For any $\Phi\in\mathsf{Assoc}$,
\[
\{J\in\CPDYUAinv{2}\;|\; (\Phi_{\beta})_{J}=\Phi_{\alpha}, J_0=1\}
=\{\gauge{u}{J_{\Phi}}\;|\; u\in\CPDYUAinv{}\}
\]
\end{theorem}
\begin{pf}
Assume $J^{(i)}=1+\sum_{k\geqslant 1}J_k^{(i)}$, $J_k^{(i)}\in\PDYUAinv{2}_k$, $i=1,2$, and
 $(\Phi_{\beta})_{J^{(i)}}=\Phi_{\alpha}$. One checks, by linearisation of \eqref{eq:reltwisteq},
 that $J_1^{(i)}$ is an element in $\PDYUAinv{2}_1$, satisfying $d_H(J_1^{(i)})=0$ and
\[
\mathsf{Alt_2}(J_1)=\mathsf{Alt_2}(r_{\beta}-r_{\alpha})=\wt{r}_{\beta}-\wt{r}_{\alpha}
\]
Up to a gauge, we may assume 
$J^{(i)}_1=\wt{r}_{\beta}-\wt{r}_{\alpha}$, $i=1,2$.
We want to show that there exists an invertible 
$u\in\CPDYUAinv{}$ such that
\begin{equation}\label{eq:twist1}
\gauge{u}{J^{(1)}}=J^{(2)}
\end{equation}
Assume that \eqref{eq:twist1} is true modulo 
$\PDYUAinv{2}_{\geqslant n}$, \ie there exists 
an invertible element $u^{(n-1)}\in\CPDYUAinv{}$ such that
\begin{equation}\label{eq:twist2}
\gaugebis{u^{(n-1)}}{J^{(1)}}= J^{(2)} \mod \PDYUAinv{2}_{\geqslant n}
\end{equation}
Let now $\wt{J}^{(1)}$ be the left--hand side of
\eqref{eq:twist2}, and $\eta\in\PDYUAinv{2}_{n}$ such that
\[
J^{(2)}=\wt{J}^{(1)}+\eta\quad\mod\PDYUAinv{2}_{\geqslant n+1}
\]
One checks that $\wt{J}^{(1)}$ satisfies $(\Phi_{\beta})_{\wt{J}^{(1)}}=\Phi_{\alpha}$ modulo 
$\PDYUAinv{2}_{\geqslant n+1}$.
Comparing with the equation $(\Phi_{\beta})_{J^{(2)}}=\Phi_{\alpha}$ 
modulo $\PDYUAinv{2}_{\geqslant n+1}$, one gets
\[
\eta_{23}+(\eta)_{1,23}-(\eta)_{12,3}-\eta_{12}=0
\]
that is $d_H(\eta)=0$. Therefore, by Proposition \ref{prop:cohoinv}, 
there exist a unique $v\in\PDYUAinv{}_{n}$ and $\mu\in\PDYextinv{2}_n$, such that
$\eta=d_H(v)+\mu$.

We claim $\mu=0$. Then, we may set $u^{(n)}=(1-v)u^{(n-1)}$, and we get
\[
\gaugebis{u^{(n)}}{J^{(1)}}=J^{(2)}\quad\mod\PDYUAinv{2}_{\geqslant n+1}
\] 
There remains to prove the claim. Set $\wt{J}^{(2)}=\gauge{u^{(n)}}{J^{(1)}}$. Then
\begin{equation*}\label{eq:twist7}
\wt{J}^{(2)}= J^{(2)}+\mu\quad\mod\PDYUAinv{2}_{\geqslant n+1}
\end{equation*}
Let $J^{(2)}_{[n+1]}, \wt{J}_{[n+1]}^{(2)}$ be the corresponding truncations. We set
\begin{align*}\label{eq:twist8}
&\xi=(J^{(2)}_{[n+1]})^{23}(J^{(2)}_{[n+1]})^{1,23}\Phi_{\beta}-\Phi_{\alpha}(J^{(2)}_{[n+1]})^{12}
(J^{(2)}_{[n+1]})^{12,3}
\quad\mod\PDYUAinv{3}_{\geqslant n+2}\\
&\wt{\xi}=(\wt{J}^{(2)}_{[n+1]})^{23}(\wt{J}^{(2)}_{[n+1]})^{1,23}\Phi_{\beta}-
\Phi_{\alpha}(\wt{J}^{(2)}_{[n+1]})^{12}
(\wt{J}^{(2)}_{[n+1]})^{12,3}
\quad\mod\PDYUAinv{3}_{\geqslant n+2}
\end{align*}
Since $\wt{J}^{(2)}$ and $J^{(2)}$ are both solutions of $(\Phi_{\beta})_{J}=\Phi_{\alpha}$, 
it follows
\begin{equation*}\label{eq:twist9}
\xi=d_H\left(J_{n+1}^{(2)}\right)\qquad\mbox{and}\qquad \wt{\xi}=d_H\left(\wt{J}_{n+1}^{(2)}
\right)
\end{equation*}
Therefore $d_H\xi=d_H\wt{\xi}=0$ and $\sfAlt{\xi}=\sfAlt{\wt{\xi}}=0$. We then observe that
\begin{equation*}
\wt{\xi}-\xi = f(\mu)
\end{equation*}
where $f(\mu)=A_r^{23}(\mu^{12}+\mu^{13})+\mu^{23}(A_r^{12}+A_r^{13})-
A_r^{12}(\mu^{13}+\mu^{23})-
\mu^{12}(A_r^{13}+A_r^{23})$
and $A_r=\wt{r}_{\beta}-\wt{r}_{\alpha}$. By straightforward computation, one checks
\begin{equation*}
\sfAlt{f(\mu)}=\llbracket \wt{r}_{\beta}-\wt{r}_{\alpha},\mu\rrbracket
\end{equation*}
where $\llbracket,\rrbracket$ is the Schouten bracket from
$\PDYext{2}\to\PDYext{3}$. Therefore $\llbracket \wt{r}_{\beta}-\wt{r}_{\alpha},\mu\rrbracket=0$. 
Since 
$\llbracket\wt{r}_{\beta}-\wt{r}_{\alpha},-\rrbracket=\llbracket\wt{r}_{\beta},-\rrbracket$ on 
$\PDYextinv{2}$, one gets $\llbracket\wt{r}_{\beta},\mu\rrbracket=0$.
It follows from \cite[Prop. 2.2]{e3} that the map $\llbracket \wt{r}_{\beta}, -\rrbracket$ 
has a trivial kernel on $\DYext{2}$ and $\PDYext{2}$. Therefore $\mu=0$, and
the theorem is proved.
\end{pf}

\noindent\remark Theorem \ref{thm:unique-univ-twist} generalises \cite[Thm. 2.1]{e3},
where it is proved for the $\PROP$ $\LBA$. In particular, the uniqueness
of the twist in $\LBA$ can be recovered by applying the forgetful functor
$\PLBA\to\LBA$. Theoren \ref{thm:unique-univ-twist} also generalises
\cite[Thm 6.1]{vtl-3}, where it is proved for a semisimple Lie
algebra.

\section{Lie bialgebras graded by semigroups}
\label{se:sgp}

In this section, we review some basic facts about partial semigroups
and Lie (co--)algebras graded by these.

\subsection{Partial semigroups}\label{ss:partial-semigroup}

A {\it partial semigroup} is a pair $(\mon, \sigma)$, where $\mon$ is a set,
and $\sigma:\mon\times\mon\to\mon$ a partial map defined on a subset
$\dmon\subseteq\mon\times\mon$ such that, for any $\alpha,\beta,\gamma\in\mon$, 
\[\sigma(\sigma(\alpha,\beta),\gamma)=\sigma(\alpha,\sigma(\beta,\gamma))\]
when both sides are defined, that is if $(\alpha,\beta)$, $(\sigma(\alpha,
\beta),\gamma)$, $(\beta,\gamma)$, $(\alpha,\sigma(\beta,\gamma))
\in\dmon$.\\

\noindent\remark 
It is common in the literature (see \eg \cite{evseev}) to assume that the
semigroup law $\sigma$ is strongly associative, \ie for any $\alpha,\beta,
\gamma\in\mon$,
\[(\alpha,\beta), (\sigma(\alpha,\beta),\gamma)\in\dmon
\quad\text{if and only if}\quad(\beta,\gamma), 
(\alpha, \sigma(\beta,\gamma))\in\dmon.\]
This definition is stronger than the one given above, and is not suited for
our purposes, since it does not hold for root systems (cf. \ref{ss:ex-simple-LA}).\footnote
{For example, in the root system of $\sl{4}$, $\alpha_2+\alpha_1$ and
$(\alpha_2+\alpha_1)+\alpha_3$ are defined, but $\alpha_1+\alpha_3$
is not.}

\subsection{Coherence}

Every partial semigroup satisfies the following coherence property.
Let $\sfBr(n)$ be the set of full bracketings on the non--associative
monomial $x_1\cdots x_{n}$. Let $\sigma_b:\mon^n \to \mon$ be the
partial map obtained by composing $\sigma$ along $b$ (\eg $\sigma
_{(\bullet\bullet)\bullet}(\alpha,\beta,\gamma)=\sigma(\sigma(\alpha,
\beta),\gamma)$). Set
\[
\mon^{(n)}=\{\ul{\alpha}\in\mon^n\;|\;\sigma_b(\alpha) \text{ is defined for any $b\in\sfBr(n)$}\}
\]

\begin{proposition}
For any $\ul{\alpha}\in\mon^{(n)}$, and $b,b'\in\sfBr(n)$, $\sigma_b(\ul{\alpha})=\sigma_{b'}(\ul{\alpha})$.
\end{proposition}

\begin{pf}
Let $b,b'\in\sfBr(n)$ two bracketings which differ by an elementary move, \ie there are $i< j
< k < l$ such that, up to a permutation $b\leftrightarrow b'$
\begin{align*}
b\;&=\cdots\big(\big((x_{i+1}\cdots x_{j})(x_{j+1}\cdots x_k)\big)(x_{k+1}\cdots x_l)\big)\cdots\\
b'&=\cdots\big((x_{i+1}\cdots x_{j})\big((x_{j+1}\cdots x_k)(x_{k+1}\cdots x_l)\big)\big)\cdots
\end{align*}
and they agree on everything else. 
Let $\ul{\alpha}\in\mon^{(n)}$, and set $\alpha=\sigma_{b_{ij}}(\alpha_{i+1},\dots,\alpha_j)$,
$\beta=\sigma_{b_{jk}}(\alpha_{j+1},\dots,\alpha_k)$, and 
$\gamma=\sigma_{b_{kl}}(\alpha_{k+1},\dots,\alpha_l)$, where $b_{rs}$ is the restriction
of $b$ and $b'$ to $(x_{r+1}\cdots x_{s})$. By associativity,
$
\sigma(\sigma(\alpha,\beta),\gamma)=
\sigma(\alpha,\sigma(\beta,\gamma))
$ 
so that $\sigma_b(\ul{\alpha})=\sigma_{b'}(\ul{\alpha})$.
Since for any $b,b'\in\sfBr(n)$, there is a sequence $b=b_0, b_1,\dots, b_r=b'$ such that
$b_i, b_{i+1}$ differ by an elementary move, $\sigma_b(\ul{\alpha})=\sigma_{b'}(\ul{\alpha})$.
\end{pf}

\subsection{Morphisms, subsemigroups and saturated subsets}\label{ss:mor sub sat}

Let $\sfS,\sfT$ be partial semigroups. A {\it morphism} $\phi:\sfS\to\sfT$ is a map
such that $(\alpha,\beta)\in\sfS^{(2)}$ if and only if $(\phi(\alpha),\phi(\beta))\in\sfT
^{(2)}$, and $\phi(\sigma_{\mon}(\alpha,\beta))=\sigma_{\sfT}(\phi(\alpha),\phi(\beta))$ 
for any $(\alpha,\beta)\in\dmon$.

Any subset $\mon'\subseteq\mon$ inherits a partial semigroup structure. Namely, 
we denote by $\sft(\mon')$ the semigroup with underlying set $\mon'$, 
\[\sft(\mon')^{(2)}=\{(\alpha,\beta)\in\mon'\times\mon'|\,(\alpha,\beta)\in\dmon\;\text{and}\;
\sigma(\alpha,\beta)\in\mon'\}\]
and semigroup law induced by that of $\mon$. The corresponding embedding
$\sft(\mon')\to\mon$ is a morphism of semigroups if and only if $\mon'$ is a {\it
subsemigroup} of $\mon$ \ie if 
$(\alpha,\beta)\in(\mon'\times\mon')\cap\dmon$ implies $\sigma(\alpha,\beta)\in
\mon'$.

For any $\alpha\in\mon$, set
\[\mon^{(2)}_{\alpha}=\{(\beta,\gamma)\in\dmon\:|\;\sigma(\beta,\gamma)=\alpha\}\]
A subset $\mon'\subseteq\mon$ is {\it saturated} if $\mon^{(2)}_{\alpha}\subseteq\mon'\times \mon'$ 
for any $\alpha\in\mon'$.

A partial semigroup is \emph{commutative} if $\dmon$ is symmetric, \ie $(\alpha,\beta)\in\dmon$ 
if and only if $(\beta,\alpha)\in\dmon$, in which case $\sigma(\alpha,\beta)=\sigma(\beta,\alpha)$.

Henceforth, by semigroup we mean a commutative partial semigroup $(\mon, +)$.

\subsection{$\mon$--graded Lie (co)algebras}\label{ss:mon Lie}

Let $\mon$ be a \commutativesg \partialsg semigroup, and $\N$
a $\sfk$--linear symmetric monoidal category $\N$.
\begin{definition}\hfill
\begin{enumerate}
\item An object $\b$ in $\N$ is {\it $\mon$--graded} if it decomposes
as ${\b=\bigoplus_{\alpha\in\mon}\b_{\alpha}}$. 
\item A morphism $\phi:\b'\to\b$ between $\mon$--graded objects in
$\N$ is {\it homogeneous} if $\phi(\b'_{\alpha})\subseteq\b_{\alpha}$
for any $\alpha\in\mon$.
\end{enumerate}
\end{definition}

\noindent
If $\b\in\N$ is $\mon$--graded, then $\b\ten\b$ is $\mon\times\mon$--graded,
and the subspace
\[\b_{\dmon}=\bigoplus_{(\beta,\gamma)\in\dmon}\b_{\beta}\ten\b_{\gamma}\subseteq\b\ten\b\] 
is naturally $\mon$--graded with
$(\b_{\dmon})_{\alpha}=\bigoplus_{(\beta,\gamma)\in\dmon_{\alpha}}\b_{\beta}\ten\b_{\gamma}$. 
Let $i_{\dmon}:\b_{\dmon}\to\b\ten\b$ and $p_{\dmon}:\b\ten\b\to\b_{\dmon}$ be the canonical
injection and projection respectively, $\dmp{\dmon}=i_{\dmon}\circ p_{\dmon}:\b\ten\b\to\b
\ten\b$ the corresponding idempotent, and set $\cdmp{\dmon}=\id_{\b\ten\b}-\dmp{\dmon}$.

\begin{definition}\hfill
\begin{enumerate}
\item A Lie algebra $(\b,[\,,\,])$ in $\N$ is {\it $\mon$--graded} if $\b$
is $\mon$--graded, $[\;,\;]\circ\cdmp{\dmon}=0$ and $[\,,\,]\circ i_{\dmon}:\b_{\dmon}\to\b$ is homogeneous.
\item A Lie coalgebra $(\b,\delta)$ in $\N$ is {\it $\mon$--graded} 
if $\b$ is $\mon$--graded, $\cdmp{\dmon}\circ\delta=0$ and 
$p_{\dmon}\circ\delta:\b\to\b_{\dmon}$ is homogeneous.
\end{enumerate}
\end{definition}

\subsection{}

Let $\b\in\N$ be an $\mon$--graded object.
For any subset $\mon'\subseteq\mon$, set $\b'=\bigoplus_{\alpha\in\mon'}\b_{\alpha}$ 
and let $i:\b'\to\b$ and $p:\b\to\b'$ be the corresponding injection and projection. The
following is straightforward.
\begin{proposition}\hfill
\begin{enumerate}
\item Assume $\b$ is an $\mon$--graded Lie algebra and set $\mu'=p\circ \mu\circ i\ten i$.  
\begin{itemize}
\item[(a)] If $\mon'$ is a subsemigroup of $\mon$, then $(\b',\mu')$ is an $\mon'$--graded
Lie algebra, and $i:\b'\to\b$ is a morphism of Lie algebras.
\item[(b)] If $\mon'$ is  a saturated subset of $\mon$, then $(\b',\mu')$ is a $\sft(\mon')$--graded
Lie algebra, and $p:\b\to\b'$ is a morphism of Lie algebras. 
\end{itemize}
\item Assume $\b$ is an $\mon$--graded Lie coalgebra, and set $\delta'=p\ten p\circ\delta\circ i$. 
\begin{itemize}
\item[(a)] If $\mon'$ is a subsemigroup of $\mon$, then $(\b',\delta')$ is an $\mon'$--graded 
Lie coalgebra, and $p:\b\to\b'$ is a morphism of Lie coalgebras.
\item[(b)] If $\mon'$ is a saturated subset of $\mon$, then $(\b',\delta')$ is a $\sft(\mon')$--graded 
Lie coalgebra, and $i:\b'\to\b$ is a morphism of Lie coalgebras.
\end{itemize}
\item In particular, if $(\b,\mu,\delta)$ is an $\mon$--graded Lie bialgebra and
$\mon'\subseteq \mon$ a saturated subsemigroup, then $(\b',\mu',\delta')$ is
an $\mon'$--graded Lie bialgebra, and $(\b,\b')$ is a split pair of Lie bialgebras
with respect to $i,p$.
\end{enumerate}
\end{proposition}

\subsection{Example}\label{ss:simple Lie root}

Let $\g$ be a complex semisimple Lie algebra with fixed Borel and
Cartan subalgebras $\g\supset\b\supset\h$ and standard Lie bialgebra
structure (see \S \ref{ss:bil on g}), and let $\sfR_+\subset\h^*$ be the
semigroup of positive roots of $\g$ relative to $\b$ . Let $\sfRo$ be
the semigroup with underlying set $\sfR_+\sqcup\{0\}$, and law extending
that of $\sfR_+$ by an element $0$ such that $\alpha+0=\alpha$ for
any $\alpha\in\sfR_+$, with $0+0$ not defined. Then, $\b$ is graded
as a Lie bialgebra by $\sfRo$, with $\b_0=\h$, and $\b_\alpha=\g_
{\alpha}$, $\alpha\in\sfR_+$.

Let $D$ be the Dynkin diagram of $\g$, $B\subseteq D$ a subdiagram,
and $\sfR_{B,+}\subseteq\sfR_+$ the subset of roots whose support
lies in $B$. $\sfR_{B,+}\sqcup\{0\}$ is a saturated subsemigroup
of $\sfRo$, and the corresponding Lie subbialgebra of $\b$ is $\b_B
=\h\oplus\bigoplus_{\alpha\in\sfR_{B,+}}\g_\alpha$.

\section{Semigroup extensions of \ul{$\mathsf{LBA}$}}\label{s:sgp-ext}

In this section, we introduce the $\PROP$ $\MLBA$ which governs Lie
bialgebras graded by a given \commutativesg \partialsg semigroup
$\mon$. 

\subsection{}\label{ss:mon-ext} 

We construct below an $\mon$--graded version of the $\PROP$s $\LA,
\LCA,\LBA$ by adding a complete family of orthogonal idempotents labeled
by the elements of the \partialsg semigroup $\mon$. We describe in details the
refinement of $\LA$, which is easily adapted to $\LCA$ and $\LBA$.
\subsubsection{}

Let $\wmLA{\mon}$ be the $\PROP$ generated by morphisms $\mu:[2]\to[1]$
and $\dmp{\alpha}:[1]\to[1]$, $\alpha\in\mon$, with relations \eqref{eq:bracket},
\begin{equation}\label{eq:ort-idem}
\dmp{\alpha}\circ\dmp{\beta}=\drc{\alpha,\beta}\cdot\dmp{\alpha}
\end{equation}
for any $\alpha,\beta\in\mon$, and
\begin{equation}\label{eq:idem-bra}
\mu\circ\dmp{\beta}\ten\dmp{\gamma}=
\left\{
\begin{array}{ccc}
\dmp{\beta+\gamma}\circ\mu\circ\dmp{\beta}\ten\dmp{\gamma} & \mbox{if} & (\beta,\gamma)\in\dmon\\[1.1ex]
0 & \mbox{if} & (\beta,\gamma)\not\in\dmon
\end{array}
\right.
\end{equation}

\subsubsection{}

In addition to the orthogonality condition \eqref{eq:ort-idem}, we wish to
impose the completeness relation
\begin{equation}\label{eq:completeness}
\sum_{\alpha\in\mon}\dmp{\alpha}=\id_{[1]}
\end{equation}
and, more generally, $\sum_{\alpha\in\mon^p}\dmp{\ul{\alpha}}=\id_{[p]}$ 
for any $p\in\IN$, where
\[\dmp{\ul{\alpha}}=\dmp{\alpha_1}\otimes\cdots\otimes\dmp{\alpha_p}\in\End([p])
\qquad
\ul{\alpha}=(\alpha_1,\ldots,\alpha_p)\in\mon^p\]
To this end, let $\sfk[\mon^p]\fin$ be the functions on $\mon^p$ with finite
support, regarded as a non--unital algebra irrespective of whether $\mon$ is
finite. Then, $\wmLA {\mon}([p],[q])$ is a $(\sfk[\mon^q]\fin,\sfk [\mon^p]\fin)
$--bimodule, with the functions $\drc{\ul{\alpha}}$, ${\alpha}\in\mon^q$,
and $\drc{\ul{\beta}}$, ${\beta}\in\mon^p$, acting as $\dmp{\ul{\alpha}}\circ-$
and $-\circ\dmp{\ul{\beta}}$ respectively. We denote by $\mLA{\mon}$ the
$\PROP$ with morphisms
\begin{eqnarray*}
\mLA{\mon}([p],[q])
&=&\sfk[\mon^q]\fin\ten_{\sfk[\mon^q]\fin}\wmLA{\mon}([p],[q])
\ctp_{\sfk[\mon^p]\fin}\sfk[\mon^p]
\label{eq:completedmor-2}
\end{eqnarray*}
where the tensor product is completed with respect to the weak topology on
$\sfk[\mon^p]$.\footnote{
A basis of neighborhoods of zero in $\sfk[\mon^p]$ is given by the subsets
\[
U_{\ul{\alpha}_1,\dots, \ul{\alpha}_{r}}\coloneqq \{f\in\sfk[\mon^p]\;|\;\forall i=1,\dots, r,\;\;f(\ul{\alpha}_i)=0\}
\]
Then $\sfk[\mon^p]=\lim_U\sfk[\mon^p]/U$ where the limit runs over the open subsets 
$U\subseteq\sfk[\mon^p]$, and, for any discrete $\sfk[\mon^p]\fin$--module $V$,
\[
\sfk[\mon^p]\ctp_{\sfk[\mon^p]\fin} V=
\lim_U\left(\sfk[\mon^p]/U\ten_{\sfk[\mon^p]\fin} V\right)
\]
} 
Explicitly, one has
\begin{eqnarray*}
\mLA{\mon}([p],[q])
&=&
{\prod_{\ul{\alpha}\in\mon^p}}\bigoplus_{\ul{\beta}\in\mon^q}
\dmp{\ul{\beta}}\circ\wmLA{\mon}([p],[q])\circ\dmp{\ul{\alpha}}\label{eq:completedmor}
\end{eqnarray*}
The composition of morphisms $\mLA{\mon}([p],[q])\ten\mLA{\mon}([q],[r])\to
\mLA{\mon}([p],[r])$ in $\mLA{\mon}$ is induced by that in $\wmLA{\mon}$,
because the multiplication of $f\in\sfk[\mon^q]$ and $g\in\sfk[\mon^q]\fin$
has finite support. The identity on $[p]$ in $\mLA{\mon}$ is precisely the
element $\sum_{\ul{\alpha}\in\mon^p}\dmp{\ul{\alpha}}$.

The following is straightforward.

\begin{lemma}\label{lem:comp mu idem}
In $\mLA{\mon}$, the compatibility condition \eqref{eq:idem-bra} 
is equivalent to
\[\dmp{\alpha}\circ\mu=
\sum_{(\beta',\gamma')\in\dmon_{\alpha}}\mu\circ\dmp{\beta'}\ten\dmp{\gamma'}\]
where $\dmon_{\alpha}=\{(\beta',\gamma')\in\dmon\;|\;\beta'+\gamma'=\alpha\}$.
\end{lemma}

Finally, we denote by $\MLA$ the closure under infinite direct sums
of the Karoubi envelope of $\mLA{\mon}$.

\subsubsection{}

The $\PROP$s $\MLCA$ and $\MLBA$ are obtained similarly. In particular,
we impose the following compatibility condition between the idempotents
$\{\dmp{\alpha}\}_{\alpha\in\mon}$ and the cobracket $\delta:[2]\to[1]$,
\begin{equation}\label{eq:idem-cobra}
\dmp{\beta}\ten\dmp{\gamma}\circ\delta=
\left\{
\begin{array}{ccc}
\dmp{\beta}\ten\dmp{\gamma}\circ\delta\circ\dmp{\beta+\gamma} & \mbox{if} & (\beta,\gamma)\in\dmon\\[1.1ex]
0 & \mbox{if} & (\beta,\gamma)\not\in\dmon
\end{array}
\right.
\end{equation}
Analogously to Proposition \ref{ss:mon-ext}, in $\mLCA{\mon}$ and $\mLBA{\mon}$
this is equivalent to the condition
\[
\delta\circ\dmp{\alpha}=\sum_{(\beta,\gamma)\in\dmon_{\alpha}}\dmp{\beta}\ten\dmp{\gamma}\circ\delta
\]

\subsubsection{}\label{sss:comm-rel-MLBA} 
We observe that, although not strictly necessary, the commutativity of
$\mon$ is a natural requirement in the case of $\PROP$s describing Lie operations.
Namely, one has, for any $\alpha,\beta,\gamma\in\mon$,
\begin{align*}
\drc{\alpha,\beta+\gamma}\cdot\mu\circ\dmp{\beta}\ten\dmp{\gamma}
&=\;\dmp{\alpha}\circ\mu\circ\dmp{\beta}\ten\dmp{\gamma}
=-\dmp{\alpha}\circ\mu\circ\dmp{\gamma}\ten\dmp{\beta}\circ(1\,2)\\
&=-\drc{\alpha,\gamma+\beta}\cdot\mu\circ\dmp{\gamma}\ten\dmp{\beta}\circ(1\,2)
=\drc{\alpha,\gamma+\beta}\cdot\mu\circ\dmp{\beta}\ten\dmp{\gamma}
\end{align*}
and $\mu\circ\dmp{\beta}\ten\dmp{\gamma}=\drc{\beta+\gamma,\gamma+\beta}\cdot\mu\circ\dmp{\beta}\ten\dmp{\gamma}$.
In particular, we see that $\mu\circ\dmp{\beta}\ten\dmp{\gamma}=0$ if either $(\beta,\gamma)$ or 
$(\gamma,\beta)$ are not in $\dmon$, or $\beta+\gamma\neq\gamma+\beta$.\\

\noindent\remark\;
Let $\N$ be a Karoubian, $\sfk$--linear symmetric monoidal category.
There is a canonical isomorphism 
\[\Fun_{\bfb}^\otimes(\mLBA{\mon},\N)\simeq\mathsf{LBA}_{\mon}(\N)\]
where $\mathsf{LBA}_{\mon}(\N)$ is the category of $\mon$--graded Lie
bialgebras in $\N$. 

\subsection{Examples}\label{ss:ex-mon}
\begin{enumerate}
\item For $\mon=\{0\}$, one has $\dmp{0}=\id$ and $\MLBA=\LBA$.
\item Let $\mon=\{0,1\}$ be the semigroup with the addition table
\[%
\begin{array}{|c||c|c|}
\hline
+ & 0 & 1\\
\hline\hline
0 & 0 & 1\\
\hline
1 & 1 & 1\\
\hline
\end{array}
\]
$\dmp{0}$ is a morphism of Lie bialgebras, $\dmp{1}=\id_{[1]}-\dmp{0}$, and
it is immediate to check that $\MLBA=\PLBA$ as described in \ref{ss:plba2}.\footnote
{Note that the equality $\MLBA=\PLBA$ holds only after taking Karoubi envelopes.
Indeed, $\mathsf{LBA}_{\mon}$ is generated by one object, while $\mathsf
{PLBA}$ is generated by the objects $[\a],[\b]$.}
\item More generally, if $\mon=\{0,\ldots,n\}$ with the tropical addition law
$p+q=\max(p,q)$, a module over $\MLBA$ consists of a Lie bialgebra 
$\b$ endowed with a sequence $\b_0\hookrightarrow\cdots\hookrightarrow\b_n=\b$ of 
split inclusions of Lie bialgebras. The Lie subbialgebra $\b_i$ is the direct
sum $\bigoplus_{p=0}^i\Im(\theta_p)$.
\end{enumerate}

\subsection{Morphisms in $\MLBA$}\label{ss:mor-mlba}

The projections
\[\dmp{\ul{\alpha}}=
\dmp{\alpha_1}\ten\cdots\ten\dmp{\alpha_N},\quad
\ul{\alpha}=(\alpha_1,\ldots,\alpha_N)\in\mon^N\]
are a complete family of orthogonal idempotents in $\MLA([N],[N])$ and $\MLCA([N],[N])$.
By construction, there is a natural right (resp. left) action of $\Gamma_{\mon,N}=\SS_N\ltimes\sfk[\mon^N]$ on
$\MLA([N],[q])$ and $\MLCA([p],[N])$, and natural identifications
\begin{equation}\label{eq:LA-fact}
\MLA([N],[q])\simeq\prod_{\ul{\alpha}\in\mon^N}\LA([N],[q])\circ\dmp{\ul{\alpha}}
\end{equation}
and
\begin{equation}
\MLCA([p],[N])\simeq\prod_{\ul{\alpha}\in\mon^N}\dmp{\ul{\alpha}}\circ\LCA([p],[N])
\end{equation}

The description of the morphisms in $\mLBA{\mon}$ is similar to those
in $\LBA$ and $\PLBA$ (cf.  \ref{ss:morplba}), but the commutativity relations
\eqref{eq:idem-bra} and \eqref{eq:idem-cobra},
\[
\forall\,(\beta,\gamma)\not\in\dmon,
\qquad\mu\circ\dmp{\beta}\ten\dmp{\gamma}
=0=
\dmp{\beta}\ten\dmp{\gamma}\circ\delta
\]
require the replacement of the free Lie algebras $\FL{N}$ with the
Lie algebras $\FL{N,\ul{\alpha}}$, $\ul{\alpha}\in\mon^N$ defined as follows. 
As we explained in \ref{ss:morlba}, it is convenient to describe the Lie 
algebras $\FL{N,\ul{\alpha}}$ in terms of labeled binary trees. 
Set $X_{\ul{\alpha}}=\{\alpha_1,\dots,\alpha_N\}$ and let $\mon_{\ul{\alpha}}\subseteq\mon$
be the \partialsg subsemigroup generated by $X_{\ul{\alpha}}$. Then
$\FL{N,\ul{\alpha}}=\T(X_{\ul{\alpha}})/J_{\ul{\alpha}}$ where 
$J_{\ul{\alpha}}$ is the ideal generated by all elements of
the form $[t,t]$, $t\in\T(X_{\ul{\alpha}})$, $[t_1,[t_2,t_3]]+[t_2,[t_3,t_1]]+[t_3,[t_2,t_1]]$,
$t_1,t_2, t_3\in\T(X_{\ul{\alpha}})$, and
$[\alpha_{i_1},[\alpha_{i_2},\cdots,[\alpha_{i_{m-1}},\alpha_{i_m}]]\cdots]$
for any $m\leqslant N$ and $\{i_1,\dots, i_m\}\subseteq\{1,\dots, N\}$ such that
$\alpha_{i_1}+(\alpha_{i_2}+(\cdots+\alpha_{i_m})\cdots)$ is not defined in $\mon_{\ul{\alpha}}$.
We observe in particular that $\FL{N,\ul{\alpha}}$ is an $\mon$--graded Lie algebra,
and
\[
\FL{N,\ul{\alpha}}=
\left\{
\begin{array}{cl}
\FL{N} & \text{if $\mon_{\ul{\alpha}}^{(2)}=\mon_{\ul{\alpha}}\times\mon_{\ul{\alpha}}$}\\
\FL{N}^{\scs\operatorname{ab}} & \text{if $\mon_{\ul{\alpha}}^{(2)}=\emptyset$}
\end{array}
\right.\]
where $\FL{N}^{\scs\operatorname{ab}}$ is
the abelian Lie algebra in $N$ generators.

By \eqref{eq:LA-fact}, there is a surjective map
from $\prod_{\ul{\alpha}\in\mon^N}\homo{\abFL{N}{\ul{\alpha}}{\ten q}}$
to $\mLA{\mon}([N],[q])$.
The injectivity follows easily by application of the realisation functor of 
$\mLA{\mon}$ on the $\mon$--graded Lie algebras $\FL{N,\ul{\alpha}}$, $\ul{\alpha}\in\mon^N$.
We then obtain an isomorphism of right $\sfk[\mon^N]\rtimes\SS_N$--modules,
\begin{equation}
\mLA{\mon}([N],[q])\simeq\prod_{\ul{\alpha}\in\mon^N}\homo{\abFL{N}{\ul{\alpha}}{\ten q}}
\end{equation}
compatible with the left action of $\SS_q$. Through the equivalence 
$\mathsf{LCA}_{\mon}\simeq\mathsf{LA}_{\mon}\op$, we then obtain the isomorphism
of left $\sfk[\mon^N]\rtimes\SS_N$--modules,
\begin{equation}
\mLCA{\mon}([p],[N])\simeq\prod_{\ul{\alpha}\in\mon^N}\homo{\abFL{N}{\ul{\alpha}}{\ten p}}
\end{equation}
compatible with the right action of $\SS_p$.
The following is clear.

\begin{proposition}\label{pr:LBAQ morphisms}\hfill
\begin{enumerate}
\item The embeddings $\MLA,\MLCA\rightarrow\MLBA$ induce an isomorphism
of $(\SS_q,\SS_p)$--bimodules
\[\MLBA([p], [q])\simeq\bigoplus_{N\geqslant0}
\MLCA([p], [N])\ten_{\Gamma_{\mon,N}}\MLA([N],[q])\]
\item There is an isomorphism of $(\SS_q,\SS_p)$--bimodules
\begin{equation*}
\MLBA([p],[q])\simeq\bigoplus_{N\geqslant0}
\left(\prod_{\ul{\alpha}\in\mon^N}\homo{\FL{N,\ul{\alpha}}^{\ten p}}\ten
\homo{\FL{N,\ul{\alpha}}^{\ten q}}\right)_{\SS_N}
\end{equation*}
where the coinvariants are taken \wrt the diagonal action of $\SS_N$.
\item
There are natural isomorphisms
\begin{align*}
\MLBATfin
&\simeq
\bigoplus_{\substack{N\geqslant0}}
\left(\prod_{\ul{\alpha}\in\mon^N}(T\FL{N,\ul{\alpha}}^{\otimes n})_{\delta_N}\ten(T\FL{N,\ul{\alpha}}^{\otimes n})_{\delta_N}\right)_{\SS_N}\\
\MLBASfin
&\simeq
\bigoplus_{\substack{N\geqslant0}}
\left(\prod_{\ul{\alpha}\in\mon^N}(S\FL{N,\ul{\alpha}}^{\otimes n})_{\delta_N}\ten(S\FL{N,\ul{\alpha}}^{\otimes n})_{\delta_N}\right)_{\SS_N}
\end{align*}
\end{enumerate}
\end{proposition}
\noindent
In particular, if $\dmon=\mon\times\mon$, we get
\[\LA_{\mon}([N],[q])\simeq\sfk[\mon^N]\ten\homo{\FL{N}^{\ten q}}
\aand
\mLCA{\mon}([p],[N])\simeq\homo{\FL{N}^{\ten p}}\ten\sfk[\mon^N]\]
This yields the following generalisation of \ref{ss:morplba}
\begin{equation*}
\MLBA([p],[q])\simeq\bigoplus_{N\geqslant0}
\left(\homo{\FL{N}^{\ten p}}\ten\sfk[\mon^N]\ten
\homo{\FL{N}^{\ten q}}\right)_{\SS_N}
\end{equation*}

\subsection{Universal Drinfeld--Yetter modules}\label{ss:CLBAcompletions}

The category $\MDY{\mon}{n}$, $n\geqslant1$, is the colored $\PROP$
generated by $n+1$ objects $\ACDY{1}$ and $\VCDY{k}$, $k=1,
\dots, n$, and morphisms
\begin{gather*}
\mu:\ACDY{2}\to\ACDY{1}\qquad\qquad\delta:\ACDY{1}\to\ACDY{2}\\
\dmp{\alpha}:\ACDY{1}\to\ACDY{1},\quad\alpha\in\mon\\
\pi_k:\ACDY{1}\ten \VCDY{k}\to \VCDY{k}\qquad\qquad
\pi_k^*:\VCDY{k}\to\ACDY{1}\ten \VCDY{k}
\end{gather*}
such that $(\ACDY{1},\mu,\delta, \{\dmp{\alpha}\})$ is an $\MLBA
$--module in $\MDY{\mon}{n}$, and every $(\VCDY{k},\pi_k,\pi_k^*)$ is
a \DYt module over $\ACDY{1}$.

Set
\[
\MUA{\mon}{n}=\pEnd{\MDY{\mon}{n}}{\bigotimes_{k=1}^n\VCDY{k}}
\]
Let $\b$ be an $\mon$--graded Lie bialgebra. Then, for any $n$--tuple 
$\{V_k,\pi_k,\pi_k^*\}_{k=1}^n$ of \DYt modules over $\b$, there is a
realisation functor 
\[\G_{(\b,V_1,\dots, V_n)}:\MDY{\mon}{n}\longrightarrow{\kvect}\]
such that $\pb\mapsto\b$, and 
$\VDY{k}\mapsto V_k$, $k=1,\ldots,n$. 

\begin{proposition}\hfill
\begin{enumerate}
\item
Let $\ff\hspace{-0.1cm}:\DrY{\b}\to{\kvect}$ be the forgetful functor,
and $\DYA{\b}{n}=\sfEnd{\ff^{\boxtimes n}}$. The functors $\G_{(\b,V_1,
\dots, V_n)}$ induce an algebra homomorphism
\[
\MDYrho{\b}{n}:\MUA{\mon}{n}\to\DYA{\b}{n}
\]
\item Set $\FA{N,\ul{\alpha}}=U\FL{N,\ul{\alpha}}$. There is a linear isomorphism
\[
\MUA{\mon}{n}
\simeq
\bigoplus_{\substack{N\geqslant0}} 
\left(\prod_{\ul{\alpha}\in\mon^N}
\homo{\FA{N,\ul{\alpha}}^{\ten n}}\ten\homo{\FA{N,\ul{\alpha}}^{\ten n}}\right)_{\SS_N}\]
In particular, if $\dmon=\mon\times\mon$,
\[\MUA{\mon}{n}
\simeq
\bigoplus_{\substack{N\geqslant0}} 
\left(\homo{\FA{N}^{\ten n}}\ten
\sfk[\mon^N]
\ten\homo{\FA{N}^{\ten n}}\right)_{\SS_N}\]
\item 
Every element in $\MUA{\mon}{n}$ is a linear combination of 
the endomorphisms of $\VDY{1}\ten\cdots\ten\VDY{n}$ given by
\begin{equation}\label{eq:colarch-KM}
\colarch{\ul{N}}{\ul{\alpha}}{\sigma}{\ul{N}'}=
\pi^{(\ul{N})}\circ\dmp{\ul{\alpha}}\ten\id^{\ten n}\circ\;\sigma\ten\id^{\ten n}\circ{\pi^*}^{(\ul{N}^\prime)}
\end{equation}
where $N\geqslant0$, $\ul{N},\ul{N}^\prime\in\IN^n$ are such that $|\ul{N}|=N=|\ul{N}
^\prime|$, $\ul{\alpha}\in\mon^N$, and $\sigma\in\SS_{N,\ul{\alpha}}$, where
\[
\SS_{N,\ul{\alpha}}=\{\sigma\in\SS_N\;|\; (\alpha_i,\alpha_j)\not\in\dmon, 
i<j \Rightarrow \sigma(i)<\sigma(j)\}
\]
\item
Let $\{b_i\}\subset\b$ be a basis and $\{b^i\}\subset\b^*$ the dual basis.
Then one has
\[
\rho^n_{\mon}(\colarch{\ul{N}}{\ul{\alpha}}{\sigma}{\ul{N}'})=
\sum_{\substack{\ul{i}=(i_1,\dots,i_N)\\i_k\in I_{\alpha_k}}} b_{\ul{N}(\ul{i})}\cdot 
b^{{\sigma}(\ul{N}')(\ul{i})}
\]
where $I_{\alpha_k}$ is the set of indices corresponding to the basis of $\b_{\alpha_k}$.
\end{enumerate}
\end{proposition}

\newcommand{\DYLA}[2]{\ul{\mathsf{LM}}_{#1}^{#2}}
\newcommand{\DYLCA}[2]{\ul{\mathsf{LCM}}_{#1}^{#2}}

\begin{pf}
$(1)$ follows by construction. 
The proof of $(2)$ is an easy generalisation of those in \ref{ss:DY action} and \ref{ss:mordy}. 
Namely, we first observe that by normal ordering there is an isomorphism
\begin{equation}\label{eq:norm-ord-DY}
\MDY{\mon}{}\left(\VDY{1},\VDY{1}\right)
\simeq\bigoplus_{N\geqslant0}
\DYLCA{\mon}{}\left(\VDY{1},[N]\ten\VDY{1}\right)
\ten_{\Gamma_{\mon,N}}
\DYLA{\mon}{}\left([N]\ten\VDY{1},\VDY{1}\right)
\end{equation}
where $\DYLA{\mon}{}$ (resp. $\DYLCA{\mon}{}$) is the $\PROP$
generated by an $\mon$--graded Lie algebra object $[1]$ and a $[1]$--module $\VDY{1}$
(resp. an $\mon$--graded Lie coalgebra object $[1]$ and a $[1]$--comodule $\VDY{1}$).  
By normal ordering in $\DYLA{\mon}{}$, we obtain a surjective map
\begin{equation}\label{eq:mor-DY}
\prod_{\ul{\alpha}\in\mon^N}\homo{U\FL{N,\ul{\alpha}}}\to
\DYLA{\mon}{}([N]\ten\VDY{1},\VDY{1})
\end{equation}
which is easily seen to be an isomorphism by considering the action of the Lie algebra
$\FL{N,\ul{\alpha}}$ on $U\FL{N,\ul{\alpha}}$ and the corresponding realisation functor.
In particular, since $\mathsf{LCM}_{\mon}\simeq\mathsf{LM}_{\mon}\op$, 
combining \eqref{eq:norm-ord-DY} and \eqref{eq:mor-DY}, we
obtain a linear isomorphism
\[\MUA{\mon}{}\simeq
\bigoplus_{\substack{N\geqslant0}} 
\left(\prod_{\ul{\alpha}\in\mon^N}
\homo{\FA{N,\ul{\alpha}}}\ten\homo{\FA{N,\ul{\alpha}}}\right)_{\SS_N}\]
The proof of the general case follows by replacing $\VDY{1}$ with $\VDY{1}\ten\cdots\ten\VDY{n}$
and $U\FL{N,\ul{\alpha}}$ with $U\FL{N,\ul{\alpha}}^{\ten n}$. 
$(iii)$ and $(iv)$ are straightforward. 
\end{pf}

\subsection{PBW theorem for $\MUA{\mon}{n}$}\label{ss:CLBAPBW}

As in the case of $\Uguniv$ and $\Ugsplit$, the tower of algebras
$\{\RDYUA{\mon}{n}\}_{n\geqslant1}$ is endowed with face maps $\Delta^n_i:\RDYUA{\mon}{n}\to\RDYUA{\mon}{n+1}$
and degeneration maps $\E_n^{i}:\RDYUA{\mon}{n}\to\RDYUA{\mon}{n-1}$ defining 
a cosimplicial structure. 

Let
\[\sfa^n:\MLBATfin\to
\MDY{\mon}{n}(\ten_{k=1}^n\VDY{k},\ten_{k=1}^n\VDY{k})\]
be the map given on $\phi_{\ul{p},\ul{q}}\in\MLBA(T^\ulp[1],T^\ulq[1])$,
by
\[
\sfa^n(\phi_{\ul{p},\ul{q}})=\pi^{(\ul{p})}\circ\phi_{\ul{p},\ul{q}}\circ{\pi^*}^{(\ul{q})}
\]

As in the case of $\LBA$ and $\PLBA$, the isomorphism $\FA{N,\ul{\alpha}}=
U\FL{N,\ul{\alpha}}\simeq S\FL{N,\ul{\alpha}}$ induces a PBW theorem for $\MUA{\mon}{n}$.

\begin{theorem}\hfill
\begin{enumerate}
\item The following diagram is commutative
\[
\xymatrix@C=0.3cm{
\mbox{$\displaystyle
\MDY{\mon}{n}(\bigotimes_{k=1}^n\VDY{k},\bigotimes_{k=1}^n\VDY{k})
$}
\ar[r] & 
\mbox{$\displaystyle
{\bigoplus_{\substack{N\geqslant0}}}
\left(
\prod_{\ul{\alpha}\in\mon^N}
\homo{\FA{N,\ul{\alpha}}^{\ten n}}\ten\homo{\FA{N,\ul{\alpha}}^{\ten n}}\right)_{\SS_N}
$}
\\
\MLBATfin\ar[u]^{\sfa^n} \ar[r] &
\mbox{$\displaystyle
{\bigoplus_{\substack{N\geqslant0}}}
\left(
\prod_{\ul{\alpha}\in\mon^N}
\homo{T{\FL{N,\ul{\alpha}}}^{\ten n}}\ten\homo{T{\FL{N,\ul{\alpha}}}^{\ten n}} \right)_{\SS_N}
$}
\ar[u]\\
\MLBASfin \ar[u]^\sfsym \ar[r] &
\mbox{$\displaystyle
{\bigoplus_{\substack{N\geqslant0}}}
\left(
\prod_{\ul{\alpha}\in\mon^N}
\homo{\SA{\FL{N,\ul{\alpha}}}^{\ten n}}\ten\homo{\SA{\FL{N,\ul{\alpha}}}^{\ten n}} \right)_{\SS_N}
$}
\ar[u]_{\sfsym\ten\sfsym}
}
\]
\item The map $\sfa^n\circ\sfsym$ is an isomorphism of cosimplicial spaces.
\end{enumerate}
\end{theorem}

\subsection{Hochschild cohomology}\label{ss:CLBAcoho}

The cosimplicial structure on $\{\RDYUA{\mon}{n}\}_{n\geqslant1}$ gives rise to the 
\emph{semigroup} universal Hochschild complex with differential 
$d^n_{\mon}=\sum_{i=0}^{n+1}(-1)^i\Delta^n_i:\RDYUA{\mon}{n}\to\RDYUA{\mon}{n+1}$
The morphisms $\{\MDYrho{\b}{n}\}_{n\geqslant1}$ defined in \ref{ss:CLBAcompletions}
define a chain map between the corresponding Hochschild complexes. As in \ref{th:Hochschild} and
\ref{ss:PLBAcoho}, we get the following

\begin{theorem}\hfill
\begin{enumerate}
\item The map $\sfa^n\circ\sfsym$ induces an isomorphism
\[H^i(\MUA{\mon}{\bullet}, d_{\mon})\cong\bigoplus_{j=0}^i\MLBA(\wedge^{j}[1],\wedge^{i-j}[1])\]
In particular, $H^{0}(\MUA{\mon}{\bullet}, d_{H,\mon})=\sfk$ and $H^{1}(\MUA{\mon}{\bullet}, d_{\mon})=0$.
\item The identification in terms of semigroup Lie algebras $\FL{N,\ul{\alpha}}$ of Proposition \ref{pr:LBAQ morphisms}
yields
\[H^i(\MUA{\mon}{\bullet}, d_H)
\cong
\bigoplus_{N\geqslant 0}
\bigoplus_{j=0}^i
\left(
\prod_{\ul{\alpha}\in\mon^N}
\left(\wedge^j \FL{N,\ul{\alpha}}\right)_{\delta_N}\otimes\left(\wedge^{i-j} \FL{N,\ul{\alpha}}\right)_{\delta_N}\right)_{\SS_N}\]
\end{enumerate}
\end{theorem}

\subsection{Gluing maps}

If $\b$ is a Lie bialgebra with Drinfeld double $\gb$, the standard multiplication
maps $U\gb^{\ten n}\to U\gb^{\ten n-1}$ cannot be lifted to the $\PROP$ic level
since this would imply, for example, that the anti--normally ordered Casimir
$m(r_{21})=\sum_i b^ib_i$ acts on any Drinfeld--Yetter module, which is not the
case if $\b$ is infinite--dimensional. However, the \emph{polarised} multiplication maps $U(\b^*)^{\ten n}\ten U\b
^{\ten n}\to U(\b^*)^{\ten n-1}\ten U\b^{\ten n-1}$ do admit a universal analogue
as maps from $\RDYUA{\mon}{n}$ to $\RDYUA{\mon}{n-1}$.
Their description in terms of associative algebras, under the identification given by
Proposition \ref{ss:CLBAcompletions}, simply corresponds to polarised multiplication
maps $\FA{N,\ul{\alpha}}^{\ten n}\ten\FA{N,\ul{\alpha}}^{\ten n}\to\FA{N,\ul{\alpha}}
^{\ten n-1}\ten\FA{N,\ul{\alpha}}^{\ten n-1}$ (cf. \cite[Prop.~A.1]{e3}).

Their intrinsic description in terms of morphisms in $\MDY{\mon}{n}$, however,
is more involved. Roughly, we consider maps $m^i_n:\RDYUA{\mon}{n}\to 
\RDYUA{\mon}{n-1}$, $i=1,\dots, n-1$, which produce an endomorphism of $\VDY{1}\ten\cdots\ten\VDY{n-1}$
from one of $\VDY{1}\ten\cdots\ten\VDY{n}$ by \emph{gluing} together the modules $\VDY{i}$ 
and $\VDY{i+1}$, as we now describe.
This definition relies on the universal Verma modules $[M]$ and $[M^{\vee}]$ constructed 
in \cite{ek-2} (see also \cite[Sec. 4]{ATL1} for more details). As objects in $\mLBA{\mon}$, $[M]= S[1]$, 
$[M^{\vee}]=\wh{S}[1]\coloneqq \prod_{N\geqslant0}S^N[1]$. They are endowed with a structure of Drinfeld--Yetter 
module over $[1]$, and they satisfy
\begin{equation*}
\pi_{[M]}\circ\id_{[1]}\ten\iota=i_{[1]}
\aand
\id_{[1]}\ten\varepsilon\circ\pi^*_{[M^{\vee}]}=p_{[1]}
\end{equation*}
where $\iota$ and $i_{[1]}$ (resp. $\varepsilon$ and $p_{[1]}$) are the canonical
injections of (resp. projections to) $[0]=S^0[1]$ and $[1]=S^1[1]$.

\subsubsection{}
Let $\DYLA{\mon}{}$ be the $\PROP$
generated by an $\mon$--graded Lie algebra object $[1]$ and $[1]$--modules 
$\VDY{k}$, $k=1,\dots, n$. Let $\pi_k:[M]\ten\VDY{k}\to\VDY{k}$ be the map obtained by
iterations of the action of $[1]$ on $\VDY{k}$. 

\begin{definition}
For any $n\geqslant 2$, $i=1,\dots, n-1$, the $n$th \emph{action--gluing} map
in position $i$
\[
(m_{\DYLA{}{}})^i_n:\DYLA{\mon}{}\left([N]\ten\bigotimes_{k=1}^n\VDY{k},\bigotimes_{k=1}^n\VDY{k}\right)
\to
\DYLA{\mon}{}\left([N]\ten\bigotimes_{k=1}^{n-1}\VDY{k},\bigotimes_{k=1}^{n-1}\VDY{k}\right)
\]
is defined by
\[
(m_{\DYLA{}{}})^i_n\left(\phi_{[N],\VDY{1},\dots,\VDY{n}}\right)=
\pi_i^{(i)}\circ
\phi_{[N],\VDY{1},\dots,\VDY{i-1},[M],\VDY{i},\dots,\VDY{n-1}}\circ
\iota^{(i)}
\]
where 
$\pi_i^{(i)}=\id_{\VDY{[1,i-1]}}\ten\pi_i\ten\id_{\VDY{[i+1,n-1]}}$, 
$\iota^{(i)}=\id_{[N]}\ten\id_{\VDY{[1,i-1]}}\ten\iota\ten\id_{\VDY{[i,n-1]}}$,
and
$\id_{\VDY{[i,j]}}\coloneqq \id_{\VDY{i}\ten\cdots\ten\VDY{j}}$, $i\leqslant j$.
\end{definition}

\subsubsection{}
Similarly, let $\DYLCA{\mon}{n}$ be the $\PROP$
generated by an $\mon$--graded Lie coalgebra object $[1]$ and $[1]$--comodules $\VDY{k}$,
$k=1,\dots, n$. Let $\pi_k^*:\VDY{k}\to[M^{\vee}]\ten\VDY{k}$ be the map obtained by
iterations of the coaction of $[1]$ on $\VDY{k}$. 

\begin{definition}
For any $n\geqslant 2$, $i=1,\dots, n-1$, the $n$th \emph{coaction--gluing} map
in position $i$
\[
(m_{\DYLCA{}{}})^i_n:\DYLCA{\mon}{}\left(\bigotimes_{k=1}^n\VDY{k},
[N]\ten\bigotimes_{k=1}^n\VDY{k}\right)
\to
\DYLCA{\mon}{}\left(\bigotimes_{k=1}^{n-1}\VDY{k},[N]\ten\bigotimes_{k=1}^{n-1}\VDY{k}\right)
\]
is defined by
\[
(m_{\DYLCA{}{}})^i_n\left(\phi_{[N],\VDY{1},\dots,\VDY{n}}\right)=
\varepsilon^{(i)}\circ
\phi_{[N],\VDY{1},\dots,\VDY{i-1},[M],\VDY{i},\dots,\VDY{n-1}}\circ
(\pi_i^*)^{(i)}
\]
where 
$(\pi_i^*)^{(i)}=\id_{\VDY{[1,i-1]}}\ten\pi_i^*\ten\id_{\VDY{[i+1,n-1]}}$, 
$\varepsilon^{(i)}=\id_{[N]}\ten\id_{\VDY{[1,i-1]}}\ten\varepsilon\ten\id_{\VDY{[i,n-1]}}$.
\end{definition}

\subsubsection{}
Set $\VDY{[1,n]}=\bigotimes_{k=1}^n\VDY{k}$.
As we observed in \eqref{eq:norm-ord-DY}, by normal
ordering, the algebra of endomorphisms 
$\MUA{\mon}{n}=\MDY{\mon}{n}\left(\VDY{[1,n]},\VDY{[1,n]}\right)$ is
isomorphic to 
\begin{equation*}
\bigoplus_{N\geqslant0}
\DYLCA{\mon}{}\left(\VDY{[1,n]},[N]\ten\VDY{[1,n]}\right)
\ten_{\Gamma_{\mon,N}}
\DYLA{\mon}{}\left([N]\ten\VDY{[1,n]},\VDY{[1,n]}\right)
\end{equation*}
where $\Gamma_{\mon,N}=\sfk[\mon^N]\rtimes\SS_N$.

\begin{definition}
For any $n\geq1$, $i=1,\dots, n$, the $n$th gluing map in position $i$, $m_n^i:\MUA{\mon}{n}\to\MUA{\mon}{n-1}$,  
is the map induced by $(m_{\DYLCA{}{}})_n^i\ten(m_{\DYLA{}{}})_n^i$.
\end{definition}

\noindent
\remark It should be clear from the description above that the maps $m_n^i$
reduce an element of $\RDYUA{\mon}{n}$ to
one of  $\RDYUA{\mon}{n-1}$ by \emph{gluing} together
$\VDY{i}$ and $\VDY{i+1}$, preserving the order of actions and coactions.
Specifically, the coactions on $\VDY{i}$ occur before (resp. any action on
$\VDY{i}$ occur after) any coaction or action on $\VDY{i+1}$.
As we anticipated, these should not be mistaken with a universal version 
of the multiplication maps $U\g^{\ten n}\to U\g^{\ten n-1}$.
For example, let $\colarch{\ul{N}}{\ul{\alpha}}{\sigma}{\ul{N}'}$ be as in \ref{eq:colarch-KM}
and set $r_{12}^{\alpha}=\colarch{(0,1)}{\alpha}{\id}{(1,0)}$, 
$r_{21}^{\alpha}=\colarch{(1,0)}{\alpha}{\id}{(0,1)}$, and
$\noc{\alpha}=\colarch{1}{\alpha}{\id}{1}$. Then, 
$m^1_2(r^{\alpha}_{12})=\noc{\alpha}=m^1_2(r_{21}^{\alpha})$
or, pictorially, \[\polarizedexample\]

\subsection{The Casimir operator of $\DYUA{}$}\label{ss:comm-rel} 

Recall that $\noc{}=\sarch{1}{\id}\in\DYUA{1}$ is the universal version
of the normally ordered Casimir $\sum_ib_ib^i$.
In $\MUA{\mon}{1}$,  we have $\wt{\rho}_{\mon}(\noc{})=\sum_{\alpha
\in\mon}\noc{\alpha}$.
For any $i=1,\dots, n$, set 
\[
{\noc{}}_i=\id_{\VDY{1}\ten\cdots\ten\VDY{i-1}}
\ten\noc{}\ten
\id_{\VDY{i+1}\ten\cdots\ten\VDY{n}}
\in\MUA{\mon}{n}
\]
\begin{proposition}\label{prop-comm-rel}
The element
$\sum_{i=1}^n(\noc{})_i$ is central in $\MUA{\mon}{n}$.
\end{proposition}

\noindent
\remark Note that in the algebra $\MUA{\mon}{}$ the notions of invariant and central
element stand in an opposite relation than they do for the algebra $U\gb$. Namely, an
invariant element is clearly central, but the opposite in not necessarily true. For example
 the Casimir element $\noc{}$ is central but not invariant.\\

\begin{pf}
The argument below is an easy generalisation of \cite[Prop.~A.1]{e3}. 
Let $\X_n\subseteq\MUA{\mon}{n}$
be the subspace of all elements $X\in\MUA{\mon}{n}$ satisfying
\[
[\noc{1}+\noc{2}+\cdots+\noc{n},X]=0
\]
The following is straightforward.
\begin{enumerate}
\item[(i)] If $X\in\X_n$, then $\Delta_i^n(X)\in\X_{n+1}$, $i=0,1,\dots, n+1$.
\item[(ii)] If $X\in\X_n$, then $X^{\sigma}\in\X_n$ for any $\sigma\in\SS_n$, where
\[X^{\sigma}=\sigma^{-1}\circ X_{\VDY{\sigma(1)}\ten\cdots\ten\VDY{\sigma(n)}}\circ\sigma\]
\item[(iii)] If $X\in\X_n$, then $m^i_n(X)\in\X_{n-1}$ for any $n\geqslant 2$ and $i=1,\dots, n-1$.
\end{enumerate}
To prove $(iii)$ it is enough to observe that, for any $X\in\MUA{\mon}{n}$, 
\[
m_n^i([X,\noc{1}+\noc{2}+\cdots+\noc{n}])=[m_n^i(X),\noc{1}+\noc{2}+\cdots+\noc{n-1}]
\]
Let $\P_N$ be the set of partitions of $\{1,\dots, 2N\}$ of the form $\{i_1,j_1\}\sqcup\cdots\{i_N,j_N\}$.
For any partition $P$ in $\P_N$ and $\ul{\alpha}\in\mon^N$, set
\[
r_{P}^{\ul{\alpha}}=\prod_{k=1}^Nr_{i_k,j_k}^{\alpha_k}
\]
where $r_{i_k,j_k}^{\alpha_k}$ denotes the composition of the coaction on $\VDY{j_k}$,
the idempotent $\dmp{\alpha_k}$, and the action on $\VDY{i_k}$. The morphisms $r_P^{\ul{\alpha}}$
are well--defined, since the elements $r_{i_k,j_k}^{\alpha_k}$, $k=1,\dots, N$, commute in 
$\MUA{\mon}{2n}$.
It follows from $(i), (ii), (iii)$ that $\X_n=\MUA{\mon}{n}$ for any $n\geqslant 1$ if and only if
\[
\{r_{P}^{\ul{\alpha}}\;|\; P\in\P_{n}, \ul{\alpha}\in\mon^n\}\subset\X_{2n}
\]
for all $n\geqslant1$. The result follows from the explicit computation
$[r^{\alpha}_{12}, \noc{1}+\noc{2}]=0$, for any $\alpha\in\mon$.

Namely, for any $\beta\in\mon$, we have
\[
[r^{\alpha}_{12}, (\noc{\beta})_1+(\noc{\beta})_2]=C_{\alpha,\beta}-\sum_{\gamma\in L(\beta)}C_{\alpha,\gamma}
\]
where $L(\beta)=\{\gamma\in\mon\;|\; \alpha+\gamma=\beta\}$ and
\[
C_{\alpha,\beta}=\colarch{(1,1)}{(\beta,\alpha)}{\id}{(2,0)}-\colarch{(1,1)}{(\beta,\alpha)}{(1\,2)}{(2,0)}
+\colarch{(0,2)}{(\alpha,\beta)}{\id}{(1,1)}-\colarch{(0,2)}{(\beta,\alpha)}{(1\,2)}{(1,1)}
\]
Set $A_{\alpha}=\{\beta\in\mon\;|\; L(\beta)\neq\emptyset\}$, 
$B_{\alpha}=\{\beta\in\mon\;|\; (\alpha,\beta)\in\dmon\}$, and 
$L(A_{\alpha})=\bigsqcup_{\beta\in A_{\alpha}}L(\beta)$. It is clear that, if
$\beta\not\in B_{\alpha}$, then $C_{\alpha,\beta}=0$.
Therefore,
\begin{align*}
[r^{\alpha}_{12}, \noc{1}+\noc{2}]
&=
\sum_{\beta\in\mon}\left(C_{\alpha,\beta}-\sum_{\gamma\in L(\beta)}C_{\alpha,\gamma}\right)
=
\sum_{\beta\in B_{\alpha}}C_{\alpha,\beta}-\sum_{\beta\in A_{\alpha}}\sum_{\gamma\in L(\beta)}C_{\alpha,\gamma}\\
&=
\sum_{\beta\in B_{\alpha}}C_{\alpha,\beta}-\sum_{\gamma\in L(A_{\alpha})}C_{\alpha,\gamma}
\end{align*}
The result follows by observing that, if $\beta\in B_{\alpha}$, then $\beta\in L(\alpha+\beta)\subset L(A_{\alpha})$.
Therefore, $B_{\alpha}=L(A_{\alpha})$ and $[r^{\alpha}_{12}, \noc{1}+\noc{2}]=0$.
\end{pf}

\section{Saturated subsemigroups and split pairs}\label{se:sub U_S}

Let $\mon$ be a \partialsg semigroup, and $\MUA{\mon}{n}$ the universal
algebras introduced in \ref{ss:CLBAcompletions}. In this section, we study
the subalgebras of $\MUA{\mon}{n}$ determined by the \partialsg saturated
subsemigroups of $\mon$.

\subsection{Subsemigroups and $\LBA$--modules in $\MLBA$}\label{ss:diagrammaticsubprop}

Recall from \ref{ss:mor sub sat} that a subsemigroup $\mon'\subseteq\mon$
is saturated if $\mon^{(2)}_{\alpha}\subseteq\mon'\times \mon'$ for any $\alpha
\in\mon'$.

\begin{proposition}
Let $\mon'\nts\subseteq\nts\mon$ be a saturated subsemigroup
in $\mon$. Then,
\begin{enumerate}
\item The idempotent $\dmp{\mon'}=\sum_{\alpha\in\mon'}\dmp{\alpha}$
satisfies
\[
\dmp{\mon'}\circ\mu=\mu\circ\dmp{\mon'}\ten\dmp{\mon'}
\aand
\delta\circ\dmp{\mon'}=\dmp{\mon'}\ten\dmp{\mon'}\circ\delta
\]
\item The object $([1]_{\mon'},\mu_{\mon'},\delta_{\mon'})$, where $[1]_{\mon'}=([1],\dmp{\mon'})$,
\[
\mu_{\mon'}=\dmp{\mon'}\circ\mu\circ\dmp{\mon'}\ten\dmp{\mon'}
\aand 
\delta_{\mon'}=\dmp{\mon'}\ten\dmp{\mon'}\circ\delta\circ\dmp{\mon'}\]
is an $\LBA$--module in $\MLBA$ and a Lie subbialgebra of $[1]$. 
\item For any saturated \partialsg subsemigroup $\mon''\subseteq \mon'\nts\subseteq\nts\mon$,
the pair $([1]_{\mon'},[1]_{\mon''})$ is a $\PLBA$--module in $\MLBA$,\ie there is a 
canonical functor
\[\rho_{(\mon',\mon'')}:\PLBA\to\MLBA\]
mapping $[\a]$ to $[1]_{\mon''}$ and $[\b]$ to $[1]_{\mon'}$.
\end{enumerate}
\end{proposition}

\begin{pf}
We prove $(1)$. Then $(2)$ and $(3)$ are obvious consequences. We
have
\[\begin{split}
\dmp{\mon'}\circ\mu
&=
\sum_{\alpha\in\mon'}\dmp{\alpha}\circ\mu
=
\sum_{\substack{\alpha\in\mon'\\(\beta,\gamma)\in\dmon_{\alpha}}}\mu\circ\dmp{\beta}\ten\dmp{\gamma}
=
\sum_{\substack{\alpha\in\mon'\\(\beta,\gamma)\in{\mon'}^{(2)}_\alpha}}\mu\circ\dmp{\beta}\ten\dmp{\gamma}\\
&=
\sum_{(\beta,\gamma)\in{\mon'}^{(2)}}\mu\circ\dmp{\beta}\ten\dmp{\gamma}
=
\mu\circ\dmp{\mon'}\otimes\dmp{\mon'}
\end{split}\]
where the second equality holds by Lemma \ref{lem:comp mu idem},
the third one by saturation of $\mon'$, and the fourth one because $\mon'$
is a \partialsg semigroup.
\end{pf}

\subsection{Semigroup subalgebras}\label{ss:subalg-sg}

Let $\mon'\subseteq\mon$ be a saturated \partialsg subsemigroup, 
$\MDY{\mon'}{n}$, $\MDY{\mon}{n}$ the corresponding $\PROP$s.

For every $n\geqslant1$ and $i=0,1,\dots, n+1$, let
$\D_{i,\mon}^n:\MDY{\mon}{n}\to\MDY{\mon}{n+1}$ and
$\E_n^{i,\mon}:\MDY{\mon}{n}\to\MDY{\mon}{n-1}$
be the $\mon$--analogues of the functors $\D_{i}^n$, $\E_n^{(i)}$
defined in \ref{ss:cosimp-DY}. The following is straightforward.

\begin{proposition}\hfill
\begin{enumerate}
\item The inclusion $\mon'\subset\mon$ induces a faithful functor
$\iota_{\mon',\mon}^n:\MDY{\mon'}{n}\to\MDY{\mon}{n}$ defined as follows: the $\mon'$--graded
Lie bialgebra object $([1]_{\mon'},\mu_{\mon'}, \delta_{\mon'}, \{\dmp{\alpha}\}_{\alpha\in\mon'})$
in $\MDY{\mon'}{n}$ maps to
\[
\left(([1]_{\mon},\dmp{\mon'}), \dmp{\mon'}\circ\mu_{\mon}\circ\dmp{\mon'}\ten\dmp{\mon'},
\dmp{\mon'}\ten\dmp{\mon'}\circ\delta_{\mon}\circ\dmp{\mon'}, \{\dmp{\alpha}\}_{\alpha\in\mon'}\right)
\]
and the Drinfeld--Yetter modules $\VDY{k}$ in $\MDY{\mon'}{n}$ map to their analogues
in $\MDY{\mon}{n}$ restricted to $([1]_{\mon},\dmp{\mon'})$, \ie
\[
(\VDY{k}, \pi_{k,\mon'}, \pi_{k,\mon'}^*)
\mapsto
(\VDY{k}, \pi_{k,\mon}\circ\dmp{\mon'}\ten\id_{\VDY{k}}, 
\dmp{\mon'}\ten\id_{\VDY{k}}\circ\pi^*_{k,\mon})
\]
\item For any $i=0,\dots, n+1$, 
\[\iota_{\mon',\mon}^{n+1}\circ\D_{i,\mon'}^n=
\D_{i,\mon}^n\circ\iota_{\mon',\mon}^n
\aand
\iota_{\mon',\mon}^{n-1}\circ\E^{i,\mon'}_n=
\E^{i,\mon}_n\circ\iota_{\mon',\mon}^{n}\]
\end{enumerate}
\end{proposition}

\noindent
It follows from Theorem \ref{ss:CLBAPBW} that the faithful functor $\iota_{\mon',\mon}^n$ 
induce an injective morphism of algebras $\MUA{\mon'}{n}\to\MUA{\mon}{n}$, which preserves 
the cosimplicial structures induced on $\MUA{\mon'}{n}$ and $\MUA{\mon}{n}$ by the functors
$\D_{i,\cdot}^n,\E^{i,\cdot}_n$. Henceforth, for every subsemigroup
$\mon'$, we will identify the algebras $\MUA{\mon'}{n}$ with their images in 
$\MUA{\mon}{n}$.

\subsection{Subsemigroup invariants}\label{sss:wtzero}

Let $\mon'\subseteq\mon$ be a \partialsg subsemigroup. If $\b$
is an $\mon$--graded Lie bialgebra, the subspace
\begin{equation}\label{eq:b_S}
\b_{\mon'}=\bigoplus_{\alpha\in\mon'}\b_{\alpha}\subseteq\b
\end{equation}
is a Lie subalgebra of $\b$, and it is a Lie subbialgebra if $\mon'$
is saturated. Denote by 
\[\DYA{\b,\bb{\mon'}}{n}=\sfEnd{\DrY{\b}\to\DrY{\bb{\mon'}}}\]
the algebra of endomorphisms of the restriction functor from the
category of \DYt modules over $\b$ to those over $\bb{\mon'}$.

Let $[\bb{\mon'}]=([1],\dmp{\mon'})\in\LBA_\mon$ be the $\LBA
$--module corresponding to $\mon'$, and $\MUA{\mon,\mon'}{n}\coloneqq(\MUA{\mon}{n})^
{[\bb{\mon'}]}\subseteq\MUA{\mon}{n}$ the subalgebra of
$[\bb{\mon'}]$--invariants, \ie the subspace of all $\phi\in\MUA{\mon}{n}$
which commute with the action and coaction of $[\bb{\mon'}]$
on $\VDY{1} \otimes\cdots\otimes\VDY{n}$ (cf. \ref{ss:univ-inv}).

\begin{proposition}
\hfill
\begin{enumerate}
\item For any $\ul{\alpha}\in(\mon\setminus\mon')^N$, the elements
$\colarch{\ul{N}}{\ul{\alpha}}{\sigma}{\ul{N}'}\in\MUA{\mon}{n}$ given
by \eqref{eq:colarch-KM} are invariant under $[\bb{\mon'}]$.
\item If $(\mon')^{(2)}=\emptyset$, then $[\bb{\mon'}]$ is an abelian Lie bialgebra and 
$\MUA{\mon,\mon'}{n}=\MUA{\mon}{n}$.
In particular, the homomorphism
$\MDYrho{\b}{n}:\MUA{\mon}{n}\to\DYA{\b}{n}$ factors through $\DYA{\b,\bb{\mon'}}{n}$.
\end{enumerate}
\end{proposition}

\begin{pf}
$(1)$ It follows from \eqref{eq:idem-bra} and \eqref{eq:idem-cobra} that,
for any $\ul{\alpha}\in\mon^N$ with $\alpha_i\not\in\mon'$, $i=1,\dots, N$,
the elements $\colarch{\ul{N}}{\ul{\alpha}}{\sigma}{\ul{N}'}$
commutes with the action and coaction of $[\bb{\mon'}]$. 
We give the proof of the invariance with respect to the action of $[\bb{\mon'}]$
for the elements $\colarch{N}{\ul{\alpha}}{\sigma}{N}$. The general case is proved
similarly. For any $\beta\in\mon'$, we have
\begin{align*}
&\xy
(0,0)*{\diagA};
\endxy
\\
\end{align*}
\begin{align*}
&
=
\nts\nts\nts\nts\nts\nts
\nts\nts\nts\nts\nts\nts
\xy
(0,0)*{\diagB};
\endxy
+\sum_{i=1}^N
\left(
\nts\nts\nts\nts\nts\nts
\nts\nts\nts\nts\nts\nts
\xy
(0,0)*{\diagC};
\endxy
\right)
\\
\end{align*}
\begin{align*}
&
=
\nts\nts\nts\nts\nts\nts
\nts\nts\nts\nts\nts\nts
\xy
(0,0)*{\diagD};
\endxy
\end{align*}
The identity is obtained by iteration of the compatibility
condition between action and coaction, and by observing that,
since $\mon'$ is saturated and $\alpha_i\not\in\mon'$, one has
\[
\dmp{\alpha_i}\ten\id\circ\delta\circ\dmp{\beta}=0
\aand
\dmp{\beta}\circ\mu\circ\dmp{\alpha_i}\ten\id=0
\]
from \eqref{eq:idem-bra} and \eqref{eq:idem-cobra}. 
Therefore $\colarch{\ul{N}}{\ul{\alpha}}{\sigma}{\ul{N}'}\in(\MUA{\mon}{n})^{[\bb{0}]}$ for any
$\ul{\alpha}\in(\mon\setminus\mon')^N$.

$(2)$ If $(\mon')^{(2)}=\emptyset$, then $[\bb{\mon'}]\subseteq[\b]$ is
an abelian Lie subbialgebra and the same proof works for any $\ul{\alpha}
\in\mon^N$. Therefore every morphisms $\MDYrho{\b}{n}(\colarch{\ul
{N}}{\ul{\alpha}}{\sigma}{\ul{N}'})$ is a morphism in the category of \DYt
$\bb{\mon'}$--modules, and $\MDYrho{\b}{n}$ factors through $\DYA
{\b,\bb{\mon'}}{n}$. 
\end{pf}

\subsection{Semigroup subalgebras of the $\mon$--universal algebra}

Let $\mon''\subseteq\mon'\nts\subseteq\nts\mon$ be saturated \partialsg subsemigroups,
$\b$ an $\mon$--graded Lie bialgebra, and $\b_{\mon''}\subseteq\b_{\mon'}
\subseteq\b_\mon=\b$ the sub Lie bialgebras defined by \eqref{eq:b_S}. For
any $n\geqslant   1$, we denote by $\rho^n_{(\mon',\mon'')}$ the morphism
$\rho^n_{(\mon',\mon'')}:\PDYUA{n}\to\DYA{\b}{n}$
corresponding to the split pair of Lie bialgebras $(\b_{\mon'},\b_{\mon''})$.
The following is clear.

\begin{proposition}\hfill
\begin{enumerate}
\item Let $\wt{\rho}^n_{\mon'}:\DYUA{n}\to\MUA{\mon}{n}$ be the linear map
which is the identity on each $\VDY{k}$, and maps the Lie bialgebra 
$[1]\in\DY{n}$ to $([1],\dmp{\mon'})\in\MDY{\mon}{n}$. In particular,
\[
\wt{\rho}^n_{\mon'}(\rarch{\ul{N}}{\sigma}{\ul{N}^\prime})=
\sum_{\ul{\alpha}\in(\mon')^N}
\colarch{\ul{N}}{\ul{\alpha}}{\sigma}{\ul{N}^\prime}
\]
Then, $\wt{\rho}_{\mon'}$ is an algebra homomorphism and satisfies
\[
\xymatrix{
\DYUA{n} \ar[r]^{\wt{\rho}^n_{\mon'}}  \ar[d]_{\rho^n_{\b_{\mon'}}} \ar@{}|{\circlearrowleft}[dr] & 
\MUA{\mon}{n} \ar[d]^{\rho^n_{\mon}}\\
\U_{\b_{\mon'}}^n\ar[r] & \DYA{\bb{\mon}}{n}
}
\]
\item Let $\wt{\rho}^n_{(\mon',\mon'')}:\PDYUA{n}\to\MUA{\mon}{n}$ be the linear map 
which is the identity on each $\VDY{k}$, and maps the split pair 
$([1], \id), ([1],\pi_0)$ in $\PDY{n}$ to the split pair $([1],\dmp{\mon'}), ([1],\dmp{\mon''})$ in 
$\MDY{\mon}{n}$. In particular,
\[
\wt{\rho}^n_{(\mon',\mon'')}(\colarch{\ul{N}}{\ul{i}}{\sigma}{\ul{N}^\prime})=
\sum_{\ul{\alpha}\in\I_{(\mon',\mon'')}^{\ul{i}}}
\colarch{\ul{N}}{\ul{\alpha}}{\sigma}{\ul{N}^\prime}
\]
where $\ul{\alpha}\in\mon^N$ belongs to $\I_{(\mon',\mon'')}^{\ul{i}}$ iff $\alpha_k\in\mon''$ 
whenever $i_k=0$ and $\alpha_k\in\mon'\setminus\mon''$ otherwise. 
Then, $\wt{\rho}^n_{(\mon',\mon'')}$ is an algebra homomorphism and 
satisfies 
\[
\xymatrix{
\PDYUA{n} \ar[r]^{\wt{\rho}^n_{(\mon',\mon'')}}  \ar[d]_{\rho^n_{(\b_{\mon'},\b_{\mon''})}} \ar@{}|{\circlearrowleft}[dr] & 
\MUA{\mon}{n} \ar[d]^{\rho^n_{\mon}}\\
\U_{\b_{\mon'}}^n\ar[r] & \DYA{\b_\mon}{n}
}
\]
\end{enumerate}
\end{proposition}

\subsection{Universal twists}

Let $\Phi\in\mathsf{Assoc}$ be a fixed associator. For any saturated
\partialsg subsemigroup $\mon'\nts\subseteq\nts\mon$, denote by $\Phi_{\mon'}=
\wt{\rho}^{3}_{\mon'}(\Phi)$ the image of $\Phi$ in $\UU_{\mon}^3$ under the
map $\wt{\rho}^{3}_{\mon'}:\TenUguniv{3}\to\MUA{\mon}{3}$.
 
Let $J\rel\in\UU\osplit^2$ be the universal relative twist constructed in \cite{ATL1}
(see Theorem \ref{th:ATL1}). For any $\mon''\subseteq\mon'\nts\subseteq\nts\mon$, 
set $J_{(\mon',\mon'')}=\wt{\rho}^{2}_{(\mon',\mon'')}(J\rel)\in\UU^2_{\mon}$.
Then, $J_{(\mon',\mon'')}\in\CRDYUA{\mon',\mon''}{2}=(\CRDYUA{\mon'}{2})^{[\bb{\mon''}]}$ and 
it satisfies the relative twist equation
\begin{equation}\label{eq:rel twist eq}
(\Phi_{\mon'})_{J_{(\mon',\mon'')}}=\Phi_{\mon''}
\end{equation}

\begin{theorem}\label{thm:twistrigid}\hfill
\begin{enumerate}
\item If $J_1,J_2\in\UU^2_{\mon',\mon''}$, are solutions of \eqref{eq:rel twist eq}, 
with $(J_i)_0=1$, there is a gauge transformation
$u\in\UU^\times_{\mon',\mon''}$, with $u_0=1$, such that $J_2=\gauge{u}{J_1}$.
\item The gauge transformation $u$ is unique.
\end{enumerate}
\end{theorem}
\begin{pf}
The proof of $(1)$ follows verbatim that of Theorem \ref{thm:unique-univ-twist}.
$(2)$ Assume that $u\in\UU^\times_{\mon',\mon''}$ is such that
\[u_1\cdot u_2\cdot J=J\cdot u_{12}\]
and $u=1$ mod $(\mathfrak{U}_{\mon',\mon''})_{\geqslant n}$. Let
$v\in(\mathfrak{U}_{\mon',\mon''})_{n}$ such that $u=1+v$ mod
$(\mathfrak{U}_{\mon',\mon''})_{\geqslant n+1}$. Taking the component
of degree $n+1$ in the above equation yields
\[
d_H(v)=v_2-v_{12}+v_1=0\]
which by Theorem \ref{ss:CLBAcoho} implies $v=0$.
\end{pf}

\subsection{Orthogonal semigroups}\label{ss:ort-sg}

By definition, two \partialsg subsemigroups $\mon',\mon''\subseteq\mon$
are \emph{orthogonal} if $(\mon'\times\mon'')\cap\mon^{(2)}=\emptyset$. 

\begin{proposition} Let $\mon',\mon''\subseteq\mon$ be
orthogonal saturated \partialsg subsemigroups.
\begin{enumerate}
\item In $\MDY{\mon}{n}$, the action and coaction of $[\bb{\mon'}]$
on each $\VDY{k}$ commute with those of $[\bb{\mon''}]$.
\item Every element in  $\MUA{\mon'}{n}$ commutes with the action and
coaction of $[\bb{\mon''}]$ on any $\VDY{k}$.
In particular, every element in $\MUA{\mon'}{n}$ is $[\bb{\mon''}]$--invariant.
\item In $\MUA{\mon}{n}$, $[\MUA{\mon'}{n},\MUA{\mon''}{n}]=0$.
\item In $\MUA{\mon}{n}$, $\MUA{\mon'\sqcup\mon''}{n}=\MUA{\mon'}{n}\cdot\MUA{\mon''}{n}$.
\end{enumerate} 
\end{proposition}

\begin{pf}
$(1)$ follows from the orthogonaliy of the subsemigroups, since it implies
that, for any $\alpha\in\mon',\beta\in\mon''$, $\mu\circ\dmp{\alpha}\ten
\dmp{\beta}=0=\dmp{\alpha}\ten\dmp{\beta}\circ\delta$. $(2)$ and $(3)$ 
are direct consequences of $(1)$ since every element in $\MUA{\mon'}{n}$
(resp. $\MUA{\mon''}{n}$) is realised as a composition of actions and
coaction of $[\bb{\mon'}]$ (resp. $[\bb{\mon''}]$). Finally, the same
argument shows that $\MUA{\mon'\sqcup\mon''}{n}=\MUA{\mon'}{n}
\cdot\MUA{\mon''}{n}$. Namely, let $\ul{\alpha}\in(\mon'\sqcup\mon'')^N$
and define a partition $I',I''$ of $\{1,\dots, N\}$ such that $\alpha_i\in\mon'$,
$i\in I'$, $\alpha_j\in\mon''$, $j\in I''$. One has
\[
\colarch{\ul{N}}{\ul{\alpha}}{\sigma}{\ul{\wt{N}}}
=
\colarch{\ul{N}'}{\ul{\alpha}'}{\sigma'}{\ul{\wt{N}}'}
\cdot
\colarch{\ul{N}''}{\ul{\alpha}''}{\sigma''}{\ul{\wt{N}}''}
\]
where $\ul{\alpha}'=(\alpha_i)_{i\in I'}$, $\ul{\alpha}''=(\alpha_j)_{j\in I''}$, $\sigma'\in\SS_{|I'|}$,
$\sigma''\in\SS_{|I''|}$ are the restrictions of $\sigma$ to $I'$ and $I''$, respectively, and
similarly for $\ul{N}',\ul{N}''$. Therefore, 
$\colarch{\ul{N}'}{\ul{\alpha}'}{\sigma'}{\ul{\wt{N}}'}\in\MUA{\mon'}{n}$,
$\colarch{\ul{N}''}{\ul{\alpha}''}{\sigma''}{\ul{\wt{N}}''}\in\MUA{\mon''}{n}$, and
$(4)$ follows.
\end{pf}

\section{Diagrams and nested sets}\label{s:diagrams}

We review in this section a number of combinatorial notions associated to a diagram
$D$, in particular the definition of  nested sets on $D$ and their relative version, 
following \cite{DCP2,vtl-4,ATL1-2}.

\subsection{Nested sets on diagrams} \label{ss:diagrams}

A {\it diagram} is an undirected graph $D$ with no multiple edges or
loops. A {\it subdiagram} $B\subseteq D$ is a full subgraph of $D$,
that is, a graph consisting of a (possibly empty) subset of vertices
of $D$, together with all edges of $D$ joining any two elements of it.

Two subdiagrams $B_1,B_2\subseteq D$ are {\it orthogonal}
if they have no vertices in common, and no two vertices $i\in B_1$, $j\in
B_2$ are joined by an edge in $D$. We denote by $B_1\sqcup B_2$ the
disjoint union of orthogonal subdiagrams. Two subdiagrams $B_1,B_2\subseteq D$ 
are {\it compatible} if either one contains the other or they are orthogonal.

A {\it nested set} on $D$ is a collection $\H$ of pairwise compatible,
connected subdiagrams of $D$ which contains the empty set and
$\cc{D}$, where $\cc{D}$ denotes the set of connected components 
of $D$.

Let $\Ns{D}$ be the partial ordered set of nested sets on $D$, ordered by reverse inclusion.
$\Ns{D}$ has a unique maximal element $\cc{D}$ and its minimal elements are the
\emph{maximal nested sets}. We denote the set of maximal nested sets on $D$ by $\Mns
{D}$. It is easy to see that the cardinality of any maximal nested set on $D$ is equal to 
$|D|+1$. Every nested set $\H$ on $D$ is uniquely determined by  a collection $\{\H_i\}_{i=1}
^r$ of nested sets on the connected components $D_i$ of $D$. We therefore obtain canonical
identifications
\[\Ns{D}=\prod_{i=1}^r \Ns{D_i}\qquad\text{and}\qquad\Mns{D}=\prod_{i=1}^r\Mns{D_i}.\]

\subsection{Relative nested sets}\label{ss:rel-ns}

If $B'\subseteq B\subseteq D$ are two subdiagrams of $D$, a nested
set on $B$ {\it relative} to $B'$ is a collection of subdiagrams of $B$, containing
$\cc{B}$ and $\cc{B'}$, in which every element is compatible with, but not properly contained in 
any of the connected components of $B'$.  We denote by $\Ns{B,B'}$ and
$\Mns{B,B'}$, respectively, the collections of nested sets and maximal 
nested sets on $B$ relative to $B'$. In particular, $\Ns{B}=\Ns{B,\emptyset}$
and $\Mns{B}=\Mns{B,\emptyset}$. Relative nested sets are endowed with
the following operations, which preserve maximal nested sets.
\begin{itemize}
\item {\bf Vertical union.}
For any $B''\subseteq B' \subseteq B$, there is an embedding
\begin{equation}\label{eq:union-nes-set}
\cup:\Ns{B,B'}\times\Ns{B',B''}\to\Ns{B,B''},
\end{equation}
given by the union of nested sets. Its image 
$\Nsr{B,B''}{B'}\subseteq\Ns{B,B''}$ is the collection of relative nested sets which
contains $\cc{B'}$.
\item {\bf Vertical decomposition.} Let $B''\subseteq B$ and $\H\in\Ns{B,B''}$. 
If $\cc{B'}\subseteq\H$ and $B''\subseteq B'$, $\H$ is in the image of \eqref{eq:union-nes-set}.
Therefore, there are uniquely defined nested sets $\H_{B''B'}\in\Ns{B',B''}$
and $\H_{B'B}\in\Ns{B,B'}$
\footnote{
More precisely, for any $\H\in\Mns{B,B'''}$ with $\cc{B'},\cc{B''}\in\H$ and $B'''\subseteq B''\subseteq B'$,
we set
\[
\H_{B''B'}=\{C\in\H\;|\; C\subseteq B', C\subsetneq B''\}
\]
}
such that
\[
\H=\H_{B''B'}\cup\H_{B'B};
\]
\item {\bf Orthogonal union.}
For any $B=B_1\sqcup B_2$ and $B'=B_1'\sqcup B_2'$
with $B_1'\subseteq B_1$, $B_2'\subseteq B_2$, there is 
a bijection
\[\Ns{B_1,B_1'}\times\Ns{B_2,B_2'}\to\Ns{B,B'},\]
mapping $(\H_1,\H_2)\mapsto\H_1\cup\H_2$.
\end{itemize}

\section{Diagrammatic semigroups and Lie bialgebras}\label{s:gr-prop}

In this section, we introduce the notion of diagrammatic semigroup.
The corresponding extension of $\mLBA{}$ allows to account for
both the diagrammatic structure of the Borel subalgebra of a complex
semisimple Lie algebra, as well as its root space decomposition.

\subsection{Lax $D$--algebras \cite[\S3]{ATL1-2}}\label{ss:D-KM}\label{ss:D algebra}

Let $D$ be a diagram. A {\it lax $D$--algebra} is the datum of
\begin{itemize}
\item for any $B\subseteq D$, a $\sfk$--algebra $A_B$
\item for any $B'\subseteq B$, a homomorphism $i_{BB'}:A_{B'}\to A_B$
\end{itemize}
such that 
\begin{itemize}
\item for any $B''\subseteq B'\subseteq B$, $i_{BB'}\circ i_{B'B''}=i_{BB''}$
\item for any $B=B'\sqcup B''$, with $B'\perp B''$, $m_B\circ i_{BB'}\otimes i_{BB''}$
is a morphism of algebras $A_{B'}\ten A_{B''}\to A_{B}$, where $m_B$ denotes the
multiplication in $A_B$.
\end{itemize}

A {\it strict morphism} $A^1\to A^2$ of lax $D$--algebras is a collection
of algebra homomorphisms $\phi_B:A^1_B\to A^2_B$ labeled by
the subdiagrams $B\subseteq D$ such that, for any $B'\subseteq
B$, $i^2_{BB'}\circ\phi_{B'}=\phi_{B}\circ i^1_{BB'}$ as morphisms
$A^1_{B'}\to A^2_B$.

\subsection{Diagrammatic Lie bialgebras \cite[\S 5]{ATL1-2}}\label{ss:diag-LBA}

A {\it diagrammatic Lie bialgebra} $\b$ is the datum of
\begin{itemize}
\item a diagram $D$
\item for any $B\subseteq D$, a Lie bialgebra $\b_{B}$
\item for any $B'\subseteq B$, a Lie bialgebra morphism $i_{BB'}
:\b_{B'}\to\b_{B}$
\end{itemize}
such that 
\begin{itemize}
\item for any $B$, $i_{BB}=\id_{\b_B}$
\item for any $B''\subseteq B'\subseteq B$, $i_{BB'}\circ i_{B'B''}=i_{BB''}$
\item for any $B=B'\sqcup B''$ with $B'\perp B''$, $i_{BB'}+i_{BB''}:\b_{B'}
\oplus\b_{B''}\to\b_{B}$ is an isomorphism of Lie bialgebras. 
\end{itemize}
The above properties imply in particular that $\b_{\emptyset}=0$,
and that $U\b$ is a lax $D$--algebra, with $(U\b)_B=U\b_B$.

\subsection{Split diagrammatic Lie bialgebras \cite[\S 5]{ATL1-2}}\label{ss:diag-sLBA}

A {\it split pair} of Lie bialgebras $(\b,\a)$
is the datum of two Lie bialgebras $\a,\b$, together with Lie bialgebra
morphisms $i:\a\to\b$ and $p:\b\to\a$ such that $p\circ i=\id_{\a}$.
These give rise to an embedding $i\oplus p^t:\ga\hookrightarrow\gb$
of the corresponding doubles, which preserves the bracket and the
inner product.

A diagrammatic Lie bialgebra $\b$ is {\it split} if there are Lie bialgebra
morphisms $p_{B'B}:\b_{B}\to\b_{B'}$ for any $B'\subseteq B$, such
that $p_{B'B}\circ i_{BB'}=\id_{\b_{B'}}$, and
\begin{itemize}
\item for any $B$, $p_{BB}=\id_{\b_B}$
\item for any $B''\subseteq B'\subseteq B$, $p_{B''B'}\circ p_{B'B}=p_{B''B}$
\item for any $B=B'\sqcup B''$ with $B'\perp B''$, $p_{B'B}\oplus p_{B''B}:
\b_{B}\to\b_{B'}\oplus\b_{B''}$ is an isomorphism of Lie bialgebras, and
is the inverse of $i_{BB'}+i_{BB''}$.
\end{itemize}

\subsection{Example}\label{ss:ex-simple-LA-diag}

Let $\g$ be a complex semisimple Lie algebra, with Borel and
Cartan subalgebras $\g\supset\b\supset\h$, Dynkin diagram $D$,
Serre generators $\{e_i,f_i,\hcor{i}\}_{i\in D}$, and standard Lie
bialgebra structure (see \ref{ss:bil on g}). Then, $\g$ is a diagrammatic
Lie bialgebra, where, for any $B\subseteq D$, $\g_B\subseteq
\g$ is the subalgebra generated by $\{e_i,f_i,\hcor{i}\}_{i\in B}$.

$\b$ is also a diagrammatic Lie bialgebra with subalgebras
$\b_B=\h_B\oplus\n_B$, where $\h_B\subseteq\h$ is the span
of $\{\hcor{i}\}_{i\in B}$ and $\n_B$ is the nilpotent subalgebra
generated by $\{e_i\}_{i\in B}$. Moreover, the diagrammatic
structure on $\b$ is split as follows. Let $\sfR_+\subset\h^*$
be the set of positive roots of $\g$ relative to $\b$ and, for
any $B\subseteq D$, let $\sfR_{B,+}\subseteq \sfR_+$ be
the subset of roots whose support lies in $B$. In particular,
$\n_B=\bigoplus_{\alpha\in\sfR_{B,+}}\g_\alpha$. Then, for
any $B'\subseteq B$, we have
\[\h_B=\h_{B'}\oplus\h_{B'}^{\perp}\aand\n_B=\n_{B'}\oplus\n_{B'}^{\perp}\]
where $\h_{B'}^\perp=\{t\in\h_B|\,\alpha_i(t)=0,\,i\in B'\}$ and 
$\n_{B'}^{\perp}=\bigoplus_{\alpha\in\sfR_{B,+}\setminus\sfR_{B',+}}\g_{\alpha}$.
The corresponding projections $p_{B'B}:\b_B=\h_B\oplus\n_B
\to\h_{B'}\oplus\n_{B'}=\b_{B'}$ are Lie bialgebra morphisms
and  give rise to a split diagrammatic Lie bialgebra structure
on $\b$.

\subsection{Diagrammatic semigroups}\label{ss:D-monoids}
A {\it diagrammatic semigroup} is a pair $\dsg=(\mon, D)$ where
$D$ is a diagram and $\mon$ a semigroup endowed with a
family of subsets $\mon(B)\subseteq\mon$ indexed by the
subdiagrams of $D$, such that
\begin{itemize}
\item $\mon(B)$ is a saturated subsemigroup in $\mon$
\item $\mon(B')\subseteq\mon(B)$ for any $B'\subseteq B$
\item for any $B'\perp B''$,
\[\mon(B'\sqcup B'')=\mon(B')\sqcup\mon(B'')
\aand
\left(\mon(B')\times\mon(B'')\right)\cap(\mon)^{(2)}=\emptyset\]
\end{itemize}
It follows in particular that $\mon(\emptyset)=\emptyset$. Moreover,
if $\b$ is an $\mon$--graded Lie bialgebra, and $B\subseteq D$ is a
subdiagram, then
\[\b_B\coloneqq \bigoplus_{\alpha\in\mon(B)}\b_{\alpha}\]
is a Lie subbialgebra of $\b$. The following is straightforward. 

\begin{proposition}\label{pr:grad diagr}
Let $\dsg=(\mon,D)$ be a diagrammatic semigroup. Then, every
$\dsg$--graded Lie bialgebra is a split diagrammatic Lie bialgebra.
\end{proposition}

\subsection{Semisimple Lie algebras}\label{ss:ex-simple-LA}

Let $\g$ be a complex semisimple Lie algebra, with Borel
subalgebra $\b$. As pointed out in \ref{ss:simple Lie root},
$\b$ is graded by $\sfRo=\sfR_+\sqcup \{0\}$, where $\sfR
_+$ is the semigroup of positive roots of $\g$ relative to $\b$.
However, the diagrammatic structure of $\b$ given in \ref{ss:ex-simple-LA-diag}
is not encoded by its $\sfRo$--grading via Proposition \ref{pr:grad diagr}.
Indeed, if $B\subsetneq D$, $\b_B=\h_B\oplus\n_B$ does
not correspond to a subset of $\sfRo$ since $\h_B\subsetneq
\h$ is not  a graded component of $\b$.

Note that $\sfR_+$ is a diagrammatic semigroup, with saturated
subsemigroups $\sfR_{B,+}$ given by the set of roots with
support in $B$, but this diagrammatic structure does not extend
to $\sfRo$ since $\sfR_{B,+}$ is not saturated in $\sfRo$. Alternatively,
one can consider the saturated subsemigroups $\sfRo(B)=
\sfR_{B,+}\sqcup\{0\}\subset\sfRo$, but the latter only detect
the Lie subbialgebras $\h\oplus\n_B\supset\b_B$ and, correspondingly,
do not satisfy the orthogonality property $(\sfRo(B)\times\sfRo
(B'))\cap\sfRo^{(2)}=\emptyset$ if $B\perp B'$.

The need to simultaneously account for the diagrammatic
and the $\sfRo$--graded structure of $\b$ motivates the
construction in the following paragraph.

\subsection{Extensions of diagrammatic semigroups and $\PROP$s}\label{ss:0-sgp-ex}

Let $\dsg=(\mon,D)$ be a diagrammatic semigroup, and $\mon_0=\mon
\sqcup\{0\}$ the semigroup which extends $\mon$ with an element $0$
such that $(0,0)\not\in\dmon_0$ and $\alpha+0=\alpha$ for any $\alpha
\in\mon$. If $\alpha\in\mon$ and $B\subseteq D$, we write $\alpha\perp B$
if $\alpha\in\mon(B')$ for some $B'\perp B$.

Let $\dLBA{\dsg}$ be the $\PROP$ generated by a module over
$\mLBA{\mon_0}$,\footnote{That is, by a complete family of orthogonal
idempotents $\dmp{\alpha}:[1]\to[1]$, $\alpha\in\mon_0$, a bracket $\mu:[2]\to[1]$, 
and a cobracket $\delta:[1]\to[2]$, with the same relations described in \ref{ss:mon-ext}.}
and a family of idempotents $\dmp{0,B}:[1]\to[1]$, $B\subseteq D$, such that
$\dmp{0,D}=\dmp{0}$, 
\begin{gather*}
\dmp{0,B'}\circ\dmp{0,B}=
\dmp{0,B'}=
\dmp{0,B}\circ\dmp{0,B'}
\quad\text{for any $B'\subseteq B$}\\[1.1ex]
\dmp{0,B'\sqcup B''}=\dmp{0,B'}+\dmp{0,B''}
\quad\text{for any $B'\perp B''$}
\end{gather*}
and the following additional relations hold
\begin{align*}
\mu\circ\dmp{0,B}\ten\dmp{\alpha}=&
\left\{
\begin{array}{cll}
0 & \text{if $\alpha\perp B$}\\[1.1ex]
\mu\circ\dmp{0}\ten\dmp{\alpha} & \mbox{if }\alpha\in\mon(B)
\end{array}
\right.
\\[1.2ex]
\dmp{0,B}\ten\dmp{\alpha}\circ\delta=&
\left\{
\begin{array}{cll}
0 &  \text{if $\alpha\perp B$}\\[1.1ex]
\dmp{0}\ten\dmp{\alpha}\circ\delta &  \mbox{if }\alpha\in\mon(B)
\end{array}
\right.
\end{align*}
The above relations imply that $\dmp{0,\emptyset}=0$, and that
$\dmp{0,B'}\circ\dmp{0,B''}=0=\dmp{0,B''}\circ\dmp{0,B'}$ for any
$B'\perp B''$ since if $p,q$ are idempotents, $p+q$ is an idempotent
if and only if $pq=0=qp$. 
\begin{proposition}\hfill
\begin{enumerate}
\item For any $B\subseteq D$, set $\dmp{B}=\dmp{0,B}+\sum
_{\alpha\in\mon(B)}\dmp{\alpha}$. Then, $\dmp{B}^2=\dmp{B}$,
\[
\dmp{B}\circ\mu=\mu\circ\dmp{B}\ten\dmp{B}
\aand
\delta\circ\dmp{B}=\dmp{B}\ten\dmp{B}\circ\delta.
\]
\item If $\N$ is a Karoubian, $\sfk$--linear symmetric monoidal category,
any module $\b\in\N$ over $\dLBA{\dsg}$ is a split diagrammatic Lie bialgebra 
with $\b_B=\dmp{B}(\b)$, $B\subseteq D$.
\end{enumerate}
\end{proposition}
\begin{pf}
(1) The relations above imply in particular that $\dmp{0,B}\circ\dmp{\alpha}
=0=\dmp{\alpha}\circ\dmp{0,B}$ for any $B\subseteq D$ and $\alpha\in\mon$,
so that $\dmp{B}^2=\dmp{B}$. Set now $\dmp{B}'=\dmp{0}
+\sum_{\alpha\in\mon(B)}\dmp{\alpha}$. Then, since $\mon(B)\cup\{0\}$
is a saturated subsemigroup in $\mon_0$, 
\[
\mu\circ\dmp{B}\ten\dmp{B}=\mu\circ\dmp{B}'\ten\dmp{B}'=\dmp{B}'\circ\mu=\dmp{B}\circ\mu,
\]
where the first equality follow from $\mu\circ\dmp{0}\otimes\dmp{0}=0
=\mu\circ\dmp{0,B}\otimes\dmp{0,B}$ and $\mu\circ(\theta_0-\theta_
{0,B})\otimes\theta_\alpha=0$ if $\alpha\in\mon(B)$, and the last one
from $\dmp{0}\circ\mu=0=\dmp{0,B}\circ\mu$. Similarly for $\delta$. 
Moreover, for any $B'\perp B''$, $\mu\circ\dmp{B'}\ten\dmp{B''}=0=\dmp{B'}\ten\dmp{B''}\circ\delta$.
It follows that $\dmp{B}(\b), B\subseteq D$, are Lie bialgebras and define a split diagrammatic 
structure on $\b$.
(2) is a direct consequence of (1).
\end{pf}

\subsection{Example}\label{ex:km-sgp}
Let $\g$ be a complex semisimple Lie algebra, and retain the notation
of \ref{ss:ex-simple-LA}. Set $\dsg=(\sfR_+,D)$, where the diagrammatic
structure on the semigroup $\sfR_+$ is given by $\sfR_+(B)=\sfR_{B,+}$.
Then, the Borel subalgebra $\b\subset\g$ is a module over $\dLBA
{\dsg}$, where the idempotents $\dmp{0,B}$  correspond to the projections 
$\h=\h_B\oplus\h_B^{\perp}\to\h_B$. In particular, both the diagrammatic
and $\sfR_0$--graded structure of $\b$ are codified by $\dLBA{\dsg}$.

Henceforth, by abuse of terminology, we say that a Lie bialgebra is 
\emph{$\dsg$--graded} if it is a $\dLBA{\dsg}$--module.

\subsection{Colimit structure of $\dLBA{\dsg}$}\label{ss:0-sgp-ex-colim}

The $\PROP$ $\dLBA{\dsg}$ is not a semigroup extension of $\LBA$
in the sense of \ref{ss:mon-ext}, since the family of idempotents
$\{\dmp{\alpha},\dmp{0,B}\}$ is not labeled by a semigroup, and is
neither complete nor orthogonal. Nevertheless, we show below that
$\dLBA{\dsg}$ is in fact a colimit of semigroup extensions of $\LBA$.

We retain the notation from Section \ref{s:diagrams}.
For any $\H\in\Ns{D}$ and $B\in\H$, denote by $\iH{B}{\H}\subset B$ 
the union of the maximal elements of $\H$ properly contained in $B$.
Let $\mon_{\H}$ be the semigroup with underlying set $\mon\sqcup\{
\zeta_B^\H\}_{B\in\H}$, which contains $\mon$ as subsemigroup and
is such that $\zeta_{B_1}^\H+\zeta_{B_2}^\H$ is undefined for any
$B_1,B_2\in\H$ and
\begin{equation*}
\alpha+\zeta^\H_B=\left\{
\begin{array}{cl}
\mbox{undefined} & \text{if $\alpha\perp B$ or $\alpha\in\mon(\iH{B}{\H})$}\\[1.1ex]
\alpha & \mbox{otherwise}  
\end{array}
\right.
\end{equation*}  

\begin{proposition}
Set $\mLBA{\mon,\H}=\mLBA{\mon_{\H}}$.
\begin{enumerate}
\item For any $\H$, there is a morphism 
$\hcmp{\H}{}:\mLBA{\mon,\H}\to\dLBA{\dsg}$ given by
$\hcmp{\H}{}(\mu)=\mu$, $\hcmp{\H}{}(\delta)=\delta$,
$\hcmp{\H}{}(\dmp{\alpha})=\dmp{\alpha}$, $\alpha\in\mon$, and
\[\hcmp{\H}{}(\theta_{\zeta_B^\H})=\dmp{0,B}-\dmp{0,\iH{B}{\H}}\]
\item For any $\H'\subseteq\H$, there is a morphism
$\hcmp{\H}{\H'}:\mLBA{\mon,\H'}\to\mLBA{\mon, \H}$ given by
$\hcmp{\H}{\H'}(\mu)=\mu$, $\hcmp{\H}{\H'}(\delta)=\delta$,
$\hcmp{\H}{\H'}(\dmp{\alpha})=\dmp{\alpha}$, $\alpha\in\mon$, and
\[\hcmp{\H}{\H'}(\theta_{\zeta_{B'}^{\H'}})=\sum_{B}\theta_{\zeta^\H_B}\]
where the sum ranges over the $B\in\H$ such that $B\subseteq B'$,
and $B\not\subset\iH{B'}{\H'}$.
\item The following holds for any $\H''\subseteq\H'\subseteq\H$,
\[\hcmp{\H}{}\circ\hcmp{\H}{\H'}=\hcmp{\H'}{}
\aand
\hcmp{\H}{\H'}\circ \hcmp{\H'}{\H''}=\hcmp{\H}{\H''}\]
as morphisms $\mLBA{\mon,\H'}\to\dLBA{\dsg}$ and
$\mLBA{\mon,\H''}\to\mLBA{\mon, \H}$ respectively.
\item The $\PROP$ $\dLBA{\dsg}$ is the colimit of
$(\mLBA{\mon,\H}, \hcmp{\H}{\H'})$.
\end{enumerate}
\end{proposition}
\begin{pf}
(1)--(3) are verified by direct inspection. (4) Let $\ul{\sfP}$ be a $\PROP$
endowed with a family of morphisms $p_{\H}:\mLBA{\mon,\H}\to\ul{\sfP}$
such that $p_{\H}\circ\hcmp{\H}{\H'}=p_{\H'}$ for any $\H'\subseteq\H$. 
Then, one can check easily that there is a unique morphism $p:\dLBA
{\dsg}\to\ul{\sfP}$ such that $p\circ\hcmp{\H}=p_{\H}$. Namely,
$p$ is determined by the assignment $p(\dmp{0,B})=p_{\{B,D\}}
(\dmp{\zeta^{\{B,D\}}_{B}})$.
\end{pf}

\subsection{Universal Drinfeld--Yetter modules}\label{ss:dyclba}

Fix henceforth a diagrammatic semigroup $\dsg=(\mon,D)$. The category
$\MDY{\dsg}{n}$, $n\geqslant1$, is the colored $\PROP$ generated by
$n+1$ objects, $\ACDY{1}$ and $\{\VCDY{k}\}_{k=1}^n$, and morphisms
\begin{itemize}
\item $\dmp{\alpha}:\ACDY{1}\to\ACDY{1}$, $\alpha\in\mon$, and 
$\dmp{0,B}:\ACDY{1}\to\ACDY{1}$, $B\subseteq D$
\item $\mu:\ACDY{2}\to\ACDY{1}$, $\delta:\ACDY{1}\to\ACDY{2}$
\item $\pi_k:\ACDY{1}\ten \VCDY{k}\to \VCDY{k}$, $\pi_k^*:\VCDY{k}\to\ACDY{1}\ten \VCDY{k}$
\end{itemize}
such that 
\begin{itemize}
\item $(\ACDY{1}, \dmp{\alpha},\dmp{0,B},\mu,\delta)$ is an $\dLBA{\dsg}$--module 
in $\MDY{\dsg}{n}$
\item every $(\VCDY{k},\pi_k,\pi_k^*)$ is a \DYt module over $\ACDY{1}$
\end{itemize}

\subsection{Universal algebras}\label{ss:colim-univ-alg-sgp}
We set $\MUA{\dsg}{n}=\pEnd{\MDY{\dsg}{n}}{\VCDY{1}\ten\VCDY{2}\ten\cdots\ten\VCDY{n}}$.
For any $\dsg$--graded Lie bialgebra $\b$ and $n$--tuple $\{V_k,\pi_k,\pi_k^*\}_{k=1}^n$ of \DYt $\b
$--modules, there is a canonical realisation functor 
\[\G_{(\b,V_1,\dots, V_n)}:\MDY{\dsg}{n}\longrightarrow{\vect_{\sfk}}\]
sending $[1]\mapsto\b$, and 
$\VDY{k}\mapsto V_k$. As usual, the functors $\G_{(\b,V_1,\dots, V_n)}$ induce 
an algebra homomorphism 
$\MDYrho{\b}{n}:\MUA{\dsg}{n}\to\DYA{\b}{n}$, where 
$\DYA{\b}{n}=\sfEnd{\ff^{\boxtimes n}}$ and $\ff\hspace{-0.1cm}:\DrY{\b}\to{\kvect}$ is 
the forgetful functor.

\subsection{Colimit structures}\label{ss:colim-univ}

As in the case of $\mLBA{\dsg}$, $\mDY{\dsg}{n}$ is a colimit of semigroup extensions
of $\mDY{}{n}$. Namely, let $\H\in\Ns{D}$ and set $\mDY{\mon,\H}{n}=\mDY{\mon_{\H}}{n}$
and $\mUA{\mon,\H}{n}=\mUA{\mon_{\H}}{n}$.
The morphism $h_{\H}:\mLBA{\mon,\H}\to\mLBA{\dsg}$ (resp. $h_{\H\H'}:
\mLBA{\mon,\H'}\to\mLBA{\mon,\H}$, $\H'\subseteq\H$) extend immediately to a 
morphism of $\PROP$s $h^n_{\H}:\mDY{\mon,\H}{n}\to\mDY{\dsg}{n}$ (resp. 
$h^n_{\H\H'}:\mDY{\mon,\H'}{n}\to\mDY{\mon,\H}{n}$) and then to a morphism 
of algebras $h^n_{\H\H'}:\mUA{\mon,\H'}{n}\to\mUA{\mon,\H}{n}$,
compatibly with the inclusion of nested sets. The following corollary of Proposition \ref{ss:0-sgp-ex-colim}
is immediate.

\begin{corollary}\hfill
\begin{enumerate}
\item
The $\PROP$ $\mDY{\dsg}{n}$ is the colimit of the system $(\mDY{\mon,\H}{n}, h^n_{\H\H'})$.
\item
The algebra $\mUA{\dsg}{n}$ is the colimit of the system $(\mUA{\mon,\H}{n}, h^n_{\H\H'})$.
\end{enumerate}
\end{corollary}

\subsection{Hochschild cohomology}\label{ss:coho-dsg}
The cosimplicial structure and the Hochschild differential on $\RDYUA{\dsg}{\bullet}$
are defined as in \ref{ss:cosimplicial}. Relying on the colimit structure described 
above, the results of Section \ref{s:sgp-ext}, in particular the PBW theorem (Thm.~\ref{ss:CLBAPBW})
and computation of the Hochschild cohomology (Thm.~\ref{ss:CLBAcoho}), extend immediately
to the case of a diagrammatic semigroup. More precisely, let $\I_{\dsg}\subseteq\mLBA{\dsg}([1],[1])$ 
be the subset containing the idempotents $\dmp{\alpha}$, $\alpha\in\mon$,
and all iterated products of $\dmp{0,B}$, corresponding to connected subdiagrams 
$B\subseteq D$. For any $\ul{\alpha}\in\I_{\dsg}^N$, the Lie algebras $\FL{N,\ul{\alpha}}$ 
are defined as in \ref{ss:mor-mlba}. We have the following
\begin{theorem}\label{thm:coho-sgp} 
	The Hochschild cohomology of $\MUA{\dsg}{\bullet}$ is described as follows:
	\begin{align*}
	H^i(\MUA{\dsg}{\bullet}, d_H)\cong&\bigoplus_{j=0}^i\mLBA{\dsg}
	(\wedge^{j}[1],\wedge^{i-j}[1])\\
	\cong&
	\bigoplus_{N\geqslant 0}
	\bigoplus_{j=0}^i
	\left(\prod_{\ul{\alpha}\in\I_{\dsg}^N}
	\left(\wedge^j \FL{N,\ul{\alpha}}\right)_{\delta_N}\otimes\left(\wedge^{i-j} \FL{N,
		\ul{\alpha}}\right)_{\delta_N}\right)_{\SS_N}
	\end{align*}
	In particular, $H^{0}(\MUA{\dsg}{\bullet}, d_H)=\sfk$ and $H^{1}(\MUA{\dsg}{\bullet}, d_H)=0$.
\end{theorem}

\subsection{Diagrammatic subalgebras}\label{ss:D-univ-algebras}
For any $B\subseteq D$, set $\dsg(B)=(\mon(B), B)$.
Let $\RDYUA{\dsg, B}{n}\coloneqq \RDYUA{\dsg(B)}{n}$ be the universal algebra in the $\PROP$ 
$\RDY{\dsg,B}{n}\coloneqq \RDY{\dsg(B)}{n}$. For any $B\subseteq B'$, there is a canonical realisation 
functor
\[\G_{\dmp{B}[1],\VCDY{1},\dots,\VCDY{n}}:\RDY{B}{n}\to\RDY{B'}{n}\]
which sends the object $[1]_{B}$ in $\RDY{B}{n}$ to 
the Lie bialgebra object $\dmp{B}[1]_{B'}=([1]_{B'},\dmp{B})$ in $\RDY{\dsg, B'}{n}$ and
induces a homomorphism $\ff_{B'B}:\RDYUA{\dsg,B}{n}\to\RDYUA{\dsg,B'}{n}$.

For any $\dsg$--graded Lie bialgebra  $\b$ 
with diagrammatic Lie subbialgebras $\b_B=\dmp{B}(\b)$ (cf. \ref{ss:0-sgp-ex}),
we set $\DYA{\b, B}{n}=\sfEnd{\ff_B^{\boxtimes n}}$, where $\ff_B:\DrY{\b_B}\to\kvect$ is the forgetful functor,
and we define $\rho^n_{\b,B}:\MUA{\dsg,B}{n}\to\DYA{\b,B}{n}$ as in \ref{ss:dyclba}.

\begin{proposition}\hfill
\begin{enumerate}
\item The subalgebras $\{\RDYUA{\dsg,B}{n}\}_{B\subseteq D}$ define a lax $D$--algebra
structure on $\RDYUA{\dsg}{n}$.
\item The subalgebras $\{\DYA{\b,B}{n}\}_{B\subseteq D}$ defines a lax $D$--algebra structure on $\DYA{\b}{n}$.
\item The collection of homomorphisms $\{\rho_{\b,B}^{n}:\MUA{\dsg,B}{n}\to\DYA{\b,B}{n}\}_{B\subseteq D}$,
defines a strict morphism of $D$--algebras $\rho^n_{\b}:\MUA{\dsg}{n}\to\DYA{\b}{n}$.
\end{enumerate}
\end{proposition}

\begin{pf}
$(1)$ follows from Proposition \ref{ss:ort-sg} $(3)$, since the algebras $\mUA{\dsg, B}{n}$ 
are colimits of semigroup universal algebras. $(2)$ is a consequence of the diagrammatic 
structure of the $\dsg$--graded Lie bialgebra $\b$, and $(3)$ is proved by direct inspection.
\end{pf}

An analogue result holds for the grading completions of $\mUA{\dsg}{n}$ and $\DYA{\b}{n}$, 
which are defined as in \ref{ss:grading} and denoted $\CRDYUA{\dsg}{n}$ and $\CDYA{\b}{n}$.

\subsection{Uniqueness of twists in $\RDYUA{\dsg}{\bullet}$}\label{ss:uniqueness-dsg}
For any subdiagrams $B'\subseteq B\subseteq D$, we denote by
$\CRDYUA{\dsg,B,B'}{n}$ the elements in $\CRDYUA{\dsg,B}{n}$
which are invariant \wrt the action and coaction of $[\b_{B'}]=([1],
\dmp{\mon(B')})$. Relying on the description of the cohomology of $\RDYUA{\dsg}{n}$
given by Theorem \ref{ss:coho-dsg}, we proceed as in Section \ref{se:sub U_S} 
and prove the analogue of Theorem \ref{thm:twistrigid}.

\begin{theorem}\label{thm:twistrigid-sgp}\hfill
	\begin{enumerate}
		\item If $J_1,J_2\in\CRDYUA{\dsg, B, B'}{2}$ are solutions of the relative twist equation $(\Phi_B)_{J_i}=\Phi_{B'}$,
		with $(J_i)_0=1$, there is a gauge transformation
		$u\in\CRDYUA{\dsg, B, B'}{\times}$, with $u_0=1$, such that $J_2=\gauge{u}{J_1}$.
		\item The gauge transformation $u$ is unique.
	\end{enumerate}
\end{theorem}

\section{Universal braided pre--Coxeter structures}
\label{s:univ-w-qC}

Let $\dsg=(\mon,D)$ be a diagrammatic semigroup, and $\MUA{\dsg}{n}$
the universal algebras arising from the $\PROP$ $\mLBA{\dsg}$. We define
in this section the notion of braided pre--Coxeter structure on $\MUA{\dsg}{\bullet}$,
and prove its rigidity. We will prove in Section \ref{se:braided Cox} that such
structures give rise to braided Coxeter categories, as defined in \cite{ATL1-2}.

\subsection{Pre--Coxeter structures on $\CRDYUA{\dsg}{\bullet}$}\label{ss:prop-weak-cox}

Recall that,
for any $B'\subseteq B\subseteq D$, $\CRDYUA{\dsg,B,B'}{n}$ 
denotes the elements in $\CRDYUA{\dsg,B}{n}$
which are invariant \wrt the action and coaction of 
$[\b_{B'}]=([1],\dmp{\mon(B')})$. 

\begin{definition}
A \emph{braided pre--Coxeter structure} $(\Phi_B,J_{\F},\DCPA
{\F}{\G})$ on $\CRDYUA{\dsg}{\bullet}$ consists of the following data. 
\begin{enumerate}
\item 
For any $B\subseteq D$, an associator $\Phi_B\in\CRDYUA{\dsg,B,B}{3}$ (cf. Definition \ref{ss:assoc}),
satisfying the following orthogonal factorisation property. For any $B=B_1\sqcup B_2$,
\[\Phi_B=\Phi_{B_1}\cdot\Phi_{B_2}\]
We set $R_{B}=\exp(\Omega_B/2)\in\CRDYUA{\dsg,B,B}{2}$, and note
that $R_{B_1\sqcup B_2}=R_{B_1}\cdot R_{B_2}$.
\item 
For any $B^\prime\subseteq B\subseteq D$, and maximal nested set
$\F\in\Mns{B,B^\prime}$, a {\it relative twist} $J_{\F}\in
\CRDYUA{\dsg,B,B'}{2}$, that is an (invertible) element such that $(J_{\F})_0=1$
and $\varepsilon_2^1(J_{\F})=1=\varepsilon_2^2(J_{\F})$, 
where $\varepsilon_2^1,\varepsilon_2^2:\CRDYUA{\dsg,B}{2}\to\CRDYUA{\dsg,B}{}$
are the degeneration homomorphisms, which is a solution of the
relative twist equation 
\[
\left(\Phi_B\right)_{J_{\F}}=\Phi_{B^\prime}
\]
where $\Phi_{J}\coloneqq J^{23}J^{1,23}\Phi (J^{12,3})^{-1}(J^{12})^{-1}$ 
(cf. Section \ref{ss:relative-twist} and equations \eqref{eq:reltwisteq}--\eqref{eq:twisted-associator}).
Moreover, the twists $J_{\F}$ satisfy the following factorisation
properties.

\begin{itemize}
\item {\bf Vertical factorisation.}
For any $B''\subseteq B'\subseteq B$,
$\F\in\Mns{B,B^\prime}$ and $\F^\prime\in\Mns{B^\prime,B''}$
\[
J_{\F\cup\F'}=
J_{\F}\cdot J_{\F'}
\]
In particular, $J_\F=1$ if $B'=B$ and $\F$ is the unique element
in $\Mns{B,B}$.
\item {\bf Orthogonal factorisation.}
For any $B=B_1\sqcup B_2$ and $B'=B_1'\sqcup B_2'$, with $B_1^\prime
\subseteq B_1$ and $B_2^\prime\subseteq B_2$, and any orthogonal pair 
$\F=(\F_1,\F_2)\in\Mns{B_1,B_1'}\times\Mns{B_2,B_2'}=\Mns{B, B'}$
\[
J_{\F}=J_{\F_1}\cdot J_{\F_2}
\]
\end{itemize}
\item 
For any $B'\subseteq B\subseteq D$, and $\F,\G\in\Mns{B,B'}$, an invertible
element $\DCPA{\G}{\F}\in\CRDYUA{\dsg,B,B'}{}$, henceforth referred to
as a \emph{De Concini--Procesi associator},
such that $(\DCPA{\G}{\F})_0=1$, $\varepsilon(\DCPA{\G}{\F})=1$, and 
\[J_{\G}=\gaugetris{\DCPA{\G}{\F}}{J_{\F}}\]
The associators $\DCPA{\G}{\F}$ satisfy the following properties.
\begin{itemize}
\item {\bf Transitivity. }
For any $\F,\G,\H\in\Mns{B,B'}$,
\[
\DCPA{\H}{\F}=\DCPA{\H}{\G}\cdot\DCPA{\G}{\F}
\]
In particular, $\Upsilon_{\F\F}=1$ and $\Upsilon_{\G\F}=\Upsilon_{\F\G}^{-1}$.
\item {\bf Vertical factorisation.}
For any $B''\subseteq B'\subseteq B$,
$\F,\G\in\Mns{B,B'}$ and $\F',\G'\in\Mns{B',B''}$,
\[\DCPA{(\G\cup\G')}{(\F\cup\F')}=
\DCPA{\G}{\F}\cdot\DCPA{\G'}{\F'}\]
\item {\bf Orthogonal factorisation.} 
For any $B=B_1\sqcup B_2$ and $B'=B_1'\sqcup B_2'$, with $B_1^\prime
\subseteq B_1$ and $B_2^\prime\subseteq B_2$, and orthogonal pairs
$\F=(\F_1,\F_2)$ and $\G=(\G_1,\G_2)$ in $\Mns{B_1,B_1'}\times
\Mns{B_2,B_2'}=\Mns{B, B'}$
\[
\DCPA{\G}{\F}=\DCPA{\G_1}{\F_1}\cdot\DCPA{\G_2}{\F_2}
\]
\end{itemize}
\end{enumerate}
\end{definition}

\subsection{Twisting of braided pre--Coxeter structures on $\CRDYUA{\dsg}{\bullet}$}\label{ss:twist-weak}

\begin{definition}\hfill
\begin{enumerate}
\item A {\it twist} in $\CRDYUA{\dsg}{\bullet}$ is a pair $(u,F)$ where
\begin{enumerate}
\item $u=\{u_{\F}\}$ is a collection of invertible elements in $\CRDYUA{\dsg,B',B}{}$,
indexed by pairs of subdiagrams $B'\subseteq B$ and a maximal nested set $\F\in\Mns{B,B'}$, 
which satisfy $\varepsilon(u_{\F})=1$, and the following factorisation properties.
\begin{itemize}
\item {\bf Vertical factorisation.}
For any $B''\subseteq B'\subseteq B$, $\F\in\Mns{B,B'}$, and
$\F'\in\Mns{B',B''}$,
\begin{equation}\label{eq:gauge-fac}
u_{\F\cup\F'}=u_{\F}\cdot u_{\F'}
\end{equation}
\item {\bf Orthogonal factorisation.}
For any $B=B_1\sqcup B_2$ and $B'=B_1'\sqcup B_2'$, with
$B_1^\prime\subseteq B_1$ and $B_2^\prime\subseteq B_2$,
and orthogonal pair $\F=(\F_1,\F_2)$ in $\Mns{B_1,B_1'}\times
\Mns{B_2,B_2'}=\Mns{B, B'}$,
\begin{equation}\label{eq:gauge-fac-2}
u_{\F}=u_{\F_1}\cdot u_{\F_2}
\end{equation}
\end{itemize}
\item $F=\{F_B\}$ is a collection of invertible elements of $\CRDYUA{\dsg,B,B}{2}$,
indexed by subdiagrams $B\subseteq D$, which satisfy $\varepsilon_2^1
(F_B)=1=\varepsilon_2^2(F_B)$, are symmetric,\ie $(F_B)_{21}=F_B$ (cf. \ref{ss:notation}), 
$d_H(F_B)_1=0$, and, 
for any $B=B_1\sqcup B_2$, 
\[
F_{B}=F_{B_1}\cdot F_{B_2}
\]
\end{enumerate}
\item The {\it twisting} of a braided pre--Coxeter structure
$\sfCox=(\Phi_B,J_{\F},\DCPA{\F}{\G})$ 
by a twist $(u,F)$  is the braided pre--Coxeter structure 
\[\sfCox_{(u,F)}=((\Phi_B)_{F_B},(J_{\F})_{(u,F)},(\DCPA{\F}{\G})_{u})\]
given by
\begin{eqnarray*}
(\Phi_B)_{F_B}&=&(F_B)_{23}\cdot (F_B)_{1,23}\cdot \Phi_B\cdot (F_B)_{12,3}^{-1}\cdot (F_B)_{12}^{-1}\\
(J_{\F})_{(u,F)}&=&F_{B'}(u_{\F})_1\cdot (u_{\F})_2\cdot J_{\F}\cdot (u_{\F})_{12}^{-1}\cdot F_{B}^{-1}\\
(\DCPA{\F}{\G})_{u}&=&u_{\F}\cdot \DCPA{\F}{\G}\cdot u_{\G}^{-1}\\
\end{eqnarray*}
\end{enumerate}
\end{definition}

\noindent
\remark The twisting of a braided pre--Coxeter structure does not affect the
$R$--matrix $R_B=\exp(\Omega_B/2)$. Specifically, set
\[(R_B)_{F_B}=(F_B)_{21}R_BF_B^{-1}\] 
Since $2\Omega_B=\Delta(\noc{B})-((\noc{B})_1+(\noc{B})_2)$, we have
\begin{align*}
(R_B)_{F_B}&=F_B^{21}\exp(\Delta(\noc{B})/2)\exp(-((\noc{B})_1+(\noc{B})_2)/2)F_B^{-1}\\
&=\exp(\Delta(\noc{B})/2)F_B^{21}\exp(-((\noc{B})_1+(\noc{B})_2)/2)F_B^{-1}\\
&=\exp(\Delta(\noc{B})/2)\exp(-((\noc{B})_1+(\noc{B})_2)/2)=R_B
\end{align*}
where the first identity follows from the invariance of $F_B$, and the second
one from the fact that $(\noc{B})_1+(\noc{B})_2$ is central in $\MUA{\dsg,B}{2}$ 
(Prop. \ref{prop-comm-rel}) and $F_B^{21}=F_B$ by assumption.

Finally, we observe that the conditions $(R_B)_{F_B}=R_B$ and $d_H(F_B)_1
=0$ guarantee that the $2$--jet of the associator is preserved, \ie $((\Phi_B)_
{F_B})_1=0$, and therefore $((\Phi_B)_{F_B})_2=[\Omega_{B,12},\Omega_{B,23}]
/24$ by \cite[Prop. 3.1]{drin-4}.

\subsection{Gauging of twists transformation}\label{ss:gauge-twist}

\begin{definition}\hfill
\begin{enumerate}
\item A {\it gauge} is a collection $a=\{a_B\}$ of
invertible elements $a_B\in\CRDYUA{\dsg,B,B}{}$ indexed 
by subdiagrams $B\subseteq D$ and satisfying, for any
$B=B_1\sqcup B_2$,
\[
a_{B}=a_{B_1}\cdot a_{B_2}
\]
\item The {\it gauging} of a twist $(u,F)$ by $a$ is the
twist $(u_a,F_a)$ given by
\begin{align*}
(u_{\F})_a&=a_{B'}\cdot u_\F\cdot a_{B}^{-1}\\
(F_B)_a		&=(a_B)_1(a_B)_2\cdot F_B\cdot(a_B)_{12}^{-1}
\end{align*}
\end{enumerate}
\end{definition}

\noindent
\begin{remark}
It is easy to see that if $(u,F)$ is a twist, and $a$
a gauge, the twist of a braided pre--Coxeter structure on
$\CRDYUA{\dsg}{\bullet}$ by $(u,F)$ is the same as that by $(u_a,F_a)$.
\end{remark}

\subsection{Uniqueness of braided pre--Coxeter structures on $\CRDYUA{\dsg}{\bullet}$}

\begin{theorem}\label{thm:proprigidity}
Let $\sfCox_k=(\Phi_B^{(k)},J^{(k)}_{\F},\DCPA{\F}{\G}^{(k)})$,
$k=1,2$, be two braided pre--Coxeter structures on $\CRDYUA{\dsg}{\bullet}$.
Then
\begin{enumerate}
\item There exists a twist $(u,F)$ such that $u_0=1, F_0=1$, and
\[
\sfCox_2=(\sfCox_1)_{(u,F)}
\]
\item The twist $(u,F)$ is unique up to a unique gauge $a$.
\end{enumerate}
\end{theorem}

\begin{pf}
We first match the associators. The proof of Drinfeld's uniqueness theorem 
\cite[Prop. 3.12]{drin-4} is easily adapted to $\CRDYUA{\dsg}{\bullet}$. Namely, given 
 $\Phi_B^{(1)},\Phi_B^{(2)}\in\CRDYUA{\dsg,B}{3}$,
there is a symmetric, invariant twist $F_B\in\CRDYUA{\dsg,B}{2}$ such that 
$(\Phi_B^{(1)})_{F_B}=\Phi^{(2)}_B$. In particular, $d_H(F_B)_1=0$. 
$F_B$ is uniquely defined up to multiplication with an element
of the form $(a_B)_1^{-1}(a_B)_2^{-1}(a_B)_{12}$, where $a_B$ belongs to the center of $\CRDYUA{\dsg,B}{}$ and
such that $(a_B)_0=1$. Further, $(R_B)_{F_B}=R_B$, since 
$R_B=\exp(\kappa_B/2)_{12}\cdot\exp(-((\kappa_B)_1+(\kappa_B)_2)/2)$, as we explain in \ref{ss:twist-weak}.

We may therefore assume that $\sfCox_k=(\Phi_B,J^{(k)}_{\F},\DCPA{\F}{\G}^{(k)})$.
We now match the twists. By \ref{thm:twistrigid-sgp}, there exists, for any $\F\in\Mns
{B,B'}$, an invertible element $u_{\F}\in\CRDYUA{\dsg,B,B'}{}$ satisfying
\[
J_{\F}^{(2)}=(u_{\F})_1(u_{\F})_2 J_{\F}^{(1)}(u_{\F})_{12}^{-1}
\]
Moreover, it follows by Theorem \ref{thm:twistrigid-sgp} that the gauge transformation is
unique, and therefore that $u=\{u_{\F}\}$ satisfies the factorisation properties \eqref{eq:gauge-fac}
and \eqref{eq:gauge-fac-2}. Therefore we can assume $\sfCox_k=(\Phi_B,J_{\F},\DCPA{\F}{\G}^{(k)})$
and
\[
\gaugetris{\DCPA{\F}{\G}^{(k)}}{J_{\G}}=J_{\F}
\]
fro $k=1,2$. Again by uniqueness, it follows that $\DCPA{\F}{\G}^{(1)}=\DCPA{\F}{\G}^{(2)}$.
$\sfCox_2$ is therefore a twist of $\sfCox_1$, and the twist
is uniquely defined up to a unique gauge.
\end{pf}

\section{Braided Coxeter categories}\label{se:braided Cox}

In this section, we review the definition of braided (pre--)Coxeter categories
given in \cite{ATL1-2}. We then show that if $\dsg$ is a diagrammatic
semigroup, and $\b$ an $\dsg$--graded Lie bialgebra, a braided pre--Coxeter
structure on the universal algebras $\CRDYUA{\dsg}{\bullet}$ 
endows \DYt modules over the diagrammatic subalgebras of $\b$ with
the structure of a braided pre--Coxeter category.

\subsection{}\label{ss:precox}

Let $D$ be a diagram. A {\it braided pre--Coxeter category} $\C$ of type
$D$ consists of the following data
\begin{itemize}
\item {\bf Diagrammatic categories.} For any subdiagram $B\subseteq
D$, a braided tensor category $\C_B$.
\item {\bf Restriction functors.} For any inclusion $B'\subseteq B$ and
relative \mns $\F\in\Mns{B,B'}$, a tensor functor $F_\F:\C_B\to\C_{B'}$.
\item {\bf \DCP associators.} For any inclusion $B'\subseteq B$ and pair
of relative \mnss $\F,G\in\Mns{B,B'}$, an isomorphism of tensor functors
$\Upsilon_{\G\F}:F_\F\Rightarrow F_\G$.
\end{itemize}

\noindent
This data is assumed to satisfy the following axioms
\begin{itemize}
\item {\bf Normalisation.} If $B\subseteq D$, and $\F$ is the unique
element in $\Mns{B,B}$, then $F_\F=\id_{C_\F}$.
\item {\bf Transitivity.} For any $B'\subseteq B$ and $\F,\G,
\H\in\Mns{B,B'}$, $\Upsilon_{\H\F}=\Upsilon_{\H\G}\circ\Upsilon_{\G\F}$
as isomorphisms $F_\F\Rightarrow F_\H$. In particular, $\Upsilon_{\F\F}
=\id_{F_\F}$ and $\Upsilon_{\G\F}=\Upsilon_{\F\G}^{-1}$.
\item {\bf Vertical factorisation.} For any $B''\subseteq B'\subseteq B$,
$\F\in\Mns{B,B'}$ and $\F'\in\Mns{B',B''}$, the tensor functor $F_{\F'
\cup\F}:\C_B\to\C_{B''}$ is equal to the composition $F_{\F'}\circ
F_\F$. Moreover, for any $\G\in\Mns{B,B'}$ and $\G'\in\Mns{B',B''}$,
the following equality holds
\[\Upsilon_{\G'\cup\G\, \F'\cup\F}=\begin{array}{l}\Upsilon_{\G\F}\\\phantom{00}\circ\\ \Upsilon_{\G'\F'}\end{array}\]
as isomorphisms $F_{\F'\cup\F}=F_{\F'}\circ F_\F\Rightarrow F_{\G'}
\circ F_\G=F_{\G'\cup\G}$.\footnote{In \cite{ATL1-2}, a more
general version of vertical factorisation is considered, where the
equalities $F_{\F'\cup\F}=F_{\F'}\circ F_\F$ and $\Upsilon_{\G'\cup
\G\, \F'\cup\F}=\begin{array}{l}\Upsilon_{\G\F}\\\Upsilon_{\G'\F'}\end{array}$
are only assumed to hold up to coherent isomorphisms. For the purposes
of the present paper, it is sufficient to assume that these isomorphisms
are equalities.}
\end{itemize}

\subsection{Morphisms}\label{ss:precox mor}

Let $\C$, $\C'$ be two braided pre--Coxeter categories of type $D$. 
 A $1$--morphism $H:\C\to\C'$ consists of the following data.
\begin{itemize}
\item For any $B\subseteq D$, a braided tensor functor $H_B: \C_B\to\C'_B$.
\item For any $B'\subseteq B\subseteq D$ and
$\F\in\Mns{B,B'}$, an isomorphism of tensor functors $\gamma_{\F}:F'
_{\F}\circ H_B\Rightarrow H_{B'}\circ F_{\F}$ such that $\DCPA{\G}{\F}
\circ\gamma_{\F}=\gamma_{\G}\circ\DCPA{\G}{\F}'$ as isomorphisms
$F'_{\F}\circ H_B\Rightarrow H_{B'}\circ F_{\G}$.
\end{itemize}
This data is assumed to satisfy the following axioms.
\begin{itemize}
\item {\bf Normalisation.} If $B\subseteq D$ and $\F$ is the unique element 
in $\Mns{B,B}$, so that $F_\F=\id_{\C_B}$ and $F'_\F=\id_{\C'_B}$, then 
$\gamma_{\F}=\id_{H_B}$.
\item {\bf Vertical factorisation.} For any $B''\subseteq B'\subseteq B$,
$\F\in\Mns{B,B'}$ and $\F'\in\Mns{B',B''}$, the following equality holds
\[\gamma_{\F'\cup\F}=\begin{array}{l}\gamma_{\F}\\\phantom{0}\circ\\ \gamma_{\F'}\end{array}\]
as isomorphisms $F'_{\F'\cup\F}\circ H_{B}=F'_{\F'}\circ F'_{\F}\circ H_{B}
\Rightarrow H_{B''}\circ F_{\F'}\circ F_{\F}=H_{B''}\circ F_{\F'\cup \F}$.
\end{itemize}
 
Let $H^1,H^2$ be two $1$--morphisms $\C\to\C'$. A $2$--morphism $v: H^1\Rightarrow H^2$ consists
of the following data.
\begin{itemize}
\item For any $B\subseteq D$, a natural transformation of braided tensor functors 
$v_B:H^1_B\Rightarrow H^2_B$
such that, for any $B'\subseteq B$ and $\F\in\Mns{B,B'}$, 
$\gamma_{\F}\circ v_B=v_{B'}\circ \gamma_{\F}$ as morphisms 
$F'_{\F}\circ H^1_B\Rightarrow H^2_{B'}\circ F_{\F}$.
\end{itemize}
 
\subsection{Generalised braid groups \cite{brieskorn}}

\begin{definition}
A {\it labeling} $\ulm$ of the diagram $D$ is the assignment of an integer
$m_{ij}\in\{2,3,\ldots,\infty\}$ to any pair $i,j$ of distinct vertices of D such
that
\[m_{ij}=m_{ji}\aand
m_{ij}=2\,\,\,\text{if $i$ and $j$ are orthogonal}\]
The \emph{generalised braid group (or Artin group)} corresponding to $D$
and a labeling $\ulm$ is the group $\BDm$ generated by $\sfS_i$, $i\in D$,
with relations\footnote{The group $\BDm$ is called an Artin group in \cite
{brieskorn-saito}. We follow here the terminology of \cite{deligne}.}

\begin{equation}\label{eq:gen-braid}
\underbrace{\sfS_i\cdot \sfS_j\cdot \sfS_i\;\cdots\;}_{m_{ij}}=
\underbrace{\sfS_j\cdot \sfS_i\cdot \sfS_j\;\cdots\;}_{m_{ij}}
\end{equation}
\end{definition}

\subsection{Braided Coxeter categories}

\newcommand {\Fi}{F_{\emptyset i}}
\newcommand {\SC}{S^\C}
\newcommand {\SCp}{S^{\C'}}

Let $\ulm$ be a labeling of $D$. A {\it braided Coxeter category} of type $(D,\ulm)$
is a braided pre--Coxeter category $\left(\C_B,F_\F,\Upsilon_{\G\F}\right)$ of
type $D$ endowed with distinguished isomorphisms $\SC_i\in\Aut(\Fi)$ for any
vertex $i$ of $D$ called {\it local monodromies}. These are assumed to satisfy
the following.

\begin{itemize}
\item {\bf Braid relations.} 
For any $B\subseteq D$, $i\neq j\in B$ and \mnss $\F,\G$ on $B$ with $\{i\}\in\F,
\{j\}\in\G$,  the following holds in $\sfAut{F_\G}$
\[
\underbrace{\mathsf{Ad}\left(\DCPA{\G}{\F}\right)(\SC_i)\cdot \SC_j
\cdot \mathsf{Ad}\left(\DCPA{\G}{\F}\right)(\SC_i)\cdots}_{m_{ij}}=
\underbrace{\SC_j\cdot\mathsf{Ad}\left(\DCPA{\G}{\F}\right)(\SC_i)\cdot \SC_j\cdots}_{m_{ij}}
\]
\item {\bf Coproduct identity.} For any $i\in D$, the following holds in 
$\Aut(\Fi\ten \Fi)$
\begin{equation}\label{eq:coxcoprod}
J_i^{-1}\circ 
\Fi(c_i)\circ\Delta(\SC_i)\circ J_i=
c_{\emptyset}\circ \SC_i\ten \SC_i
\end{equation}
where $J_i$ is the tensor structure on $\Fi$ and $c_i, c_{\emptyset}$ are
the opposite braidings in $\C_i$ and $\C_{\emptyset}$, respectively.\footnote{
In a braided monoidal category with braiding $\beta$, the opposite
braiding is $\beta^{\scs\operatorname{op}}_{X,Y}\coloneqq \beta_{Y,X}^{-1}$.
} In other words, the following diagram is commutative for any $V,W
\in\C_i$, 
\[
\xymatrix{
\Fi(V)\ten \Fi(W) \ar[d]_{J_i^{V,W}} \ar[r]^{S^V_i\ten S^W_i} & \Fi(V)\ten \Fi(W)  \ar[r]^{c_{\emptyset}} & \Fi(W)\ten \Fi(V)
\ar[d]^{J_i^{W,V}}\\
\Fi(V\ten W) \ar[r]_{S_i^{V\otimes W}} & \Fi(V\ten W)  \ar[r]_{\Fi(c_i)} & \Fi(W\ten V)  
}
\]
\end{itemize}

A 1--morphism $\C\to\C'$ of braided Coxeter categories is a 1--morphism
of the underlying braided pre--Coxeter categories which preserves the local
monodromies. That is, it consists of the data $(H_B,\gamma_{B'B})$ defined
in \ref{ss:precox mor} and such that, for any $i\in D$, $\SC_i\circ\gamma_
{\emptyset i}=\gamma_{\emptyset i}\circ\SCp_i$ as isomorphisms $F'_i\circ
H_i\Rightarrow H_\emptyset\circ F_i$.

A 2--morphism $H^1\Rightarrow H^2$ of 1--morphisms $H^1,H^2:\C\to\C'$
of braided Coxeter categories is a 2--morphism of the 1--morphisms
of braided pre--Coxeter categories.

\subsection{Braid group representations}

The axioms of a braided Coxeter category are tailored to produce
natural representations of the generalised braid group. More 
precisely, in \cite[Prop.~3.9]{ATL1-2} we show that a braided Coxeter 
category $\C$ of type $(D,\ulm)$ gives rise to a family of actions 
$\lambda_{\F}:\BDm\to\sfAut{F_\F}$ on the functors $F_\F:\C_D\to\C_\emptyset$ 
labeled by \mnss
on $D$, which are uniquely determined by the conditions
\begin{enumerate}
\item $\lambda_\F(\sfS_i)=\SC_i$ if $\{i\}\in\F$,
\item $\lambda_\G=\sfAd{\DCPA{\G}{\F}}\circ\lambda_\F$.
\end{enumerate}

 \subsection{Deformation Drinfeld--Yetter modules}\label{ss:defDY-2}
We retain the notations from Section \ref{s:gr-prop}.
Let $\dsg$ be a diagrammatic semigroup with underlying
diagram $D$ and $\b$ an $\dsg$--graded Lie bialgebra.
Let $\defDY{\b}$ be the category of Drinfeld--Yetter $\b$--modules in the category of 
topologically free $\sfk{\fml}$--modules, $\hDYA{\b}{n}$ the algebra of endomorphisms of 
the $n$--fold forgetful functor $\ff:\defDY{\b}\to{\khvect}$, and $\RDYUA{\dsg}{n}$ the 
universal algebra introduced in \ref{ss:colim-univ-alg-sgp}. Following the same procedure described
in \ref{ss:defDY}, one can rely on the category $\hDrY{\hextsub{\b}}{\adm}$ of \DYt modules 
over the Lie bialgebra $\hextsub{\b}=(\hext{\b},[\cdot,\cdot],\hbar\delta)$ whose coaction is divisible
by $\hbar$ to obtain a homomorphism $\wh{\rho}_{\b}^{n}:\RDYUA{\dsg}{n}\to\hDYA{\b}{n}$
which naturally extends to $\CRDYUA{\dsg}{n}$. 

\subsection{From universal algebras to Drinfeld--Yetter modules}\label{ss:univ-weak-cox}

\begin{proposition}\label{pr:transfer}
Let $\b$ be an $\dsg$--graded Lie bialgebra. 
\begin{enumerate}
\item A braided pre--Coxeter structure $\Cox{}$ on $\CRDYUA{\dsg}{\bullet}$ canonically
induces  a \multifibered braided pre--Coxeter structure $\Coxb$ on $\{\defDY{\b_B}\}_{B\subseteq D}$.
\item A twist $\T$ in $\CRDYUA{\dsg}{\bullet}$ canonically induces a $1$--isomorphism
$\T(\b):\C(\b)\to\C_{\T}(\b)$, where $\C_{\T}$ denotes the twisted braided pre--Coxeter structure. 
\item A gauge $g$ in $\CRDYUA{\dsg}{\bullet}$ canonically induces a $2$--isomorphism
$g(\b):\T(\b)\Rightarrow\T_{g}(\b)$, where $\T_g$ denotes the gauged twist. 
\end{enumerate}
\end{proposition}

\begin{pf}
$(1)$ Let $\Cox{}=(\Phi_B,J_{\F},\DCPA{\F}{\G})$ be a braided pre--Coxeter
structure on $\CRDYUA{\dsg}{\bullet}$. We show below that the homomorphisms 
$\wh{\rho}_{\b}^{n}$ define a \multifibered braided pre--Coxeter structure
$\Coxb$ with underlying categories $\{\defDY{\b_B}\}_{B\subseteq D}$.
\begin{itemize}
\item {\em Diagrammatic categories.} For any $B\subseteq D$, set
$\Coxb_B=\hDrY{\b_B}{\Phi_B}$, the braided monoidal category of 
topologically free \DYt $\b_B$--modules, with associativity and commutativity 
constraints given by $\wh{\rho}_{B}^{3}(\Phi_B)$ and $\wh{\rho}_{B}^{2}(R_B)$
respectively.
\item {\em Restriction functors.} For any $B'\subseteq B$ and
$\F\in\Mns{B,B'}$, let $F_{\F}^{\C}$ be the standard restriction functor 
$\Res_{B'B}=\Res_{\b_{B'},\b_{B}}$ with tensor structure $\wh{\rho}_{BB'}^{2}(J_{\F})$. 
\item {\em De Concini--Procesi associators.} For any $\F,\G\in\Mns{B,B'}$, let $\Upsilon_{\G\F}^{\C}:
F_{\F}^{\C}\Rightarrow F_{\G}^{\C}$ be the tensor isomorphism defined by $\wh{\rho}_{BB'}(\DCPA{\G}{\F})$.
\end{itemize}

We now show that this datum satisfies the properties required in \ref{ss:precox}.

\begin{itemize}
\item {\em Normalisation.} If $B\subseteq D$, and $\F$ is the unique
element in $\Mns{B,B}$, then, by the vertical factorisation property of the relative twists in 
$\CRDYUA{\dsg}{\bullet}$ (cf. Definition \ref{ss:prop-weak-cox}), $J_{\F}=J_{\F}\cdot J_{\F}$.
In particular, $J_{\F}=1$ and $F_{\F}=\id_{\C(\b)_B}$.
\item {\em Transitivity.} This follows from the horizontal factorisation of
De Concini--Procesi associators in $\CRDYUA{\dsg}{\bullet}$.
\item {\em Vertical factorisation.} This follows from the vertical factorisation
of the relative twists and De Concini--Procesi associators in $\CRDYUA{\dsg}{\bullet}$.
\end{itemize}

$(2)$ Let $\T=(u,F)$ be a twist in $\CRDYUA{\dsg}{\bullet}$ and $\C'=\C_{(u,F)}$ the twisted
pre--Coxeter structure (cf. \ref{ss:twist-weak}). Define a $1$--isomorphism 
$\T(\b)=(H_{B}^{\T}, \gamma_{\F}^{\T}): \C(\b)\to\C'(\b)$ as follows.
\begin{itemize}
\item For any $B\subseteq D$, we denote by $H_B$ the identity
functor on $\C(\b)_B$ endowed with the tensor structure $\wh{\rho}_B^{2}(F_{B})$. In particular, 
it follows immediately from Definition \ref{ss:twist-weak} that $H_B$ is a braided tensor equivalence
$\C(\b)_B\to\C'(\b)_B$. 
\item For any $B'\subseteq B\subseteq D$ and $\F\in\Mns{B,B'}$, 
we denote by $\gamma_{\F}^{\T}$ the natural isomorphism $F^{\C'}_{\F}\circ H^{\T}_B\Rightarrow 
H^{\T}_{B'}\circ F^{\C}_{\F}\circ H$ induced by $\wh{\rho}_{BB'}(u_{\F})$. Therefore, by definition of $u$, 
$\gamma_{\F}^{\T}$ is a well--defined isomorphism of tensor functors satisfying the vertical 
factorisation property.
\end{itemize}

$(3)$ Finally, let $g$ be a gauge in $\CRDYUA{\dsg}{\bullet}$ and $\T'=\T_g$ the gauged twist
(cf. \ref{ss:gauge-twist}). Then, we define a $2$--isomorphism $g(\b):\T(\b)\Rightarrow\T'(\b)$ 
as follows. For any $B\subseteq D$, we denote by $v^g_B$ the isomorphism of braided tensor functors
$H^{\T}_{B}\Rightarrow H^{\T'}_{B}$ given by $\wh{\rho}_B(g_B)$. Then, it follows from the definition of $g$ 
that $\gamma^{\T'}_{\F}\circ v^g_{B}=v^g_{B'}\circ\gamma^{\T}_{\F}$.
\end{pf}

\subsection{Universal braided pre--Coxeter structures}

\begin{definition}
Let $\b$ be an $\dsg$--graded  Lie bialgebra. A braided pre--Coxeter
structure (resp. $1$--morphism, $2$--morphism) on $\{\defDY{\b_B}\}_{B\subseteq D}$ 
is \emph{universal} if it is induced by a braided pre--Coxeter
structure (resp. twist, gauge) on $\CRDYUA{\dsg}{\bullet}$
via Proposition \ref{pr:transfer}.
\end{definition}

The following is a direct consequence of Theorem \ref{thm:proprigidity}.

\begin{theorem}\label{th:rigidity}
Let $\b$ be an $\dsg$--graded Lie bialgebra, and $\C_1,\C_2$
two universal braided pre--Coxeter structures with diagrammatic categories 
$\{\defDY{\b_B}\}_{B\subseteq D}$. Then, there is a universal $1$--isomorphism 
$\C_1\to\C_2$, which is unique up to a unique universal $2$--isomorphism.
\end{theorem}

In \cite[Thm. 9.1]{ATL1-2}, we show that, for any $\dsg$--graded
Lie bialgebra $\b$, there is a canonical universal braided pre--Coxeter
structure on Drinfeld--Yetter modules. Although we do not need this result
for the purposes of this paper, we observe that, combined with
the uniqueness result above, this implies the following.

\begin{corollary}
Let $\b$ be an $\dsg$--graded Lie bialgebra. Then, there exists
an essentially unique universal braided pre--Coxeter structure on the
categories of deformation \DYt modules.
\end{corollary}

\section{\nqcs structures and Kac--Moody algebras}\label{s:km-rigidity}

In this section, we consider the diagrammatic semigroup of
positive roots of a symmetrisable \KM algebra $\g$ with negative
Borel subalgebra $\b$. We then rely on the results of Sections
\ref{s:univ-w-qC}--\ref{se:braided Cox} to prove the uniqueness
of braided \nqcs structures on integrable \DYt modules over
$\b$, and category $\O$ modules over $\g$.

\subsection{Kac--Moody algebras \cite{K}}\label{ss:km-recap}

Throughout this section, we fix a finite set $\bfI$, a matrix $\GCM{A}
=(a_{ij})_{i,j\in\bfI}$ with entries in $\sfk$, and a realisation $(\h,\Pi,\Pi^{\vee})$
of $\GCM{A}$. 
Thus, $\h$ is a $\sfk$--vector space of dimension $2|\bfI|-\rk(\GCM{A})$,
and $\Pi=\{\alpha_i\}_{i\in\bfI}\subset\h^*$, $\Pi^{\vee}=\{\hcor{i}\}_{i\in\bfI}
\subset\h$ are linearly independent subsets such that $\alpha_i(\hcor{j})=
a_{ji}$.

Let $\wt{\g}=\wt{\g}(\GCM{A})$ be the Lie algebra generated by $\h$,
$\{e_i, f_i\}_{i\in\bfI}$ with relations $[h,h']=0$, for any $h,h'\in\h$, and
\[
[h,e_i]=\alpha_i(h)e_i
\qquad
[h,f_i]=-\alpha_i(h)f_i
\qquad
[e_i,f_j]=\drc{ij}\hcor{i}
\]
The Kac--Moody algebra corresponding to $\GCM{A}$ is the Lie algebra
$\g=\g(\GCM{A})=\wt{\g}/I$, where $I$ is the sum of all two--sided ideals in $\wt{\g}$
having trivial intersection with $\h\subset\wt{\g}$. 
If $\GCM{A}$
is a generalised Cartan matrix (\ie $a_{ii}=2$, $a_{ij}\in\IZ_{\leqslant 0}$ if $i\neq j$,
and $a_{ij}=0$ implies $a_{ji}=0$), the ideal $I$ is generated by the Serre relations 
$\mathsf{ad}(e_i)^{1-a_{ij}}(e_j)=0=\mathsf{ad}(f_i)^{1-a_{ij}}(f_j)$ for any $i\neq j$.

Set ${\mathsf{Q}}_+=\bigoplus_{i\in\bfI}\IZ_{\geqslant0}{\alpha}_i\subseteq{\h}^*$,
so that $\g$ has the root space decomposition $\g=\n_-\oplus\h\oplus\n+$, where
$\n_\pm=\bigoplus_{\alpha\in\mathsf{Q}_+\setminus\{0\}}\g_{\pm\alpha}$, and
${\g}_{{\alpha}}=\{X\in{\g}\;|\;[h,X]=\alpha(h)X,\;\forall h\in{\h}\}$. Denote by
${\mathsf{R}}_+=\{\alpha\in\mathsf{Q}_+\;|\; {\g}_{\alpha}\neq0\}$ the set
of positive roots of ${\g}$.

\subsection{Extended \KM algebras \cite{ATL1-2}}\label{ss:ext-KM}

Let $\wh{\g}=\wh{\g}(\GCM{A})$ be the Lie algebra generated by 
$\{e_i,f_i, \hcor{i}, \fcw_i\}_{i\in\bfI}$ with relations $[\hcor{i},
\hcor{j}]=[\fcw_i,\fcw_j]=[\hcor{i},\fcw_j]=0$ for any $i,j\in\bfI$, and 
\[[\hcor{i},e_j]=a_{ij}e_j
\qquad [\hcor{i},f_j]=-a_{ij}f_j
\qquad
[\fcw_i,e_j]=\drc{ij}e_j 
\qquad [\fcw_i,f_j]=-\drc{ij}f_j\]

\begin{definition}
The \emph{extended Kac--Moody algebra} corresponding to $\GCM{A}$
is the $\sfk$--Lie algebra $\ekm{\g}=\wh{\g}/I$, where $I$ is the sum
of all two--sided ideals in $\wh{\g}$ having trivial intersection with the
abelian subalgebra $\ekm{\h}\subset\wh{\g}$ spanned by $\{\hcor{i},\fcw
_i\}_{i\in\bfI}$.
\end{definition}

Let $D$ be the Dynkin diagram of $\g$ and, for any $B\subseteq D$,
let $\ekm{\g}_B\subseteq\ekm{\g}$ be the Lie subalgebra generated
by $\{e_i,f_i,\hcor{i},\fcw_i\}_{i\in B}$ if $B\neq\emptyset$, and $\ekm{\g}
_{\emptyset}=\{0\}$ otherwise.

\begin{proposition}\label{pr:d-alg-Ug}
The extended Kac--Moody algebra
$\ekm{\g}$ is a diagrammatic Lie algebra with 
Lie subalgebras $\ekm{\g}_B$, $B\subseteq D$.
\end{proposition}

\begin{pf}
Clearly, for any $B_1\subseteq B_2$, $\ekm{\g}_{B_1}\subseteq \ekm{\g}_{B_2}$. 
If $B_3\perp B_4$, then for any $i\in B_3, j\in B_4$, $e_i, f_i$ commute with $e_j, f_j$ 
\cite[Lemma 1.6]{K}, and, since $[\hcor{i},e_j]=0=[\hcor{i},f_j]$ and $[\fcw_i,e_j]=0=[\fcw_i,f_j]$,
$[\ekm{\g}_{B_3}, \ekm{\g}_{B_4}]=0$. Finally, if $B=B_1\sqcup B_2$, 
$\ekm{\g}_B=\ekm{\g}_{B_1}\oplus \ekm{\g}_{B_2}$.
\end{pf}

\noindent
\remark The definition of $\ekm{\g}$ takes its cue from \cite{FZ}. Its
use is prompted by the fact that not all (symmetrisable) \KM algebras
are diagrammatic \cite[11]{ATL1-2}.

\subsection{Relation between $\g$ and $\ekm{\g}$}\label{ss:h''}

We show in \cite[\S11.6]{ATL1-2} that $\ekm{\g}$ is non--canonically
a split central extension of $\g$, with a $\rk(\GCM{A})$--dimensional kernel.
Namely, set $r=\rk{\GCM{A}}$, $\ell=|\bfI|$, and assume for simplicity that the first $r$ rows
of $\GCM{A}$ are linearly independent. Let $\h'\subset\h$ be the ($\ell$--dimensional)
span of $\{\hcor{i}\}_{i\in\bfI}$, and $\h''\subset\h$ a subspace with basis
$\{d_j\}_{j=r+1}^\ell$ such that $\alpha_i(d_j)=\delta_{ij}$, $1\leqslant   i\leq\ell$, $r+1\leqslant   j
\leq\ell$. Let $\varpi^{\vee}_i=\sum_{j=1}^r c_{ij}\hcor{j}$ be the fundamental coweights
corresponding to $\{\alpha_1,\dots, \alpha_r\}$. Then, the elements $\gamma_i
\coloneqq \varpi^{\vee}_i-\fcw_i\in\ekm{\h}$, $i=1,\dots, r$, are central in $\ekm{\g}$.
Denote by $\c$ the subspace spanned by $\{\gamma_i\}_{i=1}^r$.

\begin{proposition}[\cite{ATL1-2}]
The choice of the complementary subspace $\h''\subset\h$ determines 
\begin{enumerate}
\item An embedding $\g\subset\ekm{\g}$, mapping $e_i,f_i,\hcor{i}\mapsto
e_i,f_i,\hcor{i}$, and $d_j\mapsto \fcw_j$ for any $i\in\bfI,j=r+1,\ldots,\ell$.
\item A projection $\ekm{\g}\to\ekm{\g}/\c\to\g$
mapping $e_i,f_i,\hcor{i}\to e_i,f_i,\hcor{i}$, $i\in\bfI$, $\fcw_j\mapsto\varpi^\vee_j$, if  
$j=1,\dots, r$, and $\fcw_j\mapsto d_j$, if $j=r+1,\dots, \ell$. 
\item A Lie algebra isomorphism $\ekm{\g}\simeq\g\oplus\c$.
\end{enumerate}
\end{proposition}

\subsection{Root space decomposition of $\ekm{\g}$}

Let $\{\ol{\alpha}_i\}_{i\in\bfI}$ be the linear forms on $\ekm{\h}$ defined
by
\[\ol{\alpha}_i(\hcor{j})=a_{ji}
\aand
\ol{\alpha}_i(\fcw_j)=\delta_{ij}\]
so that, for any ${h}\in\ekm{\h}$, $[{h}, e_i]=\ol{\alpha}_i({h})e_i$ and 
$[{h}, f_i]=-\ol{\alpha}_i({h})f_i$.
Set $\ol{\mathsf{Q}}_+=\bigoplus_{i\in\bfI}\IZ_{\geqslant0}\ol{\alpha}_i\subseteq\ekm{\h}^*$.
Then, $\ekm{\g}$ has the root space decomposition
\[
\ekm{\g}=
\bigoplus_{\substack{{\alpha}\in\ol{\mathsf{Q}}_+\\{\alpha}\neq0}}\ekm{\g}_{{\alpha}}
\oplus\ekm{\h}\oplus
\bigoplus_{\substack{{\alpha}\in\ol{\mathsf{Q}}_+\\{\alpha}\neq0}}\ekm{\g}_{-{\alpha}}
\]
where $\ekm{\g}_{{\alpha}}=\{X\in\ekm{\g}\;|\;[h,X]=\alpha(h)X\;\forall\,h\in\ekm{\h}\}$.
Let $\ol{\cdot}:\mathsf{Q}_+\to\ol{\mathsf{Q}}_+$ be the $\IZ_{\geqslant0}$--linear map 
sending $\alpha_i$ to $\ol{\alpha}_i$, $i\in\bfI$. It follows from the proposition above
that, for any $\alpha\in\mathsf{Q}_+$, $\alpha\neq0$, $\ekm{\g}_{\ol{\alpha}}$ identifies
canonically with $\g_{\alpha}$, and 
\[
\ekm{\g}=
\bigoplus_{\alpha\in\sfR_+}\g_{\alpha}
\oplus\ekm{\h}\oplus
\bigoplus_{\alpha\in\sfR_+}\g_{-\alpha}
\]

\subsection{Symmetrisable Kac--Moody algebras}\label{sss:sym-ext-km}\label{ss:bil on g}

Assume henceforth that the matrix $\GCM{A}$ is symmetrisable, and fix an invertible
diagonal $\GCM{D}=\Diag(d_i)_{i\in\bfI}$ such that $\GCM{AD}$ is symmetric.


Let $\h'\subset\h$ be the span of $\{\hcor{i}\}_{i\in\bfI}$, and $\h''\subset\h$ a complementary
subspace. Then, there is a symmetric, non--degenerate bilinear form $\<\cdot,\cdot\>$
on $\h$, which is uniquely determined by $\<\hcor{i},\cdot\>=d_i\alpha_i(\cdot)$ and
$\<\h'',\h''\>=0$. The form $\iip{\cdot}{\cdot}$ uniquely extends to an invariant symmetric
bilinear form on $\wt{\g}$, and $\iip{e_i}{f_j}=\drc{ij}d_i$. The kernel of this form is
precisely $I$, so that $\iip{\cdot}{\cdot}$ descends to a nondegenerate form on
$\g$.

Set $\bb{\pm}=\h\oplus\bigoplus_{\alpha\in\sfR_+}\g_{\pm\alpha}\subset\g$. 
The bilinear form induces a canonical isomorphism of graded vector spaces 
$\bb{+}\simeq\bb{-}^{\star}$, where 
$\bb{-}^{\star}=\h^*\oplus\bigoplus_{\alpha\in\sfR_+}\g_{-\alpha}^*$,
and determines on $\g$ a natural structure of Lie bialgebra with cobracket 
$\delta:\g\to\g\wedge\g$ given by
\[
\delta|_{\h}=0
\qquad
\delta(e_i)=d_i^{-1} \hcor{i}\wedge e_i
\qquad
\delta(f_i)=d_i^{-1} \hcor{i}\wedge f_i
\]

\subsection{Extended symmetrisable Kac--Moody algebras}\label{sss:dyn-ext-km}


The extended Lie 
algebra ${\ekm{\g}}={\ekm{\g}}(\GCM{A})$ is endowed with an invariant, symmetric 
and non--degenerate bilinear form $\iip{\cdot}{\cdot}:\ekm{\g}\ten\ekm{\g}\to\sfk$ uniquely
determined by the table
\[
\begin{array}{|c||c|c|c|c|}
\hline
\iip{\cdot}{\cdot} & e_j & \hcor{j} & \fcw_j & f_j \\
\hline \hline
e_i & 0 & 0 & 0 & \drc{ij}d_i  \\
\hline
\hcor{i} & 0 & d_ja_{ij} & \drc{ij}d_i & 0 \\
\hline
\fcw_i & 0 & \drc{ij}d_j & 0 & 0 \\
\hline 
f_i & \drc{ij}d_j & 0 & 0 & 0 \\
\hline
\end{array}
\]
If the bilinear form $\<\cdot,\cdot\>$ on $\g$ given in \ref{ss:bil on g} is obtained
from a subspace $\h''\subset\h$ satisfying the requirements of \ref{ss:h''}, the
embedding $\g\subset\ekm{\g}$ corresponding to $\h''$ is compatible with the
bilinear forms and the Lie bialgebra structures.

There is a natural structure of Lie bialgebra on $\ekm{\g}$ with cobracket 
$\delta:\ekm{\g}\to\ekm{\g}\wedge\ekm{\g}$
\[
\delta|_{\ekm{\h}}=0
\qquad
\delta(e_i)=d_i^{-1} \hcor{i}\wedge e_i
\qquad
\delta(f_i)=d_i^{-1} \hcor{i}\wedge f_i
\]

\subsection{Drinfeld double realisation}\label{ss:drinfeld-double-ext}

It is well--known that any symmetrisable Kac--Moody algebra is a central
quotient of the restricted Drinfeld double of its Borel subalgebra.
An analogous result holds for a symmetrisable extended \KMA 
$\olg$ \cite[10.7]{ATL1-2}.
Specifically, consider the Lie algebra $\olg^{(2)}=\olg\oplus\zolh$, where
$\zolh=\olh$ is central in $\olg^{(2)}$, and endow it with the 
inner product $\iip{\cdot}{\cdot}\oplus-\left\iip{\cdot}{\cdot}
\right|_{\zolh\times \zolh}$. Let $\pi_0:\olg\to\olh$ be the projection, and
$\olbtwo_\pm\subset\olgtwo$ the subalgebra
\[\olbtwo_\pm=\{(X,h)\in\olb_\pm\oplus \zolh|\,\pi_0(X)=\pm h\}\]
The projection $\olgtwo\to\olg$ onto the first component restricts
to an isomorphism $\olbtwo_\pm\to\olb_\pm$ with inverse $\olb_\pm\ni X\to
(X,\pm\pi_0(X))\in\olbtwo_\pm$. Then, it is easy to see that $\olgtwo=\olg\oplus\zolh$ 
is the restricted Drinfeld double of $\olbtwo_-\simeq\olb_-$. 

\subsection{Diagrammatic semigroup structure}

Let $\dsg=(\mon,D)$ be the diagrammatic semigroup introduced in \ref
{ex:km-sgp}. Thus, $D$ is the Dynkin diagram of $\g$, $\mon=\sfR_+$
its partial semigroup of positive roots, and $\mon(B)=\sfR_{B,+}$ for
any $B\subseteq D$.

For any $B\subseteq D$, let $\ekm{\b}_{B,-}$ (resp. $\ekm{\b}_{B,+}$)
be the Lie subbialgebra of $\olg_{B}$ generated by $\{\hcor{i},\fcw_{i}, 
f_i\}_{i\in B}$ (resp. $\{\hcor{i},\fcw_i,e_i\}_{i\in B}$). Then,
\[\ekm{\b}_{B,\pm}=
\ekm{\h}_{B}\oplus\ekm{\n}_{B,\pm}\]
where $\ekm{\h}_B\subseteq\ekm{\h}$ is spanned by $\{\hcor{i},\fcw_{i}\}_{i\in
B}$, and $\ekm{\n}_{B,\pm}=\bigoplus_{\alpha\in\sfR_{B,+}}\olg_{\pm\alpha}$,
with $\sfR_{B,+}=\sfR_+\cap\bigoplus_{i\in B}\IZ_{\geqslant 0}\alpha_i$. For any 
$B'\subseteq B$, set
\begin{equation}\label{eq:split}
\ol{\h}_B=\ol{\h}_{B'}\oplus\ol{\h}_{B'}^{\perp}
\aand
\ol{\n}_{B,\pm}=\ol{\n}_{B',\pm}\oplus\ol{\n}_{B',\pm}^{\perp}
\end{equation}
where $\ol{\h}_{B'}^\perp\subseteq\{t\in\ol{\h}_B\,|\,\alpha_i(t)=0,\,i\in B'\}$ 
is a chosen complement to $\ol{\h}_{B'}$ in $\ol{\h}_{B}$, and $\ol{\n}_{B',
\pm}^{\perp}=\bigoplus_{\alpha\in\sfR_{B,+}\setminus\sfR_{B',+}}\olg_{\pm
\alpha}$. For any $B'\subseteq B$, let 
\begin{align*}
i_{0,BB',\pm}:\ol{\h}_{B'}\to\ol{\b}_{B,\pm}
\qquad\qquad
p_{0,B'B,\pm}:\ol{\b}_{B,\pm}\to\ol{\h}_{B'}
\end{align*}
and
\begin{align*}
i_{BB',\pm}:\ol{\b}_{B',\pm}\to\ol{\b}_{B,\pm}
\qquad\qquad
p_{B'B,\pm}:\ol{\b}_{B,\pm}\to\ol{\b}_{B',\pm}
\end{align*}
be the canonical injections and the projections corresponding to the splitting
\eqref{eq:split}. Set  $\dmp{B,\pm}=i_{DB,\pm}\circ p_{BD,\pm}$
and $\dmp{0,B,\pm}=i_{0,DB,\pm}\circ p_{0,BD,\pm}$ in $\End(\ol{\b}_{\pm})$.
Then,
\[
\dmp{B,\pm}=\dmp{0,B,\pm}+\sum_{\alpha\in\sfR_B^+}\dmp{\pm\alpha}
\]
where $\dmp{\pm\alpha}$ denotes the standard idempotent projecting over 
the root space $\olg_{\pm\alpha}$. The following is straightforward.

\begin{proposition}\label{pr:split diagr olb}
The data $(\dmp{\pm\alpha}, \dmp{0,B,\pm})$ induce an $\dsg$--graded
Lie bialgebra structure on $\ol{\b}_{\pm}$ (cf. \ref{ss:0-sgp-ex}). In particular,
the maps $\dmp{B,\pm}$ are morphisms of Lie bialgebras and $\ol{\b}_{\pm}$
is split diagrammatic with Lie subbialgebras $\olb_{B,\pm}$, $B\subseteq D$.
\end{proposition}

By Proposition \ref{pr:transfer}, a braided pre--Coxeter structure $\Cox$ on 
$\CRDYUA{\dsg}{\bullet}$ induces a universal braided pre--Coxeter 
structure $\Cox(\ol{\b}_{\pm})$ on Drinfeld--Yetter modules.
Theorem \ref{th:rigidity} then yields the following

\begin{corollary}
Let $\C_1,\C_2$ be two universal braided pre--Coxeter structures on
$\{\defDY{\b_B}\}_{B\subseteq D}$. Then, there is a universal $1$--isomorphism
$\C_1\to\C_2$, which is unique up to a unique universal $2$--isomorphism.
\end{corollary}

\subsection{The category $\O_{\olg}$}\label{ss:cat-O-lba}

A $\olg$--module $V$ is in category $\O_{\olg}$ if the following holds.
\begin{enumerate}
\item[($\O1$)] $V=\bigoplus_{\lambda\in\olh^*}V_{\lambda}$, where
$V_{\lambda}=\{v\in V|\,h\,v=\lambda(h)v,\,h\in\olh\}$\\[.05ex]
\item[($\O2$)]  $\dim V_{\lambda}<\infty$ for any $\lambda\in\sfP(V)
=\{\lambda\in\olh^*|\,V_{\lambda}\neq0\}$\\[.05ex]
\item[($\O3$)]  $\sfP(V)\subseteq D(\lambda_1)\cup\cdots\cup D
(\lambda_m)$, for some $\lambda_1,\ldots,\lambda_m\in\olh^*$
\end{enumerate}
where $D(\lambda)=\{\mu\in\olh^*\;|\;\mu\leqslant\lambda\}$, with
$\mu\leqslant\lambda$ iff $\lambda-\mu\in\sfQ_+$.
The category $\O_{\olg}$ has a natural symmetric tensor structure
inherited from $\Rep\olg$.

We observed in \ref{ss:drinfeld-double-ext} that the restricted Drinfeld
double of the negative Borel subalgebra $\olb_-$ of $\olg$ is isomorphic
to the trivial central extension $\olgtwo=\olg\oplus\zolh$ of $\g$ by $\zolh
=\olh$. It follows by \ref{ss:pre-DY}--\ref{ss:res double} that the category
of \DYt modules over $\olb_-$ is equivalent to the category $\E_{\olgtwo}
$ of $\olgtwo$--modules, where $\olgtwo=\olg\oplus\zolh$, which carry a
locally finite action of $\olbtwo_+\subset\olgtwo$. This implies the following.

\begin{proposition}
\hfill
\begin{enumerate}
\item
The category $\O_{\olg}$ is isomorphic to the full tensor subcategory
of $\E_{\olgtwo}$ consisting of those modules carrying a trivial action
of $\zolh$ and satisfying, as a module over $\olh\subset\olg\subset
\olgtwo$, the conditions $(\O1)$--$(\O3)$ above. 
\item
Under the equivalence $\E_{\olgtwo}\simeq\DrY{\olb_-}$, 
$\O_{\olg}$ is isomorphic to the full tensor subcategory of $\DrY{\olb_-}$
consisting of those modules $V$ such that the action and the coaction
of $\olh$ on $V$ coincide under $\iip{\cdot}{\cdot}_{\olh}$, \ie
\begin{equation}\label{eq:DY-ss-0}
\pi_V\circ i_0\ten\id_V=\iip{\cdot}{\cdot}_{\olh}\ten\id_V\circ \id_{\olh}\ten 
p_0\ten\id_V\ten\id\circ\id_{\olh}\ten\pi_V^*
\end{equation}
and, as a module over $\olh\subset\olb_-$, $V$ satisfies the conditions
$(\O1)$--$(\O3)$ above.
\end{enumerate}
\end{proposition}

\subsection{Pre--Coxeter structures and category $\O_{\olg}$}\label{ss:cat-O-infty}

Condition ($\O2$)  on the finite--dimensionality of weight spaces in
\ref{ss:cat-O-lba} is not stable under restriction from $\olg=\olg_D$
to $\olg_B$ if $B\subsetneq D$, which makes category $\O_{\olg}$
unsuitable to the axiomatic of braided pre--Coxeter structures. We
therefore omit it, and denote  by $\O_{\infty,\olg}$ the category of
$\olg$--modules satisfying conditions ($\O1$) and ($\O3$). Proposition
\ref{ss:cat-O-lba} shows that $\O_{\infty,\olg}$ is a full subcategory
of $\DrY{\olb_-}$. Moreover, any universal braided pre--Coxeter
structure on $\{\defDY{\olb_{B,-}}\}_{B\subseteq D}$ restricts to
one on $\{\O_{\infty,\olg_B}^{\hbar}\}_{B\subseteq D}$.

\subsection{Braid group actions}\label{ss:qC-LBA}

Assume now that $\sfA$ is a symmetrisable generalised Cartan matrix,
let $W$ be the corresponding Weyl group with set of simple reflections
$\{s_i\}_{i\in\bfI}$, and set $\ulm=(m_{ij})$, where $m_{ij}$ is the order
of $s_is_j$ in $W$. 

Let $\C_{\olg}^{\sint}$ be the category of integrable $\olg$--modules,
 \ie $\olh$--semisimple
modules endowed with a locally nilpotent action of the elements $\{e_i,f_i\}_{i\in\bfI}$.
Let  $\wh{\C}_{\olg}^{\sint}$ be the algebra $\sfEnd{\C_{\olg}^{\sint}\to\vect}$ and,
for any $i\in D$, denote by $\wt{s}_i\in\wh{\C}_{\olg}^{\sint}$ the triple exponential
\[\wt{s}_i=\exp(e_i)\cdot\exp(-f_i)\cdot\exp(e_i).\]
It is well--known (cf. \cite{tits}) that these satisfy the generalised braid relations
\eqref{eq:gen-braid}.

Let $\hDrY{\olb_-}{\sint,0}$ be the category of integrable Drinfeld--Yetter $\olb_-
$--modules in $\DrY{\olb_-}$, \ie $\olh$--diagonalisable, endowed with a locally
nilpotent action of the elements $\{f_i\}_{i\in D}\subseteq\olb_-$, and satisfying
\eqref{eq:DY-ss-0}, so as to give rise to integrable modules over $\olg$.
In particular, the triple exponential $\wt{s}_i$ acts on the objects in $\hDrY{\olb_-}{\sint,0}$ and
the subcategory of integrable modules in $\O_{\infty,\olg}$, denoted $\O^{\sint}_{\infty,\olg}$, 
is isomorphic to a braided tensor subcategory of $\hDrY{\olb_-}{\sint,0}$.

\subsection{Universal braided \nqcs structures on Kac--Moody algebras}
Set $\ekm{\b}=\ekm{\b}_-$.
Let $\intDY{\ekm{\b}}$ be the category of integrable deformation Drinfeld--Yetter
$\ekm{\b}$--modules. As usual, we denote by $\hDYA{B}{n}$ (resp. $\hDYA{B,0}{n}$) 
the algebra of endomorphisms of the forgetful functor $\ff_B^{\boxtimes n}:(\defDY{\ekm{\b}_B})^{n}\to{\khvect}$
(resp. $\ff_{B,0}^{\boxtimes n}:(\defDY{\ekm{\b}_B})^{ n}\to\defDY{\ekm{\h}}{}$).
For any $X\in\hDYA{B}{n}$, we denote by $\sfp(X)$ the induced endomorphism
of the forgetful functor $(\intDY{\ekm{\b}_B})^{n}\to{\khvect}$.

\begin{definition}
A braided Coxeter structure of type $(D,\ulm)$ with diagrammatic
categories $\{\intDY{\olb_B}\}_{B\subseteq D}$ is
\emph{good} (resp. \emph{universal}) if the underlying braided pre--Coxeter structure is induced 
by a weight--zero\footnote{A braided pre--Coxeter structure on $\{\defDY{\ekm{\b}_B}\}_{B\subseteq D}$ 
is \emph{weight--zero} if it is defined over the category $\defDY{\ekm{\h}}$.} (resp. universal) 
pre--Coxeter structure on $\defDY{\ekm{\b}}$, and its local monodromies have the form
\begin{equation}\label{eq:locmon}
S_i=\wt{s}_i\cdot\sfp(\ul{S}_i)
\end{equation}
where $\ul{S}_i\in\hDYA{\vrtx{i},0}{1}$, $\ul{S}_i=1\mod\hbar$, and $\wt{s}_i=\exp(e_i)\exp(-f_i)\exp(e_i)$.
\end{definition}

\noindent\remark 
It follows from Proposition \ref{sss:wtzero} (2) that any universal braided Coxeter structure on $\{\intDY{\olb_B}\}_{B\subseteq D}$
is good. Moreover, it is important to observe that, since $\defDY{\ekm{\b}_i}{}\simeq\Rep\hext{U\olgtwo_i}$ 
with $\g_i=\mathfrak{sl}_2^{\alpha_i}$, we have $\hDYA{\vrtx{i}}{n}=\hext{(U\olgtwo_i)^{\ten n}}$. 
In particular, $\sfp(\ul{S}_i)$ is an element in $\hext{(U\olg_i)^{\olh_i}}$.\\

The twisting of a braided pre--Coxeter structure on $\{\defDY{\olb_B}\}_{B\subseteq D}$
extends to a twisting of a braided \nqcs structure on $\{\intDY{\olb_B}\}_{B\subseteq D}$.
Namely, if $(\sfCox, S_i)$ is a good braided \nqcs structure on 
the latter, where $\sfCox$ is the corresponding braided pre--Coxeter 
structure  on the former, and $(u, F)$ is a weight--zero twist of $\sfCox$, 
then
\[
(\sfCox, S_i)_{(u,F)}\coloneqq (\sfCox_{(u,F)}, S_i^u\coloneqq u_{\vrtx{i}}\cdot S_i\cdot u^{-1}_{\vrtx{i}})
\] 
is a good braided \nqcs structure on $\{\intDY{\olb_B}\}_{B\subseteq D}$. Moreover,
the representations of $B_D^{\ul{m}}$ corresponding to $(\sfCox, S_i)$
and $(\sfCox, S_i)_{(u,F)}$ are equivalent, \ie $(\lambda_{(u,F)})_{\F}=\sfAd{u_{\F}}\circ\lambda_{\F}$.

\subsection{The local monodromies $S_i$}\label{ss:local-monodromies}

Let $i\in D$ be a fixed vertex, and set $\ekm{\h}_i=\sfk \hcor{i}\oplus\sfk\fcw_i$. 

\begin{lemma}\label{le:monod}
Let $S_{i}^{(1)},S_{i}^{(2)}$ be two elements of  the form \eqref{eq:locmon}
which satisfy the relation
	\begin{equation}\label{eq:coprod-2}
	J(S_{i}^{(k)})_{12}J^{-1}=J\cdot R_{21}\cdot J_{21}^{-1}(S_{i}^{(k)})_1(S_{i}^{(k)})_2
	\end{equation}
for some $R,J\in1^{\ten 2}+\hbar\cdot\hext{((U\olg_i)^{\ten 2})^{\olh_i}}$.
Then, there exist unique $u,v\in \hext{\sfk}$ such that
	\[
	S_{i}^{(2)}=\sfAd{e^{u\hcor{i}+v\fcw_i}}(S_{i}^{(1)})
	\]
\end{lemma}

\begin{pf}
	Let $S_{i}^{(k)}=\wt{s}\cdot\sfp(\ul{S}_{i}^{(k)})$, where
	\[
	\sfp(\ul{S}_{i}^{(k)})=1+\sum_{N\geqslant0}\hbar^ns_n^{(k)}\qquad s_n^{(k)}\in (U\olg_i)^{\olh}
	\]
	The identity above reads
	\begin{equation}\label{eq:coxcoprod2}
	\sfp(\ul{S}_{i}^{(k)})_{12}\cdot J^{-1}=R_{21}^{\theta}\cdot (J^{-1}_{21})^{\theta}\cdot
	\sfp(\ul{S}_{i}^{(k)})_1\cdot\sfp(\ul{S}_{i}^{(k)})_2
	\end{equation}
	where $\theta$ is the Chevalley involution, acting as $-1$ on $\ekm{\h}_i$.
	
	We construct two sequences
	\[
	u_n=\sum_{k=0}^n a_k\hbar^n
	\aand
	v_n=\sum_{k=0}^n b_k\hbar^n\qquad a_k,b_k\in\sfk
	\]
	such that
	\begin{equation}\label{eq:indmonod}
	{S}_{i}^{(2)}=e^{u_n\hcor{i}+v_n\fcw_i}{S}_{i}^{(1)}e^{-u_n\hcor{i}-v_n\fcw_i}\quad\mod\hbar^{n+1}
	\end{equation}
	Since $S_2=S_1=\wt{s}$ modulo $\hbar$, we may assume $a_0=0=b_0$.
	Assume therefore $a_k, b_k$ defined for $k=0,1,\dots, n$ for some $n\geqslant0$.
	Let $(S_i^{(1)})^\prime$ be given by the right--hand side of \eqref{eq:indmonod}, so that
	\[
	\sfp(\ul{S}_i^{(2)})=\sfp(\ul{S}_i^{(1)})^\prime+\hbar^{n+1}\eta\quad\mod\hbar^{n+2}
	\]
	for some $\eta\in(U\olg_i)^{\olh_i}$. One readily checks that $\sfp(\ul{S}_i^{(1)})^\prime$ satisfies
	\eqref{eq:coxcoprod2}, since $e^{u_n\hcor{i}+v_n\fcw_i}$ is group--like element in $\hext{U\olh_i}$.
	Subtracting from this the coproduct identity for $\sfp(\ul{S}_i^{(2)})$, and computing modulo
	$\hbar^{n+2}$, we find that
	\[
	d_H(\eta)=\eta_2-(\eta)_{12}+\eta_1=0
	\]
	Therefore, $\eta$ is a primitive element in $(U\olg_i)^{\olh_i}$.
	It follows $\eta=c\cdot \hcor{i}+d\cdot\fcw_i$, for some $c,d\in\sfk$. Then for $a_{n+1}=-c/2$,
	$b_{n+1}=-d/2$, we get 
	\[
	e^{(a_{n+1}\hcor{i} + b_{n+1}\fcw_i)\hbar^{n+1}}({S}_i^{(1)})^\prime 
	e^{-(a_{n+1}\hcor{i} + b_{n+1}\fcw_i)\hbar^{n+1}}=
	({S}_i^{(1)})^\prime e^{(c\hcor{i}+d\fcw_i)\hbar^{n+1}}=\ul{S}_i^{(2)}
	\]
	modulo $\hbar^{n+2}$.
	By induction, one gets $u,v\in\hext{\sfk}$ such that
	\[
	S_i^{(2)}=\sfAd{e^{u\hcor{i}+v\fcw_i}}(S_i^{(1)})
	\]
\end{pf}

Since the coproduct identity \eqref{eq:coxcoprod} has the form \eqref{eq:coprod-2},
where $R=R_i\in\hDYA{\vrtx{i},0}{2}$ is an $R$--matrix and $J=J_i\in\hDYA{\vrtx{i},0}{2}$ 
is a twist, we get the following

\begin{corollary}
	Up to gauge transformation, a good (resp. universal) braided pre--Coxeter structure 
	on $\{\intDY{\ekm{\b}_B}\}_{B\subseteq D}$ can be completed to at most one universal braided \nqcs structure.
\end{corollary}
\subsection{Coxeter structures on extended Kac--Moody algebras}

Let $\ekm{\g}$ be an extended symmetrisable Kac--Moody algebra with negative
Borel subalgebra $\ekm{\b}$ and Dynkin diagram $D$, and $\intDY{\ekm{\b}}$ the 
deformation category of integrable Drinfeld--Yetter $\ekm{\b}$--modules.

The following is the main result of this paper.

\begin{theorem}\label{thm:main}
Let $k=1,2$, and
\[\sfCox_k=(\Phi_B^{(k)},J^{(k)}_{\F},\DCPA{\F}{\G}^{(k)},S_i^{(k)})\]
two universal braided \nqcs structures on $\{\intDY{\ekm{\b}_B}\}_{B\subseteq D}$
corresponding to a fixed labeling $\ulm$ on $D$. 
Then,
\begin{enumerate}
\item There is a twist $(u,F)$ such that
$\sfCox_2=(\sfCox_1)_{(u,F)}$.
\item The twist $(u,F)$ is unique up to a unique gauge $a$.
\end{enumerate}
\end{theorem}

\begin{pf}
Let $(\sfCox_k, \{S_{i}^{(k)}\})$, $k=1,2$, be two universal Coxeter
structures on $\DrY{\ekm{\b}}^{\intm}$. By \ref{thm:proprigidity},
there is a universal twisting $(u,F)$ such that
\[\sfCox_2=(\sfCox_1)_{(u,F)}\]
where $u$ is uniquely determined, and $F$ is uniquely determined
up to multiplication with elements of the form $(a_B)_1^{-1}(a_B)_2
^{-1}(a_B)_{12}$, where $a_B$ belongs to the center of $\hDYA{\ekm
{\b}}{n}$. Therefore, $S_i^{(2)}$ and $(S_i^{(1)})_{a}$ are two Coxeter
extensions of $\sfCox_2$. By Lemma \ref{le:monod}, there is a unique
tuple $\ul{v}=(v_1,\dots,v_n, v'_1,\dots, v'_n)$, $v_i, v'_i\in\hext{\sfk}$,
such that
\[\sfAd{e^{v_i\hcor{i}+v_i'\fcw_i}}(S_i^{(1)})_u=S_i^{(2)}\]
and
\[(\sfCox_2, \{S_i^{(2)}\})=(\sfCox_1, \{S_i^{(1)}\})_{(\ul{v}\circ u, F)}\]
The theorem is proved.
\end{pf}

Let $\O_{\infty,\olg}^{\hbar,\sint}$ be the category of deformation, integrable, category $\O_{\infty}$
$\olg$--modules. From \ref{ss:cat-O-infty} and \ref{ss:qC-LBA}, we get the following

\begin{corollary}
Any two universal braided Coxeter structures on $\{\O_{\infty,\olg_B}^{\hbar,\sint}\}_{B\subseteq D}$
are twist equivalent, with respect to a universal twist, which is unique up to a 
unique gauge.
\end{corollary}

\noindent
\remark Since the labeling of the diagram $D$ plays no role in the
proof of the rigidity of braided \nqcs structures, the latter
yields the following strenghtening of Theorem \ref{thm:main}. If 
$\sfCox_1,\sfCox_2$ are two universal braided
\nqcs structures on $\intDY{\ekm{\b}}$ corresponding
to the labelings $\{m^1_{ij}\}$, $\{m^2_{ij}\}$, then there is a
twist $(u,F)$ such that $\sfCox_2=(\sfCox_1)_{(u,F)}$,
which is unique up to a unique gauge. In particular, the local
monodromies of $\sfCox_1,\sfCox_2$ satisfy the
braid relations \wrt the labeling $\{\min(m^1_{ij},m^2_{ij})\}$.

\subsection{Coxeter structures on diagrammatic Kac--Moody algebras}

We mention in Remark \ref{ss:ext-KM} that the definition of extended
Kac--Moody algebra is prompted by the fact that not all Kac--Moody 
algebras are diagrammatic and, more specifically, not all symmetrisable
Kac--Moody algebras are graded, as Lie bialgebras, over the diagrammatic 
semigroup $\dsg$ associated to their root system (cf.~\ref{ex:km-sgp}). 
Nonetheless, one observes easily that a large class of (non--extended) 
symmetrisable Kac--Moody algebras are $\dsg$--graded, including those 
of finite, affine, and hyperbolic type. In \cite[11]{ATL1-2}, we refer to these 
as \emph{Cartan diagrammatic} symmetrisable Kac--Moody algebras. It is 
evident that the results described above hold verbatim for these Lie 
bialgebras. Therefore, we get the following

\begin{theorem}
Let $\g$ be a Cartan diagrammatic symmetrisable Kac--Moody algebra
(in particular, of finite, affine, or hyperbolic type). Any two universal braided 
Coxeter structures on $\{\O_{\infty,\g_B}^{\hbar,\sint}\}_{B\subseteq D}$ are 
twist equivalent, with respect to a universal twist, which is unique up to a 
unique gauge.
\end{theorem}



\begin{thebibliography}{99}
%
\bibitem{ATL1}
A. Appel, V. Toledano Laredo, 
\emph{A $2$--categorical extension of Etingof--Kazhdan quantisation},
Selecta Math. (N.S.) {\bf 24} (2018), 3529--3617.
%
\bibitem{ATL1-2}
A. Appel, V. Toledano Laredo,
\emph{Coxeter categories and quantum groups}, Selecta
Math. (N.S.) (2019), in press. \href{http://arxiv.org/abs/1610.09741}{\sf arXiv:1610.09741},
63 pp.
%
\bibitem{ATL3}
A. Appel, V. Toledano Laredo,  
\emph{Monodromy of the Casimir connection of a symmetrisable Kac--Moody algebra}.
\href{http://arxiv.org/abs/1512.03041}{\sf arXiv:1512.03041}, 48 pp.

\bibitem{baez}
J. Baez, T. Trimble, \emph{Schur functors I-II}, 
\href{https://ncatlab.org/nlab/show/Schur+function}{https://ncatlab.org/nlab/show/Schur+function},
\href{https://ncatlab.org/nlab/show/Schur+functor}{https://ncatlab.org/nlab/show/Schur+functor},
nLab articles.

\bibitem{brieskorn}
E. Brieskorn, \emph{Die Fundamentalgruppe des Raumes der regul\"{a}ren Orbits einer endlichen komplexen Spiegelungsgruppe}, 
Invent. Math. {\bf 12} (1971), 57--61.

\bibitem{brieskorn-saito}
E. Brieskorn, K. Saito, \emph{Artin--Gruppen und Coxeter--Gruppen},
Invent. Math. {\bf 17} (1972), 245--271.

\bibitem{DCP2}
C. De Concini, C. Procesi, 
\emph{Hyperplane arrangements and holonomy equations}, 
Selecta Math. (N.S.) {\bf 1} (1995), 495--535.

\bibitem{deligne}
P. Deligne, \emph{Les immeubles des groupes de tresses g\'{e}n\'{e}ralis\'{e}s}, Invent. Math. {\bf 17} (1972), 273--302.

\bibitem{drin-4}
V. G. Drinfeld, \emph{Quasi--Hopf algebras}. 
Leningrad Math. J. \textbf{1} (1990), 1419--1457.

\bibitem{e1}
B. Enriquez, 
\emph{On some universal algebras associated to the category of Lie bialgebras}. 
Adv. Math. \textbf{164} (2001), 1--23.
%
\bibitem{e2}
B. Enriquez, 
\emph{Quantization of Lie bialgebras and shuffle algebras of Lie algebras}. 
Selecta Math. (N.\,S.) \textbf{7} (2001), 321--407.
%
\bibitem{e3}
B. Enriquez, 
\emph{A cohomological construction of quantization functors of Lie bialgebras}. 
Adv. Math. \textbf{197} (2005), 430--479.
%
\bibitem{ee}
B. Enriquez, P. Etingof,
\emph{On the invertibility of quantization functors}, 
J. Algebra \textbf{289} (2005), 321--345.
%
\bibitem{ek-1}
P. Etingof, D. Kazhdan, 
\emph{Quantization of Lie bialgebras. I}. 
Selecta Math. (N.\,S.) \textbf{2} (1996), 1--41.
%
\bibitem{ek-2}
P. Etingof, D. Kazhdan, 
\emph{Quantization of Lie bialgebras. II}. 
Selecta Math. (N.\,S.) \textbf{4} (1998), 213--231.
%
\bibitem{evseev}
A. E. Evseev, \emph{A survey of partial groupoids}. In Ben Silver, \emph{Nineteen Papers on Algebraic Semigroups}. 
American Mathematical Society, 1988. 
%
\bibitem{FZ}
B. Feigin, A. Zelevinsky, 
\emph{Representations of contragradient Lie algebras and the Kac-Macdonald identities}.
Representations of Lie groups and Lie algebras (Budapest, 1971), 25--77, Akad. 
Kiado, Budapest, 1985.
%
\bibitem{FMTV}
G. Felder, Y. Markov, V. Tarasov, A. Varchenko,
\emph{Differential Equations Compatible with KZ Equations},
Math. Phys. Anal. Geom. \textbf{3} (2000), 139--177.
%
\bibitem{K}
V. Kac, 
\emph{Infinite--dimensional Lie algebras}. 
Cambridge University Press (1991), 3rd edition.
%
\bibitem{La}
F. W. Lawvere, {\it Functorial semantics of algebraic theories},
Proc. Nat. Acad. Sci. U.S.A. {\bf 50} (1963), 869--872.
%
\bibitem{mac}
S. MacLane,
{\it Categorical algebra}, Bull. Amer. Math. Soc. {\bf 71}
(1965), 40--106.
%
\bibitem{pol}
L. Positselsky,
\emph{Letter to M. Finkelberg and R. Bezrukavnikov}. 1995.
%
\bibitem{reut}
C. Reutenauer, 
\emph{Free Lie algebras}. London Mathematical Society Monographs.
New Series, 7. Oxford Science Publications. The Clarendon Press,
Oxford University Press, 1993.

\bibitem{tits}
J. Tits,
\emph{Normalisateurs de tores I. Groupes de Coxeter \'{e}tendus},
J. Algebra \textbf{4} (1966), 96--116.

\bibitem{vtl-2}
V. Toledano Laredo, 
\emph{A Kohno--Drinfeld theorem for quantum Weyl groups}. 
Duke Math. J. \textbf{112} (2002), 421--451.
%
\bibitem{vtl-3}
V. Toledano Laredo, 
\emph{Cohomological construction of relative twists}. 
Adv. Math. \textbf{210} (2007), 375--403.
%
\bibitem{vtl-4}
V. Toledano Laredo, 
\emph{Quasi--Coxeter algebras, Dynkin diagram cohomology and quantum Weyl groups}. 
Int. Math. Res. Pap. IMRP 2008, Art. ID rpn009, 167 pp.
%
\bibitem{vtl-6}
V. Toledano Laredo,  
\emph{Quasi--Coxeter quasitriangular quasibialgebras and the Casimir connection}. 
\href{http://arxiv.org/abs/1601.04076}{\sf arXiv:1601.04076}, 55 pp.
\end{thebibliography}
\end{document}